\documentclass{memo-l} 
\usepackage{amsxtra, amssymb}
\usepackage[mathscr]{eucal}

\input xy
\xyoption{matrix}\xyoption{frame}

\theoremstyle{plain}
\newtheorem{theorem}{Theorem}[chapter]
\newtheorem{lemma}[theorem]{Lemma}
\newtheorem{definition-theorem}[theorem]{Definition-Theorem}
\newtheorem{proposition}[theorem]{Proposition}
\newtheorem{corollary}[theorem]{Corollary}
\newtheorem{convention}[theorem]{Convention}
\newtheorem*{mainth1gen*}{Main Theorem I: General Form}
\newtheorem*{thA*}{Theorem A}
\newtheorem*{mainth2spec*}{Main Theorem II: Specific Form}
\newtheorem*{mainth3add*}{Main Theorem III: Additional Clusters}
\newtheorem*{step1*}{Step 1}
\newtheorem*{step2*}{Step 2}
\newtheorem*{step3*}{Step 3}
\newtheorem*{step4*}{Step 4}

\theoremstyle{definition}
\newtheorem{definition}[theorem]{Definition}
\newtheorem{example}[theorem]{Example}
\newtheorem{remark}[theorem]{Remark}
\newtheorem*{defA*}{Definition A}

\newcommand \bth[1] { \begin{theorem}\label{t#1} }
\newcommand \ble[1] { \begin{lemma}\label{l#1} }

\newcommand \bpr[1] { \begin{proposition}\label{p#1} }
\newcommand \bco[1] { \begin{corollary}\label{c#1} }
\newcommand \bde[1] { \begin{definition}\label{d#1}\rm }
\newcommand \bex[1] { \begin{example}\label{e#1}\rm }
\newcommand \bre[1] { \begin{remark}\label{r#1}\rm }
\newcommand \bcj[1] { \begin{conjecture}\label{j#1}\rm }

\newcommand \bnota[1] { \begin{notation}\label{n#1}\rm }
\renewcommand {\eth} { \end{theorem} }
\newcommand {\ele} { \end{lemma} }

\newcommand {\epr} { \end{proposition} }
\newcommand {\eco} { \end{corollary} }
\newcommand {\ede} { \end{definition} }
\newcommand {\eex} { \end{example} }
\newcommand {\ere} { \end{remark} }
\newcommand {\ecj} { \end{conjecture} }

\newcommand {\enota} { \end{notation} }
\newcommand \thref[1]{Theorem \ref{t#1}}
\newcommand \leref[1]{Lemma \ref{l#1}}
\newcommand \prref[1]{Proposition \ref{p#1}}
\newcommand \coref[1]{Corollary \ref{c#1}}

\newcommand \deref[1]{Definition \ref{d#1}}
\newcommand \exref[1]{Example \ref{e#1}}
\newcommand \reref[1]{Remark \ref{r#1}}

\numberwithin{section}{chapter}
\numberwithin{equation}{chapter}
\newcommand{\nn}{\hfill\nonumber}
\newcommand \circled[1]{\xymatrix{#1 \save**\frm{o} \restore}}

\newcommand \Rset {{\mathbb R}}         

\newcommand \KK {{\mathbb K}}
\newcommand \Ffield {{\mathbb F}}
\newcommand \Zset {{\mathbb Z}}

\newcommand \Qset {{\mathbb Q}}

\newcommand \B  {{\mathcal{B}}}               
\newcommand \Abb {{\mathcal{A}}}

\newcommand \FF {{\mathcal{F}}}

\newcommand \Tbb {{\mathcal{T}}}

\newcommand \CC {{\mathcal{C}}}

\newcommand \QQ {{\mathcal{Q}}}
\newcommand \OO {{\mathcal{O}}}
\newcommand \PP {{\mathcal{P}}}

\newcommand \UU {{\mathcal{U}}}
\newcommand \RR {{\mathcal{R}}}
\newcommand \Sbb {{\mathcal{S}}}

\newcommand \qb {{\bf{q}}}
\newcommand \rbf {{\bf{r}}}

\newcommand \tb {{\bf{t}}}
\newcommand \ex {{\bf{ex}}}
\newcommand \inv {{\bf{inv}}}
\newcommand \Scr {{\mathscr{S}}}
\newcommand \al {\alpha}
\newcommand \be {\beta}
\newcommand \vpi {\varpi}
\newcommand \la {\lambda}

\newcommand \om {\omega}
\newcommand \Om {\Omega}
\newcommand \ga {\gamma}
\newcommand \Ga {\Gamma}
\newcommand \Sig {\Sigma}

\newcommand \vp {\varphi}
\newcommand \ep {\epsilon}
\newcommand \De {\Delta}

\newcommand \mt  {\mapsto}

\newcommand \hra {\hookrightarrow}

\newcommand \sy  {\ast}                         
\newcommand \bu  {\bullet}                 
\newcommand \ci  {\circ}
\newcommand \rcor {\rangle}
\newcommand \lcor {\langle}
\newcommand \ol {\overline}
\newcommand \wt {\widetilde}
\newcommand \wh {\widehat}

\newcommand \id { {\mathrm{id}} }

\newcommand \sign { {\mathrm{sign}} }

\DeclareMathOperator \rk { {\mathrm{rk}} }

\newcommand \g  {\mathfrak{g}}   

\newcommand \slfrak {\mathfrak{sl}}

\newcommand \n  {\mathfrak{n}}

\DeclareMathOperator \Span { {\mathrm{Span}} }
\DeclareMathOperator \Aut { {\mathrm{Aut}} }

\DeclareMathOperator \charr { {\mathrm{char}} }

\DeclareMathOperator \ad { {\mathrm{ad}} }

\DeclareMathOperator \Ker { {\mathrm{Ker}} }

\DeclareMathOperator \GKdim {{\mathrm{GK \, dim}}}

\DeclareMathOperator \lt  { {\mathrm{lt}} }

\DeclareMathOperator \supp { {\mathrm{supp}} }
\DeclareMathOperator \range { {\mathrm{range}} }
\DeclareMathOperator \Fract { {\mathrm{Fract}} }

\renewcommand \max { {\mathrm{max}} }
\newcommand\kx{\KK^*}

\newcommand\HH{{\mathcal{H}}}
\newcommand\xh{X(\HH)}

\DeclareMathOperator \chr {char}
\newcommand \Znn {\Zset_{\ge 0}}

\newcommand \Hmax {\HH_{\max}}
\newcommand\ismi{[i,s^m(i)]}

\newcommand\lab{{\boldsymbol \lambda}}
\newcommand \gab {{\boldsymbol \gamma}}
\newcommand\nub{{\boldsymbol \nu}}
\newcommand{\Det}{\widetilde{\De}}
\newcommand{\OqM}{\OO_q(M_{m,n}(\KK))}

\makeindex

\begin{document}

\frontmatter

\title[Quantum cluster algebras]
{Quantum cluster algebra structures on \\ quantum nilpotent algebras}

\author[K. R. Goodearl]{K. R. Goodearl}
\address{
Department of Mathematics \\
University of California\\
Santa Barbara, CA 93106 \\
U.S.A.
}
\email{goodearl@math.ucsb.edu}
\thanks{The research of K.R.G. was partially supported by NSF grant DMS-0800948.}

\author[M. T. Yakimov]{M. T. Yakimov}
\address{
Department of Mathematics \\
Louisiana State University \\
Baton Rouge, LA 70803 \\
U.S.A.
}
\email{yakimov@math.lsu.edu}
\thanks{The research of M.T.Y. was partially supported  by NSF grants DMS-1001632 and DMS-1303038.}

\date{23 December 2013}



\dedicatory{To the memory of Andrei Zelevinsky}

\begin{abstract} 
All algebras in a very large, axiomatically defined class of quantum nilpotent 
algebras are proved to possess quantum cluster algebra structures under mild conditions. 
Furthermore, it is shown that these quantum cluster algebras always equal the corresponding 
upper quantum cluster algebras. Previous approaches to these problems for the 
construction of (quantum) cluster algebra structures on (quantized) 
coordinate rings arising in Lie theory were done on a case by case 
basis relying on the combinatorics of each concrete family. 
The results of the paper have a broad range of applications to these problems, 
including the construction of quantum cluster algebra structures on
quantum unipotent groups and quantum double Bruhat cells (the Berenstein--Zelevinsky 
conjecture), and treat these problems from a unified perspective. 
All such applications also establish equality 
between the constructed quantum cluster algebras and their upper counterparts. 
The proofs rely on Chatters' notion of noncommutative unique 
factorization domains. Toric frames are constructed 
by considering sequences of homogeneous prime elements 
of chains of noncommutative UFDs (a generalization 
of the construction of Gelfand--Tsetlin subalgebras) and mutations are obtained 
by altering chains of noncommutative UFDs.
Along the way, an intricate (and unified) 
combinatorial model for the homogeneous prime elements
in chains of noncommutative UFDs and their alterations is developed.
When applied to special families, this recovers the combinatorics 
of Weyl groups and double Weyl groups previously used in 
the construction and categorification of cluster algebras. 
It is expected that this combinatorial model of sequences of homogeneous prime elements will 
have applications to the unified categorification of quantum nilpotent algebras.
\end{abstract}

\maketitle

\tableofcontents

\mainmatter


\chapter{Introduction}
\label{intro}
\section{Quantum cluster algebras and a general formulation of the main theorem}
\label{1.1} Cluster algebras were invented by Fomin and Zelevinsky in \cite{FZ}
based on a novel construction of producing infinite generating sets 
via a process of mutation. The initial goal was to set up 
a combinatorial framework for studying canonical bases and total positivity \cite{Fo}.
Remarkably, for the past twelve years cluster algebras and the procedure of mutation 
have played an important role in a large number of diverse areas of mathematics, including 
representation theory of finite dimensional algebras,
combinatorial and geometric Lie theory, Poisson geometry,
integrable systems, topology,
commutative and noncommutative algebraic geometry, and mathematical physics.
The quantum counterparts of cluster algebras were introduced by Berenstein 
and Zelevinsky in \cite{BZ}. We refer the reader to the recent surveys
\cite{GLSrev,Ke,R,W} and the book \cite{GSVb} for more information on some 
of the abovementioned aspects of this theory.

A major direction in the theory of cluster algebras is to prove 
that important (quantized) coordinate rings of algebraic varieties arising from Lie theory 
admit (quantum) cluster algebra structures or upper (quantum) cluster algebra structures.
For example, the Berenstein--Zelevinsky conjecture \cite{BZ} states that the quantized 
coordinate rings of double Bruhat cells in all finite dimensional simple Lie groups admit 
explicit quantum cluster algebra structures. The motivation for this type 
of problem is that once (quantum) cluster algebra structures are constructed on 
families of (quantized) coordinate rings, they can then be used in the study of canonical 
bases of those rings. In the classical case, a cluster algebra structure on the 
coordinate ring of a variety can be used to investigate its totally positive 
part.
 
Going back to the general problem, a second part asks if the constructed 
upper (quantum) cluster algebra equals the corresponding (quantum) cluster algebra.
For example, ten years ago Berenstein, Fomin and Zelevinsky proved \cite{BFZ} that 
the coordinate rings of double Bruhat cells in all simple algebraic groups 
admit upper cluster algebra structures. Yet it was unknown if these 
upper cluster algebras equal the corresponding cluster algebras, i.e., 
if the coordinate rings of double Bruhat cells are actually cluster algebras.

Previous approaches to the above problems relied on a construction 
of an initial seed and some adjacent seeds in terms of (quantum) minors and
related regular functions
\cite{BFZ,BZ,GLSh0,GLSh}. After that point, two different approaches were followed.
The first one, due to Berenstein, Fomin and Zelevinsky \cite{BFZ}, 
used the methods of unique factorization domains to prove that the coordinate rings 
under consideration are upper cluster algebras. It was first 
applied to coordinate rings of double Bruhat cells \cite{BFZ}. 
This approach was developed further in \cite{GSVb,GSV-is} and \cite{GLSh2}.
The second approach was via the 
construction of a categorification based on concrete combinatorial data from Weyl 
groups and then to prove that the corresponding (quantum) cluster algebra
equals the (quantized) coordinate ring under consideration.
This approach is due to Gei\ss--Leclerc--Schr\"oer \cite{GLSh0,GLSh}, who applied it 
to the coordinate rings of the unipotent groups $U_+ \cap w(U_-)$ 
and the quantum Schubert cell algebras 
$\UU_q(\n_+ \cap w(\n_-))$ (also called quantum unipotent groups)
for symmetric Kac--Moody groups $G$, where $w$ is a Weyl group 
element.

In both of the above approaches, one relied on specific data in terms of 
Weyl group combinatorics
for the concrete family of coordinate rings.
Moreover, the initial (quantum) seeds were built via a direct
construction by considering (quantum) minors. 

The goal of this paper is to present a new algebraic approach to quantum cluster 
algebras based on noncommutative ring theory. We produce a general construction of quantum cluster algebra 
structures on a broad class of algebras and construct 
initial clusters and mutations in a uniform and intrinsic way, without ad hoc constructions 
with quantum minors.
We first state the main theorem of the paper in a general form. 
The following sections contain a precise formulation of it.

\begin{mainth1gen*}
Each algebra in a very large, axiomatically defined 
class of quantum nilpotent algebras admits a quantum cluster algebra 
structure. Furthermore, for all such algebras, the latter equals the corresponding 
upper quantum cluster algebra.
\end{mainth1gen*}

The theorem has a broad range of applications because many important families of algebras
fall within this axiomatic class. In particular, the previously mentioned families 
are special subfamilies of this class of algebras. Furthermore, the required axioms
are easy to verify for additional families of algebras. 
The proof of the theorem 
is constructive, so one obtains an explicit quantum cluster algebra 
structure in each case. Initial clusters are constructed
intrinsically as finite sequences of (homogeneous) elements in chains 
of noncommutative unique factorization domains. Another key feature of the result is that 
it holds for arbitrary base fields: there are no restrictions 
on their characteristic and they do not need to be algebraically closed. 
The proof of the theorem is based on purely ring theoretic arguments 
which are independent of the characteristic of the field and 
do not use specialization. Finally, when the methods are applied to 
algebras arising from quantum groups, the deformation parameter 
$q$ only needs to be a non-root of unity while the previous 
methods needed $q$ to be transcendental over $\Qset$.

In this paper, we apply the theorem to construct explicit quantum cluster 
algebra structures on the quantum Schubert cell algebras  
$\UU_q(\n_+ \cap w(\n_-))$ for all finite dimensional 
simple Lie algebras $\g$. (The technique works for all 
Kac--Moody Lie algebras $\g$, but the general case requires technicalities 
which would only increase the size of the paper and obstruct the main idea.
Because of this, the minor additional details for infinite dimensional Lie algebras 
$\g$ will appear elsewhere.) If $\g$ is symmetric, this result is due to 
Gei\ss, Leclerc and Schr\"oer \cite{GLSh}. In this case we obtain the same quantum cluster algebra
structure, but under milder assumptions on the base field and the 
deformation parameter.

In a forthcoming publication \cite{GY-BZconj},
we will give a proof of the Berenstein--Zelevinsky 
conjecture \cite{BZ} using the above theorem.
We will show that the quantized coordinate rings of all double 
Bruhat cells are localizations of quantum nilpotent algebras and 
that applying the above theorem produces precisely the conjectured 
explicit (upper) quantum cluster algebra structure of Berenstein and Zelevinsky \cite{BZ}.
In fact, the final result is stronger than the conjecture in that we prove that in each case 
the upper quantum cluster algebra coincides with the quantum cluster algebra.

The main theorem of the paper has an analog for a certain (very large, 
axiomatically defined) 
class of Poisson structures on 
polynomial algebras which can be considered as the semiclassical limits 
of quantum nilpotent algebras. A direct application of that theorem produces 
explicit cluster algebra structures on unipotent groups of the form 
$U_+ \cap w(U_-)$ (for any Kac--Moody group $G$) and on all 
double Bruhat cells. The latter result also proves that the upper cluster 
algebras of Berenstein, Fomin and Zelevinsky \cite{BFZ} 
for double Bruhat cells are equal to the corresponding 
cluster algebras, thus solving one of the abovementioned problems.
This will also appear in a forthcoming publication.  

\section{Definition of quantum nilpotent algebras}
\label{1.2}
Next, we proceed with the definition of quantum nilpotent algebras.

Throughout, $\KK$ \index{K@$\KK$} will denote a base field. Its characteristic can  be arbitrary, 
and it need not be algebraically closed. All automorphisms and skew-derivations of $\KK$-algebras 
will be $\KK$-linear. For two integers $j$ and $k$ we will denote 
$[j,k] := \{j, j+1, \ldots, k \}$ if $j \leq k$ and $[j,k] := \varnothing$
otherwise.  \index{jk@$[j,k]$}

Consider an iterated Ore extension (or iterated skew polynomial algebra)
\begin{equation}
\label{Ore0}
R := \KK[x_1][x_2; \sigma_2, \delta_2] \cdots [x_N; \sigma_N, \delta_N], 
\end{equation}
and denote the intermediate algebras $R_k:=\KK[x_1][x_2; \sigma_2, \delta_2] \cdots [x_k; \sigma_k, \delta_k]$  \index{Rk@$R_k$}  for
$k \in [0,N]$. (Thus, $R_0 = \KK$ and $R_N = R$.) Our conventions for Ore extensions are detailed in Convention \ref{skewpoly}.

\begin{defA*}
 An iterated Ore extension $R$ as in \eqref{Ore0}  
is called a \emph{Cauchon--Goodearl--Letzter} (\emph{CGL}) \emph{extension} \cite{LLR}
 \index{Cauchon--Goodearl--Letzter extension} \index{CGL extension}
if it is equipped with a rational action of a $\KK$-torus $\HH$  \index{H@$\HH$} 
by $\KK$-algebra automorphisms satisfying the following conditions:
\begin{enumerate}
\item[(i)] The elements $x_1, \ldots, x_N$ are $\HH$-eigenvectors.
\item[(ii)] For every $k \in [2,N]$, $\delta_k$ is a locally nilpotent 
$\sigma_k$-derivation of $R_{k-1}$. 
\item[(iii)] For every $k \in [1,N]$, there exists $h_k \in \HH$ such that 
$\sigma_k = (h_k \cdot)$ and the $h_k$-eigenvalue of $x_k$, to be denoted by $\la_k$,  \index{lambdak@$\la_k$}  is not a root of unity.  \index{hk@$h_k$}
\end{enumerate}
\end{defA*}

The quantum Schubert cell algebras
$\UU_q(\n_+ \cap w(\n_-))$ mentioned earlier, for non-roots of unity $q \in \kx$, are examples of CGL extensions.
We refer the reader to Chapter \ref{q-gr} for details. Particular cases include generic 
quantized coordinate rings of affine spaces, matrix varieties, symplectic and euclidean spaces, 
and generic quantized Weyl algebras. 

Every nilpotent Lie algebra of dimension $m$ 
contains a chain of ideals
$$
\n = \n_m  \vartriangleright \n_{m-1} \vartriangleright \ldots \vartriangleright \n_1 \vartriangleright
\n_0 = \{ 0 \}
$$
with the properties that $\dim (\n_k /\n_{k-1}) =1$ and $[\n,\n_k] \subseteq \n_{k-1}$, 
for all $1 \leq k \leq m$. For $k \in [1,m]$, let $x_k \in \n_k \setminus \n_{k-1}$. This set of elements gives rise to the following 
iterated Ore extension presentation of $\UU(\n)$:
$$
\UU(\n) = \KK [x_1][x_2; \id, \delta_2] \ldots [x_m; \id, \delta_m].
$$
The derivations $\delta_k = \ad_{x_k}$ are locally nilpotent.
Thus, all universal enveloping algebras of finite dimensional nilpotent 
Lie algebras can be presented as iterated Ore extensions \eqref{Ore0} 
with all $\sigma_k = \id$ and $\delta_k$ locally nilpotent. 
One can consider $\HH = \{ 1 \}$ acting trivially on 
them and then all conditions in Definition A will be satisfied 
except for the second part of (iii) -- in this case all eigenvalues $\la_k$ will be 
equal to $1$.

{\em{We consider the class of CGL extensions to be the best current definition of 
quantum nilpotent algebras from a ring theoretic perspective}}.
\index{quantum nilpotent algebra}
On the level of presentations, the local nilpotency of $\delta_k$ from 
the classical situation is kept in the definition. The torus $\HH$
is used to get hold of the 
eigenvalues $\la_k$. The condition that they are non-roots of unity 
is the main feature of the quantum case. On the level of deeper ring theoretic 
properties, these algebras exhibit the common property of having
only finitely many $\HH$- invariant prime ideals \cite{Ca,GL}. This property, first derived from the pioneering
works \cite{HL1, HL2, J} on the spectra of quantum groups, has played a key role in the study of spectra of ``quantum algebras".
From another ring theoretic perspective, all CGL extensions are conjectured 
to be catenary just as Gabber's theorem for universal enveloping algebras 
of solvable Lie algebras (though this is currently established only 
for quantized Weyl algebras
\cite{GLen-cate} and quantum Schubert cell algebras \cite{Y-cate}).

The class of CGL extensions appears to be very large. We do not know 
any classification results except for very low dimensions \cite{Ri}.
As mentioned in Section \ref{1.1},
in addition to the quantum Schubert cell algebras, in \cite{GY-BZconj} we prove 
that all quantized coordinate rings of double Bruhat cells are localizations 
of CGL extensions. When the main theorem is applied to them, these are precisely 
the localizations with respect to all frozen variables.

The quantum Schubert cell algebras are special members in the class of CGL extensions.
The former have the property that 
their gradings by the character lattices of the torus $\HH$ can be specialized to 
$\Zset_{\geq 0}$-gradings that are connected. We provide an example of a CGL extension 
that does not have this property in \exref{non-N-grad}. This example can be easily generalized 
to families of higher dimensional examples.

CGL extensions are very rarely Hopf algebras. Since Drinfeld's definition of 
quantized universal enveloping algebras \cite{Dr} imposes the condition 
that they are topological Hopf algebras, we use the term {\em{quantum 
nilpotent algebras}} as opposed to {\em{quantized universal enveloping algebras 
of nilpotent Lie algebras}}. Along these lines one should also note 
that it is currently unknown whether 
every CGL extension is a deformation of the universal enveloping algebra 
of a nilpotent Lie algebra.

All important CGL extensions that we are aware of are {\em{symmetric}} in the sense that they are also CGL extensions when the 
generators $x_1, \ldots, x_N$ are adjoined in the reverse order
$$
R= \KK [x_N] [x_{N-1}; \sigma^\sy_{N-1}, \delta^\sy_{N-1}] \ldots [x_1; \sigma^\sy_1, \delta^\sy_1].
$$
(This symmetry condition may only hold for certain orderings of the initial generators.)
We refer the reader to Section \ref{2.5} for a detailed discussion of this condition.
Here we note that  symmetricity of a CGL extension is equivalent to imposing 
a very mild condition on the action of $\HH$ and the following 
abstract Levendorskii--Soibelman type straightening law:

{\em{For all $j <k$, the element $x_k x_j - \la_{kj} x_j x_k$
belongs to the unital subalgebra of $R$ generated by $x_{j+1}, \ldots, x_{k-1}$ for some 
scalar $\la_{kj} \in \kx$.}}

One easily shows that in the setting of Definition A the scalar $\la_{kj}$ equals 
the $h_k$-eigenvalue of $x_j$, i.e., $h_k \cdot x_j = \la_{kj} x_j$.

{\em{The class of algebras covered by the main theorem is the class of symmetric 
CGL extensions satisfying two very mild conditions on the scalars involved.}}

\section[A precise statement of the main result]{Constructing an initial quantum seed and a precise statement of the main result}
\label{1.3}
The main technique that we use to generate quantum clusters in 
CGL extensions is that of 
{\em{noncommutative unique factorization domains}} (UFDs), as defined by Chatters \cite{Cha}.
We briefly recall some background.
A nonzero element $p$ of an integral domain $R$ is called prime if it is normal (meaning that $Rp=pR$)
and the factor $R/pR$ is a domain. A noetherian 
domain $R$ is called a UFD if every nonzero prime ideal of $R$ contains a 
prime element. It is well known that, in the commutative noetherian situation, this definition 
is equivalent to the more common one. If $R$ is equipped with a rational
action of a $\KK$-torus $\HH$ by algebra automorphisms, 
then $R$ has a canonical grading by the rational character lattice 
$\xh$ of $\HH$. Such an algebra $R$ is called an $\HH$-UFD if every 
nonzero $\HH$-prime ideal of $R$ contains a homogeneous prime 
element. As in the commutative situation, two prime elements $p_1, p_2 \in R$ 
are associates if and only if there exists a unit $a \in R$ such that 
$p_2 = a p_1$.

Chatters proved \cite{Cha} that the universal enveloping algebra of every solvable Lie 
algebra over the field of complex numbers is a UFD. Launois, Lenagan and Rigal proved \cite{LLR} that 
every CGL extension is an $\HH$-UFD. (One can view these results as 
another reason for thinking of CGL extensions as quantum nilpotent algebras.) 
As noted in Section \ref{1.2}, CGL extensions have only finitely many 
$\HH$-prime ideals \cite{Ca,GL}, so they have only finitely many homogeneous 
prime elements up to associates.

Our key method for generating quantum clusters, and the first step towards 
proving the main theorem, is as follows:

\begin{step1*}
 For an $\HH$-UFD $R$ and a chain  of $\HH$-UFDs 
$$
\CC : \quad \{0\} \subset R_1 \subsetneq \ldots \subsetneq R_N = R,
$$
define the subset
$$
P(\CC):= \bigcup_{k=1}^N \{ \mbox{homogeneous prime elements of $R_k$ up to associates} \}
$$
of $R$ and prove that {\rm(}under some general assumptions\/{\rm)} $P(\CC)$ is a quantum cluster in $R$.
\end{step1*}

From general properties of prime elements, one gets that under mild assumptions 
the elements of $P(\CC)$ quasi-commute, i.e., for every $y, y' \in P(\CC)$, 
$y y' = \xi y' y$ for some $\xi \in \KK^*$. So one of the major points in 
realizing Step 1 in a particular situation is to ensure 
that $P(\CC)$ has the right size, i.e., that for each $k$, 
$R_k$ has precisely one homogeneous prime element (up to taking associates) 
that is not a prime element of some $R_i$ for $i<k$. 
Furthermore, in order for this to work out the chain $\CC$ should 
be sufficiently fine, e.g., $\GKdim R_k = k$ for all $k$.
This condition is obviously satisfied for the canonical chain 
of subalgebras associated to the CGL extension \eqref{Ore0}. 

For arbitrary CGL extensions, Step 1 was carried out in \cite{GY} in full generality. 
In the statement of this 
result, for a function $\eta : [1,N] \to \Zset$, we will use the  
canonical predecessor and successor functions 
$p:= p_\eta : [1,N] \to [1,N] \sqcup \{ - \infty \}$ 
and $s:= s_\eta : [1,N] \to [1,N] \sqcup \{ + \infty \}$ 
for the level sets of $\eta$ (see \eqref{p}--\eqref{s}).

\begin{thA*}
{\rm \cite[Theorem 4.3]{GY}}
Let $R$ be an arbitrary CGL extension of length $N$. 
Then there exist a function $\eta : [1,N] \to \Zset$ and elements
$$
c_k \in R_{k-1} \; \; \mbox{for all} \; \; k \in [2,N] \; \; 
\mbox{with} \; \; p(k) \neq - \infty
$$
such that the elements $y_1, \ldots, y_N \in R$, recursively defined by 
\begin{equation}
\label{y0}
y_k := 
\begin{cases}
y_{p(k)} x_k - c_k, &\mbox{if} \; \;  p(k) \neq - \infty \\
x_k, & \mbox{if} \; \; p(k) = - \infty,  
\end{cases}
\end{equation}
\index{yk@$y_k$}
are homogeneous and have the property that for every $k \in [1,N]$,
$$
\{y_j \mid j \in [1,k] , \, s(j) > k \}
$$
is a list of the homogeneous prime elements of $R_k$ up to scalar multiples.
The elements $y_1, \ldots, y_N \in R$ with these properties are unique.
The function $\eta$ with these properties is not unique but its level sets are. 
Furthermore, the function $\eta$ has the property that 
$p(k) = - \infty$ if and only if $\delta_k=0$.
\end{thA*}

In different words the theorem states that, if $\CC_k$ denotes 
the truncated chain $\{0\} \subset R_1 \subset \ldots \subset R_k$, 
then in the setting of the theorem
$$
P(\CC_k) = P(\CC_{k-1}) \sqcup \{ y_k \}
$$
and the prime element $y_k$ of $R_k$ is determined by the linear 
expression \eqref{y0} with leading term equal to 1 or to a prime element of the previous algebra
$R_{k-1}$. The function $\eta$ (or more precisely its level sets) 
keeps track of the leading terms of this recursive formula.

To construct exchange matrices that go with the quantum clusters from 
Step 1 (Theorem A), we need to introduce some more notation.
The quantum clusters of an algebra generate quantum tori. Recall 
that a matrix $\qb= (q_{kj}) \in M_N(\kx)$ is called {\em{multiplicatively 
skew-symmetric}}  \index{multiplicatively skew-symmetric} if $q_{kk} =1$ for all $k$ and $q_{kj} q_{jk} =1$ for all 
$k \neq j$. Such a matrix gives rise to the {\em{quantum torus}}  \index{quantum torus}
\begin{equation}
\label{qTor}
\Tbb_\qb := \frac{ \KK \lcor Y_1^{\pm 1}, \ldots, Y_N^{\pm 1} \rcor}
{ ( Y_k Y_j - q_{kj} Y_j Y_k, \; \forall k \neq j ) }\cdot
\end{equation}
\index{Tq@$\Tbb_\qb$}

Given a CGL extension $R$ of length $N$, we consider the unique 
multiplicatively skew-symmetric matrix $\lab \in M_N(\kx)$  \index{lambda@$\lab$} 
whose entries $\la_{kj}$ for $1 \leq j < k \leq N$ 
are the $h_k$-eigenvalues of $x_j$
in the notation of Definition A.
Define the (multiplicatively skew-symmetric) bicharacter
$$
\Om_\lab : \Zset^N \times \Zset^N \to \kx
\quad \mbox{by} 
\quad 
\Om_\lab(e_k, e_j) := \la_{kj}, \; \; 
\forall k,j \in [1,N],
$$
\index{Omegalambda@$\Om_\lab$}
where $\{ e_1, \ldots, e_N\}$ is the standard basis of $\Zset^N$.
In the setting of Theorem A, define
$$
\ol{e}_k : = \sum_{n \in \Zset_{\geq 0}, \, \, p^n(k) \neq - \infty} e_{p^n(k)} 
\in \Zset^N, \; \; 
\forall k \in [1,N].
$$
\index{ekbar@$\ol{e}_k$}
The set $\{ \ol{e}_1, \ldots, \ol{e}_N \}$ is another basis 
of $\Zset^N$. We consider the (multiplicatively skew-symmetric) bicharacter
$$
\Om : \Zset^N \times \Zset^N \to \kx
\quad \mbox{such that} 
\quad 
\Om(e_k, e_j):= \Om_\lab(\ol{e}_k, \ol{e}_j), \; \; 
\forall k,j \in [1,N].
$$
\index{Omega@$\Om$}
The sequence of prime elements in Theorem A produces a quantum cluster of the CGL extension 
$R$ in the following sense:

{\em{For any CGL extension $R$, the elements $y_1, \ldots, y_N \in R$ 
from Theorem A and their inverses 
generate a copy of the quantum torus $\Tbb_\qb$ inside 
$\Fract(R)$, where
$$
q_{kj} = \Om(e_k, e_j) = \Om_\lab(\ol{e}_k, \ol{e}_j), \; \; 
\forall k, j \in [1,N] 
$$
\index{qkj@$q_{kj}$}
and 
$$
R \subseteq \Tbb_\qb \subset \Fract(R).
$$
}}
Consider the subset
$$
\ex:=\{ k \in [1,N] \mid s(k) \neq + \infty \}
$$
\index{ex@$\ex$}
of $[1,N]$. The columns of all $N \times |\ex|$ matrices will be indexed by $\ex$.
The set of such matrices with integral entries will be denoted by $M_{N \times \ex}(\Zset)$.  \index{MNXZ@$M_{N \times X} (\Zset)$}

\begin{mainth2spec*}
For every symmetric 
CGL extension $R$, the following hold under very mild assumptions
on the base field $\KK$ and the scalars 
$\la_1, \ldots, \la_N \in \kx$:

{\rm (a)} For each $k \in \ex$, there exists a unique vector $b^k = \sum_{l=1}^N b_{lk} e_l \in \Zset^N$ 
such that 
$$
\Om(b^k, e_j)^2 = \la_{s(k)}^{ - \delta_{kj}}, \; \; 
\forall j \in [1,N]
$$
and $y_1^{b_{1k}} \ldots y_N^{b_{Nk}}$ is fixed under the $\HH$-action.
Denote by $\wt{B} \in M_{N \times \ex}(\Zset)$ the matrix with columns
$b^k \in \Zset^N$, $k \in \ex$.

{\rm (b)} The matrix $\wt{B}$ is compatible with the 
quantum cluster $(y_1, \ldots, y_N)$ of $R$ in the sense that they define
a quantum seed, call it $\Sig$, see \cite{BZ} and \S{\rm\ref{3.2}} for definitions.
Let $\Abb(\Sig)$ and $\, \UU(\Sig)$ be
the quantum cluster algebra and upper quantum cluster algebra 
associated to this quantum seed for the set of exchangeable indices
$\ex$ where none of the frozen cluster variables are inverted.

{\rm (c)} After an explicit rescaling of 
$y_1, \ldots, y_N$, we have the equality 
$$
R = \Abb(\Sig) = \UU(\Sig).
$$ 
Furthermore, the quantum cluster algebra $\Abb(\Sig)$ is generated by finitely many cluster 
variables.
\end{mainth2spec*}

We refer the reader to \thref{cluster} for a full statement of the main theorem.
The homogeneity condition in part (a) of the theorem can be written in an explicit form 
taking into account that for each $k \in [1,N]$, 
$y_k$ is an $\HH$-eigenvector with eigenvalue equal to the product 
of the $\HH$-eigenvalues of $x_{p^n(k)}$ for those $n \in \Zset_{\geq 0}$ 
such that $p^n(k) \neq - \infty$ (see \thref{cluster} (a)).

The quantum cluster algebras in question are built from general (multiparameter) quantum tori 
by an extension of the Berenstein--Zelevinsky construction \cite{BZ}. 
The general setting for such quantum cluster algebras is developed 
in Chapter \ref{qClust}.

We finish this section with the statements of the two mild conditions 
which are imposed in the Main Theorem. In \prref{la-equal}, we prove that 
$$ 
\la_k=\la_l \; \; \mbox{for all $k, l \in [1,N]$ such that 
$\eta(k) = \eta(l)$ and $p(k) \neq - \infty$, $p(l) \neq - \infty$.} 
$$
Recall also that $\la_{kl}$ denotes the $h_k$-eigenvalue of $x_l$ for $1 \leq l < k \leq N$.
In this setting, the two conditions imposed in the main theorem are:

{\bf{Condition (A).}} The base field $\KK$ contains square roots $\sqrt{\la_{kl}}$ of the scalars
$\la_{kl}$ for $1 \leq l < k \leq N$ such that the subgroup of $\kx$ generated 
by all of them contains no elements of order $2$.

{\bf{Condition (B).}} There exist positive integers $d_n$, $n \in \range(\eta)$,
such that 
$$
\la_k^{d_{\eta(l)}} = \la_l^{d_{\eta(k)}} 
$$
for all $k, l \in [1,N]$, $p(k) \neq - \infty$, $p(l) \neq - \infty$.

All symmetric CGL extensions that we are aware of (with the exception of certain very specific 
two-cocycle twists) satisfy these conditions after a field extension of the base field $\KK$. 
We give a detailed account of this in Remarks \ref{rtorsion-free} and \ref{rd-pr}. The first 
condition is needed when dealing with toric frames as opposed to individual cluster variables.
The second condition is needed to ensure that the principal parts of certain exchange matrices 
that are constructed are skewsymmetrizable. Both issues are minor and the main constructions of 
quantum clusters and mutations work without imposing those conditions.

\section{Additional clusters}
\label{1.4}
From now on, $R$ will denote a symmetric CGL extension as in \eqref{Ore0}.

Consider the subset $\Xi_N$  \index{XiN@$\Xi_N$}  of the symmetric group $S_N$ consisting of all 
permutations $\tau \in S_N$ such that 
$$
\tau(k) = \max \, \tau( [1,k-1]) +1 \; \;
\mbox{or} 
\; \; 
\tau(k) = \min \, \tau( [1,k-1]) - 1, 
\; \; \forall k \in [2,N].
$$
It is easy to see that $\Xi_N$ can also be defined as the set of all $\tau \in S_N$ such that 
$\tau([1,k])$ is an interval for all $k \in [2,N]$. Each symmetric CGL extension has a different 
CGL extension presentation associated to every element $\tau \in \Xi_N$, of the form
\begin{equation}
\label{taupres}
R= \KK[x_{\tau(1)}][x_{\tau(2)}; \sigma''_{\tau(2)}, \delta''_{\tau(2)}] \ldots 
[x_{\tau(N)}; \sigma''_{\tau(N)}, \delta''_{\tau(N)}].
\end{equation}
The proof of this fact and a description of the automorphisms $\sigma''_{\tau(k)}$ and the 
skew-derivations $\delta''_{\tau(k)}$ are given in \prref{tauOre}. 
This presentation gives rise to a chain of subalgebras of $R$,
$$
\CC_\tau : \quad \{0\} \subset R_{\tau, 1} \subset \ldots \subset R_{\tau, N} = R,
$$
\index{Ctau@$\CC_\tau$}
where $R_{\tau,k}$  \index{Rtauk@$R_{\tau,k}$}  is the unital subalgebra generated by $x_{\tau(1)}, \ldots, x_{\tau(k)}$.
Theorem A associates a 
quantum cluster
$$
(y_{\tau, 1}, \ldots, y_{\tau, N}):= P(\CC_\tau)
$$
to $R$ for each $\tau \in \Xi_N$.
Clearly, $\id \in \Xi_N$ and  
$$
(y_{\id, 1}, \ldots, y_{\id, N}) = P(\CC_\id) = P(\CC)=(y_1, \ldots, y_N).
$$

\begin{mainth3add*}
Let $R$ be a symmetric
CGL extension satisfying the conditions {\rm (A)} and {\rm (B)}. For each $\tau \in \Xi_N$, 
the quantum cluster algebra $\Abb(\Sig)$ from Main Theorem {\rm II} has a quantum seed 
with a set of cluster variables obtained from $(y_{\tau, 1}, \ldots, y_{\tau, N})$
by a precise rescaling and permutation {\rm(}see Theorem {\rm\ref{tcluster} (a)--(b))}. The exchange 
matrix $\wt{B}_\tau$  \index{Btautilde@$\wt{B}_\tau$}  for this seed is the unique solution of the equations in part {\rm (a)} 
of Main Theorem {\rm II} for the bicharacter associated to the 
CGL extension presentation \eqref{taupres}. 

In particular, an appropriate scalar multiple of each generator $x_k$ of $R$ is a cluster variable.
\end{mainth3add*}

The last part of the theorem follows from the first because for each $k \in [1,N]$ there exists an
element $\tau \in \Xi_N$ such that $\tau(1) =k$.

\section{Strategy of the proof}
\label{1.5}
The second step of the proof of the main theorem is to connect 
the quantum clusters indexed by the elements of $\Xi_N$ by a chain of mutations:

\begin{step2*}
Let $\tau \neq \tau'$ in $\Xi_N$ be such that for some $k \in [1,N-1]$,
$\tau(j) = \tau'(j)$ for all $j \neq k, k+1$ {\rm(}i.e., the subalgebras in the 
chains $\CC_\tau$ and $\CC_{\tau'}$ only differ in position $k$\/{\rm)}. We prove that 
in this setting the following hold for every symmetric CGL extension $R$:

{\rm (a)} If $\tau(k)$ and $\tau(k+1)$ have different images under the $\eta$-function for the 
presentation \eqref{taupres}, then the quantum clusters $P(\CC_\tau)$ and $P(\CC_{\tau'})$
are equal. More precisely,
$$
y_{\tau',j} = y_{\tau, j} \; \; \mbox{for $j \neq k,k+1$ and} \; \; 
y_{\tau',k} = y_{\tau, k+1}, \; \; y_{\tau',k+1} = y_{\tau, k}. 
$$ 

{\rm (b)} If $\tau(k)$ and $\tau(k+1)$ have the same images under the $\eta$-function for the 
presentation \eqref{taupres}, then the quantum cluster $P(\CC_{\tau'})$ of $R$ 
is a one-step mutation of the quantum cluster $P(\CC_\tau)$. More precisely,
\begin{equation}
\label{one-step}
y_{\tau', k} = y_{\tau, k}^{-1} 
\left(\xi_1 \prod_{j \in [1,N], \, c_j>0} y_{\tau,j}^{c_j} + \xi_2 \prod_{j \in [1,N], \, c_j < 0} y_{\tau,j}^{-c_j} \right) 
\end{equation}
for some $\xi_1, \xi_2 \in \kx$ and an integral vector $(c_1, \ldots, c_N)$ such that $c_k=0$, and
$$
y_{\tau',j} \; \; \mbox{is a scalar multiple of} \; \; y_{\tau,j} \; \; 
\mbox{for all} \; \;  j \neq k.
$$
\end{step2*}

This result has a simple proof given in Theorems \ref{t1}, \ref{t2}. It is derived from the uniqueness part of Theorem A 
and the uniqueness of the decomposition of a normal element of a noncommutative UFD 
as a product of prime elements. {\em{Thus, in a nutshell, cluster mutation is forced by the uniqueness 
in Theorem A and the unique decomposition in noncommutative UFDs.}}

The simultaneous normalization (rescaling) of all
cluster variables $y_{\tau, k}$, $\tau \in \Xi_N$, $k \in [1,N]$
for which all scalars $\xi_1$ and $\xi_2$ in \eqref{one-step} become equal to 1, 
turns out to be a rather delicate issue. It is 
carried out in Chapters \ref{prime-sym} and \ref{mut-sym}
in an explicit form.

The next step is to extend the clusters $P(\CC_\tau)$ to quantum seeds, which amounts 
to the construction of exchange matrices that are compatible with these quantum clusters.
We also need to connect those seeds with cluster mutations.
Two elements $\tau, \tau' \in \Xi_N$ will be called adjacent if they satisfy the condition in Step 2, 
in other words, if there exists $k \in [1,N-1]$ such that
\begin{equation}
\label{kk+1}
\tau' = \tau (k,k+1).
\end{equation}

\begin{step3*}
 Construct exchange matrices $\wt{B}_\tau$ 
for quantum seeds going with the quantum clusters $P(\CC_{\tau})$, $\tau \in \Xi_N$.

Show that the seeds corresponding to adjacent pairs $\tau, \tau' \in \Xi_N$ 
are obtained by one-step mutations from each other.
\end{step3*}

Let $\tau \in \Xi_N$. Step 2 constructs one column of the desired matrix $\wt{B}_\tau$ for each element $\tau' \in \Xi_N$ 
which is adjacent to $\tau$. More precisely, this is the $k$-th column of $\wt{B}_\tau$ where 
$k \in [1,N-1]$ is the index from \eqref{kk+1}. Step 3 involves constructing the rest of the matrix $\wt{B}_\tau$
and then showing that the exchange matrices for adjacent elements of $\Xi_N$ are related by mutation.
One should note that there are not enough adjacent elements to construct 
the full matrix $\wt{B}_\tau$ directly from Step 2.

The realization of Step 3 and the proof of the remaining  parts of the main theorem  
are based on using a chain of successively adjacent elements of $\Xi_N$ of the form
\begin{equation}
\label{chainS}
\tau_0 = \id \to \tau_1 \to \ldots \to \tau_{M-1} \to \tau_M=w_\ci ,
\end{equation}
where $w_\ci$ is the longest element of $S_N$ and $M=N(N-1)/2$ is the length of $w_\ci$.
The arrows are drawn for the convenience of the description of the procedure. 

We first construct $\wt{B}=\wt{B}_\id$. 
Each chain \eqref{chainS} has the property that
\begin{equation}
\label{every-mut}
\forall k \in [1,N-1], \; \; 
\exists i \in [1,M-1] \; \; 
\mbox{\em{such that}} \; \; \tau_{i+1} = \tau_i (k, k+1). 
\end{equation}
{\em{This implies that every exchangeable index in}} $\ex = \{ k \in [1,N] \mid s(k) \neq + \infty \}$
{\em{gets mutated when Step {\rm2} is applied to the arrows in the chain}}.

Let $k \in \ex$. Leaving aside certain technical details concerning a reenumeration of cluster variables,
to define the $k$-th column of $\wt{B}$ we choose an 
index $i \in [1,M-1]$ satisfying \eqref{every-mut}. Step 2 defines the 
$k$-th columns of $\wt{B}_{\tau_i}$ and $\wt{B}_{\tau_{i+1}}$. 
We construct the $k$-th column of $\wt{B}$ by inverse mutation 
from the $k$-th column of $\wt{B}_{\tau_i}$ along the chain \eqref{chainS}. 
At this point, one has to prove that this is independent of the choice of $i$. 
There two problems here. The first is that the integers $c_j$ in \eqref{one-step} 
are not determined from easily accessible ring theoretic data. They are powers 
in the decomposition of some normal elements as products of primes, but 
the setting is very general to get hold of them. Secondly, the inverse mutation from such data 
gets combinatorially overwhelming and one cannot keep track of it in general ring theoretic terms.

We present a very general (ring theoretic) solution to this problem. Firstly, we establish that there is 
at most one matrix $\wt{B}$ that satisfies the conditions in Theorem II (a), using 
a strong rationality result for the $\HH$-primes 
of CGL extensions \cite[Theorem II.6.4]{BG}, a special case of which is that 
$$
Z(\Fract(R))^\HH = \KK ,
$$
where $Z(.)$  \index{Z@$Z(.)$}  stands for the center of an algebra.
We obtain the uniqueness statement in Theorem II (a) by phrasing the conditions for the columns 
of $\wt{B}$ in terms of commutation relations 
and using the above result. 
We then use the idea for inverse mutation along the chain \eqref{chainS} to prove that a 
matrix $\wt{B}$ with the properties in Theorem II (a) exists. 

All other exchange matrices $\wt{B}_\tau$ for $\tau \in \Xi_N$ are constructed by mutating 
$\wt{B}$ along a chain which is entirely within $\Xi_N$, starts with $\id$, ends at $\tau$, and
for which each two consecutive elements are adjacent. It is not difficult to see that such 
a chain exists for every $\tau \in \Xi_N$. 
Typically, there are many chains with these properties and the independence of $\wt{B}_\tau$ from the choice of a chain is 
proved using one more time the strong rationality result. The last part of Step 3, 
that the exchange matrices for adjacent elements of $\Xi_N$ are related by mutation, follows from 
their construction.

At this point, we have a construction of quantum seeds $\Sig_\tau$ indexed by the elements 
of $\Xi_N$ and mutations between them. What is left is to establish 
part (c) of Main Theorem II. This is the last step of the proof of the theorem:

\begin{step4*}
Show that 
$$
R = \Abb(\Sig) = \UU(\Sig)
$$
for $\Sig := \Sig_\id$, i.e., that $R$ coincides with the 
quantum cluster algebra and the upper quantum cluster algebra 
for the quantum seed $\Sig$.
\end{step4*}

In the remaining part of this section we sketch how this step is carried out.
As previously mentioned, for every $k \in [1,N]$ there is an element $\tau \in \Xi_N$ 
with $\tau(1) =k$. Since all quantum seeds $\Sig_\tau$ are mutation equivalent to each other,  
a scalar multiple of each generator $x_k$ is a cluster variable of $\Abb(\Sig)$. Therefore, 
$$
R \subseteq \Abb(\Sig).
$$
By the quantum Laurent phenomenon,
$$
\Abb(\Sig) \subseteq \UU(\Sig).
$$
For $\tau \in \Xi_N$, denote by $E_\tau$  \index{Etau@$E_\tau$}  the multiplicative subset of $R$ generated by all
exchangeable cluster variables in $\Sig_\tau$. We first prove that $E_\tau$ is an Ore 
set in $R$. Consider a chain of elements of $\Xi_N$ as in \eqref{chainS}. 
We have 
$$
\UU(\Sig) \subseteq \bigcap_{i=0}^M R[E_{\tau_i}^{-1}].
$$
Combining all of the above embeddings, we obtain 
\begin{equation}
\label{chain-incl}
R \subseteq \Abb(\Sig) \subseteq \UU(\Sig) \subseteq \bigcap_{i=0}^M R[E_{\tau_i}^{-1}].
\end{equation}
At this point, we use techniques from noncommutative UFDs and iterated Ore extensions to describe the intersection 
of $R[E_\tau^{-1}]$ and $R[E_{\tau'}^{-1}]$ for two adjacent elements of $\tau$ and $\tau'$ of $\Xi_N$. Based on it 
we prove that 
\begin{equation}
\label{RRloc}
\bigcap_{i=0}^M R[E_{\tau_i}^{-1}] = R.
\end{equation}
This forces all inclusions in \eqref{chain-incl} to be equalities.

\section{Organization of the paper and reading suggestions} 
\label{1.6}
The paper is organized as follows: 

I. Chapter \ref{qClust} provides the general framework for quantum cluster algebras 
generalizing the Berenstein--Zelevinsky construction \cite{BZ} to mutations between 
multiparameter quantum tori (i.e., quantum tori \eqref{qTor} for which the scalars 
$q_{kl}$ are not necessarily integral powers of a single element $q \in \kx$). 
This chapter is self contained and can be read independently 
of \cite{BZ}. Readers who are familiar with \cite{BZ} might look only at 
the main definitions since the treatment is similar to \cite{BZ}. 
Chapter \ref{review} reviews the theory of noncommutative UFDs and the 
material from \cite{GY} that goes with Theorem A.

II. Chapter \ref{mCGL} carries out Step 2. For readers who are interested in 
understanding the intrinsic reason for the appearance of cluster mutations 
in the general framework of quantum nilpotent algebras, we suggest reading 
the proofs of Theorems \ref{t1}, \ref{t2}. 

Chapters \ref{prime-sym} and \ref{mut-sym} carry out 
in an explicit form the simultaneous normalization of all cluster variables $y_{\tau,k}$ 
so that all scalars $\xi_1$, $\xi_2$ in part (b) of Step 2 are equal to 1. 
These chapters are technical and we suggest that readers restrict to the definitions 
of $y_{[.,.]}$, $u_{[.,.]}$ and $\pi_{[.,.]}$ in \thref{y-int} (c), \coref{u-elem} and Eq. \eqref{pi-f},
as well as the statement of \prref{resc} and the condition \eqref{pi-cond}. 
This information will be sufficient for understanding the statement of the main result.

Chapter \ref{Integr} proves Eq. \eqref{RRloc} which essentially carries out Step 4. 
Readers who are only interested in understanding the main result could skip it.

III. The main result appears in \thref{cluster}. Chapter \ref{main} 
also carries out Step 3 and the general part of Step 4.

IV. Chapter \ref{q-Schu} illustrates how \thref{cluster} is applied. It proves that 
the quantum Schubert cell algebras $\UU_q(\n_{\pm} \cap w(\n_{\mp}))$, for all simple Lie algebras $\g$ 
and Weyl group elements $w$, are quantum cluster algebras which 
in addition coincide with the corresponding upper quantum cluster algebras (\thref{Uw-main}). 
Chapter \ref{q-gr} contains the needed details on quantum groups and their relations to 
quantum nilpotent algebras.

We suggest the following route for readers who are interested in understanding 
the final result (\thref{cluster}) without the details of its proof:

1. Basics of CGL extensions Sections \ref{2.2}--\ref{2.5} 
and the base change in Section \ref{4.2} to include $\sqrt{\la_{lj}}$ in $\KK$.

2. The setting of Section \ref{6.1} and the statement of \thref{cluster}. We impose conditions (A) and (B), but the unique rescaling is a feature of the algebra. One can read the statement of the theorem knowing only that such a rescaling exists. If one needs to know its precise form, one can return  to the definitions of $y_{[.,.]}$, $u_{[.,.]}$ and $\pi_{[.,.]}$ and the condition \eqref{pi-cond} 
as indicated in II above.

Readers could also consult \cite{GYann} for an abridged version of \thref{cluster}, 
and its setting and applications.

\section{Notation}
\label{1.7}
We finish the introduction with some notation to be used throughout
the paper. For $m, n \in \Zset_{>0}$ and a commutative ring $D$, 
we will denote by $M_{m,n} (D)$  \index{MmnD@$M_{m,n} (D)$}  the set of matrices of size $m \times n$ 
with entries from $D$. We will use a multiplicative (exponential) version of the standard 
matrix multiplication: For $A = (a_{ij}) \in M_{m,n}(\Zset)$, 
$B = (b_{jk}) \in M_{n,p}(D)$, and $C = (c_{kl}) \in M_{p,s}(\Zset)$, 
denote the matrix 
$$
{}^A \! B^C := \Big( \prod_{j,k} b_{jk}^{a_{ij} c_{kl}} \Big) \in M_{m,s}(D).
$$
\index{ABC@${}^AB^C$}
The two obvious special cases of this operation are given by
$$
{}^A \! B = \Big( \prod_{j} b_{jk}^{a_{ij}} \Big) \in M_{m,p}(D), \quad
B^C = \Big( \prod_{k} b_{jk}^{c_{kl}} \Big) \in M_{n,s}(D).
$$ 
Clearly, 
\begin{equation}
\label{mult}
{}^{A_1} \big( {}^A \! B^C \big) = {}^{A_1 A} \! B^C \quad
\mbox{and} \quad \big( {}^A \! B^C \big)^{C_1} = {}^{A} \! B^{C C_1}
\end{equation} 
for all matrices $A_1$ and $C_1$ with integer entries, having the 
appropriate sizes for which the operations are defined. The transpose
map on matrices will be denoted by $A \mt A^T$. 

Elements of
$\Zset^N$ will be thought of as column vectors, and the standard basis 
of $\Zset^N$ will be denoted by $\{ e_k \mid k \in [1,N] \}$.  \index{ek@$e_k$}
For $g = \sum_k m_k e_k \in \Zset^N$, denote 
\begin{equation}
\label{g+-}
[g]_+ := \sum_k \max (m_k, 0) e_k, \quad
[g]_- := \sum_k \min (m_k, 0) e_k. 
\end{equation}
\index{zz@$[\;]_{\pm}$}
\noindent
{\bf Acknowledgements.} 
We would like to thank A. Berenstein, Ph.~Di Fran\-ces\-co, S. Fomin, R. Kedem, B. Keller, 
A. Knutson, 
B. Leclerc, N. Reshetikhin, D. Rupel and 
A. Zelevinsky for valuable discussions and comments. We are also thankful to 
MSRI for its hospitality during the programs in ``Cluster Algebras'' (2012)
and ``Noncommutative Algebraic Geometry and Representation 
Theory'' (2013) when parts of this project were completed.

\chapter{Quantum cluster algebras}
\label{qClust}

In this chapter we give a definition of quantum cluster algebras 
and upper quantum cluster algebras (over arbitrary 
commutative rings)
that extends the one of Berenstein and Zelevinsky \cite{BZ}. The construction 
in \cite{BZ} used uniparameter quantum tori (quantum tori for which the commutation 
scalars are integral powers of a single element $q \in \kx$). 
We extend the construction of \cite{BZ} to mutation of cluster 
variables coming from general quantum tori. Furthermore, the 
construction allows for arbitrary frozen variables to be 
inverted or not. Following the argument of Fomin--Zelevinsky
a quantum Laurent phenomenon is established in this generality.

\section{General quantum tori}
\label{new3.1}
Recall that a matrix $\qb = (q_{kj}) \in M_N(\kx)$ is called \emph{mul\-ti\-plicatively 
skew-symmetric} \index{multiplicatively skew-symmetric} if $q_{kk} = 1$ and $q_{jk} = q_{kj}^{-1}$ for all $j \neq k$.
Such a matrix gives rise to the skew-symmetric bicharacter
\begin{equation}
\label{La}
\Om_\qb : \Zset^N \times \Zset^N \to \kx \; \; 
\mbox{given by} \; \; \Om_\qb(e_k,e_j) := q_{kj}, \; \forall k,j \in [1,N],
\end{equation}
\index{Omegaq@$\Om_\qb$} i.e.,
\begin{equation}
\label{La2}
\Om_\qb(f,g) = {}^{f^T} \! \qb^g, \quad \forall f, g \in \Zset^N.
\end{equation}
(Recall that $e_1, \ldots, e_N$ denote the standard basis vectors of $\Zset$ 
and that all elements of $\Zset^N$ are thought of as column vectors.)

A multiplicatively skew-symmetric matrix $\qb \in M_N(\kx)$
gives rise to the quantum torus \index{quantum torus}
\begin{equation}
\label{q-torus}
\Tbb_\qb := \frac{ \KK \lcor Y_1^{\pm 1}, \ldots, Y_N^{\pm 1} \rcor}
{ ( Y_k Y_j - q_{kj} Y_j Y_k, \; \forall k \neq j ) } \cdot
\end{equation}
\index{Tq@$\Tbb_\qb$} The subalgebra of $\Tbb_\qb$
\begin{equation}
\label{q-aff}
\Abb_\qb := \KK \lcor Y_1, \ldots, Y_N \rcor \subset \Tbb_\qb 
\end{equation}
\index{Aq@$\Abb_\qb$} is called a \emph{quantum affine space algebra}. \index{quantum affine space algebra} The quantum torus $\Tbb_\qb$ has 
the $\KK$-basis
\begin{equation}
\label{ord-basis}
\{ Y^f:=Y_1^{m_1} \cdots Y_N^{m_N} \mid f=(m_1, \ldots, m_N)^T \in \Zset^N \}.
\end{equation}
\index{Yf@$Y^f$}

Given any subring $D$ of $\KK$ that contains the subgroup of $\kx$ generated by all entries 
of $\qb$, one defines versions of $\Tbb_\qb$ and $\Abb_\qb$ over $D$. One can either 
use \eqref{q-torus} and \eqref{q-aff} with free algebras over $D$ instead of $\KK$, or 
consider the $D$-subalgebras of $\Tbb_\qb$ and $\Abb_\qb$ generated by $Y_k^{\pm 1}$, 
respectively $Y_k$ $(k \in [1,N])$.

\section{Based quantum tori}
\label{3.1}
For a multiplicatively skew-symmetric matrix $\rbf := (r_{kj}) \in M_N(\kx)$ 
denote
$$
\rbf^{\cdot 2}:=(r_{kj}^2) \in M_N(\kx).
$$
\index{rdot2@$\rbf^{\cdot 2}$}
The quantum torus $\Tbb_{\rbf^{\cdot 2}}$ has a $\KK$-basis which is obtained by a 
type of symmetrization of the monomials in \eqref{ord-basis}.
This basis consists of the elements
\begin{equation}
\label{basis}
Y^{(f)} := \Scr_\rbf(f) Y^f = 
\Scr_\rbf(f) Y_1^{m_1} \ldots Y_N^{m_N}
\; \; 
\mbox{for} \; \; f= (m_1, \ldots, m_N)^T \in \Zset^N,
\end{equation}   
\index{Y(f)@$Y^{(f)}$}
where the symmetrization scalar $\Scr_\rbf(f) \in \kx$ is given by
\index{Sr@$\Scr_\rbf$}
\begin{equation}
\label{Scr}
\Scr_\rbf(f):= \prod_{j <k} r_{jk}^{- m_j m_k}.
\end{equation}
We have $Y^{(e_k)}= Y_k$, $\forall k \in [1,N]$. 
Furthermore, one easily verifies that 
\begin{equation}
\label{mult-id}
Y^{(f)} Y^{(g)} = \Om_\rbf(f,g) Y^{(f+g)}, \quad
\forall f,g \in \Zset^N,
\end{equation}
recall \eqref{La}.
We will call the torus $\Tbb_{\rbf^{\cdot 2}}$ with the $\KK$-basis 
$\{Y^{(f)} \mid f \in \Zset^N \}$ the {\em{based quantum torus}} \index{based quantum torus}
associated to the matrix $\rbf \in M_N(\kx)$ (or to the corresponding 
bicharacter $\Om_\rbf : \Zset^N \times \Zset^N \to \kx$).

Each $\sigma \in GL_N(\Zset)$ gives rise to another generating 
set of $\Tbb_{\rbf^{\cdot 2}}$, consisting of the elements
$$
Y_{\sigma, k} : = Y^{(\sigma(e_k))}, \; \; k \in [1,N].
$$ 
\index{Ysigmak@$Y_{\sigma, k}$}
Define the multiplicatively skew-symmetric matrix
$$
\rbf_\sigma : = {}^{\sigma^T} \! \rbf^{\sigma} \in M_N(\kx).
$$
\index{rsigma@$\rbf_\sigma$}
It follows from \eqref{mult-id} that there is a $\KK$-algebra isomorphism
\begin{equation}
\label{psisig}
\psi_\sigma : \Tbb_{\rbf_\sigma^{\cdot 2}} \to \Tbb_{\rbf^{\cdot 2}}, 
\; \; 
\mbox{given by} \; \; 
\psi_\sigma(Y_k) = Y_{\sigma, k}, \; k \in [1,N].
\end{equation}
\index{psisigma@$\psi_\sigma$}
The elements 
$$
Y_\sigma^{(f)} := \psi_\sigma ( Y^{(f)}) = 
\Big(\prod_{j <k} (\rbf_\sigma)_{jk}^{- m_j m_k} \Big) Y_{\sigma,1}^{m_1} \ldots Y_{\sigma,N}^{m_N}
\; \; 
\mbox{for} \; \; f= (m_1, \ldots, m_N)^T \in \Zset^N
$$
\index{Yfsigma@$Y_\sigma^{(f)}$}
also form a basis of $\Tbb_{\rbf^{\cdot 2}}$. The next proposition shows that this basis coincides with the 
original basis \eqref{basis} of $\Tbb_{\rbf^{\cdot 2}}$. In other words, the based quantum torus 
associated with $\rbf \in M_N(\kx)$ is 
invariant under the canonical action of $GL_N(\Zset)$.

\bpr{glN} For all $\sigma \in GL_N(\Zset)$ and all multiplicatively skew-symmetric matrices $\rbf \in M_N(\kx)$, we have 
\begin{equation}
\label{sig-inv}
Y^{(f)} = Y_\sigma^{(\sigma^{-1}(f))}, \; \; 
\forall f \in \Zset^N.
\end{equation}
\epr

\begin{proof} For all $f, g \in \Zset^N$, 
$$
\Om_\rbf( \sigma(f), \sigma(g)) = {}^{(\sigma f)^T} \! \rbf^{\sigma g} = 
{}^{f^T} \! \! ({}^{\sigma^T} \! \rbf^\sigma)^g = \Om_{\rbf_\sigma}(f,g).   
$$
Because of \eqref{mult-id}, to show 
$$
Y^{(\sigma(f))} = Y_\sigma^{(f)}, \; \; 
\forall f \in \Zset^N,
$$
it is sufficient to verify the identity  
for $f = e_k$, $k \in [1,N]$. This is straightforward:
\begin{equation*}
Y^{(\sigma(e_k))} = Y_{\sigma,k} = Y_\sigma^{(e_k)}.  \qedhere
\end{equation*} 
\end{proof}

\bde{toric-fra} Let $\FF$ be a division algebra over $\KK$. 
We say that a map 
$$
M : \Zset^N \to \FF
$$
is a {\em{toric frame}} \index{toric frame} if there exists a multiplicatively skew-symmetric 
matrix $\rbf \in M_N(\kx)$ such that the following conditions 
are satisfied:
\begin{enumerate}
\item[(i)] There is an algebra embedding $\vp : \Tbb_{\rbf^{\cdot 2}} \hra \FF$ 
given by $\vp(Y_i) = M(e_i)$, $\forall i \in [1,N]$ 
and $\FF = \Fract(\vp(\Tbb_{\rbf^{\cdot 2}}))$. 
\item[(ii)] For all $f \in \Zset^N$, $M(f) = \vp( Y^{(f)})$.
\end{enumerate}
\ede

In other words, a toric frame is an embedding of a based quantum torus 
into a division algebra $\FF$ such that $\FF$ equals the classical 
quotient ring of the image.

For a toric frame $M : \Zset^N \to \FF$, there exists
a unique matrix $\rbf = (r_{kj}) \in M_N(\kx)$ with the properties 
of the definition, and it is determined by $M$ by
$$
r_{kj} = M(e_k) M(e_j) M(e_k + e_j)^{-1}, \quad \forall 1 \le k < j \le N.
$$
It will be called {\em{the matrix of the toric frame}} and will be denoted by $\rbf(M) := \rbf$.
\index{matrix of a toric frame} \index{rbM@$\rbf(M)$}

Suppose that $M : \Zset^N \to \FF$ is a toric frame, with matrix $\rbf$ and algebra embedding $\vp: \Tbb_{\rbf^{\cdot 2}} \hra \FF$. Let $\sigma \in GL_N(\Zset)$, and let $\psi_\sigma : \Tbb_{\rbf_\sigma^{\cdot 2}} \to \Tbb_{\rbf^{\cdot 2}}$ be the isomorphism \eqref{psisig}. Then $\vp \psi_\sigma : \Tbb_{\rbf_\sigma^{\cdot 2}} \to \FF$ is an algebra embedding such that $\vp \psi_\sigma (Y_i) = M\sigma(e_i)$ for all $i \in [1,N]$. Moreover, using \eqref{sig-inv} we see that $\vp \psi_\sigma (Y^{(f)}) = M\sigma(f)$ for all $f \in \Zset^N$. Thus,
\begin{equation}
\label{Msig}
M\sigma : \Zset^N \to \FF \; \; \mbox{is a toric frame with matrix} \; \; \rbf_\sigma = {}^{\sigma^T}\rbf^\sigma.
\end{equation}

\section{Compatible pairs}
\label{3.2}
We fix in addition a positive integer $n$ such that $n \leq N$ and a 
subset $\ex \subseteq [1,N]$ of cardinality $n$. The indices in
$\ex$ will be called \emph{exchangeable}, and those in $[1,N] \backslash \ex$ 
\emph{frozen}. 
\index{ex@$\ex$} \index{exchangeable indices} \index{frozen indices}

{\em{The entries of all $N \times n$, $n \times N$, and  $n \times n$ matrices 
will be indexed by the sets $N \times \ex$, $\ex \times N$, and $\ex \times \ex$, 
respectively.}} For a matrix 
$\wt{B} = (b_{kj}) \in M_{N \times n} (\Zset)$, its $\ex \times \ex$ 
submatrix will be denoted by $B$ and will be called 
the \emph{principal part} \index{principal part} of $\wt{B}$. The columns of 
$\wt{B}$ will be denoted by $b^j \in \Zset^N$, $j \in \ex$. Recall that $B$ is said to be \emph{skew-sym\-met\-riz\-able} \index{skew-sym\-met\-riz\-able} provided there exist positive integers $d_k$, $k \in \ex$ such that $d_k b_{kj} = - d_j b_{jk}$ for all $k,j \in \ex$.

Let $\rbf \in M_N(\kx)$ be a multiplicatively skew-symmetric matrix and
$\wt{B} = (b_{kj}) \in M_{N \times n} (\Zset)$. 
Define the matrix $\wt{\tb} := {}^{\wt{B}^T} \rbf \in M_{n \times N} (\kx)$.  \index{t@$\wt{\tb}$}
Its entries are given by 
\begin{equation}
\label{t-matr}
t_{kj} := \Om_\rbf(b^k, e_j)= \big( {}^{\wt{B}^T} \! \rbf \big)_{kj} = 
{}^{(b^k)^T} \! \rbf^{e_j} =
\prod_{l=1}^N r_{lj}^{b_{lk}},
\; \; \forall k \in \ex, \; j \in [1,N].
\end{equation}
The pair $(\rbf, \wt{B})$ is {\em{compatible}} \index{compatible pair}
if the following two conditions are satisfied:
\begin{align}
&
t_{kj} = 1, \; \; \forall k \in \ex, \; j \in [1,N], \; k \neq j
\quad \mbox{and}
\label{comp1}
\\
&t_{kk} \; \; 
\mbox{are not roots of unity}, \; \; \forall k \in \ex.
\label{comp2}
\end{align}

Denote the $\ex \times \ex$ submatrix of $\wt{\tb}$ by $\tb$. 
The condition \eqref{comp1} implies that
\begin{equation}
\label{tt}
\tb^B = \wt{\tb}^{\wt{B}}.
\end{equation}

\bpr{fullrank} If the pair $(\rbf, \wt{B})$ is compatible,
then $\wt{B}$ has full rank and the matrix
$$
\tb^B=(t_{kk}^{b_{kj}})_{k,j \in \ex} \in M_{n \times n} (\kx)
$$
is multiplicatively skew-symmetric.
\epr

\begin{proof} First we show that $\wt{B}$ has full rank.
If this is not the case, then $\Ker(\wt{B}) \cap \Zset^n \neq 0$.
Let $\sum_{j \in \ex} a_j e_j \in \Ker(\wt{B}) \cap \Zset^n$ be a nonzero vector.
Applying \eqref{comp1}, we obtain
$$
t_{jj}^{a_j} = \prod_k t_{kj}^{a_k}  = \prod_{k,l} r_{lj}^{b_{lk} a_k} =1, \; \; \forall j \in \ex.  
$$
The condition \eqref{comp2} implies $a_j = 0$ for all $j \in \ex$, which is a contradiction.

Since $\tb^B = \wt{\tb}^{\wt{B}}= {}^{\wt{B}^T} \! \rbf^{\wt{B}}$, 
cf. \eqref{tt}, the second statement follows at once from the fact that 
$\rbf$ is multiplicatively skew-symmetric.
\end{proof}

Unlike the uniparameter case \cite[Proposition 3.3]{BZ}, compatibility of $(\rbf, \wt{B})$ does not in general imply that $B$ is skew-symmetrizable. The next lemma describes an instance 
when this condition appears naturally.

\ble{alternative} Assume that $(\rbf, \wt{B})$ is a compatible pair. If there exist positive integers 
$d_k$, $k \in \ex$ such that   
\begin{equation}
\label{td}
t_{kk}^{d_j} = t_{jj}^{d_k}, \; \;
\forall j,k \in \ex,
\end{equation}
then the principal part of $\wt{B}$ is skew-symmetrizable via these $d_k$, i.e., $d_k b_{kj} = - d_j b_{jk}$ for all $k,j \in \ex$.
\ele

\begin{proof} By \prref{fullrank},
$t_{kk}^{b_{kj}} = t_{jj}^{- b_{jk}}$ for all $j, k \in \ex$.
Combining this with \eqref{td} leads to 
$$
t_{jj}^{d_k b_{kj}}=  t_{kk}^{b_{kj} d_j} = t_{jj}^{ - d_j b_{jk}}.
$$ 
Thus $d_k b_{kj} = - d_j b_{jk}$ because of \eqref{comp2}.
\end{proof}

\bre{powers} The condition \eqref{td} is satisfied if there exists 
a non-root of unity $q \in \kx$ such that each $t_{kk} = q^{m_k}$ for some 
$m_k \in \Zset_{> 0}$. Then one can set $d_k = m_k$. 
\ere

\section{Mutation of compatible pairs}
\label{3.3}
We proceed with the definition of mutation of a compatible pair $(\rbf, \wt{B})$ 
in direction $k \in \ex$. The entries $(b'_{ij})$ 
of the mutation $\mu_k(\wt{B})$ \index{mukB@$\mu_k(\wt{B})$} of the second matrix are given 
by the Fomin--Zelevinsky formula \cite{FZ}:
\begin{equation}
\label{Bm}
b'_{ij} := 
\begin{cases} 
- b_{ij}, &\mbox{if} \; \; 
i=k \; \; \mbox{or} \; \; j=k
\\
b_{ij} + \frac{ |b_{ik}| b_{kj} + b_{ik} | b_{kj}|}{2}, 
&\mbox{otherwise}.
\end{cases}
\end{equation}
By \cite[Eq. (3.2)]{BFZ}, the matrix $\mu_k(\wt{B})$ 
is also given by 
$$
\mu_k ( \wt{B}) = E_\ep \wt{B} F_\ep 
$$ 
for both choices of signs $\ep = \pm$, where $E_\ep = E_\ep^{\wt{B}}$ and $F_\ep = F_\ep^{\wt{B}}$ are 
the $N \times N$ and $n \times n$ matrices with entries given by \index{Eepsilon@$E_\ep$} \index{Fepsilon@$F_\ep$}
\begin{equation}
\label{Eep}
\begin{aligned}
(E_\ep)_{ij} &= 
\begin{cases}
\delta_{ij}, & \mbox{if} \; j \neq k \\
-1, & \mbox{if} \; i=j=k \\
\max(0, - \ep b_{ik}), & \mbox{if} \; i \neq j = k
\end{cases}  \\
(F_\ep)_{ij} &= 
\begin{cases}
\delta_{ij}, & \mbox{if} \; i \neq k \\
-1, & \mbox{if} \; i=j=k \\
\max(0, \ep b_{kj}), & \mbox{if} \; i =k \neq j.
\end{cases}
\end{aligned}
\end{equation}
We refer the reader to \cite{FZ,BZ} for a detailed 
discussion of the properties of the mutation of the matrix $\wt{B}$. 
Finally, we define the mutation in direction $k$ 
of the first matrix of the compatible pair $(\rbf, \wt{B})$ by
\begin{equation}
\label{r-mut-k}
\mu_k(\rbf) := {}^{E_\ep^T} \! \rbf^{E_\ep}.
\end{equation}
\index{mukr@$\mu_k(\rbf)$}
This is a multiplicative version of \cite[Eq. (3.4)]{BZ}.

\bpr{pair-mut} Let $(\rbf, \wt{B})$ be a compatible pair.

{\rm(a)} The matrix $\mu_k(\rbf)$, defined in \eqref{r-mut-k}, 
does not depend on the choice of sign $\ep= \pm$. 

{\rm(b)} Assume also that the principal part of $\wt{B}$ is skew-symmetrizable. Then the pair $(\mu_k(\rbf),\mu_k(\wt{B}))$ is compatible, the principal part of $\mu_k(\wt{B})$ is skew-symmetrizable, and the 
$\wt{\tb}$-matrix of this pair coincides with the one 
of the pair $(\rbf, \wt{B})$.
\epr

We define the \emph{mutation in direction $k \in \ex$} of the compatible 
pair $(\rbf, \wt{B})$ to be the compatible pair $(\mu_k(\rbf),\mu_k(\wt{B}))$.
\index{mutation of a compatible pair}

\begin{proof}[Proof of \prref{pair-mut}] We follow the argument of the proof of 
\cite[Prop\-osition 3.4]{BZ} phrased in terms of multiplicative expressions.

(a) As in \cite[Eqs. (3.5) and (3.7)]{BZ}, $E_\ep^2 = I_N$, $F_\ep^2 = I_n$ 
and (a correction) $E_+ E_-$ equals the $N \times N$ matrix $G$ with columns
$$ 
g_j = e_j + \delta_{jk} b^k, \; \; j \in [1,N].
$$
Because of \eqref{mult}, part (a) is equivalent to saying that 
${}^{G^T} \! \rbf^G = \rbf$, which is verified as follows
(using \eqref{comp1}--\eqref{comp2}):
\begin{align*}
\big( {}^{G^T} \! \rbf^G \big)_{ij} &= 
{}^{e_i^T + \delta_{ik}(b^k)^T} \! 
\rbf^{e_j + \delta_{jk} b^k} 
= \Om_\rbf (e_i + \delta_{ik} b^k, e_j + \delta_{jk} b^k) 
\\
 &= \Om_\rbf (e_i,e_j) t_{kk}^{\delta_{ik} \delta_{jk}} t_{kk}^{-\delta_{ik} \delta_{jk}}
= \Om_{\rbf}(e_i,e_j) = r_{ij}, \; \; \forall i,j \in [1,N].
\end{align*}

(b) We have
$$
{}^{\mu_k(\wt{B})^T} \! (\mu_k(\rbf)) = 
{}^{F_\ep^T \wt{B}^T E_\ep^T E_\ep^T} \! \rbf^{E_\ep} = 
{}^{F_\ep^T \wt{B}^T} \! \rbf^{E_\ep} = {}^{F_\ep^T} \! {\wt{\tb}}^{E_\ep}.
$$
The second statement in \prref{fullrank} and the fact that for all $i, j \in \ex$,
$b_{ij}$ and $-b_{ji}$ have the same signs (which follows from the skew-symmetrizability assumption) imply 
$$
\wt{\tb}^{E_\ep} = {}^{F_\ep^T} \! \wt{\tb}. 
$$
Therefore 
$$
{}^{\mu_k(\wt{B})^T} \! (\mu_k(\rbf)) = {}^{F_\ep^T} \! \wt{\tb}^{E_\ep} = \wt{\tb}
$$
and the pair $(\mu_k(\rbf), \mu_k( \wt{B}))$ satisfies the 
conditions \eqref{comp1}--\eqref{comp2}. Furthermore, the $\wt{\tb}$-matrices of the compatible pairs 
$(\rbf, \wt{B})$ and $(\mu_k(\rbf), \mu_k(\wt{B}))$ are equal. The principal part of $\mu_k(\wt{B})$, namely $\mu_k(B)$, is skew-symmetrizable (for the same choice of positive integers $d_j$, $j \in \ex$ as for $B$), by the observations in \cite[Proposition 4.5]{FZ}.
\end{proof}

\section{Quantum seeds and mutations}
\label{3.4}
\bde{Q-seed} A {\em{quantum seed}} \index{quantum seed} of a division algebra $\FF$ is a pair
$(M, \wt{B})$ consisting of a toric frame $M : \Zset^N \to \FF$ and 
an integer matrix $\wt{B} \in M_{N \times n} (\Zset)$
satisfying the following two conditions:
\begin{enumerate}
\item[(i)] The pair $(\rbf(M), \wt{B})$ is compatible, where $\rbf(M)$ 
is the matrix of $M$. 
\item[(ii)] The principal part of $\wt{B}$ is skew-symmetrizable.
\end{enumerate}

The elements $M(e_1), \ldots, M(e_N) \in \FF$ will be called 
{\em{cluster variables}} of the seed $(M, \wt{B})$. 
The cluster variables $M(e_k)$, $k \in \ex$ will be called 
{\em{exchangeable}} and the $M(e_k)$ for $k \in [1,N] \backslash \ex$ 
{\em{frozen}}. \index{cluster variables} \index{exchangeable variables} \index{frozen variables}
\ede

The following lemma is an analog of \cite[Proposition 4.2]{BZ}. It will be used 
in defining mutations of quantum seeds.

\ble{special-aut} Assume that $\rbf \in M_N(\kx)$ is a multiplicatively skew-symmetric matrix, and that $k \in [1,N]$
and $g \in \Zset^N$ are such that $\Om_\rbf(g, e_j) = 1$ for $j \neq k$ 
and $\Om_\rbf(g, e_k)$ is not a root of unity. For each $\ep= \pm$, there exists 
a unique automorphism $\rho_{g,\ep} = \rho^\rbf_{g,\ep}$  \index{rhorgepsilon@$\rho^\rbf_{g,\ep}$}  of $\Fract(\Tbb_{\rbf^{\cdot 2}})$ such that
\begin{equation}
\label{rho}
\rho_{g,\ep}( Y^{(e_j)} ) = \begin{cases}
Y^{(e_k)} + Y^{(e_k + \ep g)}, &\mbox{if} \; \; j=k  \\
Y^{(e_j)}, &\mbox{if} \; \; j \ne k. \end{cases}
\end{equation}
\ele

\begin{proof} Since $\Om_\rbf(e_k + \ep g, e_j) = \Om_\rbf(e_k, e_j)$ for $j \ne k$, we have
$$\bigl( Y^{(e_k)} + Y^{(e_k + \ep g)} \bigr) Y^{(e_j)} = r_{kj}^2 Y^{(e_j)} \bigl( Y^{(e_k)} + Y^{(e_k + \ep g)} \bigr),  \; \; \forall j \ne k.$$
Hence, there exists a $\KK$-algebra homomorphism $\rho_{g,\ep} : \Tbb_{\rbf^{\cdot 2}} \rightarrow \Fract(\Tbb_{\rbf^{\cdot 2}})$ satisfying \eqref{rho}. 
Denote $g = \sum_j n_j e_j$. We have $\Om_\rbf(g, e_k)^{n_k} = \Om_\rbf(g,g) =1$.
Thus $n_k =0$ because $\Om_\rbf(g, e_k)$ is assumed to be a non-root of unity.
This implies that $\rho_{g, \ep}(Y^{(g)}) = Y^{(g)}$. Consequently, $\rho_{g,\ep}$ is the identity on the set $G := \KK[Y^{(g)}]\backslash \{0\}$. Now $G$ is an Ore set in $\Tbb_{\rbf^{\cdot 2}}$ (observe that $Y^{(f)} G = G Y^{(f)}$ for all $f \in \Zset^N$), and so $\rho_{g,\ep}$ extends to a $\KK$-algebra endomorphism of $\Tbb_{\rbf^{\cdot 2}} G^{-1}$.

Similarly, there exists a $\KK$-algebra endomorphism $\rho'$ of $\Tbb_{\rbf^{\cdot 2}} G^{-1}$ such that
$$\rho'( Y^{(e_j)} ) = \begin{cases}
\bigl( 1 + \Om_\rbf(g,e_k)^{-\ep} Y^{(\ep g)} \bigr)^{-1} Y^{(e_k)}, &\mbox{if} \; \; j=k  \\
Y^{(e_j)}, &\mbox{if} \; \; j \ne k.  \end{cases}$$
Now $\rho_{g,\ep} \rho' (Y^{(e_j)}) = \rho' \rho_{g,\ep} (Y^{(e_j)}) = Y^{(e_j)}$ for all $j \in [1,N]$. Therefore $\rho_{g,\ep}$ and $\rho'$ are inverse automorphisms of $\Tbb_{\rbf^{\cdot 2}} G^{-1}$, so they extend to inverse automorphisms of $\Fract(\Tbb_{\rbf^{\cdot 2}})$. The uniqueness of $\rho_{g,\ep}$ is clear, because the elements $Y^{(e_j)}$, $j \in [1,N]$ generate $\Fract(\Tbb_{\rbf^{\cdot 2}})$ as a division algebra.
\end{proof}

Straightforward calculations yield that
\begin{equation}
\label{fullrho}
\rho^\rbf_{g,\ep} (Y^{(f)}) = \begin{cases}
P^{\rbf,m_k}_{g,\ep,+} Y^{(f)}, &\mbox{if} \; \; m_k \ge 0  \\
\bigl( P^{\rbf,-m_k}_{g,\ep,-} \bigr)^{-1} Y^{(f)}, &\mbox{if} \; \; m_k < 0  \end{cases}
\end{equation}
for $f = (m_1,\dots,m_N)^T \in \Zset^N$, where
$$P^{\rbf,m}_{g,\ep,\pm} := \prod_{i=1}^m \bigl( 1 + \Om_\rbf(g,e_k)^{\mp \ep (2i-1)} Y^{(\ep g)} \bigr), \; \; \forall m \in \Zset_{\ge 0}.$$

Unlike the uniparameter case \cite[Eq. (4.9)]{BZ}, here $\rho_{-g,-\ep}$ equals $\rho_{g,\ep}$ rather than $\rho_{g,\ep}^{-1}$. This stems from the fact that in the general situation, there is no analog of the minimal positive values $d(-)$ used in \cite[Subsec. 4.2]{BZ}.

For the rest of this section, we assume that $(M, \wt{B})$ is a quantum seed of $\FF$
and denote $\rbf := \rbf(M)$. We also fix $k \in \ex$. 
Since $E_\ep \in GL_N(\Zset)$, \eqref{Msig} shows that
\begin{equation}
\label{tfra}
M E_\ep : \Zset^N \to \FF \; \; 
\mbox{is a toric frame with matrix} \; \; 
\mu_k(\rbf) = {}^{E_\ep^T} \! \rbf^{E_\ep}. 
\end{equation}
By abuse of notation, we will identify $\FF$ with $\Fract(\Tbb_{\mu_k(\rbf)^{\cdot 2}})$
via the embedding $\Tbb_{\mu_k(\rbf)^{\cdot 2}} \hra \FF$ associated 
to this toric frame.    From \eqref{comp1} and the fact that $E_\ep b^k = b^k$, it follows that
$$
\Om_{\mu_k(\rbf)} (b^k, e_j) = 
\Om_\rbf(b^k, E_\ep e_j ) = t_{jj}^{- \delta_{jk}}, \; \; \forall  j \in [1,N].
$$
Applying the condition \eqref{comp2} and \leref{special-aut}, 
we obtain the automorphisms $\rho_{b^k, \ep} = \rho_{b^k,\ep}^{\mu_k(\rbf)} \in \Aut( \FF)$. When applying these automorphisms, any term $Y^{(f)} = M(f)$ in an expression such as \eqref{rho} or \eqref{fullrho} must be replaced by $M E_\ep (f)$.

\bpr{seed-mut} Assume that $(M, \wt{B})$ is a quantum seed and $k \in \ex$.

{\rm(a)} The map $\rho_{b^k, \ep} M E_\ep : \Zset^N \to \FF$ is 
a toric frame of $\FF$ and is independent of the choice of sign $\ep = \pm$.
The matrix of this toric frame equals $\mu_k (\rbf(M))$.

{\rm(b)} The pair $(\rho_{b^k, \ep} M E_\ep, \mu_k(\wt{B}))$ is a quantum seed of $\FF$. 
\epr

\begin{proof}
(a) The fact that $\rho_{b^k, \ep} M E_\ep : \Zset^N \to \FF$ is a toric frame 
of $\FF$ is immediate because $E_\ep \in GL_N(\Zset)$ and 
$\rho_{b^k, \ep} \in \Aut( \FF)$.
The last fact and \eqref{tfra} imply that 
its matrix is $\mu_k(\rbf)$. By \eqref{rho},
\begin{equation}
\label{m1}
\rho_{b^k, \ep} M E_\ep(e_j) = M E_\ep(e_j) = M(e_j), \; \; \forall j \neq k
\end{equation}
and
\begin{equation}
\label{m2}
\begin{aligned}
\rho_{b^k, \ep} M E_\ep(e_k) &= ME_\ep (e_k) + ME_\ep (e_k + \ep b^k) = M(E_\ep e_k) + M(E_\ep e_k + \ep b_k)  \\
 &= M(-e_k + [-\ep b^k]_+ ) + M(-e_k + [-\ep b^k]_+ + b^k)  \\
 &=  M(-e_k + [b^k]_+) + M(- e_k - [b^k]_-),
\end{aligned}
\end{equation}
cf. \eqref{g+-}.
Hence, the toric frame $\rho_{b^k, \ep} M E_\ep$ 
is independent of the choice of sign $\ep = \pm.$

(b) By \prref{pair-mut} (b), the pair $(\rbf(\rho_{b^k, \ep} M E_\ep), \mu_k(\wt{B})) =(\mu_k(\rbf), \mu_k( \wt{B}))$ 
is compatible, and the principal part of $\mu_k(\wt{B})$ 
is skew-symmetrizable.
\end{proof}

\bde{q-mut} Define the \emph{mutation in direction $k \in \ex$}
of the quantum seed $(M, \wt{B})$ of a division algebra $\FF$ 
to be the 
quantum seed $(\mu_k(M), \mu_k(\wt{B}))$, where $\mu_k(M) := \rho_{b^k, \ep} M E_\ep = \rho_{b^k, \ep}^{\mu_k(\rbf)} M E_\ep^{\wt{B}}$.
\index{mutation of a quantum seed} \index{mukM@$\mu_k(M)$}
\ede

\bco{inv} For all quantum seeds $(M, \wt{B})$ of $\FF$ and $k \in \ex$,
the following hold:

{\rm(a)} The matrix of the toric frame $\mu_k(M)$ equals $\mu_k(\rbf(M))$. This toric 
frame satisfies
\begin{equation}
\label{clust-mut}
\begin{aligned}
\mu_k(M)(e_j) &= M(e_j), \; \; \forall j \neq k,  \\
\mu_k(M)(e_k) &= M(-e_k + [b^k]_+) + M(- e_k - [b^k]_-),
\end{aligned}
\end{equation}
cf. \eqref{g+-}.

{\rm(b)} Mutation is involutive:
$$
\mu_k^2 (M, \wt{B}) = (M, \wt{B}).
$$
\eco

\begin{proof}
Part (a) follows from \prref{seed-mut} (a) and \eqref{m1}--\eqref{m2}. Turning to part (b), we use the independence of signs (\prref{seed-mut}(a)) to write
$$\mu_k^2(M) = \rho_{\mu_k(\wt{B})^k,-\ep}^{\mu_k^2(\rbf)} \bigl( \rho_{b^k, \ep}^{\mu_k(\rbf)} M E_\ep^{\wt{B}} \bigr) E_{-\ep}^{\mu_k(\wt{B})} = \rho_{-b^k,-\ep}^\rbf \rho_{b^k, \ep}^{\mu_k(\rbf)} M E_\ep^{\wt{B}} E_\ep^{\wt{B}} ,$$
and we apply \eqref{fullrho} to see that $\mu_k^2(M) = M$ (it suffices to check this on the $Y^{(e_j)}$, $j \in [1,N]$). Finally, we recall the well known fact \cite[p. 510]{FZ} 
that $\mu_k^2 (\wt{B}) = \wt{B}$. 
\end{proof}

Two quantum seeds $(M, \wt{B})$ and $(M', \wt{B}')$ of a division algebra 
$\FF$ will be called {\em{mutation-equivalent}} if they can be obtained from 
each other by a sequence of mutations. \deref{q-mut} implies at once:
\index{mutation-equivalent quantum seeds}

\bco{qmut-aut} If $(M, \wt{B})$ and $(M', \wt{B}')$ are mutation-equivalent quantum seeds of $\FF$, 
then there exists 
$\sigma \in GL_N(\Zset)$ and $\psi \in \Aut(\FF)$ such that 
$$
M' = \psi M \sigma : \Zset^N \to \FF.
$$
In particular, the following subgroups of $\kx$ are equal:
$$
\langle \Om_{\rbf(M)}(f,g) \mid f, g \in \Zset^N \rangle = 
\langle \Om_{\rbf(M')}(f,g) \mid f, g \in \Zset^N \rangle.
$$
\eco

The first part of \coref{qmut-aut} is one of the conditions
in the Berenstein--Zelevinsky definition \cite{BZ} of toric frame in the 
case of uniparameter quantum tori. 

The following well known equivariance of mutation of pairs will be needed 
in Chapter \ref{main}. For a group $G$ we denote by $X(G)$ its 
character lattice. For a $G$-eigenvector $u$, $\chi_u$ will
denote its $G$-eigenvalue.

\ble{1eig} Let $\FF$ be a division algebra and $G$ a group acting 
on $\FF$ by algebra automorphisms. Assume that $(M,\wt{B})$ is a 
quantum seed for $\FF$ such that $M(f)$ is a $G$-eigenvector 
for $f \in \Zset^N$ {\rm(}or, which is the same, for $f = e_1, \ldots, e_N${\rm)}
and $\chi_{M(b^j)} =1$ for all columns $b^j$ of $\wt{B}$. Then 
all mutations $(\mu_k(M), \mu_k(\wt{B}))$, $k \in \ex$ 
have the same properties. 
\ele 

\begin{proof} Fix $k \in \ex$ and denote for brevity
$M':=\mu_k(M)$ and $\wt{B}':= \mu_k( \wt{B})$. 
The entries and columns of $\wt{B}'$ will be 
accordingly denoted by $b'_{lj}$ and $(b')^j$. 
Since $\chi_{M(b^k)}=1$, we have
$\chi_{ M( - e_k + [b^k]_+)} = \chi_{M(- e_k - [b^k]_-)}$. 
Hence, $M'(e_k)$ is a $G$-eigenvector and 
\begin{equation}
\label{eigenvalues}
\chi_{M'(e_k)} = \chi_{ M( - e_k + [b^k]_+)} = \chi_{M(- e_k - [b^k]_-)}.
\end{equation}
Of course, $M'(e_j) = M(e_j)$ is a $G$-eigenvector for all $j \ne k$.

For the second condition, we have $\chi_{M'((b')^k)} = \chi_{M(-b^k)}=1$ 
since $b'_{kk}=0$. Now, let $j \neq k$. If $b_{kj}=0$, then 
$\chi_{M'((b')^j)}= \chi_{M(b^j)}=1$. If $b_{kj} \neq 0$,
denote $\ep = \sign(b_{kj})$. By \eqref{Bm},
$(b')^j = b^j + \ep b_{kj}[b^k]_\ep - 2 b_{kj} e_k$,
and using \eqref{eigenvalues} we obtain
\begin{align*}
\chi_{M'((b')^j)} &= \chi_{M(b^j + \ep b_{kj}[b^k]_\ep - b_{kj} e_k)} 
\chi_{M'( - b_{kj} e_k)}  \\
&= \chi_{M(b^j + \ep b_{kj}[b^k]_\ep - b_{kj} e_k)} 
\chi_{M(b_{kj} e_k -\ep b_{kj} [b^k]_\ep )} = \chi_{M(b^j)} = 1.  \qedhere
\end{align*}
\end{proof}

\section{Quantum cluster algebras and the Laurent phenomenon}
\label{3.5}
Let $(M, \wt{B})$ be a quantum seed of a division algebra $\FF$ over $\KK$. 
Fix a subset $\inv \subseteq [1,N] \backslash \ex$. \index{inv@$\inv$} This set will be used 
to determine which frozen variables 
(indexed by $[1,N] \backslash \ex$) will be inverted 
in the definition of quantum cluster algebras. 
Let $D$ be a unital subring of $\KK$ containing the subgroup 
\begin{equation}
\label{La-subgr}
\langle \Om_{\rbf}(f,g) \mid f,g \in \Zset^N \rangle
\end{equation}
of $\kx$, where $\rbf = \rbf(M)$. The ring $D$ will be used as a base ring.

\bde{q-cl} Define the {\em{quantum cluster algebra}} $\Abb(M, \wt{B}, \inv)_D$ associated to the above 
data to be the unital $D$-subalgebra of $\FF$ 
generated by all cluster variables $M'(e_k)$, $k \in [1,N]$ and the 
inverses $M'(e_l)^{-1}$, $l \in \inv$ for all quantum seeds $(M', \wt{B}')$ of $\FF$ 
which are mutation-equivalent to $(M, \wt{B})$.
\index{quantum cluster algebra}  \index{AMBinvD@$\Abb(M, \wt{B}, \inv)_D$}

Define the {\em{upper quantum cluster algebra}} $\UU(M, \wt{B}, \inv)_D$ to be the intersection 
of all $D$-sub\-al\-ge\-bras of $\FF$ of the form
\begin{equation}
\label{tor}
\Tbb\Abb_{(M', \wt{B}'), D}:=
D \lcor M'(e_l)^{\pm 1}, \, M'(e_k) \mid l \in \ex \sqcup \inv, \, k \in [1,N] \backslash
(\ex \sqcup \inv) \rcor
\end{equation}
taken over all quantum seeds $(M', \wt{B}')$ for $\FF$ 
mutation-equivalent to $(M, \wt{B})$.
\index{upper quantum cluster algebra}  \index{UMBinvD@$\UU(M, \wt{B}, \inv)_D$}  \index{TAMBD@$\Tbb\Abb_{(M', \wt{B}'), D}$}
\ede

Since $M'(e_l)$ is a frozen variable for $l \in [1,N] \backslash \ex$, we have that 
$M'(e_l) = M(e_l)$ for all $l \in [1,N] \backslash \ex$ and quantum seeds 
$(M', \wt{B}')$ that are mutation equivalent to $(M, \wt{B})$. 
Because of this, in the definition of quantum cluster algebras 
it is sufficient to consider the generators $M'(e_j)$, 
$j \in \ex$, $M(e_l)^{\pm 1}$, $l \in \inv$ and 
$M(e_l)$, $l \in [1,N] \backslash (\ex \sqcup \inv)$.

The subalgebra $\Tbb\Abb_{(M', \wt{B}'),D}$ is a mixed quantum torus--quantum affine 
space $D$-subalgebra of $\FF$. The $D$-subalgebra of $\FF$ generated by 
$M'(e_j)$, $j \in [1,N]$ is isomorphic to a quantum affine space algebra.
The algebra $\Tbb\Abb_{(M', \wt{B}'),D}$ is equal to the localization of it 
by the multiplicative subset generated by $M'(e_l)$, $l \in \ex \sqcup \inv$. 
If one localizes instead by all $M'(e_l)$, $l \in [1,N]$, 
then the resulting algebra will be the corresponding quantum torus.
The mixed nature of $\Tbb\Abb_{(M', \wt{B}'),D}$ has to do with this partial localization.   

The proof of the Laurent phenomenon of Fomin and Zelevinsky \cite{FZ, FZ-l} 
carries over to the quantum case. We state this result and provide all needed details 
for the quantum case in the rest of this section.
In the uniparameter case and when all frozen variables are inverted, 
a quantum Laurent phenomenon was established by  
Berenstein and Zelevinsky \cite{BZ} following their previous work 
with Fomin \cite{BFZ} using upper bounds.

\bth{q-l} For all quantum seeds $(M, \wt{B})$ of a division algebra 
$\FF$ over a field $\KK$, subrings $D$ of $\KK$ containing the subgroup 
\eqref{La-subgr} of $\kx$, and subsets $\inv \subseteq [1,N] \backslash \ex$,
we have the inclusion
$$
\Abb(M, \wt{B}, \inv)_D \subseteq \UU(M, \wt{B}, \inv)_D.
$$
\eth

We will need the following lemma.

\ble{prim-tor} Let $\rbf \in M_N(\kx)$ be a multiplicatively 
skew-symmetric matrix and $D$ a subring of $\KK$ which contains the 
subgroup \eqref{La-subgr} of $\kx$. Consider elements
$$
V := \prod_{i=1}^m \bigl( Y^{(f)} + \xi_i \bigr) \qquad \mbox{and} \qquad 
V' := \prod_{i=1}^{m'} \bigl( Y^{(f')} + \xi'_i \bigr)
$$
in the $D$-quantum torus $D \lcor Y_1^{\pm 1}, \ldots, Y_N^{\pm 1} \rcor \subseteq \Tbb_{\rbf^{\cdot 2}}$, 
where $f,f' \in \Zset^N$ are such that
$$\Om_\rbf(f, e_k)^2 = \Om_\rbf(f', e_k)^2 = 1, \; \; \forall k \in [1,N]$$
and $\xi_1,\dots, \xi_m, \xi'_1, \dots, \xi'_{m'} \in D$.

If $f$ and $f'$ are $\Zset$-linearly independent, then $V$ and $V'$ 
are relatively prime in the center 
of $D \lcor Y_1^{\pm 1}, \ldots, Y_N^{\pm 1} \rcor$. 
\ele

\begin{proof} We may assume that 
$\xi_1,\dots, \xi_m, \xi'_1, \dots, \xi'_{m'}$ are nonzero elements of $D$,
since the factors with vanishing $\xi$'s are central units. Similarly, we may assume that $f$ and $f'$ are nonzero.

The center of $D \lcor Y_1^{\pm 1}, \ldots, Y_N^{\pm 1} \rcor$
is given by 
\begin{equation}
\label{Z}
Z = \sum_{g \in L} D Y^{(g)} 
\end{equation}
for the sublattice $L \subseteq \Zset^N$ defined by
$$
L = \{ g \in \Zset^N \mid \Om_\rbf(g, e_k)^2 = 1, \; \forall k \in [1,N] \}.
$$

After reducing to the center of $D \lcor Y_1^{\pm 1}, \ldots, Y_N^{\pm 1} \rcor$, replacing $D$ by the algebraic closure of its quotient field, and changing notation appropriately, we may assume that $D$ is an algebraically closed field and $D \lcor Y_1^{\pm 1}, \ldots, Y_N^{\pm 1} \rcor$ is a commutative Laurent polynomial ring (hence, a UFD).

Write $f = ng$ and $f' = n'g'$ where $n,n' \in \Zset_{>0}$ and the gcd of the entries of $g$, resp.~of $g'$, is $1$. By hypothesis, $g \ne \pm g'$. Now the irreducible factors of $V$, up to associates, are $Y^{(g)} - \zeta$ for $n\,$th roots $\zeta$ of $-\xi_i$, and similarly for $V'$. No $Y^{(g)} - \zeta$ is an associate of any $Y^{(g')} - \zeta'$, so $V$ and $V'$ have no irreducible common factors. Thus, they are relatively prime.
\end{proof}

The proof of \thref{q-l} is analogous to the proof of the Caterpillar Lemma
of Fomin and Zelevinsky \cite[Theorem 2.1]{FZ-l}. We will indicate the 
needed ingredients for the inductive step of the proof in the noncommutative 
situation.

From now on we will assume that 
$$
(M_s, \wt{B}_s), \; \; s \in [0,3] \quad
\mbox{are four quantum seeds of $\FF$} 
$$
and
$$
(M_1, \wt{B}_1) = \mu_{i} (M_0, \wt{B}_0),
\; \; 
(M_2, \wt{B}_2) = \mu_{j} (M_1, \wt{B}_1),
\; \;
(M_3, \wt{B}_3) = \mu_{i} (M_2, \wt{B}_2)
$$
for some $i \neq  j$ in $\ex$. Denote
\begin{align*}
x &:= M_0(e_i),
\\
y &:= M_0(e_j) = M_1(e_j),
\\
z &:= M_1(e_i) = M_2(e_i),
\\
u &:= M_2(e_j) = M_3(e_j),
\\
v &:= M_3(e_i).
\end{align*}
Set for brevity $\Tbb\Abb_s:=\Tbb\Abb_{(M_s, \wt{B}_s),D}$. 
There exist unique elements $P \in \Tbb\Abb_0$, $Q \in \Tbb\Abb_1$ and $R \in \Tbb\Abb_2$ such that
\begin{equation}
\label{zuv}
z= x^{-1} P, \; \; u = y^{-1} Q \; \; \mbox{and} \; \; 
v = z^{-1} R.
\end{equation}
Each of $P$, $Q$, $R$ is a sum of two nonzero monomials in the cluster variables for the 
seeds $(M_s, \wt{B}_s)$, $c \in [0,2]$. We will denote by 
$P(y)$, $Q(z)$, $R(u)$ the elements $P$, $Q$, $R$ written as {\em{polynomials}} 
in $y$, $u$, $ z$ with left hand coefficients that are {\em{polynomials}} in 
the other cluster variables in the seeds $(M_0, \wt{B}_0)$, 
$(M_1, \wt{B}_1)$, $(M_2, \wt{B}_2)$, respectively. 
(This uniquely defines $P(y)$, $Q(z)$, $R(u)$.) In this setting 
one can substitute any element of $\FF$ for $y$, $z$, or $u$ and the result will 
be an element of $\FF$. 

Clearly, 
\begin{align}
\label{uTA}
z&= x^{-1} P, \; \;
u= y^{-1} Q = y^{-1} Q(x^{-1} P(y)) \in \Tbb\Abb_0 \quad
\mbox{and}
\\
v&= z^{-1} R = z^{-1}R(y^{-1}Q(z)) \in \Tbb\Abb_1.
\label{vTA}
\end{align}

Recall that an element $p$ of a (noncommutative) domain $A$ is called prime if it is normal
(i.e., $Ap = pA$) and $Ap$ is a completely prime ideal (i.e., $A/Ap$ is a domain).
We refer the reader to Section \ref{2.1} for details on this notion. For $c \in A$, we write 
$p \mid c $ if $c \in Rp$. For $b,c \in A$, we 
write $b \mid_l c$  \index{zzz1@$\mid_l$}  if $c \in b A$.  

Denote by 
$\Tbb\Abb^{\vee}_0$ the $D$-subalgebra of $\Tbb\Abb_0$ generated by 
the list of cluster variables and their inverses which appear in \eqref{tor} 
except for $x^{\pm 1} = M_0(e_i)^{\pm 1}$, i.e.,  
$$
\Tbb\Abb^{\vee}_0 :=
D \lcor M_0(e_l)^{\pm 1}, \, M_0(e_k) \mid l \in (\ex \sqcup \inv)\backslash \{ i \}, \, k \in
[1,N] \backslash (\ex \sqcup \inv) \rcor.
$$
This is another mixed quantum torus--quantum affine space algebra, but this time
of GK dimension $N-1$. Clearly, 
\begin{equation}
\label{TA}
\Tbb\Abb_0 = \Tbb\Abb^{\vee}_0[x^{\pm 1}; \sigma],
\end{equation}
where $\sigma \in \Aut(\Tbb\Abb^{\vee}_0)$ is given by 
$$
\sigma(M_0(e_l)) = \rbf(M_0)_{il}^2 M_0(e_l), \; \; \forall l \in [1,N] \backslash \{ i \}.
$$
We will say that $b \in \Tbb\Abb_0\backslash \{0\}$ has \emph{leading term} $c x^m$ 
where $m \in \Zset$ and $c \in \Tbb\Abb^{\vee}_0$ if 
$$
b - c x^m \in \bigoplus_{m' < m} \Tbb\Abb^{\vee}_0 x^{m'}.
$$

The localizations
\begin{align*}
\Tbb_0^{\vee} &:= \Tbb\Abb_0^{\vee} [M_0(e_k)^{-1} \mid k \in [1,N] \backslash (\ex \sqcup \inv)]  \\ 
\Tbb_0 &:= \Tbb\Abb_0 [M_0(e_k)^{-1} \mid k \in [1,N] \backslash (\ex \sqcup \inv)]
\end{align*}
are quantum tori and 
$$
\Tbb_0 \cong \Tbb_0^{\vee} [x^{\pm 1}; \sigma],  
$$
where by abuse of notation we denote by the same letter the canonical 
extension of $\sigma \in \Aut(\Tbb\Abb_0^{\vee})$ to an automorphism 
of $\Tbb_0^{\vee}$.

\ble{quant-Cater} The following hold in the above setting:

{\rm(a)} The element 
$$
L:= z^{-1} R(y^{-1}Q(0)) z P^{-1}
$$
is a monomial in the elements
$M_0(e_k)$, $k \in [1,N] \backslash \{i,j\}$ and $y^{-1}= M_0(e_j)^{-1}$. 
In par\-tic\-u\-lar, it belongs to $\Tbb\Abb_0^{\vee}$. 

{\rm(b)} The cluster variables $u$ and $v$ belong to $\Tbb\Abb_0$ and their leading terms
are equal to $y^{-1} Q(0)$ and $Lx$, respectively.

{\rm(c)} If 
$$
z^m F = G u^{m_1} v^{m_2}
$$
for some $F, G \in \Tbb\Abb_0 \backslash \{ 0 \}$ and
$m,m_1, m_2 \in \Zset_{\geq 0}$, then $z^m \mid_l G$ 
as elements of $\Tbb\Abb_0$.
\ele
\begin{proof} 
We will denote by $b_{s, kl}$ the $kl$-th entry of the exchange matrix 
$\wt{B}_s$ and by $b^l_s$ the $l$-th column of $\wt{B}_s$.

(a) Set $\ep:= \sign(b_{1,ij})$. By a direct computation one obtains that 
$$
L = \xi M_2(E_\ep [b_2^i]_\ep - [b_1^i]_\ep ),
$$
where $\xi$ is an element of the subgroup \eqref{La-subgr} of $\kx$ and
$E_\ep$ denotes the first matrix in \eqref{Eep} for the mutation 
$(M_2, \wt{B}_2) = \mu_j(M_1, \wt{B}_1)$.
The $i$-th coordinate of the vector 
in the right hand side equals $0$, so 
$$
L = \xi M_1(E_\ep [b_2^i]_\ep - [b_1^i]_\ep ).
$$
Part (a) now follows by computing the signs of the coordinates of the vector 
in the right hand side of the last formula.

(b) The statement for $u$ follows from \eqref{uTA}. By part (a), 
$z^{-1} R(y^{-1}Q(0)) = LP z^{-1} = Lx$. Analogously to \cite{FZ-l},
write
\begin{align*}
v &= z^{-1} \left( R(y^{-1}Q(z)) - R(y^{-1}Q(0)) \right) + z^{-1} R(y^{-1}Q(0))
\\
&= 
z^{-1} \left( R(y^{-1}Q(z)) - R(y^{-1}Q(0)) \right) + Lx.
\end{align*}
Since 
$$
R(y^{-1}Q(z)) - R(y^{-1}Q(0)) \in \bigoplus_{m' >0} \Tbb\Abb_0^{\vee} z^{m'},
$$
we have 
$$
z^{-1} \left( R(y^{-1}Q(z)) - R(y^{-1}Q(0)) \right) \in \bigoplus_{m' \geq 0}\Tbb\Abb_0^{\vee} z^{m'}
\subseteq \bigoplus_{m' \geq 0} \Tbb\Abb_0^{\vee} x^{-m'}.
$$
This proves both statements for the cluster variable $v$.

(c) It is straightforward to see that $M_0(e_k)$, $k \in [1,N] \backslash (\ex \sqcup \inv)$ are 
prime elements of $\Tbb\Abb_0$ and that $M_0(e_k) \nmid z$. Since a prime element $p$ 
has the property that $p\mid ab  \Rightarrow$ $p \mid a$ or $p \mid b$, 
it is sufficient to prove that 
$$
z^m \mid_l G \quad 
\mbox{in} \quad
\Tbb_0.
$$
By part (a), $L$ is a unit of $\Tbb_0^{\vee}$ (in particular, a $\sigma$-eigenvector).
Noting that $P \in \Tbb \Abb_0^\vee$, using 
$$
z^m = (\sigma^{-1}(P) \ldots \sigma^{-m}(P)) x^{-m},
$$ 
and taking into account part (b) leads to 
\begin{equation}
\label{PFGident}
(\sigma^{-1}(P) \ldots \sigma^{-m}(P)) F' = G \bigl( (y^{-1}Q(0))^{m_1} 
x^{m_2} + \mbox{lower order terms} \bigr)
\end{equation}
for some $F'  \in \Tbb_0$. 

Denote $G = \sum_{m' \in \Zset}g_{m'} x^{m'}$, where $g_{m'} \in \Tbb\Abb_0^{\vee}$. 
We claim that
$$
(\sigma^{-1}(P) \ldots \sigma^{-m}(P)) \mid_l g_{m'} \quad 
\mbox{in} \quad
\Tbb_0^{\vee}, \; \; \forall m' \in \Zset.
$$
This will imply that $z^m \mid_l G$ in $\Tbb_0$ and will complete the proof of part (c).
Assume that this is not the case. Let $\ol{m}$ be the largest integer such that 
\begin{equation}
\label{contrad}
(\sigma^{-1}(P) \ldots \sigma^{-m}(P)) \nmid_l g_{\ol{m}}
\quad \mbox{in} \quad 
\Tbb_0^{\vee}.
\end{equation}
Comparing the components in 
$\Tbb_0^{\vee} x^{\ol{m}+m_2}$ in \eqref{PFGident} 
(with respect to the decomposition $\Tbb_0 = \oplus_{m' \in \Zset} \Tbb_0^{\vee} x^{m'}$)
we obtain that
\begin{equation}
\label{divv}
(\sigma^{-1}(P) \ldots \sigma^{-m}(P)) \mid_l g_{\ol{m}} \sigma^{\ol{m}} ((y^{-1} Q(0))^{m_1})
\quad
\mbox{in} 
\quad
\Tbb_0^{\vee}.
\end{equation}

If $b_{1,ij} \neq 0$, then $y^{-1} Q(0)$ is a unit of $\Tbb_0^{\vee}$ (in particular
a $\sigma$-eigenvector) and \eqref{divv} implies that 
$(\sigma^{-1}(P) \ldots \sigma^{-m}(P)) \mid_l g_{\ol{m}}$ 
in $\Tbb_0^{\vee}$ which contradicts \eqref{contrad}.

Now consider the case $b_{1,ij} = 0$. Then $Q(0)=Q$
and
\begin{equation}
\label{bij}
b_1^j, b_1^i \in \sum_{k \neq i,j} \Zset e_k,
\end{equation}
where the second property follows from the fact that 
the principal part of $\wt{B}_1$ is skew-symmetrizable.
Thus,
\begin{align*}
\sigma^{\ol{m}} ((y^{-1} Q(0))^{m_1}) &= (M_1(b_1^j) + \xi'_1) \ldots (M_1(b_1^j)+ \xi'_{m_1}) \theta'
\\&= (M_0(b_1^j) + \xi_1) \ldots (M_0(b_1^j)+ \xi_{m_1}) \theta'
\end{align*}
for some $\xi'_1, \ldots, \xi'_{m_1},\xi_1, \ldots, \xi_{m_1} \in D$ and
$\theta' \in \Tbb_0^{\vee}$ such that $\theta'$ is a Laurent monomial in the generators 
of $\Tbb_0^{\vee}$ and its coefficients and $\xi_1, \ldots, \xi_{m_1}$ belong to 
the subgroup \eqref{La-subgr} of the units of $D$. 
We also have 
\begin{align*}
\sigma^{-1}(P) \ldots \sigma^{-m}(P) &= \theta (M_0(b_0^i) + \xi_1) \ldots (M_0(b_0^i)+ \xi_m)
\\&= \theta (M_0(-b_1^i) + \xi_1) \ldots (M_0(-b_1^i)+ \xi_m),
\end{align*}
where $\xi_1, \ldots, \xi_m \in D$ and $\theta \in \Tbb\Abb_0^{\vee}$ 
have the same properties.
In particular, $\theta$ and $\theta'$ are units of $\Tbb_0^{\vee}$.
Denote the polynomials
$$
h(t):= (t+ \xi_1) \ldots ( t+ \xi_m), \; \; 
h'(t):=(t+ \xi'_1) \ldots (t+\xi'_{m_1}) \in D [t].
$$

Substituting the above two expressions in \eqref{divv} and clearing the units gives
\begin{equation}
\label{divv2}
h(M_0(-b_1^i)) \mid_l g_{\ol{m}} h'(M_0(b_1^j)) \quad
\mbox{in} \quad \Tbb_0^{\vee}.
\end{equation}
To get to the setting of \leref{prim-tor}, consider the subtorus
$$
\Tbb_0^{\vee\vee}:= D \lcor M_0(e_k)^{\pm 1}, k \neq i,j \rcor 
\subset \Tbb_0^{\vee}.
$$
It satisfies
$\Tbb_0^{\vee} = \Tbb_0^{\vee\vee}[y^{\pm 1}; \vp]$, where 
$\vp \in \Aut(\Tbb_0^{\vee\vee})$ is given by $\vp(M_0(e_k)) = \rbf(M_0)_{jk}^2 M_0(e_k)$,
$k \neq i,j$.
By \eqref{t-matr}--\eqref{comp1} and \eqref{bij}, $M_0(-b_1^i)$ and $M_0(b_1^j)$ belong to the 
center of $\Tbb^{\vee\vee}_0$. 
Denote $q:= \vp(M_0(b_1^j))$ and write
$$
g_{\ol{m}} = \sum_{l \in \Zset} g_{\ol{m},l} y^l, \quad 
g_{\ol{m}, l} \in \Tbb_0^{\vee\vee}.
$$
Eq. \eqref{divv2} is equivalent to 
$$
h(M_0(-b_1^i)) \mid_l g_{\ol{m},l} h'(q^lM_0(b_1^j))
\quad \mbox{in} \quad \Tbb_0^{\vee\vee}, \quad 
\forall l \in \Zset.
$$
Since $\wt{B}_1$ has full rank,
the vectors $-b_1^i, b_1^j \in \sum_{k \neq i,j} \Zset e_k$ are $\Zset$-linearly independent.
\leref{prim-tor} implies that $h(M_0(-b_1^i))$ and $h'(q^lM_0(b_1^j))$ 
are relatively prime in the center of $\Tbb^{\vee\vee}_0$, $\forall l \in \Zset$.

It follows from \eqref{Z} that $\Tbb^{\vee\vee}_0$ is a 
free module over its center. Applying this and the 
relative primeness of $h(M_0(-b_1^i))$ and $h'(q^lM_0(b_1^j))$ 
in the center of $\Tbb^{\vee\vee}_0$, for all $l \in \Zset$ gives
$$
h(M_0(-b_1^i)) \mid_l g_{\ol{m},l}
\quad \mbox{in} \quad \Tbb_0^{\vee\vee}, \quad 
\forall l \in \Zset.
$$
Hence, $h(M_0(-b_1^i)) \mid_l g_{\ol{m}}$ and thus
$(\sigma^{-1}(P) \ldots \sigma^{-m}(P)) \mid_l g_{\ol{m}}$ 
in $\Tbb_0^{\vee\vee}$ which again contradicts \eqref{contrad}.
This completes the proof of part (c).
\end{proof}

\thref{q-l} follows from \leref{quant-Cater} (a), (c) by induction using the 
argument of Fomin and Zelevinsky in the first paragraph of p.~127 in \cite{FZ-l}.

\chapter{Iterated skew polynomial algebras and noncommutative UFDs}
\label{review}

In this chapter, we gather some facts concerning iterated skew polynomial rings (Ore extensions) 
and noncommutative unique factorization domains 
which will be used in the paper.

\section{Equivariant noncommutative unique factorization domains}
\label{2.1}

Recall that a {\em{prime element}} \index{prime element} of a domain $R$ is any nonzero normal element $p\in R$ (\emph{normality} meaning that $Rp = pR$) \index{normal element}
such that $Rp$ is a completely prime ideal, i.e., $R/Rp$ is a domain.
Assume that in addition $R$ is a $\KK$-algebra and $\HH$ a group acting on $R$ 
by $\KK$-algebra automorphisms. An {\em{$\HH$-prime ideal}} \index{Hprimeideal@$\HH$-prime ideal} of $R$ is any proper 
$\HH$-stable ideal $P$ of $R$ such that $(IJ\subseteq P \implies I\subseteq P$ 
or $J\subseteq P)$ for all $\HH$-stable ideals $I$ and $J$ of $R$. 
One says that $R$ is an {\em{$\HH$-UFD}} \index{HUFD@$\HH$-UFD} if each nonzero $\HH$-prime ideal 
of $R$ contains a prime $\HH$-eigenvector. This is an equivariant 
version of Chatters' notion \cite{Cha} of noncommutative unique 
factorization domain given in \cite[Definition 2.7]{LLR}.

The following fact is an equivariant version of results of 
Chatters and Jordan \cite[Proposition 2.1]{Cha},
\cite[p. 24]{ChJo}, see \cite[Proposition 2.2]{GY} and
\cite[Proposition 6.18 (ii)]{Y-sqg}. 

\bpr{factorHUFD}
Let $R$ be a noetherian $\HH$-UFD. Every normal $\HH$-eigenvector in $R$ 
is either a unit or a product of prime $\HH$-eigenvectors. The factors are unique 
up to reordering and taking associates.
\epr

\section{CGL extensions}
\label{2.2}
\begin{convention}

\label{skewpoly}
We use the standard notation $S[x;\sigma,\delta]$ \index{Sxsigmadelta@$S[x;\sigma,\delta]$} for a \emph{skew polynomial ring}, \index{skew polynomial ring} or \emph{Ore extension}; \index{Ore extension} it denotes a ring generated by a unital subring $S$ and an element $x$ satisfying $xs = \sigma(s) x + \delta(s)$ for all $s \in S$, where $\sigma$ is a ring endomorphism of $S$ and $\delta$ is a {\rm(}left\/{\rm)} $\sigma$-derivation of $S$. The ring $S[x;\sigma,\delta]$ is a free left $S$-module, with the nonnegative powers of $x$ forming a basis. 

For every Ore extension  $S[x;\sigma,\delta]$ appearing in this paper, $S$ is a $\KK$-algebra, $\sigma$ is a $\KK$-linear automorphism of $S$, and $\delta$ is a $\KK$-linear $\sigma$-derivation. As a result,  $S[x;\sigma,\delta]$ is a $\KK$-algebra.
\end{convention}

For convenient reference, we reiterate Definition A:

\bde{CGL} An iterated skew polynomial extension 
\begin{equation} 
\label{itOre}
R := \KK[x_1][x_2; \sigma_2, \delta_2] \cdots [x_N; \sigma_N, \delta_N]
\end{equation}
is called a \emph{CGL extension} \index{CGL extension} 
\cite[Definition 3.1]{LLR} if it is equipped with a rational action of a $\KK$-torus $\HH$ 
\index{H@$\HH$}
by $\KK$-algebra automorphisms satisfying the following conditions:
\begin{enumerate}
\item[(i)] The elements $x_1, \ldots, x_N$ are $\HH$-eigenvectors.
\item[(ii)] For every $k \in [2,N]$, $\delta_k$ is a locally nilpotent 
$\sigma_k$-derivation of the algebra $R_{k-1} = \KK[x_1][x_2; \sigma_2, \delta_2] \cdots [x_{k-1}; \sigma_{k-1}, \delta_{k-1}]$.  \index{Rk@$R_k$}
\item[(iii)] For every $k \in [1,N]$, there exists $h_k \in \HH$  \index{hk@$h_k$}  such that 
$\sigma_k = (h_k \cdot)|_{R_{k-1}}$ and the $h_k$-eigenvalue of $x_k$, to be denoted by $\la_k$,  \index{lambdak@$\la_k$}  is not a root of unity.
\end{enumerate}
\index{lambdak@$\la_k$}

Conditions (i) and (iii) imply that 
$$
\sigma_k(x_j) = \la_{kj} x_j \; \; \mbox{for some} \; \la_{kj} \in \kx, \; \; \forall 1 \le j < k \le N.
$$
We then set $\la_{kk} :=1$ and $\la_{jk} := \la_{kj}^{-1}$ for $j< k$.
This gives rise to a multiplicatively skew-symmetric 
matrix $\lab := (\la_{kj}) \in M_N(\kx)$   \index{lambda@$\lab$}  and a corresponding skew-symmetric bicharacter $\Om_\lab$,  \index{Omegalambda@$\Om_\lab$}  recall \eqref{La}, \eqref{La2}.

The CGL extension $R$ is called {\em{torsionfree}} \index{torsionfree CGL extension} if the subgroup $\langle \la_{kj} \mid k,j \in [1,N] \rangle$
of $\kx$  is torsionfree.
Define the \emph{length} of $R$ to be $N$ and the  \emph{rank} of $R$ by 
\index{length of a CGL extension} \index{rank of a CGL extension}
\begin{equation}
\label{rk}
\rk(R) := \{ k \in [1,N] \mid \delta_k = 0 \} \in \Zset_{ > 0}
\end{equation}
(cf.~\cite[Eq. (4.3)]{GY}).  \index{rkR@$\rk(R)$}
Denote the character group of the torus $\HH$ by $\xh$. \index{XH@$\xh$} Up through Chapter \ref{main}, we view $\xh$ as a multiplicative group, with identity $1$.  The action of 
$\HH$ on $R$ gives rise to an $\xh$-grading of $R$. The $\HH$-eigenvectors 
are precisely the nonzero homogeneous elements with respect to this grading. 
The $\HH$-eigenvalue of a nonzero homogeneous element $u \in R$ will be denoted by $\chi_u$. \index{chiu@$\chi_u$} In other words, $\chi_u = \xh\mbox{-deg}(u)$ in terms of the $\xh$-grading.
\ede 

\bex{OqMmn}
For positive integers $m$ and $n$ and a non-root of unity $q \in \kx$, let $\OqM$ \index{OqMmnK@$\OqM$} be the standard single parameter quantized coordinate ring of the matrix variety $M_{m,n}(\KK)$. \index{quantum matrix algebra} This is a $\KK$-algebra with generators $t_{ij}$ for $i \in [1,m]$ and $j \in [1,n]$, and relations
\begin{align*}
t_{ij} t_{kj} &= q t_{kj} t_{ij},  &&\quad \text{for} \; i < k,  \\
t_{ij} t_{il} &= q t_{il} t_{ij},  &&\quad \text{for} \; j < l,  \\
t_{ij} t_{kl} &=  t_{kl} t_{ij},  &&\quad \text{for} \; i < k, \; j > l,  \\
t_{ij} t_{kl} - t_{kl} t_{ij} &= (q-q^{-1}) t_{il} t_{kj},  &&\quad \text{for} \; i < k, \; j < l.
\end{align*}
It is well known that $\OqM$ is an iterated skew polynomial extension
$$
\OqM = \KK[x_1][x_2; \sigma_2, \delta_2] \cdots [x_N; \sigma_N, \delta_N],
$$
where $N = mn$ and $x_{(r-1)n+c} = t_{rc}$ for $r \in [1,m]$ and $c \in [1,n]$. It is easy to determine the elements $\sigma_k(x_j)$ and $\delta_k(x_j)$ for $N \ge k > j \ge 1$ and then to see that the skew derivations $\delta_k$ are locally nilpotent. For later reference, we note that the scalars $\la_{kj}$ are given by the following powers of $q$:
\begin{equation}
\label{OqMlakj}
\la_{(r-1)n+c, (r'-1)n+c'} = \begin{cases}
q^{\, \sign(r'-r)}, &\quad \text{if} \; c = c',  \\
q^{\, \sign(c'-c)}, &\quad \text{if} \; r= r',  \\
1,  &\quad \text{otherwise}, \; 
\end{cases}
\; \; \begin{matrix}
\forall r,r' \in [1,m],  \\ 
c,c' \in [1,n].
\end{matrix}
\end{equation}

There is a rational action of the torus $\HH := (\kx)^{m+n}$ on $\OqM$ by $\KK$-algebra automorphisms such that
$$
(\xi_1, \dots, \xi_{m+n}) \cdot t_{rc} = \xi_r \xi_{m+c}^{-1} t_{rc}
$$
for all $(\xi_1, \dots, \xi_{m+n}) \in \HH$, $r \in [1,m]$, $c \in [1,n]$. Define
$$
h_{rc} := (1, \dots, 1, q^{-1}, 1, \dots, 1, q, 1, \dots, 1) \in \HH,
$$
where the entries $q^{-1}$ and $q$ reside in positions $r$ and $m+c$, respectively. Then $\sigma_{(r-1)n+c} = (h_{rc} \cdot)$ and $h_{rc} \cdot t_{rc} = q^{-2} t_{rc}$. Thus, $\OqM$ is a torsionfree CGL extension of length $N$. (This is well known; see, e.g., \cite[Corollary 3.8]{LLR}.)

The algebra $\OqM$ and the accompanying $\HH$-action are, of course, defined for any nonzero scalar $q$. The resulting algebra is an iterated skew polynomial ring which satisfies all but one of the conditions for a CGL extension. However, if $m,n \ge 2$ and $q$ is a root of unity, there is no element $h \in \HH$ such that $\sigma_{n+2} = (h \cdot)$ and $h \cdot x_{n+2} = \la x_{n+2}$ with $\la$ not a root of unity.
\eex

The algebras of quantum matrices are special cases of the family of quantum Schubert cell algebras.
The latter algebras are treated in detail in Chapters \ref{q-gr} and \ref{q-Schu}. 
On the other hand, the quantum Schubert 
cell algebras are special members in the class of CGL extensions. They have the property that 
their gradings by the character lattices of the torus $\HH$ can be specialized to 
$\Zset_{\geq 0}$-gradings that are connected, meaning that the degree $0$ component of the algebra is just $\KK\cdot1$. For example, in the case of quantum matrices we can 
define a connected $\Zset_{\geq 0}$-grading with $\deg t_{ij} =1$ for all $i,j$.
The next example contains a CGL extension that does not have this property.

\bex{non-N-grad} Equip the quantized Weyl algebra \index{quantized Weyl algebra}
\[
R = A_1^q(\KK) := \frac{\KK \lcor x_1, x_2 \rcor}{\langle x_2 x_1 - q x_1 x_2 -1 \rangle}
\]
\index{Aq1K@$A_1^q(\KK)$}
with the action of $\HH = \kx$ by $\KK$-algebra automorphisms such that $a \cdot x_1 = a x_1$ and $a \cdot x_2 = a^{-1} x_2$ for $a \in \HH$.  It is straightforward to verify (and well known)
that as long as $q$ is not a root of unity, $R$ is a CGL extension $\KK[x_1] [x_2; \sigma_2, \delta_2]$ with $h_1 = h_2 := q$.

The grading of $R$ by $\xh$ cannot be specialized to any 
nontrivial $\Znn$-grading. Indeed, if $\deg x_1 = n$, then $\deg x_2 = -n$, hence, 
$n=0$.
\eex

By \cite[Proposition 3.2, Theorem 3.7]{LLR}, every CGL extension is an 
$\HH$-UFD. A recursive description of the sets of homogeneous prime elements 
of the intermediate algebras $R_k$ of a CGL extension $R$ was obtained in  
\cite{GY}. It is expressed using the following functions.

For a function $\eta : [1,N] \to \Zset$, we 
denote the predecessor and successor functions
for its level sets by
\index{predecessor function} \index{successor function}
$$
p = p_\eta : [1,N] \to [1,N] \sqcup \{ - \infty \}, \qquad
s = s_\eta : [1,N] \to [1,N] \sqcup \{ + \infty \}.
$$
\index{peta@$p_\eta$} \index{seta@$s_\eta$}
They are given by
\begin{equation}
\label{p}
p(k) = 
\begin{cases}
\max \{ j <k \mid \eta(j) = \eta(k) \}, 
&\mbox{if $\exists j < k$ such that $\eta(j) = \eta(k)$}, 
\\
- \infty, \; & \mbox{otherwise} 
\end{cases}
\end{equation}
and
\begin{equation}
\label{s}
s(k) = 
\begin{cases}
\min \{ j > k \mid \eta(j) = \eta(k) \}, 
&\mbox{if $\exists j > k$ such that $\eta(j) = \eta(k)$}, 
\\
+ \infty, & \mbox{otherwise}. 
\end{cases}
\end{equation}

\bth{CGL} \cite[Theorem 4.3]{GY} Let $R$ be a CGL extension of length $N$ and rank $\rk(R)$ as 
in \eqref{itOre}. There exist a function $\eta : [1,N] \to \Zset$  \index{eta@$\eta$}
whose range has cardinality $\rk(R)$ and elements
$$
c_k \in R_{k-1} \; \; \mbox{for all} \; \; k \in [2,N] \; \; 
\mbox{with} \; \; p(k) \neq - \infty
$$
such that the elements $y_1, \ldots, y_N \in R$, recursively defined by  \index{yk@$y_k$}
\begin{equation}
\label{y}
y_k := 
\begin{cases}
y_{p(k)} x_k - c_k, &\mbox{if} \; \;  p(k) \neq - \infty \\
x_k, & \mbox{if} \; \; p(k) = - \infty,  
\end{cases}
\end{equation}
are homogeneous and have the property that for every $k \in [1,N]$,
\begin{equation}
\label{prime-elem}
\{y_j \mid j \in [1,k] , \, s(j) > k \}
\end{equation}
is a list of the homogeneous prime elements of $R_k$ up to scalar multiples.

The elements $y_1, \ldots, y_N \in R$ with these properties are unique.
The function $\eta$ satisfying the above 
conditions is not unique, but the partition of $[1,N]$ into a disjoint 
union of the level sets of $\eta$ is uniquely determined by $R$, as are the predecessor and successor functions $p$ and $s$.
The function $p$ has the property that $p(k) = - \infty$
if and only if $\delta_k =0$.
\eth

The uniqueness of the level sets of $\eta$ was not stated in \cite[Theorem 4.3]{GY}, 
but it follows at once from \cite[Theorem 4.2]{GY}. This uniqueness immediately implies the uniqueness of $p$ and $s$. Although one can 
state the theorem using a partition of $[1,N]$ into $\rk(R)$ subsets, 
the use of a function will be better suited for combinatorial purposes.
In the setting of the theorem, the rank of $R$ is also given by
\begin{equation}
\label{rankR}
\rk(R) = |\{ j \in [1,N] \mid s(j) > N \}|
\end{equation}
\cite[Eq. (4.3)]{GY}.

\bex{OqMmn2}
Let $R = \OqM$ be the CGL extension in \exref{OqMmn}. For any subsets $I \subseteq [1,m]$ and $J \subseteq [1,n]$ of the same cardinality $d$, let $\De_{I,J}$ denote the $d\times d$ quantum minor with row index set $I$ and column index set $J$. \index{quantum minor} \index{DeltaIJ@$\De_{I,J}$} Namely, if $I = \{ i_1 < \cdots < i_d \}$ and $J = \{ j_1 < \cdots < j_d \}$, then
$$
\De_{I,J} := \sum_{\tau \in S_d} (-q)^{\ell(\tau)} t_{i_1, j_{\tau(1)}} t_{i_2, j_{\tau(2)}} \cdots t_{i_d, j_{\tau(d)}} ,
$$
where $\ell(\tau)$ \index{ltau@$\ell(\tau)$} is the length of $\tau$ as a minimal length product of simple transpositions. The homogeneous prime elements $y_k$ from \thref{CGL} are ``solid" quantum minors of the following form: \index{solid quantum minor}
$$
y_{(r-1)n+c} = \De_{ [r-\min(r,c)+1,r], [c-\min(r,c)+1, c]} , \; \; \forall r \in [1,m], \; c \in [1,n].
$$
The function $\eta: [1,N] \rightarrow \Zset$ can be chosen as
$$
\eta((r-1)n+c) := c-r , \; \; \forall r \in [1,m], \; c \in [1,n].
$$
It is easily checked that for the CGL extension presentation of $R$ in \exref{OqMmn}, we have $\delta_k = 0$ if and only if $k \in [1,n]$ or $k= (r-1)n+1$ for some $r \in [2,m]$. Hence, $\rk(R) = m+n-1$.
\eex

There is a universal maximal choice for the torus $\HH$ in \deref{CGL} 
given by the following theorem. Consider the rational action of the torus $(\kx)^N$ 
on $R$ by invertible linear transformations given by
$$
(\ga_1, \dots, \ga_N) \cdot (x_1^{m_1} \cdots x_N^{m_N}) = 
\ga_1^{m_1} \cdots \ga_N^{m_N} x_1^{m_1} \cdots x_N^{m_N}.
$$
Denote by $\HH_\max(R)$ \index{HmaxR@$\HH_\max(R)$} the closed subgroup of $(\kx)^N$ consisting 
of those $\psi \in (\kx)^N$ whose actions give $\KK$-algebra automorphisms of $R$.

\bth{univ-torus} \cite[Theorems 5.3, 5.5]{GY} For every CGL extension 
$R$, $\HH_\max(R)$ is a $\KK$-torus of rank $\rk(R)$, and the pair $(R, \HH_\max(R))$ 
is a CGL extension for the Ore extension presentation \eqref{itOre}. 
\eth

\thref{univ-torus} implies that the torus $\HH_\max(R)$ is universal 
(and maximal) in the sense 
that, if $\HH$ is any torus
acting rationally on $R$ by algebra automorphisms such that 
$(R, \HH)$ is a CGL extension for the presentation \eqref{itOre},
then the action of $\HH$ on $R$ factors through an algebraic 
group morphism $\HH \to \HH_\max(R)$. The torus $\HH_\max(R)$ is 
explicitly described in \cite[Theorem 5.5]{GY} on the basis 
of the sequence $y_1, \ldots, y_N$ from \thref{CGL}.

\bex{OqMmn3} In the case of $R= \OqM$ (\exref{OqMmn}), $\Hmax(R)$ can be computed from \cite[Theorem 5.5]{GY}. It can be expressed as follows:
$$
\Hmax(R) = \{ \psi \in (\kx)^{mn} \mid \psi_{(r-1)n+c} = \psi_1^{-1} \psi_c \psi_{(r-1)n+1} \;\; \forall r \in [2,m], \; c \in [2,n] \}.
$$
This torus is related to the torus $\HH = (\kx)^{m+n}$ of \exref{OqMmn} via a surjective morphism of algebraic groups $\pi : \HH \rightarrow \Hmax(R)$, given by
$$
\pi(\xi)_{(r-1)n+c} = \xi_r \xi_{m+c}^{-1} , \; \; \forall r \in [1,m], \; c \in [1,n],
$$
and $\pi$ transports the $\HH$-action on $R$ to the $\Hmax(R)$-action.
\eex

The next result provides a constructive method for finding the 
sequence of prime elements $y_1, \ldots, y_N$ for a given 
CGL extension $R$. Note that the elements $y'_k$ are not a priori assumed to be prime, only normal.

\bpr{constr-y} Let $R$ be a CGL extension of length $N$. Assume that 
$y'_1, \ldots, y'_N$ and $c'_1, \ldots, c'_N$ are two sequences of 
elements of $R$ such that 
\begin{enumerate}
\item[(i)] $y'_1, \ldots, y'_N$ are homogeneous normal elements of $R_1, \ldots, R_N$, respectively.
\item[(ii)] $c'_k \in R_{k-1}$, $\forall k \in [1,N]$.
\item[(iii)] For every $k \in [1,N]$, either $y'_k = x_k - c'_k$ or
there exists $j \in [1,k-1]$ such that $y'_k = y'_j x_k - c'_k$.
\item[(iv)] If $p(k) = - \infty$, then the first equality in {\rm(iii)} holds.
\end{enumerate}

Then $y'_1, \ldots, y'_N$ is precisely the sequence of homogeneous prime 
elements from Theorem {\rm\ref{tCGL}}, and the function $p$ 
satisfies $p(k) := j$ if the second equality in {\rm(iii)}
holds and $p(k) := - \infty$ otherwise. 
\epr

\begin{proof} The given assumptions imply that $c'_1 \in \KK$ and $y'_1 = x_1 - c'_1$. Since $y'_1$ is homogeneous, we must have $y'_1 = x_1 = y_1$.

Now let $k \in [2,N]$. We will prove that, if $y'_i = y_i$, $\forall i \in [1,k-1]$, then 
$y'_k = y_k$. This implies the first statement of the proposition by induction.
By \prref{factorHUFD} and \thref{CGL},
$$
y'_k = \xi \prod \{ y_i^{m_i} \mid i \in [1,k], \, s(i)>k  \} 
$$
for some $\xi \in \kx$ and $m_i \in \Zset_{\geq 0}$, where the factors $y_i^{m_i}$ are taken in ascending order with respect to $i$. Comparing the coefficients 
of $x_k$ and using the form of $y_k$ from \thref{CGL},
we obtain that $m_k = 1$. One of the following two situations holds:
\begin{enumerate}
\item[(a)] $m_i=0$, $\forall i \in [1,k-1]$ with $s(i) > k$.
\item[(b)] $m_{i_0} > 0$ for some $i_0 \in [1,k-1]$ with $s(i_0)>k$. 
\end{enumerate}

First we rule out (b). The $x_k$-coefficient of $y'_k$ is either $1$ or $y'_j$, hence either a unit or a prime element, by induction.  If (b) holds, then $y_k= x_k$ and  
$$
y'_k = \xi y_{i_0} x_k ,
$$
which contradicts the condition (iv) because $y_k = x_k$ only occurs when $p(k) = -\infty$.

In the situation (a), we have $y'_k = \xi y_k$. So, either 
\begin{align*}
&x_k - c'_k = \xi x_k, \; \; \mbox{or}
\\
&y'_j x_k - c'_k = \xi(y_{p(k)} x_k - c_k).
\end{align*}
In the first case, we have $\xi = 1$, $c'_k=0$ and $y'_k = y_k$. In the second case, 
using the fact that $y'_i= y_i$, $\forall i \in [1,k-1]$, we obtain 
$\xi =1$, $j= p(k)$, $c'_k = c_k$ and $y'_k= y_k$. 
This argument also proves the second 
statement of the proposition.
\end{proof}

Although this will not be used later, we also note that 
the conclusion of the proposition implies
$$
c'_k = 
\begin{cases}
c_k, & \mbox{if} \; \; p(k) \neq - \infty
\\ 
0, & \mbox{if} \; \; p(k) = - \infty.
\end{cases}
$$

The following fact from \cite{GY} will be extensively used in the paper to construct quantum 
clusters.

\bpr{CGLcluster} \cite[Eq. (4.17) and Theorem 4.6]{GY} For each CGL extension $R$, the elements $y_k$ of Theorem {\rm \ref{tCGL}} quasi-commute: there are scalars $q_{kj} \in \kx$, given in {\rm\eqref{q}} below, such that  \index{quasi-commuting elements}
\begin{equation}
\label{ycomm}
y_k y_j = q_{kj} y_j y_k, 
\quad \forall j, k \in [1,N].
\end{equation}
The quantum torus $\Tbb_\qb$ embeds in $\Fract(R)$ via the $\KK$-algebra homomorphism $\vp : \Tbb_\qb \hra \Fract(R)$ given 
by $\vp(Y_i)= y_i$, $\forall i \in [1,N]$, and this embedding gives rise to inclusions 
$$
\vp(\Abb_\qb) \subseteq R \subset \vp(\Tbb_\qb) \subset \Fract(R),
$$
recall {\rm\eqref{q-aff}}.
\epr

\section{Symmetric CGL extensions}
\label{2.5}
The CGL extensions for which we will establish quantum cluster algebra structures are those with suitably many different CGL extension presentations, in which the variables $x_1,\dots,x_N$ are permuted in various ways. A symmetry condition, which we introduce now, is sufficient to guarantee this.

Given an iterated Ore extension $R$ as in \eqref{itOre}, 
for $j,k \in [1,N]$  
denote by $R_{[j,k]}$ the unital subalgebra of
$R$ generated by $\{ x_i \mid j \le i \le k \}$. So, $R_{[j,k]} = \KK$ if $j \nleq k$. \index{Rjk@$R_{[j,k]}$}

\bde{symmetric} We call a CGL extension $R$ of length $N$ as in 
\deref{CGL} {\em{symmetric}} if the following two conditions hold: \index{symmetric CGL extension}
\begin{enumerate}
\item[(i)] For all $1 \leq j < k \leq N$,
$$
\delta_k(x_j) \in R_{[j+1, k-1]}.
$$
\item[(ii)] For all $j \in [1,N]$, there exists $h^\sy_j \in \HH$ 
such that 
$$
h^\sy_j \cdot x_k = \la_{kj}^{-1} x_k = \la_{jk} x_k, \; \; \forall 
k \in [j+1, N]
$$
and $h^\sy_j \cdot x_j = \la^\sy_j x_j$ for some $\la^\sy_j \in \kx$ which is not a root of unity.
\end{enumerate}
\index{lambdastarj@$\la^\sy_j$}
\ede

For such an algebra $R$, set 
$$
\sigma^\sy_j := (h^\sy_j \cdot) \in \Aut (R), \; \; \forall j \in [1,N-1].
$$
Then for all $j \in [1,N-1]$, 
the inner $\sigma^\sy_j$-derivation on $R$ 
given by $a \mapsto x_j a - \sigma^\sy_j(a) x_j$ restricts to a $\sigma^\sy_j$-derivation 
$\delta^\sy_j$ of $R_{[j+1, N]}$. It is given by
$$
\delta^\sy_j(x_k): = x_j x_k - \la_{jk} x_k x_j = - \la_{jk} \delta_k(x_j), \; \; 
\forall k \in [j+1, N].
$$ 
For all $1 \leq j < k \leq N$, $\sigma_k$ and $\delta_k$ preserve $R_{[j,k-1]}$
and $\sigma^\sy_j$ and $\delta^\sy_j$ preserve $R_{[j+1,k]}$. This gives 
rise to the skew polynomial extensions
\begin{equation}
\label{step1}
R_{[j,k]} = R_{[j,k-1]}[x_k; \sigma_k, \delta_k] \; \; 
\mbox{and} \; \; 
R_{[j,k]} = R_{[j+1,k]} [ x_j; \sigma^\sy_j, \delta^\sy_j].
\end{equation}
In particular, it follows that $R$ has a skew polynomial presentation with the variables $x_k$ in descending order:
$$R = \KK[x_N] [x_{N-1}; \sigma^*_{N-1}, \delta^*_{N-1}] \cdots [x_1; \sigma^*_1, \delta^*_1].$$
This is the reason for the name ``symmetric".

\bex{OqMmn4}
For example, the CGL extension $\OqM$ of \exref{OqMmn} is symmetric. It is clear that condition (i) of \deref{symmetric} holds in this example. Condition (ii) can be verified for the elements
$$
h^*_{(r-1)n+c} := (1, \dots, 1, q, 1, \dots, 1, q^{-1}, 1, \dots, 1) \in \HH,
$$
where the entries $q$ and $q^{-1}$ reside in positions $r$ and $m+c$, respectively. Then $h^* _j \cdot x_j = q^2 x_j$ for all $j$, and so all $\la^*_j = q^2$ here.
\eex

Denote the following subset of the symmetric group $S_N$:
\begin{multline}
\label{tau}
\Xi_N := \{ \tau \in S_N \mid 
\tau(k) = \max \, \tau( [1,k-1]) +1 \; \;
\mbox{or} 
\\
\tau(k) = \min \, \tau( [1,k-1]) - 1, 
\; \; \forall k \in [2,N] \}.
\end{multline}
\index{XiN@$\Xi_N$}
In other words, $\Xi_N$ consists of those $\tau \in S_N$ 
such that $\tau([1,k])$ is an interval for all $k \in [2,N]$. 
For each $\tau \in \Xi_N$, we have the iterated 
Ore extension presentation
\begin{equation}
\label{tauOre}
R = \KK [x_{\tau(1)}] [x_{\tau(2)}; \sigma''_{\tau(2)}, \delta''_{\tau(2)}] 
\cdots [x_{\tau(N)}; \sigma''_{\tau(N)}, \delta''_{\tau(N)}],
\end{equation}
where $\sigma''_{\tau(k)} := \sigma_{\tau(k)}$ and 
$\delta''_{\tau(k)} := \delta_{\tau(k)}$ if 
$\tau(k) = \max \, \tau( [1,k-1]) +1$, while 
$\sigma''_{\tau(k)} := \sigma^\sy_{\tau(k)}$ and 
$\delta''_{\tau(k)} := \delta^\sy_{\tau(k)}$ if 
$\tau(k) = \min \, \tau( [1,k-1]) -1$.

\bpr{tauOre} 
\cite[Remark 6.5]{GY} 
For every symmetric CGL extension $R$ and $\tau \in \Xi_N$,
the iterated Ore extension presentation \eqref{tauOre} of $R$ 
is a CGL extension presentation for the 
same choice of $\KK$-torus $\HH$, and the associated elements 
$h''_{\tau(1)}, \ldots, h''_{\tau(N)} \in \HH$ required by Definition {\rm\ref{dCGL}(iv)} 
are given by $h''_{\tau(k)} = h_{\tau(k)}$ if $\tau(k) = \max \, \tau( [1,k-1]) +1$  
and $h''_{\tau(k)} = h^\sy_{\tau(k)}$ if $\tau(k) = \min \, \tau( [1,k-1]) -1$.
\epr

\thref{univ-torus} implies that the rank of a symmetric CGL extension does
not depend on the choice of CGL extension presentation \eqref{tauOre}.

When describing permutations $\tau \in S_N$ as functions, we will use the one-line notation, \index{one-line notation for permutations}
\begin{equation}
\label{one-line}
\tau = [\tau(1),\tau(2),\dots,\tau(N)] := \begin{bmatrix} 1 &2 &\cdots &N\\ \tau(1) &\tau(2) &\cdots &\tau(N) \end{bmatrix} .
\end{equation}
A special role is played by the longest element of $S_N$,  
\index{wcirc@$w_\circ$}
\begin{equation}
\label{tau-ci}
w_\circ := [N, N-1, \ldots, 1].
\end{equation}
The corresponding CGL extension presentation from \prref{tauOre} is symmetric, 
while the ones for the other elements of $\Xi_N$ do not possess 
this property in general.

\section{Further CGL details}
\label{detailCGL}
In order to work closely with CGL extensions, we need some further notation and relations, which we give in this section. Throughout, $R$ is a CGL extension as in \eqref{itOre}; the elements $y_k$, $k \in [1,N]$ and the functions $p$ and $s$ are as in \thref{CGL}.

As is easily checked (e.g., see the proof of \cite[Lemma 2.6]{LLR}),
\begin{equation}
\label{hdelta}
(h\cdot) \delta_k = \chi_{x_k}(h) \delta_k (h\cdot), \; \; \forall h \in \HH, \; k \in [2,N].
\end{equation}

We will use the following convention:
\begin{equation}
\label{conven}
\begin{aligned}
\mbox{All} \; &\mbox{products} 
\; \; 
\prod_{j \in P} y_j^{m_j} \; \; 
\mbox{with} \; \; P \subseteq [1,N] \; \; \mbox{and} \; \;
m_j \in \Znn \; \; 
\mbox{are taken}  \\
 &\mbox{with the terms in increasing order with respect to} \; \;j.  
\end{aligned}
\end{equation}
Set also 
\begin{equation}
\label{inf}
y_{-\infty} := 1.
\end{equation}

Define the order functions $O_\pm : [1,N] \to \Znn$ by \index{order function} \index{Oplus@$O_+$} \index{Ominus@$O_-$}
\begin{equation}  \label{Opm}
\begin{aligned}
O_+(k) &:= \max \{ m \in \Znn \mid s^m(k) \neq + \infty \}  \\
O_-(k) &:= \max \{ m \in \Znn \mid p^m(k) \neq - \infty \},
\end{aligned}
\end{equation}
where as usual $p^0 = s^0 = \id$. For $k \in [1,N]$, define 
\begin{equation}
\label{ol-e}
\ol{e}_k := \sum_{m=0}^{O_-(k)} e_{p^m(k)} \in \Zset^N.
\end{equation}
\index{ekbar@$\ol{e}_k$}

The algebra $R$ has the $\KK$-basis 
$$
\{ x^f := x_1^{m_1} \cdots x_N^{m_N} \mid 
f = (m_1, \ldots, m_N)^T \in \Znn^N \}.
$$
Denote by $\prec$ \index{zzz@$\prec$} the 
reverse lexicographic order on $\Znn^N$: 
\begin{multline*}
(m'_1, \ldots, m'_N )^T \prec (m_1, \ldots, m_N)^T \;\; \text{iff} \;\; \exists j \in [1,N]\;\; \text{such that} \\
  m'_j < m_j \;\; \text{and} \;\; m'_k = m_k, \; \forall k \in [j+1, N] .
\end{multline*}
We say that 
$b \in R \setminus \{ 0 \}$ 
has \emph{leading term} \index{leading term} $\xi x^f$ where $\xi \in \kx$ and $f \in \Znn^N$ 
if 
$$
b = \xi x^f + \sum_{g \in \Znn^N,\; g \prec f} \xi_g x^g
$$
for some $\xi_g \in \KK$.
Set $\lt(b) := \xi x^f$. \index{ltb@$\lt(b)$} Then  
\begin{multline}
\label{prod-lt}
\lt( x^f x^{f'} ) = 
\biggl( \prod_{k>j} \la_{kj}^{m_k m'_j} \biggr) x^{f + f'}, \\
\forall \; f= (m_1, \ldots, m_N)^T, \;\; f'=(m'_1, \ldots, m'_N)^T \in \Znn^N. 
\end{multline}
We have \cite[Eq. (4.13)]{GY}
\begin{equation}
\label{tty}
\lt(y_k) = x^{\ol{e}_k}, \; \; \forall k \in [1,N],
\end{equation}
cf. \eqref{ol-e}.

For $k,j \in [1,N]$, set
\begin{align}
\al_{kj} &:= \Om_\lab (e_k, \ol{e}_j) = \prod_{m=0}^{O_-(j)} \la_{k, p^m(j)}  \in \kx  \label{al} \\
q_{kj} &:= \Om_\lab (\ol{e}_k, \ol{e}_j) = \prod_{m=0}^{O_-(k)} \prod_{l=0}^{O_-(j)}
\la_{p^m(k), p^l(j)}  = \prod_{m=0}^{O_-(k)} \al_{p^m(k), j} \in \kx,  \label{q}
\end{align}
\index{alphakj@$\al_{kj}$} \index{qkj@$q_{kj}$}
recall \eqref{La} and \eqref{ol-e}. 
Since $\lab =(\la_{kj})$ is a multiplicatively skew-symmetric matrix, so is
$\qb := (q_{kj}) \in M_N(\kx)$. In addition to \eqref{ycomm}, we have \cite[Corollary 4.8;  p.~21;  Proposition 4.7(b)]{GY}
\begin{align}
y_j x_k &= \al_{kj}^{-1} x_k y_j, \; \; \forall j,k \in [1,N] \; \; \mbox{such that} \; \; s(j) > k  \label{yxcomm}  \\
\sigma_k(y_j) &= \al_{kj} y_j, \; \; \mbox{for} \; 1 \le j < k \le N  \label{sigkyj}  \\
\delta_k(y_{p(k)}) &= \al_{kp(k)} (\la_k - 1) c_k \ne 0, \; \; \forall k \in [2,N] \; \; \mbox{such that} \; \; p(k) \ne -\infty.  \label{delkypk}
\end{align}
Set also, for $k \in [1,N]$,
\begin{equation}
\label{alqinf}
\al_{k,-\infty} = q_{k,-\infty} := 1.
\end{equation}

For $k \in [0,N]$, denote
\begin{equation}
\label{P(k)}
P(k) := \{ j \in [1,k] \mid s(j) > k \}.
\end{equation} 
\index{Pk@$P(k)$}
Then $\{y_j \mid j \in P(k) \}$ is a list of the homogeneous prime elements 
of $R_k$ up to scalar multiples, and $|P(N)| = \rk(R)$ \eqref{rankR}.

\chapter{One-step mutations in CGL extensions}
\label{mCGL}
In this chapter we obtain a very general way of 
constructing mutations of toric frames in CGL extensions.
The key idea is that, if an algebra $R$ has 
two different CGL extension presentations
obtained by reversing the order in which two adjacent variables 
$x_k$ and $x_{k+1}$ are adjoined, then 
the corresponding sequences of prime elements 
from \thref{CGL} are obtained by a type of 
mutation formula. This is realized in Section \ref{4.1}. 
In Section \ref{4.2}--{4.4} we construct toric frames 
from \thref{CGL} and treat their mutations. 
One problem arises along the way: In the analog of 
the mutation formula \eqref{clust-mut} for the current situation, 
the last term has a nonzero coefficient which 
does not equal one in general. 
For a one-step mutation such a coefficient 
can be always made 1 after rescaling, but for the purposes of constructing quantum cluster 
algebras one needs to be able to synchronize those rescalings to obtain a chain of mutations.
This delicate issue is 
resolved in the next two chapters.

We investigate a CGL extension
\begin{equation}
\label{firstCGL}
R := \KK[x_1][x_2; \sigma_2, \delta_2] \cdots 
[x_k; \sigma_k, \delta_k]
[x_{k+1}; \sigma_{k+1}, \delta_{k+1}] \cdots
[x_N; \sigma_N, \delta_N]
\end{equation}
of length $N$ as in Definition {\rm\ref{dCGL}} such that, for some $k \in [1,N-1]$, $R$ has a second CGL extension presentation of the form 
\begin{multline}
\label{secondCGL}
R := \KK[x_1][x_2; \sigma_2, \delta_2] \cdots 
[x_{k-1}; \sigma_{k-1}, \delta_{k-1}]
[x_{k+1}; \sigma'_k, \delta'_k]
[x_k; \sigma'_{k+1}, \delta'_{k+1}] \\ 
[x_{k+2}; \sigma_{k+2}, \delta_{k+2}]
\cdots
[x_N; \sigma_N, \delta_N].
\end{multline}

\section{A general mutation formula}
\label{4.1}
\ble{CGLtranspose}
Assume that $R$
is a CGL extension as in \eqref{firstCGL}, and that 
$R$ has a second CGL extension presentation of the form \eqref{secondCGL}.
Then $\delta_{k+1}$ and $\delta'_{k+1}$ map $R_{k-1}$ to itself, and
\begin{equation}
\label{sig'de'}
\begin{aligned}
\sigma'_k &= \sigma_{k+1}|_{R_{k-1}}  &\delta'_k &= \delta_{k+1}|_{R_{k-1}}  \\
\sigma'_{k+1}|_{R_{k-1}} &= \sigma_k  &\delta'_{k+1}|_{R_{k-1}} &= \delta_k .
\end{aligned}
\end{equation}
Moreover, $\delta_{k+1}(x_k) \in \KK$ and
\begin{align}
\label{moresig'de'}
\sigma'_{k+1}(x_{k+1}) &= \la_{k,k+1} x_{k+1}  &\delta'_{k+1}(x_{k+1}) &= - \la_{k,k+1} \delta_{k+1}(x_k).
\end{align}
\ele

\begin{proof} Note first that $R_{k-1}$ is stable under $\sigma_{k+1}$ and $\sigma'_{k+1}$.
 For $a \in R_{k-1}$, we have $x_{k+1} a = \sigma'_k(a) x_{k+1} + \delta'_k(a)$ with $\sigma_k'(a), \delta'_k(a) \in R_{k-1}$. Comparing this relation with $x_{k+1} a = \sigma_{k+1}(a) x_{k+1} + \delta_{k+1}(a)$, and using the fact that $1$, $x_{k+1}$ are left linearly independent over $R_k$, we conclude that $\sigma'_k(a) = \sigma_{k+1}(a)$ and $\delta'_k(a) = \delta_{k+1}(a)$. Thus, $R_{k-1}$ is stable under $\delta_{k+1}$, and the first line of \eqref{sig'de'} holds. By symmetry (since we may view \eqref{secondCGL} as the initial CGL extension presentation of $R$ and \eqref{firstCGL} as the second one), $R_{k-1}$ is stable under $\delta'_{k+1}$, and the second line of \eqref{sig'de'} holds.

Now $x_k x_{k+1} = \la'_{k+1,k} x_{k+1} x_k + \delta'_{k+1}(x_{k+1})$, and so we have
$$\la'_{k,k+1} x_k x_{k+1} - \la'_{k,k+1} \delta'_{k+1}(x_{k+1}) = x_{k+1} x_k = \la_{k+1,k} x_k x_{k+1} + \delta_{k+1}(x_k),$$
with $\delta'_{k+1}(x_{k+1}) \in R'_k = \bigoplus_{l=0}^\infty R_{k-1} x_{k+1}^l$ and $\delta_{k+1}(x_k) \in R_k = \bigoplus_{l=0}^\infty R_{k-1} x_k^l$. Moreover,
\begin{equation}
\label{Rk+1basis}
R_{k+1} \; \; \mbox{is a free left} \; \; R_{k-1}\mbox{-module with basis} \; \; \{x_k^{l_k} x_{k+1}^{l_{k+1}} \mid l_k, l_{k+1} \in \Zset_{\geq 0 } \}.
\end{equation}
Hence, we conclude that $\la'_{k,k+1} = \la_{k+1,k}$ and
$$- \la'_{k,k+1} \delta'_{k+1}(x_{k+1}) = \delta_{k+1}(x_k) \in \KK,$$
from which \eqref{moresig'de'} follows.
\end{proof}

\bth{1} Assume that $R$ is a CGL extension of length $N$ and rank $\rk(R)$ as in {\rm\eqref{firstCGL}},
and $k \in [1,N-1]$.
Denote by $y_1, \ldots, y_N$ and $\eta : [1,N] \to \Zset$ 
the sequence and function from Theorem {\rm\ref{tCGL}}. Assume that 
$R$ has a second CGL extension presentation of the form {\rm\eqref{secondCGL}}, and let $y'_1, \ldots, y'_N$ 
be the corresponding sequence from Theorem {\rm\ref{tCGL}}.

{\rm(a)} If $\eta(k) \neq \eta(k+1)$, then $y'_j = y_j$ for $j \neq k, k+1$ 
and $y'_k = y_{k+1}$, $y'_{k+1} = y_k$. 

{\rm(b)} If $\eta(k) = \eta(k+1)$, then 
$$
y_k y'_k - \al_{k p(k)} y_{p(k)} y_{k+1}
$$ 
is a homogeneous normal element of $R_{k-1}$, recall \eqref{inf}, \eqref{al}, \eqref{alqinf}.
It normalizes the elements of $R_{k-1}$ in exactly the same way as $y_{p(k)} y_{k+1}$.
Furthermore,
\begin{equation}
\label{y-equal}
y'_j = 
\begin{cases} 
\la_{k+1,k} y_j, &\mbox{if} 
\; \; j = s^l(k+1) \; \; \mbox{for some} \; \; l \in \Znn 
\\
y_j, &\mbox{if} 
\; \; j < k, \; \mbox{or} \; j > k+1 \; 
\mbox{and} \; j \neq s^l(k+1) \; \forall l \in \Znn.
\end{cases}
\end{equation}

In both cases {\rm(a)} and {\rm(b)}, the function $\eta' : [1,N] \to \Zset$ from
Theorem {\rm\ref{tCGL}} associated to the second presentation can be chosen 
to be $\eta' = \eta (k, k+1)$, 
where the last term denotes a transposition in $S_N$. In particular, 
the ranges of $\eta$ and $\eta'$ coincide 
and the rank of $R$ is the same 
for both CGL extension presentations.
\eth

\begin{proof} Let $R'_j$ be the $j$-th algebra in the 
chain \eqref{secondCGL} for $j \in [0,N]$. Obviously 
$R'_j = R_j$ for $j \neq k, k+1$ and $y'_j = y_j$
for $j \leq k-1$, and we may choose $\eta'(j) = \eta(j)$ for $j \le k-1$. Since $R_{k+1}$ is a free 
left $R_{k-1}$-module with basis \eqref{Rk+1basis},
and $R_k$ and $R'_k$ equal the $R_{k-1}$-submodules with bases 
$\{x_k^{l_k} \mid l_k \in \Zset_{\geq 0 } \}$ 
and 
$\{x_{k+1}^{l_{k+1}} \mid l_{k+1} \in \Zset_{\geq 0 } \}$
respectively, we have
\begin{equation}
\label{Rkk'intersect}
R_k \cap R'_k = R_{k-1}.
\end{equation}
In particular, $y'_k \notin R_k$.

Denote by $L$ the number of homogeneous prime elements of 
$R_{k+1}$ up to taking associates that do not belong to $R_{k-1}$.

(a) The condition $\eta(k) \neq \eta(k+1)$ implies $L =2$ and thus 
$\eta'(k) \neq \eta'(k+1)$. Moreover, $y'_k$ and $y'_{k+1}$ 
should be scalar multiples of either $y_k$ and $y_{k+1}$ or $y_{k+1}$ and $y_k$. 
Since $y'_k \notin R_k$, we must have $y'_k = \xi_{k+1} y_{k+1}$ and $y'_{k+1} = \xi_k y_k$
for some $\xi_k, \xi_{k+1} \in \kx$. Invoking \thref{CGL} 
and looking at leading terms gives $\xi_k = \xi_{k+1} = 1$ and allows us to choose $\eta'(k) = \eta(k+1)$ and $\eta'(k+1) = \eta(k)$. It follows
from $R'_j = R_j$ for $j>k+1$ that $y'_j = y_j$ for such $j$, and that we can choose $\eta'(j) = \eta(j)$ 
for all $j > k+1$. With these choices, $\eta' = \eta (k,k+1)$.

(b) In this case, $L=1$ and $y'_{k+1}$ should be a scalar multiple
of $y_{k+1}$. Furthermore, $\eta'(k)$ and $\eta'(k+1)$ must agree, and they need to equal $\eta'(p(k))$ if $p(k) \ne -\infty$, so we can choose $\eta'(k) = \eta'(k+1) = \eta(k) = \eta(k+1)$. 
By \thref{CGL},
\begin{equation}
\label{ykyk'}
\begin{aligned}
y_k &= \begin{cases} y_{p(k)} x_k - c_k, &\quad\mbox{if} \; \; p(k) \ne -\infty\\  x_k, &\quad\mbox{if} \; \; p(k) = -\infty \end{cases}  &\qquad y_{k+1} &= y_k x_{k+1} - c_{k+1}  \\
y'_k &= \begin{cases} y_{p(k)} x_{k+1} - c'_k, &\mbox{if} \; \; p(k) \ne -\infty\\  x_{k+1}, &\mbox{if} \; \; p(k) = -\infty \end{cases}  &y'_{k+1} &= y'_k x_k - c'_{k+1}
\end{aligned}
\end{equation}
for some $c_k, c'_k \in R_{k-1}$, $c_{k+1} \in R_k$, and $c'_{k+1} \in R'_k$. Write the above elements in terms of the basis \eqref{Rk+1basis}. The coefficients of $x_k x_{k+1}$ in $y_{k+1}$ and $y'_{k+1}$ are $y_{p(k)}$ and $\la_{k+1,k} y_{p(k)}$, recall \eqref{inf}. Since $y'_{k+1}$ is a scalar multiple of $y_{k+1}$, we thus see that
\begin{equation}
\label{y'k+1}
y'_{k+1} = \la_{k+1,k} y_{k+1}.
\end{equation}
Eq. \eqref{y-equal} and the fact that we may choose $\eta'(j) = \eta(j)$ for $j > k+1$ now follow easily from \thref{CGL}. In particular, $\eta' = \eta = \eta (k, k+1)$.

Next, we verify that 
\begin{equation}
\label{in}
y_k y'_k - \al_{k p(k)} y_{p(k)} y_{k+1} \in R_{k-1}.
\end{equation}
By \eqref{delkypk}, $c_{k+1}$ is a nonzero scalar multiple of $\delta_{k+1}(y_k)$.
If $p(k) = -\infty$, then $y_k = x_k$, and so $c_{k+1} \in \kx$ by \leref{CGLtranspose}. From \eqref{ykyk'} and \eqref{inf}, \eqref{alqinf}, we then obtain
\begin{equation}
\label{eltp(k)-inf}
y_k y'_k - \al_{kp(k)} y_{p(k)} y_{k+1} = x_k x_{k+1} - (x_k x_{k+1} - c_{k+1}) = c_{k+1} \in \kx .
\end{equation}
 This verifies \eqref{in} (and also shows that $y_k y'_k - \al_{kp(k)} y_{p(k)} y_{k+1}$ is a homogeneous normal element of $R_{k-1}$) in the present case.

Now assume that $p(k) \ne -\infty$. Invoking \eqref{ycomm} and the observation that $q_{kp(k)} = \al_{kp(k)}$, we have $y_k y_{p(k)} = \al_{kp(k)} y_{p(k)} y_k$. 
First, we obtain
\begin{equation}
\label{theelt}
\begin{aligned}
y_k y'_k - \al_{k p(k)} y_{p(k)} y_{k+1} 
 &= y_k ( y_{p(k)} x_{k+1} - c'_k) - \al_{k p(k)} 
y_{p(k)} (y_k x_{k+1} - c_{k+1}) \\
 &= - y_k c'_k + \al_{k p(k)} y_{p(k)} c_{k+1} \in R_k .
\end{aligned} 
\end{equation}
Using \eqref{y'k+1} together with the fact that $\sigma'_{k+1}(y'_k) = \al'_{k+1,k} y'_k$ \eqref{sigkyj} and the observation that $\al'_{k+1,k} = \la_{k+1,k}^{-1} \al_{kp(k)}$, we obtain
\begin{align*}
y_k y'_k - \al_{kp(k)} y_{p(k)} y_{k+1} &= y_{p(k)} ( \al'_{k+1,k} y'_k x_k + \delta'_{k+1}(y'_k) ) - c_k y'_k  \\
 &\qquad\qquad\qquad\qquad - \al_{kp(k)} y_{p(k)} \la_{k+1,k}^{-1} ( y'_k x_k - c'_{k+1} )  \\
 &= y_{p(k)} \delta'_{k+1}(y'_k) - c_k y'_k + \al_{kp(k)} \la_{k+1,k}^{-1} y_{p(k)} c'_{k+1} \in R'_k .
\end{align*}
This equation, combined with \eqref{theelt} and \eqref{Rkk'intersect}, yields \eqref{in}.

We now use \eqref{theelt} to verify that $y_k y'_k - \al_{kp(k)} y_{p(k)} y_{k+1}$ is homogeneous. Note first that $y_k c'_k$ and $y_{p(k)} c_{k+1}$ are homogeneous. By \eqref{delkypk}, $c'_k$ and $c_{k+1}$ are scalar multiples of $\delta'_k(y_{p(k)}) = \delta_{k+1}(y_{p(k)})$ and $\delta_{k+1}(y_k)$, respectively. Hence, it follows from \eqref{hdelta} that
$$\xh\mbox{-deg}(y_k c'_k) = \chi_{y_k} + \chi_{x_{k+1}} + \chi_{y_{p(k)}} = \xh\mbox{-deg} (y_{p(k)} c_{k+1} ).$$
Thus, $-y_k c'_k + \al_{kp(k)} y_{p(k)} c_{k+1}$ is homogeneous, as desired.

Finally, whether $p(k)= -\infty$ or not, it follows from \eqref{yxcomm} that 
$$
(y_k y'_k) x_j = \be_j x_j (y_k y'_k) \quad \mbox{and} 
\quad (y_{p(k)} y_{k+1}) x_j = \ga_j x_j (y_{p(k)} y_{k+1}), \quad
\forall j \in [1,k-1],
$$
where 
\begin{align*}
\beta_j &= (\al'_{jk})^{-1} \al^{-1}_{jk} = \biggl( \la_{j,k+1} \prod_{l=1}^{O_-(k)} \la_{j,p^l(k)} \prod_{m=0}^{O_-(k)} \la_{j,p^m(k)} \biggr)^{-1}  \\
&= \bigl( \al_{j,p(k)} \al_{j,k+1} \bigr)^{-1} = \ga_j,
\end{align*}
and so $\bigl( y_k y'_k - \al_{k p(k)} y_{p(k)} y_{k+1} \bigr) x_j = \ga_j x_j \bigl( y_k y'_k - \al_{k p(k)} y_{p(k)} y_{k+1} \bigr)$. This shows that $y_k y'_k - \al_{k p(k)} y_{p(k)} y_{k+1}$ is 
a normal element of $R_{k-1}$ which normalizes the elements of $R_{k-1}$ in exactly the same way as $y_{p(k)} y_{k+1}$, and completes the proof of the theorem.
\end{proof}

Our next result turns the conclusion of \thref{1} (b) into a cluster 
mutation statement.

\bth{2} In the setting of Theorem {\rm\ref{t1} (b)}, there exist
a collection of nonnegative integers
$\{ m_i \mid i \in P(k-1), \, i \neq p(k)\}$ and 
$\kappa \in \kx$ such that 
\begin{equation}
\label{thm2-disp}
y'_k = y_k^{-1} \Big( \al_{k p(k)} y_{p(k)} y_{k+1} + \kappa 
\prod_{i \in P(k-1), \, i \neq p(k)} y_i^{m_i} \Big),
\end{equation}
in the conventions from \eqref{conven} and \eqref{inf}, recall also \eqref{al}. If $p(k) = -\infty$, then all $m_i = 0$.
\eth

\begin{proof} By \thref{1} (b), $y_k y'_k - \al_{k p(k)} y_{p(k)} y_{k+1}$
is a homogeneous normal element of $R_{k-1}$. Applying \prref{factorHUFD} and \thref{CGL} we obtain 
$$
y_k y'_k - \al_{k p(k)} y_{p(k)} y_{k+1} =
\kappa \prod_{i \in P(k-1)} y_i^{m_i}
$$
for some $\kappa \in \KK$ and a collection of nonnegative integers
$\{ m_i \mid i \in P(k-1) \}$. Recall from \eqref{eltp(k)-inf} that if $p(k) = -\infty$, then $y_k y'_k - \al_{k p(k)} y_{p(k)} y_{k+1}$ is a nonzero scalar. Hence, $m_i = 0$ for all $i \in P(k-1)$ in this case.
We need to prove that $\kappa \neq 0$, 
and that $m_{p(k)} = 0$ if $p(k) \neq - \infty$.

Suppose that $\kappa = 0$. Then  
$$
y_k y'_k - \al_{k p(k)} y_{p(k)} y_{k+1} = 0
$$ 
which is a contradiction since 
$y_{k+1}$ is a prime element of $R_{k+1}$ which does not divide either $y_k$ or $y'_k$.

Now suppose that $p(k) \neq - \infty$ and $m_{p(k)} \neq 0$. Then
$y_{p(k)}$ is a prime element of $R_{k-1}$ and 
$$
y_k y'_k - \al_{k p(k)} y_{p(k)} y_{k+1}  \in y_{p(k)} R_{k-1}.
$$
Hence,
$$
y_k y'_k \in y_{p(k)} R_{k+1}.
$$
Furthermore, by \thref{CGL}, $y_k = y_{p(k)} x_k - c_k$ 
and $y'_k = y_{p(k)} x_{k+1} - c'_k$ for some 
$c_k,  c'_k \in R_{k-1}$ such that $y_{p(k)}$ does not divide $c_k$ or 
$c'_k$. Taking into account the normality 
of $y_{p(k)}$ in $R_{k-1}$ gives
$$
c_k c'_k = y_k y'_k - y_{p(k)} x_k y'_k + c_k y_{p(k)} x_{k+1} 
\in y_{p(k)} R_{k+1} \cap R_{k-1} = y_{p(k)} R_{k-1}.
$$
This contradicts the fact that $y_{p(k)}$ 
is a prime element of $R_{k-1}$ which does not divide $c_k$ or 
$c'_k$.
\end{proof}

\section[A base change and normalization of the elements $y_j$]{A base change and normalization of the elements $y_j$\\ from \thref{CGL}}
\label{4.2}
If $R$ is a $\KK$-algebra with a rational action of 
a torus $(\kx)^r$ by algebra automorphisms and $\Ffield/\KK$ 
is a field extension, then there is a canonical rational 
action of $(\Ffield^*)^r$ on $R \otimes_\KK \Ffield$ by $\Ffield$-algebra 
automorphisms. This easily implies the validity of the 
following lemma.

\ble{Fext} Assume that $R$ is a CGL extension as in \eqref{firstCGL} 
with an action of the torus $\HH=(\kx)^r$, 
and $\Ffield/\KK$ is a field extension. Then, as $\Ffield$-algebras,
$$
R \otimes_\KK \Ffield \cong \Ffield [x_1] [x_2; \sigma^\circ_2, \delta^\circ_2] 
\ldots [x_N; \sigma^\circ_N, \delta^\circ_N]
$$
where $\sigma^\circ_k$ and $\delta^\circ_k$ are the automorphisms and skew derivations
obtained from the ones in \eqref{itOre} by base change. Furthermore, 
$R \otimes_\KK \Ffield$ is a CGL extension with the canonical 
induced rational action of $(\Ffield^*)^r$ and the same choice of 
elements $h_k$.
\ele

\bre{sqroots} After a field extension, every (symmetric) CGL extension $R$ 
as in Definitions \ref{dCGL} and \ref{dsymmetric} 
can be brought to one that satisfies 
$\sqrt{\la_{lj}} \in \KK$ for all $l,j \in [1,N]$
by applying \leref{Fext} with
$$
\Ffield:= \KK[ \sqrt{\la_{lj}} \; | \; 1 \leq j < l \leq N].  
$$
\ere

Let $R$ be a CGL extension as in \eqref{firstCGL} that satisfies
$\sqrt{\la_{lj}} \in \KK$ for all $l,j \in [1,N]$. Fix
\begin{equation}
\label{nu1}
\nu_{lj} \in \KK \; \; 
\mbox{for} \; \;
1 \leq j < l \leq N 
\quad \mbox{such that} \; \; 
\nu^2_{lj} = \la_{lj}.
\end{equation}
Set
\begin{equation}
\label{nu2}
\nu_{ll} := 1 \; \; 
\mbox{for} \; \; l \in [1,N]
\quad \mbox{and} \quad
\nu_{jl} := \nu_{lj}^{-1} \; \; 
\mbox{for} \; \; 1 \leq j < l \leq N.
\end{equation}
Then $\nub:=(\nu_{lj}) \in M_N(\kx)$ is a 
multiplicatively skew-symmetric matrix 
and $\nub^{\cdot2} = \lab$. \index{nu@$\nub$}

Analogously to \eqref{q}, for 
$l,j  \in [1,N]$ set
\begin{equation}
\label{r}
r_{lj} := \Om_\nub(\ol{e}_l, \ol{e}_j)= \prod_{m=0}^{O_-(l)} \prod_{n=0}^{O_-(j)}
\nu_{p^m(l), p^n(j)} \in \kx,
\end{equation}
recall \eqref{ol-e}.
The matrix $\rbf = (r_{lj}) \in M_N(\kx)$  \index{r@$\rbf$}
is multiplicatively skew-symmetric, 
$\rbf^{\cdot 2} = \qb$, and 
\begin{equation}
\label{La-r}
\Om_\rbf(e_l, e_j) := \Om_\nub(\ol{e}_l, \ol{e}_j), \quad
\forall l, j \in [1,N].
\end{equation}

We normalize the elements $y_j$ from \thref{CGL} by 
$$
\ol{y}_j = \Big( 
\prod_{0 \le n < m \le O_-(j)} 
\nu_{p^m(j), p^n(j)}^{-1} \Big)
y_j, \; \; \forall j \in [1,N].
$$
\index{yj@$\ol{y}_j$}
The meaning of this normalization is as follows. The generators $x_j$ of $R$ skew commute 
up to lower terms of $R$ in the partial order from Section \ref{detailCGL}: 
$$
x_l x_j - \la_{lj} x_j x_l \in \Span \{ x^g \mid g \in \Zset_{\geq 0}^N, \; g \prec e_l + e_j \}.
$$
Ignoring these 
lower order terms, the above normalization is precisely the normalization 
from \eqref{basis} applied to the leading term $x^{\ol{e}_j}$ of $y_j$ (cf. 
\eqref{tty}), thought of as an element of the abstract based quantum 
torus for the multiplicatively skew-symmetric matrix $\nub$
(recall Section \ref{3.1}); that is
\begin{equation}
\label{new-y}
\ol{y}_j = \Scr_{\nub}(\ol{e}_j) y_j, 
\; \; \forall j \in [1,N].
\end{equation}

\prref{CGLcluster} immediately implies:

\bpr{tframe} Consider a CGL extension presentation of an algebra $R$  
as in \eqref{itOre} for which $\sqrt{\la_{lj}} \in \KK$ for all $l,j \in [1,N]$,
and define the matrices $\nub, \rbf \in M_N(\kx)$ by 
\eqref{nu1}, \eqref{nu2}, and \eqref{r}. There exists a unique 
toric frame $M : \Zset^N \to \Fract(R)$ having matrix $\rbf$ and such that 
$M(e_j) = \ol{y}_j$, $\forall j \in [1,N]$. 
\epr 

\bex{OqMmn5} 
In the case of $R = \OqM$, the scalars $\la_{lj}$ for $l > j$ are either $1$ or $q^{-1}$, so the only square root needed in $\KK$ is  $\sqrt q$. (We fix a choice of this scalar.) However, $\sqrt q$ does not appear in the normalizations $\ol{y}_j$. In fact, $\nu_{p^s(j), p^t(j)} = 1$ for $0 \le t < s \le O_-(j)$ (recall the choice of $\eta$ from \exref{OqMmn2} and the description of the $\la_{kl}$ in \eqref{OqMlakj}), and so in this example, $\ol{y}_j = y_j$ for all $j \in [1,N]$.
\eex

\section{Almost cluster mutations between CGL extension presentations}
\label{4.3}
Let $R$ be a $\KK$-algebra with a CGL extension presentation
as in \eqref{firstCGL}. We will assume that 
the torus $\HH$ acting on $R$ is the universal maximal torus 
from \thref{univ-torus}.
It is easy to see that in this setting, 
the assumption that \eqref{secondCGL} is a second CGL extension 
presentation of $R$ is equivalent to the following condition:
\begin{enumerate}
\item[(i)] $\delta_{k+1} (x_j) \in R_{k-1}$, $\forall j \in [1,k]$ and 
$\exists h'_{k+1} \in \HH$ such that $h'_{k+1} \cdot x_k = \la'_{k+1} x_k$ for some $\la'_{k+1} \in \kx$ which is not a root of unity, and
$h'_{k+1} \cdot x_j = \la_{kj} x_j$, $\forall j \in [1,k-1]\cup \{ k+1 \}$. 
\end{enumerate}
The rest of the data for the second CGL extension presentation of $R$
($\KK$-torus $\HH'$, scalars $\la'_{kj} \in \kx$, and elements $h'_k \in \HH'$) 
is given by \leref{CGLtranspose} and the following:  
\begin{enumerate}
\item[(ii)] The torus $\HH'$ acting on the second CGL extension presentation 
can be taken as the original universal maximal torus $\HH$
from \thref{univ-torus}. The corresponding elements $h'_j \in \HH$ 
are given by (i) for $j = k+1$ and $h'_j = h_j$ for $j \neq k+1$.
\item[(iii)] $\sigma'_{k+1} = (h'_{k+1} \cdot)|_{R'_k}$, 
where $R'_k$ is the unital $\KK$-subalgebra of $R$ 
generated by $x_1, \ldots, x_{k-1}$, and $x_{k+1}$.
\item[(iv)] $\la'_{lj} = \la_{(k,k+1)(l), (k, k+1) (j)}$ for $l,j \in [1,N]$.
\end{enumerate}

{\em{Here and below, $(k,k+1)$ denotes a transposition 
in the symmetric group $S_N$, and $S_N$ is embedded in $GL_N(\Zset)$ 
via permutation matrices.}}

Assume that the condition (i) is satisfied and 
the base field $\KK$ contains all square roots $\sqrt{\la_{lj}}$ for  
$j,l \in [1,N]$, recall \reref{sqroots}. Fix 
$\nu_{lj} \in \kx$ for $1 \leq j < l \leq N$ such that
$\nu_{lj}^2 = \la_{lj}$ and define a multiplicatively 
skew-symmetric matrix $\nub:=(\nu_{lj}) \in M_N(\kx)$ by 
\eqref{nu2}. Let $\nub' := (\nu'_{lj}) \in M_N(\kx)$
be the multiplicatively skew-symmetric matrix 
\begin{equation}
\label{nuprime}
\nub' := (k,k+1) \nub (k,k+1) \in M_N(\kx).
\end{equation}
By property (iv), $(\nu'_{lj})^2 = \la'_{lj}$, $\forall l,j \in [1,N]$. 
Denote by $\rbf' := (r'_{lj})\in M_N(\kx)$ 
the matrix obtained from $\nub'$ by \eqref{r}, where $p$ is replaced by the predecessor function for $\eta'=\eta (k,k+1)$.

\prref{tframe} implies that $((\ol{y}_1, \ldots, \ol{y}_N), \rbf)$ 
and $((\ol{y}'_1, \ldots, \ol{y}'_N), \rbf')$ define 
two toric frames $M, M' : \Zset^N \to \Fract(R)$. Their matrices 
equal $\rbf$ and $\rbf'$, respectively, 
and $M(e_j) = \ol{y}_j$, $M'(e_j) = \ol{y}'_j$, $\forall j \in [1,N]$.
The next theorem describes the relationship between these 
two toric frames on the basis of \thref{2}. Let 
$E_+^\ci \in M_N(\Zset)$ be the matrix with entries
$$
(E_+^\ci)_{lj} = 
\begin{cases}
\delta_{lj}, & \mbox{if} \; j \neq k \\
-1, & \mbox{if} \; l=j=k \\
1, & \mbox{if} \; j = k, \; \mbox{and $l=p(k)$ or $k+1$}, \\
0, & \mbox{if} \; j = k, \; \mbox{and $l \ne p(k)$, $k$, or $k+1$},
\end{cases}
$$
cf. \eqref{Eep}.

For the following theorems, we adopt the convention that
\begin{equation}
\label{e-inf}
e_{-\infty} = \ol{e}_{-\infty} := 0.
\end{equation}

\bth{3} Assume the setting of Theorem {\rm\ref{t1}} and that 
$\sqrt{\la_{lj}} \in \KK$ for all $j,l \in [1,N]$.

{\rm(a)} If $\eta(k) \neq \eta(k+1)$, then $\rbf' = (k,k+1) \rbf (k,k+1)$ and $\ol{y}'_j = \ol{y}_{(k, k+1) j}$, 
for all $j \in [1,N]$, i.e., 
$M'=M (k,k+1)$.

{\rm(b)} If $\eta(k) = \eta(k+1)$, then $M'(e_j) = \ol{y}'_j = \ol{y}_j= M(e_j)$, 
$\forall j \neq k$ and 
\begin{equation}
\label{MM'}
M'(e_k)=
\ol{y}'_k = M ( - e_k + e_{p(k)} + e_{k+1}) + \zeta
M \Big( - e_k + \sum_{i \in P(k-1), \, i \neq p(k)} m_i e_i \Big)
\end{equation}
for the collection of nonnegative integers $\{m_i \mid i \in P(k-1), \; i \neq p(k)\}$
from Theorem {\rm\ref{t2}} and some $\zeta \in \kx$. Furthermore,
\begin{equation}
\label{rr'}
\rbf' = {}^{E^{\ci{}T}_+} \! \rbf^{E^\ci_+}.
\end{equation}
\eth

\begin{proof} Part (a) of the theorem and the equality $\ol{y}'_j = \ol{y}_j$, 
$\forall j \neq k$ in (b) follow from \thref{1}, once one verifies that
$$
\Scr_{\nub'} (\ol{e}_j) = \begin{cases} 
\la_{k,k+1} \Scr_{\nub}(\ol{e}_j), &\mbox{if} 
\; \; j = s^l(k+1) \; \; \mbox{for some} \; \; l \in \Znn 
\\
\Scr_{\nub}(\ol{e}_j), &\mbox{if} 
\; \; j < k, \; \mbox{or} \; j > k+1 \; 
\mbox{and} \; j \neq s^l(k+1) \; \forall l \in \Znn.
\end{cases}
$$

Next we verify \eqref{MM'}. From \eqref{ycomm} and the observation that $q_{kp(k)} = \al_{kp(k)}$, we have 
$$
\al_{k p(k)} y_k^{-1} y_{p(k)} y_{k+1} = y_{p(k)} y_k^{-1} y_{k+1}.
$$ 
In view of \eqref{new-y}, \eqref{basis}, and \eqref{Scr}, \thref{2} implies that \eqref{MM'} 
is equivalent to the identity
\begin{multline*}
r_{p(k),k} r_{k, k+1} r^{-1}_{p(k), k+1}
\Scr_\nub(\ol{e}_{p(k)}) \Scr_\nub(\ol{e}_k)^{-1} \Scr_\nub(\ol{e}_{k+1}) = 
\Scr_{\nub'}(\ol{e}_k) =
\\
\prod \{ \nu^{-1}_{ij} \mid i, j = p^{O_-(k)}(k), \ldots, p(k), k+1; \; \; 
i<j \},
\end{multline*}
where $r_{-\infty,j} := 1$ and $\ol{e}_{-\infty} := 0$.
This identity follows from \eqref{r}.

To show \eqref{rr'}, first note that Eq.~\eqref{La-r} (applied to the presentation \eqref{secondCGL}) implies $\Om_{\rbf'}(e_l,e_j) = \Om_{\nub}((k,k+1)\ol{e}_l, (k,k+1)\ol{e}_j)$ for all $l,j \in [1,N]$ (recall \eqref{nuprime}), whence 
\begin{align*}
&\Om_{\rbf'}(e_l, e_j) = \Om_\nub(\ol{e}_l, \ol{e}_j), \quad \forall l, j \neq k,  
\\
&\Om_{\rbf'}(e_k, e_j) = \Om_\nub(\ol{e}_{p(k)} + e_{k+1}, \ol{e}_j), \quad \forall j \neq k.
\end{align*}
Since $\ol{e}_{p(k)} + \ol{e}_{k+1} - \ol{e}_k = \ol{e}_{p(k)} + e_{k+1}$, 
taking into account \eqref{La-r} applied to the presentation \eqref{firstCGL}, 
we obtain
\begin{align*}
&\Om_{\rbf'}(e_l, e_j) = \Om_\rbf(e_l, e_j), \quad \forall l,j \neq k,
\\
&\Om_{\rbf'}(e_k, e_j) = \Om_\rbf(e_{p(k)} + e_{k+1} - e_k, e_j), \quad \forall j \neq k.
\end{align*}
These two equalities are equivalent to \eqref{rr'}. 
\end{proof}

\section{Scalars associated to mutation of prime elements}
\label{4.4}
Next, we derive formulas for certain scalars which play the role of the 
entries of the matrix $\wt{\tb}$ from Section \ref{3.2} for the quantum seeds which we 
construct in Chapter \ref{main}. Assume that there exist square roots $\nu_{lj} = \sqrt{\la_{lj}} \in \KK$ for $1 \le j < l \le N$.

In the next theorem and Chapters \ref{mut-sym} and \ref{main}, 
we will consider the following mild condition:
\begin{equation}
\label{-1}
\begin{aligned}
{\mbox{The subgroup of $\kx$ generated by\;}} &\{\nu_{lj} \mid 1 \leq j < l \leq N \}  \\ 
 &{\mbox{contains no elements of order $2$}}.
 \end{aligned}
\end{equation}
(This condition automatically holds if $\chr \KK = 2$; otherwise, it just means that the given subgroup of $\kx$ does not contain $-1$.)
For all torsionfree CGL extensions $R$, cf. Section \ref{2.2}, one can always choose 
the $\nu_{lj}$ 
in such a way that \eqref{-1} is satisfied; in fact, it suffices to assume that there are no elements of order $2$ in the subgroup of $\kx$ generated by $\{ \la_{lj} \mid 1 \le j < l \le N \}$. Condition \eqref{-1} implies that the bicharacter $\Om_\nub : \Zset^N \times \Zset^N \to \kx$ 
associated to the multiplicatively skew-symmetric matrix $\nub$ does not take the value
$-1$ if $\chr \KK \ne 2$.

Recall from Section \ref{2.2} that for an $\HH$-eigenvector $u \in R$, $\chi_u \in \xh$ denotes
its eigenvalue.

\bth{scalars} Assume the setting of Theorem {\rm\ref{t2}} and that there exists $\nu_{lj} = \sqrt{\la_{lj}} \in \KK$ for $1 \le j < l \le N$.  
Then
\begin{equation}
\label{Lar1}
\Om_\rbf \Big(e_{p(k)} + e_{k+1} - \sum_{i \in P(k-1),\, i \neq p(k)} m_i e_i, \,
e_j \Big)^2 = 1, \quad \forall j \neq k
\end{equation}
and
\begin{equation}
\label{Lar2}
\begin{aligned}
\Om_\rbf \Big(e_{p(k)} + e_{k+1} - \sum_{i \in P(k-1),\, i \neq p(k)} m_i e_i, \,
e_k \Big)^2 &= \chi_{x_{k+1}}(h_{k+1})^{-1}  \\
 &= \chi_{x_k}(h'_{k+1}) 
 \end{aligned}
\end{equation}
for the collection of nonnegative integers $\{m_i \mid i \in P(k-1), \, i \neq p(k)\}$
from Theorems {\rm\ref{t2}} and {\rm\ref{t3} (b)}. {\rm(}Recall also condition {\rm(i)} in Section {\rm\ref{4.3}.)}

If in addition \eqref{-1} is satisfied, then 
\begin{equation}
\label{Lar3}
\Om_\rbf \Big(e_{p(k)} + e_{k+1} - \sum_{i \in P(k-1),\, i \neq p(k)} m_i e_i, \,
e_j \Big) = 1, \quad \forall j \neq k.
\end{equation}
\eth

\begin{proof} Denote for brevity the elements
$$
g := \sum_{i \in P(k-1), \, i \neq p(k)} m_i e_i \quad \text{and} \quad
\ol{g} := \sum_{i \in P(k-1), \, i \neq p(k)} m_i \ol{e}_i
$$
in $\Zset^N$.
Using the fact that $\ol{y}_j$, $j \neq k$ belongs to the images of both 
toric frames $M$ and $M'$, we obtain
\begin{align*}
\xi_{kj} \ol{y}_j \ol{y}'_k &= \ol{y}'_k \ol{y}_j = 
\big( M(- e_k + e_{p(k)} +e_{k+1}) + \zeta M( - e_k + g) \big) \ol{y}_j  \\
&= \ol{y}_j \big( \Om_\rbf (- e_k + e_{p(k)} + e_{k+1}, e_j)^2 M(- e_k + e_{p(k)} + e_{k+1}) \\
&\qquad\qquad\qquad \qquad\qquad\qquad + \zeta \Om_\rbf(- e_k + g, e_j)^2 M( -e_k + g) \big),
\end{align*}
for some $\xi_{kj} \in \kx$. Hence, 
$$
\Om_\rbf (- e_k + e_{p(k)} + e_{k+1}, e_j)^2= \Om_\rbf(- e_k + g, e_j)^2, \quad
\forall j \neq k, 
$$
which implies \eqref{Lar1} and \eqref{Lar3}. Applying \eqref{Lar1} for 
$j = k+1$ and \eqref{La-r} leads to 
\begin{align*} 
\Om_\rbf (e_{p(k)} + e_{k+1} - g, e_k)^2  &= 
\Om_\rbf (e_{p(k)} + e_{k+1} - g, e_k)^2 \Om_\rbf (e_{p(k)} + e_{k+1} - g, e_{k+1})^{-2}  \\
&= \Om_\nub (\ol{e}_{p(k)} + \ol{e}_{k+1} - \ol{g}, \ol{e}_k)^2 
\Om_\nub (\ol{e}_{p(k)} + \ol{e}_{k+1} - \ol{g}, \ol{e}_{k+1})^{-2}  \\
&= \Om_\nub (\ol{e}_{p(k)} + \ol{e}_{k+1} - \ol{g}, e_{k+1})^{-2}  \\
&= \Om_\nub (\ol{e}_{p(k)} + \ol{e}_k - \ol{g}, e_{k+1})^{-2}.  
\end{align*}
Since $\ol{y}'_k$ is an $\HH$-eigenvector, it follows from \eqref{MM'} that $\chi_{M(e_{p(k)} + e_{k+1} -g)} = 1$. 
Using that $\chi_{M(e_{p(k)} + e_{k+1} - g)} = \chi_{M(e_{p(k)} + e_k -g)} + \chi_{x_{k+1}}$,
we obtain
\begin{align*}
\Om_\nub (\ol{e}_{p(k)} + \ol{e}_k - \ol{g}, e_{k+1})^{-2} &= \chi_{M(e_{p(k)} + e_k -g)}(h_{k+1})  \\
&= \chi_{M(e_{p(k)} + e_{k+1} -g)}(h_{k+1}) \chi_{x_{k+1}}(h_{k+1})^{-1}  \\
 &= \chi_{x_{k+1}}(h_{k+1})^{-1},
\end{align*}
which proves the first equality in \eqref{Lar2}.

It follows from \eqref{Lar1} that $M(e_{p(k)} + e_{k+1} - g)$ commutes with $M(g)$, and hence also with $M(e_{p(k)} + e_{k+1})$. By \thref{3} (b), $\ol{y}_k \ol{y}_k'$ is a linear combination of $M(e_{p(k)} + e_{k+1})$ and $M(g)$, so $M(e_{p(k)} + e_{k+1} - g)$ also commutes with $\ol{y}_k \ol{y}'_k$. By the first equality in \eqref{Lar2},
$$
M(e_{p(k)} + e_{k+1} -g) \ol{y}_k = \chi_{x_{k+1}}(h_{k+1})^{-1} 
\ol{y}_k M(e_{p(k)} + e_{k+1} -g). 
$$
Interchanging the roles of $x_k$ and $x_{k+1}$ and using the symmetric nature 
of the assumption of \thref{1} (b) shows that 
$$
M'(e_{p(k)} + e_{k+1} -g) \ol{y}'_k = \chi_{x_k}(h'_{k+1})^{-1} 
\ol{y}'_k M'(e_{p(k)} + e_{k+1} -g ).
$$ 
In view of \thref{3} (b), the element $M'(e_{p(k)} + e_{k+1} - g)$ is a scalar multiple of $M(e_{p(k)} + e_{k+1} - g)$, so we conclude that
$$M(e_{p(k)} + e_{k+1} - g) \ol{y}'_k = \chi_{x_k}(h'_{k+1})^{-1} \ol{y}'_k M(e_{p(k)} + e_{k+1} - g).$$
Therefore, $\chi_{x_{k+1}}(h_{k+1}) \chi_{x_k}(h'_{k+1}) = 1$, 
which proves the second equality in \eqref{Lar1}.
\end{proof}

\chapter[Homogeneous prime elements for subalgebras of CGL extensions]{Homogeneous prime elements for subalgebras of symmetric CGL extensions}
\label{prime-sym}
Each symmetric CGL extension $R$ of length $N$ has many different 
CGL extension presentations given by \eqref{tauOre}. They are parametrized by 
the elements of the subset $\Xi_N$ of $S_N$, cf. \eqref{tau}.
In order to phrase \thref{3} into a mutation statement 
between toric frames for $\Fract(R)$ associated to the elements of $\Xi_N$ and 
to make the scalars $\zeta$ from \thref{3} equal to one, 
we need a good picture of the sequences of homogeneous prime 
elements $y_1, \ldots, y_N$ from \thref{CGL} 
associated to each presentation \eqref{tauOre}. 
This is obtained in \thref{sym-prime}. \thref{y-int} contains
a description of the homogeneous prime elements that enter
into this result. Those prime elements (up to rescaling)
comprise the cluster variables that 
will be used in Chapter \ref{main} to construct 
quantum cluster algebra structures on symmetric 
CGL extensions. Along the way, we explicitly describe the elements 
of $\Xi_N$ and prove an invariance property of the scalars 
$\la_l$ and $\la_l^*$ from Definitions \ref{dCGL} and \ref{dsymmetric}.
\thref{u-prod}, which appears at the end of the chapter, contains
a key result used in the next chapter to normalize the generators $x_j$ 
of symmetric CGL extensions so that all scalars $\zeta$ in \thref{3} 
become equal to one.

Throughout the chapter, $\eta$ will denote a function $[1,N] \rightarrow \Zset$ satisfying the conditions of \thref{CGL}, with respect to the original CGL extension presentation \eqref{itOre} of $R$, and $p$ and $s$ will denote the corresponding predecessor and successor functions. We will repeatedly use the one-line notation \eqref{one-line} for permutations.
\section{The elements $y_{[i,s^m(i)]}$}
\label{4a.1}
Recall from Section \ref{2.5} that for a symmetric CGL extension $R$ of rank $N$
and $1 \leq j \leq k \leq N$, $R_{[j,k]}$ denotes the unital subalgebra 
of $R$ generated by $x_j , \ldots, x_k$. It is an Ore extension of both 
$R_{[j,k-1]}$ and $R_{[j+1,k]}$. All such subalgebras 
are (symmetric) CGL extensions and \thref{CGL} applies 
to them. 

For $i \in [1,N]$ and $0 \le m \le O_+(i)$, recall \eqref{Opm} (i.e., 
$s^m(i) \in [1,N]$), set 
\begin{equation}
\label{e-int}
e_{[i, s^m(i)]} = e_i + e_{s(i)} + \cdots + e_{s^m(i)} \in \Zset^N.
\end{equation}
\index{eismi@$e_{[i, s^m(i)]}$}
The vectors \eqref{ol-e} are special cases of those:
$$
\ol{e}_k = e_{[p^{O_-(k)}(k), k]}, \quad
\forall k \in [1,N]. 
$$
We also set $e_\varnothing = 0$.
The next theorem treats the prime elements that will appear as 
cluster variables for symmetric CGL extensions. It will be proved in Section \ref{4a.4}.

\bth{y-int} Assume that $R$ is a symmetric CGL extension of length $N$, and $i \in [1,N]$ and 
$m \in \Zset_{\geq 0}$ are such that $s^m(i) \in [1,N]$, i.e., 
$s^m(i) \neq + \infty$. Then the following hold:

{\rm(a)} All homogeneous prime elements of $R_{[i, s^m(i)]}$ 
that do not lie in $R_{[i, s^m(i)-1]}$ are associates of each other.

{\rm(b)} All homogeneous prime elements of $R_{[i, s^m(i)]}$ 
that do not lie in $R_{[i+1, s^m(i)]}$ are associates of each other.
In addition, the set of these homogeneous prime elements coincides 
with the set of homogeneous prime elements in part {\rm(a)}.

{\rm(c)} The homogeneous prime elements in {\rm(a)} and {\rm(b)} have leading terms of the form 
$$
\xi x^{e_{[i,s^m(i)]}} = \xi x_{i} x_{s(i)} \ldots x_{s^m(i)}
$$ 
for some $\xi \in \kx$, see \S{\rm\ref{detailCGL}}. For each $\xi \in \kx$, there 
is a unique homogeneous prime element of $R_{[i,s^m(i)]}$ with such a leading term.
Denote by $y_{[i,s^m(i)]}$ the 
prime element with leading term $x_{i} x_{s(i)} \ldots x_{s^m(i)}$.
Let $y_\varnothing :=1$.
\index{yismi@$y_{[i,s^m(i)]}$}

{\rm(d)} We have 
$$
y_{[i,s^m(i)]} = y_{[i, s^{m-1}(i)]} x_{s^m(i)} - c_{[i,s^m(i)-1]}
= x_i y_{[s(i), s^m(i)]} - c'_{[i+1, s^m(i)]}
$$
for some $c_{[i, s^m(i)-1]} \in R_{[i,s^m(i)-1]}$ 
and $c'_{[i+1, s^m(i)]} \in R_{[i+1, s^m(i)]}$.

{\rm(e)} For all $k \in [1,N]$ such that $p(i)< k< s^{m+1}(i)$, we have
$$
y_{[i,s^m(i)]} x_k = \Om_\lab( e_{[i,s^m(i)]}, e_k ) x_k y_{[i,s^m(i)]} .$$
\eth

The case $m=0$ of this theorem is easy to verify.
In that case $y_{[i,i]}= x_i$, $\forall i \in [1,N]$.

\bex{OqMmn6}
In the case of $R= \OqM$, the elements $y_{[i, s^l(i)]}$ of \thref{y-int} are solid quantum minors, just as in \exref{OqMmn2}. More precisely, if $i \in [1,N]$ and $l \in \Znn$ with $s^l(i) \ne +\infty$, then $i = (r-1)n+c$ with $r, r+l \in [1,m]$ and $c,c+l \in [1,n]$, and $y_{[i, s^l(i)]} = \De_{[r, r+l], [c, c+l]}$.
\eex

The following theorem describes the $y$-sequences from \thref{CGL} associated to the 
CGL extension presentations \eqref{tauOre} in terms of the prime elements 
from \thref{y-int}. It will be proved in Section \ref{4a.4}. Recall that 
for every $\tau \in \Xi_N$ and $k \in [1,N]$, $\tau([1,k])$ is an 
interval.

\bth{sym-prime} Assume that $R$ is a symmetric CGL extension of length $N$
and $\tau$ an element of the subset $\Xi_N$ of $S_N$, cf. \eqref{tau}.
Let $y_{\tau, 1}, \ldots, y_{\tau, N}$ be the sequence 
in $R$ from Theorem {\rm\ref{tCGL}} applied to the CGL extension presentation 
\eqref{tauOre} of $R$ corresponding to $\tau$. Let $k \in [1,N]$. 

If $\tau(k) \geq \tau(1)$, then $y_{\tau, k}$ is a scalar multiple of 
$y_{[p^m(\tau(k)), \tau(k)]}$, where 
\begin{equation}
\label{m-p}
m = \max \{ n \in \Znn \mid p^n( \tau(k)) \in \tau([1, k]) \}.
\end{equation}

If $\tau(k) \leq \tau(1)$, then $y_{\tau, k}$ is a scalar multiple of 
$y_{[\tau(k), s^m(\tau(k))]}$, where 
\begin{equation}
\label{m-s}
m = \max \{ n \in \Znn \mid s^n( \tau(k)) \in \tau([1, k]) \}.
\end{equation}
In both cases, the predecessor and successor functions are with 
respect to the original CGL extension presentation 
\eqref{itOre} of $R$.
\eth

\section{The elements of $\Xi_N$}
\label{4a.3}
In this and the next section, we investigate the 
elements of the subset $\Xi_N$ of $S_N$ defined 
in \eqref{tau}. It follows from \eqref{tau} that every element $\tau \in \Xi_N$ 
has the property that either $\tau(N) = 1$ or $\tau(N) = N$. This implies the following recursive 
description of $\Xi_N$.

\ble{SigN-ind} For each $\tau \in \Xi_N$, there exists $\tau' \in \Xi_{N-1}$ 
such that either
$$
\tau(i) = \tau'(i), \; \; \forall i \in [1,N-1] \quad 
\mbox{and} \quad \tau(N) = N
$$
or 
$$
\tau(i) = \tau'(i) + 1, \; \;  \forall i \in [1,N-1]
\quad \mbox{and} \quad \tau(N) =1.
$$
For all $\tau' \in \Xi_{N-1}$, the above define elements of $\Xi_N$.
\ele

Given $k \in [1,N]$ and a sequence $k \leq j_k \leq \cdots \leq j_1 \leq N$, 
define
\begin{equation}
\label{tau-seq}
\tau_{(j_k, \ldots, j_1)} 
:= (k (k+1) \ldots j_k) 
\ldots
(2 3 \ldots j_2)
(1 2 \ldots j_1) \in S_N,
\end{equation}
\index{taujkj1@$\tau_{(j_k, \ldots, j_1)}$}
where in the right hand side we use the standard notation 
for cycles in $S_N$. \leref{SigN-ind} implies by induction the 
following characterization of the elements of $\Xi_N$. We 
leave its proof to the reader.

\ble{charSigN} The subset $\Xi_N\subset S_N$ consists 
of the elements of the form $\tau_{(j_k, \ldots, j_1)}$, 
where $k \in [1,N]$ and $k \leq j_k \leq \cdots \leq j_1 \leq N$.  
\ele

The representation of an element of $\Xi_N$ 
in the form \eqref{tau-seq} is not unique. One way to visualize 
$\tau_{(j_k, \ldots, j_1)}$ is that the sequence $\tau(1), \ldots, \tau(N)$ 
is obtained from the sequence $1, \ldots, N$ by the following procedure:

$(*)$ {\em{The number $1$ is pulled to the right to position $j_1$ {\rm(}preserving the order of the other numbers\/{\rm)}, 
then the number $2$ is pulled to the right to position 
$j_2-1$, ..., at the end the number $k$ is pulled to the right to position $j_k-k+1$.}}

For example, for $k=2$ the following illustrates how $\tau_{(j_2, j_1)}$ is obtained from the 
identity permutation:
\begin{align*}
&[\circled{1}, \circled{2},3, 4, \ldots, j_2, j_2+1, \ldots, j_1, j_1+1, \ldots, N] \mt
\\
&[\circled{2},3, 4, \ldots, j_2, j_2+1, \ldots, j_1, \circled{1}, j_1+1, \ldots, N] \mt
\\
&[3, 4, \ldots, j_2, \circled{2}, j_2+1, \ldots, j_1, \circled{1}, j_1+1, \ldots, N],
\end{align*}
where the numbers that are pulled ($1$ and $2$) are circled.

If we perform the above procedure one step at a time, so in each step we only interchange the positions 
of two adjacent numbers, then the elements of $S_N$ from all intermediate steps will 
belong to $\Xi_N$. For example, this requires factoring the cycle $(1,2,\dots,j_1)$ as $(1,j_1)(1,j_1-1) \cdots (1,2)$ rather than as $(1,2) (2,3) \cdots (j_1-1,j_1)$. This implies at once the first part of the next corollary.
 
\bco{steps} Let $R$ be a symmetric CGL extension of length $N$ and $\tau \in \Xi_N$. 

{\rm(a)} There exists a sequence $\tau_0 = \id, \tau_1, \ldots, \tau_n = \tau$ in $\Xi_N$ such that 
for all $l \in [1,n]$,
$$
\tau_l = ( \tau_{l-1}(k_l), \tau_{l-1} (k_l+1)) \tau_{l-1} = \tau_{l-1}(k_l, k_l+1)
$$
for some $k_l \in [1,N-1]$ such that $\tau_{l-1}(k_l) < \tau_{l-1}(k_l +1)$.

{\rm(b)} If $\eta : [1,N] \to \Zset$ is a function satisfying the conditions of 
Theorem {\rm\ref{tCGL}} for the original CGL presentation of $R$, then 
\begin{equation}
\label{etatau}
\eta_\tau:=\eta \tau : [1,N] \to \Zset 
\end{equation}
\index{etatau@$\eta_\tau$}
satisfies the conditions of Theorem {\rm\ref{tCGL}} for the $\eta$-function of the 
CGL extension presentation \eqref{tauOre} of $R$ corresponding to $\tau$.
\eco

The above described sequence for the first part of the corollary for the element
$\tau_{(j_k, \ldots, j_1)} \in \Xi_N$ has length $j_1 + \cdots + j_k - k(k+1)/2$. 
The second part of the corollary follows by recursively applying \thref{1} to the 
CGL extension presentations \eqref{tauOre} for the elements
$\tau_{l-1}$ and $\tau_l$. \coref{steps} (b) gives a second proof of the 
fact that the rank of a symmetric CGL extension does not depend on the choice 
of CGL extension presentation of the form \eqref{tauOre}, see Section \ref{2.5}.

\bco{restr-eta} Assume that $R$ is a symmetric CGL extension of length $N$ 
and $\eta : [1,N] \to \Zset$ is a function satisfying the conditions of 
Theorem {\rm\ref{tCGL}}. Then for all $1 \leq j \leq k \leq N$, 
the function $\eta_{[j,k]} \colon [1,k-j+1] \to \Zset$ given by 
$$
\eta_{[j,k]}(l) = \eta(j+l-1), \; \; \forall l \in [1, k-j+1]
$$
satisfies
the conditions of Theorem {\rm\ref{tCGL}} for the symmetric CGL extension 
$R_{[j,k]}$.
\eco

The meaning of \coref{restr-eta} is that the $\eta$-function for ``interval 
subalgebras'' of symmetric CGL extensions can be chosen to be
the restriction of the original $\eta$-function up to a shift. This fact
makes it possible to run induction on the rank of $R$ in various situations. 
\coref{restr-eta} follows by applying \coref{steps} (b) to 
$$
\tau = [j, \ldots, k, j-1, \ldots, 1, k+1, \ldots, N] \in \Xi_N 
$$
and considering the $(k-j+1)$-st subalgebra of the corresponding Ore extension,
which equals $R_{[j,k]}$.

\section{A subset of $\Xi_N$}
\label{4a.2}
 We investigate a subset of $\Xi_N$ which will 
play an important role in Chapter \ref{main}. Write elements $\tau$ of the symmetric group $S_N$ in the form $[\tau(1), \tau(2), \dots, \tau(N)]$.
Recall \eqref{tau-ci} that $w_\ci = [N, N-1, \ldots, 1]$ denotes the longest element of $S_N$.
For $1 \leq i \leq j \leq N$, define the following elements of $\Xi_N$:
\begin{equation}
\label{tauij}
\tau_{i,j} := [i+1, \ldots, j, i, j+1, \ldots, N, i-1, i-2, \ldots, 1] \in \Xi_N.
\end{equation}
\index{tauij@$\tau_{i,j}$}
They satisfy
$$
\tau_{1,1} = \id, \quad \tau_{i,N} = \tau_{i+1, i+1}, \; \forall i \in [1,N-1], \quad
\tau_{N,N} = w_\ci.
$$
Denote by $\Ga_N$ \index{GammaN@$\Ga_N$} the subset of $\Xi_N$ consisting 
of all $\tau_{i,j}$'s and consider the 
following linear ordering on it
\begin{multline}
\label{sequence}
\Ga_N := \{
\id = \tau_{1,1} \prec \ldots \prec \tau_{1,N} = \tau_{2,2} \prec \ldots \prec \tau_{2, N} = \tau_{3,3} \prec \ldots \prec \\ 
\prec \tau_{3,N}= \tau_{4,4} \prec \ldots \prec
\tau_{N-2, N} = \tau_{N-1,N-1} \prec \tau_{N-1,N} = \tau_{N,N} = w_\ci \}.
\end{multline}
\index{zzz@$\prec$}
In the notation of \eqref{tau-seq}, the elements of $\Ga_N$ are given by
$$
\tau_{i,j} := \tau_{(j, N, \ldots, N)}, \, \, \text{where} \; N \; \text{is repeated} \; i-1 \; \text{times}, \, \, \forall 1\le i \le j \le N.
$$
The sequence of elements \eqref{sequence} is nothing but the sequence of the intermediate steps 
of the procedure $(*)$ from Section \ref{4a.3} applied to the longest element $w_\ci \in S_N$ (in which 
case $k=N-1$ and $j_1= j_2 = \cdots = j_{N-1} = N$).

Assume that $R$ is a symmetric CGL extension.  
To each element of $\Ga_N$, \prref{tauOre} associates a CGL extension presentation of $R$.  
Each two consecutive presentations are associated to a pair $\tau_{i,j}, \tau_{i,j+1} \in \Xi_N$ 
for some $1 \leq i \leq j < N$. They have the forms
\begin{multline}
\label{present1}
R = 
\KK[x_{i+1}] \ldots [x_j; \sigma_j, \delta_j][x_i; \sigma^\sy_i, \delta^\sy_i] 
[x_{j+1}; \sigma_{j+1}, \delta_{j+1}] \ldots 
\\
[x_N; \sigma_N, \delta_N] 
[x_{i-1}; \sigma^\sy_{i-1}, \delta^\sy_{i-1}] \ldots [x_1; \sigma^\sy_1, \delta^\sy_1] 
\end{multline}
and
\begin{multline}
\label{present2}
R= 
\KK[x_{i+1}] \ldots [x_j; \sigma_j, \delta_j] 
[x_{j+1}; \sigma_{j+1}, \delta_{j+1}] [x_i; \sigma^\sy_i, \delta^\sy_i] [x_{j+2}; \sigma_{j+2}, \delta_{j+2}]
\ldots 
\\
[x_N; \sigma_N, \delta_N] 
[x_{i-1}; \sigma^\sy_{i-1}, \delta^\sy_{i-1}] \ldots [x_1; \sigma^\sy_1, \delta^\sy_1], 
\end{multline}
respectively. These two presentations satisfy the assumptions of \thref{1}.
Recall from Definitions \ref{dCGL} and \ref{dsymmetric}
that $\sigma_l = (h_l \cdot)$, $\sigma^\sy_l = (h^\sy_l \cdot)$ 
and $h_l \cdot x_l = \la_l x_l$, $h^\sy_l \cdot x_l = \la^\sy_l x_l$
with $h_l, h^\sy_l \in \HH$, $\la_l, \la^\sy_l \in \kx$. By \coref{steps}(b), the $\eta$-function for the CGL extension presentation \eqref{present1} can be taken to be $\eta \tau_{i,j}$. In particular, the values of this function on $j+1-i$ and $j+2-i$ are $\eta(i)$ and $\eta(j+1)$, respectively.
If $\eta(i)=\eta(j+1)$, then Eq. \eqref{Lar2} of \thref{scalars} applied to the 
presentations \eqref{present1} and \eqref{present2} implies 
\begin{equation}
\label{lala*}
\la^\sy_i \la_{j+1} =
\chi_{x_i}(h^\sy_i) \chi_{x_{j+1}}(h_{j+1}) = 1.
\end{equation}

\bpr{la-equal} Let $R$ be a symmetric CGL extension of length $N$
and $a \in \Zset$ be such that $|\eta^{-1}(a)|>1$. Denote
$$
\eta^{-1}(a) = \{l, s(l), \ldots, s^m(l) \} 
$$
where $l \in [1,N]$ and $m = O_+(l) \in \Zset_{> 0}$. Then 
\begin{equation}
\label{la-eq}
\la^*_l = \la^*_{s(l)} = \ldots = \la^*_{s^{m-1}(l)} = 
\la^{-1}_{s(l)} = \la^{-1}_{s^2(l)}= \ldots = \la^{-1}_{s^m(l)}.
\end{equation}
\epr

\begin{proof} It follows from \eqref{lala*} that 
$$
\la^\sy_{s^{m_1}(l)} = \la_{s^{m_2}(l)}^{-1}, \quad \forall 0 \leq m_1 < m_2 \leq m
$$
which is equivalent to the statement of the proposition. 
\end{proof}

\section{Sequences of homogeneous prime elements} 
\label{4a.4}
We proceed with the proofs of the two theorems formulated in Section \ref{4a.1}.

\begin{proof}[Proof of \thref{y-int}] As already noted, the case $m=0$ is easily verified. Assume now that $m > 0$.

(a) Consider the CGL extension 
presentation of $R$ associated to 
\begin{multline}
\label{y-int1}
\tau_{(s^m(i)-1, s^m(i), \ldots, s^m(i))} = \\ 
[i +1 , \ldots, s^m(i)-1, i, s^m(i), i-1, \ldots, 1, s^m(i)+1, \ldots, N] \in \Xi_N,
\end{multline}
where $s^m(i)$ is repeated $i-1$ times.
The $(s^m(i)-i)$-th and $(s^m(i)-i+1)$-st algebras in the chain 
are precisely $R_{[i, s^m(i)-1]}$ and $R_{[i,s^m(i)]}$.
\thref{CGL} implies that the homogeneous prime elements of $R_{[i,s^m(i)]}$ 
that do not belong to $R_{[i, s^m(i)-1]}$ are associates of each other.
Denote by $z_{[i, s^m(i)]}$ one such element. Set $z_\varnothing = 1$.
(One can prove part (a) using 
the simpler presentation of $R$ associated to the permutation 
$$[i, i + 1 , \ldots, s^m(i), 
i-1, \ldots, 1, s^m(i)+1, \ldots, N] \in \Xi_N$$
but 
the first presentation will also play a role in the proof of part (b).) The equivariance 
property of the $\eta$-function from \coref{steps} (b) and \thref{CGL} imply that 
\begin{equation}
\label{y-int2}
z_{[i,s^m(i)]} = \xi_{i;m}( z_{[i,s^{m-1}(i)]} x_{s^m(i)} - c_{[i,s^m(i)-1]})
\end{equation}
for some $c_{[i,s^m(i)-1]} \in R_{[i,s^m(i)-1]}$ and $\xi_{i;m} \in \kx$.

(b) Now consider the CGL extension 
presentation of $R$ associated to 
\begin{equation}
\label{y-int3}
\tau_{(s^m(i), \ldots, s^m(i))} =
[i +1 , \ldots, s^m(i), i, \ldots, 1, s^m(i)+1, \ldots, N] \in \Xi_N,
\end{equation}
where $s^m(i)$ is repeated $i$ times.
The $(s^m(i)-i)$-th and $(s^m(i)-i+1)$-st algebras in the chain 
are $R_{[i+1, s^m(i)]}$ and $R_{[i,s^m(i)]}$.
\thref{CGL} implies that the homogeneous prime elements of $R_{[i,s^m(i)]}$ 
that do not belong to $R_{[i+1, s^m(i)]}$ are associates of each other.
They are associates of $z_{[i,s^m(i)]}$. This follows from 
\thref{1} (b) applied to the CGL extension presentations of $R$ associated to 
the elements \eqref{y-int1} and \eqref{y-int3}. The fact that 
these two CGL extension presentations satisfy the assumptions of 
\thref{1} (b) follows from  
the equivariance of the $\eta$-function from \coref{steps} (b).
This equivariance and \thref{CGL} applied to the presentation for \eqref{y-int3} 
also imply that 
\begin{equation}
\label{y-int4}
z_{[i,s^m(i)]} = \xi'_{i;m}( x_{i} z_{[s(i),s^m(i)]} - c'_{[i+1,s^m(i)]})
\end{equation}
for some $c'_{[i+1,s^m(i)]} \in R_{[i+1,s^m(i)]}$ and $\xi'_{i;m} \in \kx$.

The rest of part (b) and parts (c)--(d) follow at once by comparing the leading terms in 
Eqs. \eqref{y-int2} and \eqref{y-int4}, and the fact that the group of units 
of an Ore extension is reduced to scalars.

For part (e), we apply \eqref{yxcomm} in $R_{[\min\{k,i\},\max\{k,s^m(i)\}]}$ with $j = s^m(i)$.
\end{proof}

\begin{proof}[Proof of  \thref{sym-prime}] For $k \in [0,N]$, 
denote by $R_{\tau, k}$ the $k$-th algebra in the chain
\eqref{tauOre}. We consider the case when $\tau(k) \geq \tau(1)$, 
leaving the analogous case $\tau(k) \leq \tau(1)$ to the 
reader. Since $\tau([1,j])$ is an interval for all $j \leq k$,
$$
\tau([1, k]) = [\tau(i), \tau(k)] \quad
\mbox{for some} \; \; i \in [1,k].
$$
Therefore $R_{\tau, k} = R_{[\tau(i), \tau(k)]}$ and 
$R_{\tau, k-1} = R_{[\tau(i), \tau(k)-1]}$.
For $m \in \Zset_{\geq 0}$ given by \eqref{m-p} 
we have 
$$
\tau(i) \leq p^m(\tau(k)) \leq \tau(k).
$$
\thref{CGL} implies that $y_{[p^m(\tau(k)), \tau(k)]}$
is a homogeneous prime element of the algebra $R_{[\tau(i),\tau(k)]} = R_{\tau, k}$. 
It does not belong to $R_{\tau, k-1} = R_{[\tau(i), \tau(k)-1]}$
because of \thref{y-int} (c). It follows from \thref{CGL} that 
the homogeneous prime elements of $R_{\tau, k}$ that do not belong to $R_{\tau, k-1}$
are associates of  $y_{\tau, k}$.  
Hence, $y_{\tau, k}$ is a scalar multiple of 
$y_{[p^m(\tau(k)), \tau(k)]}$.
\end{proof}

\bco{int-commut} If $R$ is a symmetric CGL extension of length 
$N$ and $i, j \in [1,N]$, $m, n \in \Znn$ are such that
$$
i \leq j \leq s^n(j) \leq s^m(i) \leq N,
$$
then 
$$
y_{[i,s^m(i)]} y_{[j, s^n(j)]} = 
\Om_\lab(e_{[i,s^m(i)]}, e_{[j,s^n(j)]})
y_{[j, s^n(j)]} y_{[i,s^m(i)]},
$$
recall \eqref{La}.
\eco

The corollary follows by applying \thref{sym-prime} to the CGL extension 
presentation of $R$ associated to 
$$
\tau':= [j, \ldots, s^m(i),j-1, \ldots, 1, s^m(i) +1, \ldots, N] \in \Xi_N 
$$
and then using \eqref{ycomm}, \eqref{q}. For this permutation, \thref{sym-prime} 
implies that $y_{\tau', s^n(j)-j+1}$ and $y_{\tau', s^m(i)-i+1}$ are 
scalar multiples of $y_{[j, s^n(j)]}$ and $y_{[i,s^m(i)]}$, 
respectively.

\section{An identity for normal elements}
\label{4a.5}
For $1 \le j \le l \le N$, set 
\begin{equation}
\label{genO}
O_-^j(l) = \max \{ m \in \Zset_{\geq 0} \mid p^m(l) \geq j \}.
\end{equation}
\index{Ojminus@$O_-^j$}
The following fact follows directly from Theorems \ref{tCGL} and \ref{tsym-prime}. 

\bco{primeRjk} Let $R$ be a symmetric CGL extension 
of length $N$ and $1 \leq j \leq k \leq N$. 
Then
$$
\Big\{ y_{[p^{O_-^j(i)}(i), i]} \Bigm| i \in [j,k] , \; s(i)>k \Big\}
$$
is a list of the homogeneous prime elements of $R_{[j,k]}$ up to scalar multiples.
\eco

Applying Theorems \ref{t1}, \ref{t2} to the CGL extension presentations \eqref{tauOre} of a 
symmetric CGL extension of length $N$ associated to the elements 
$$
\tau_{i, s^m(i)-1} = [i+1, \ldots, s^{m}(i)-1, i, s^m(i), \ldots, N, i-1, \ldots, 1] \in \Xi_N
$$
and 
$$
\tau_{i, s^m(i)} = [i+1, \ldots, s^{m}(i)-1, s^m(i), i, s^m(i)+1, \ldots, N, i-1, \ldots, 1] \in \Xi_N
$$ 
for $i \in [1,N]$ and $m \in \Zset_{>0}$ such that $s^m(i) \in [1,N]$
and using \thref{sym-prime}
yields the following result.

\bco{u-elem} Assume that $R$ is a symmetric CGL extension of length $N$,
and $i \in [1,N]$ and $m \in \Zset_{>0}$ are such that $s^m(i) \in [1,N]$.
Then
\begin{align}
u_{[i,s^m(i)]} &:= y_{[i, s^{m-1}(i)]} y_{[s(i), s^m(i)]} \,-   \label{uuu}  \\
&\qquad\qquad\qquad  \Om_\lab(e_i, e_{[s(i), s^{m-1}(i)]}) 
y_{[s(i), s^{m-1}(i)]} y_{[i,s^m(i)]}
\nn  \\
&= y_{[i, s^{m-1}(i)]} y_{[s(i), s^m(i)]} \,-  \label{uuuu}  \\
&\qquad\qquad\qquad \Om_\lab(e_{s^m(i)}, e_{[s(i), s^{m-1}(i)]})^{-1}
y_{[i,s^m(i)]}  y_{[s(i), s^{m-1}(i)]}  \nn 
\end{align}
\index{uismi@$u_{[i,s^m(i)]}$}
is a nonzero homogeneous normal element of $R_{[i+1,s^m(i)-1]}$ 
which is not a multiple of $y_{[s(i), s^{m-1}(i)]}$ if $m \geq 2$, 
recall \eqref{La}. It normalizes the elements of $R_{[i+1,s^m(i)-1]}$
in exactly the same way as $y_{[s(i), s^{m-1}(i)]} y_{[i,s^m(i)]}$ does.  Moreover,
\begin{equation}  \label{uformula}
u_{[i,s^m(i)]} = \psi \prod_{k\in P} y^{m_k}_{[p^{O_-^{i+1}(k)}(k),k]}
\end{equation}
where $\psi \in \kx$,
\begin{equation}  \label{Pset}
P = P_{[i,s^m(i)]} := \{ k\in [i,s^m(i)] \setminus \{i,s(i),\dots,s^m(i)\} \mid s(k) > s^m(i) \},
\end{equation}
and the integers $m_k$ are those from Theorems {\rm\ref{t2}, \ref{t3}}. Consequently, the leading term of $u_{[i,s^m(i)]}$ 
has the form 
$$
\xi x^{n_{i+1} e_{i+1} + \cdots + n_{s^m(i)-1} e_{s^m(i)-1}}
$$
for some $\xi \in \kx$ and $n_{i+1}, \ldots, n_{s^m(i)-1} \in \Znn$ 
such that $n_{s(i)}= \ldots = n_{s^{m-1}(i)} =0$
and $n_j = n_l$ for all $i+1 \leq j \leq l \leq s^m(i)-1$ with 
$\eta(j) = \eta(l)$.
\eco

The equality \eqref{uuuu} follows from \coref{int-commut}.
The scalar in the right hand side of Eq. \eqref{uuu} is 
the one that cancels out the leading terms of the 
two products, cf. \thref{y-int} (c).  
The last statement of \coref{u-elem} follows from Eqs. \eqref{uformula} and \eqref{prod-lt}.

Since $y_\varnothing=1$ 
and $y_{[i,i]}=x_i$, the $m=1$ case of the corollary states that
\begin{equation}
\label{ucc'}
u_{[i,s(i)]}:= x_i x_{s(i)} - y_{[i,s(i)]}= c_{[i,s(i)-1]}= c'_{[i+1, s(i)]}
\end{equation}
is a nonzero normal element of $R_{[i+1,s(i)-1]}$, cf. \thref{y-int} (d).
We set 
$$
u_{[i,i]} :=1.
$$

\bex{OqMmn7}
In the case of $R= \OqM$, \exref{OqMmn2} shows how to express each $u_{[i, s^l(i)]}$ as a difference of products of solid quantum minors. In particular, if $i = (r-1)n+c$ with $s(i) \ne +\infty$, then $r < m$, $c < n$, and
\begin{equation}
\label{OqMuisi}
\begin{aligned}
u_{[i,s(i)]} &:= x_i x_{s(i)} - y_{[i,s(i)]}= t_{rc} t_{r+1, c+1} - \Delta_{[r, r+1], [c, c+1]}  \\
 &= q t_{r, c+1} t_{r+1, c} = q x_{i+1} x_{i+n} .
\end{aligned}
\end{equation}
More generally, if $i = (r-1)n+c$, $l \in \Znn$, and $s^l(i) \ne +\infty$, then
$$
u_{[i, s^l(i)]} = q \De_{[r, r+l-1], [c+1, c+l]} \De_{[r+1, r+l], [c, c+l-1]} = q y_{[i+1, s^{l-1}(i+1)]} y_{[i+n, s^{l-1}(i+n]} .
$$
This follows from an identity for quantum minors:
\begin{multline*}
\De_{[r, r+l-1], [c, c+l-1]} \De_{[r+1, r+l], [c+1, c+l]} - q \De_{[r, r+l-1], [c+1, c+l]} \De_{[r+1, r+l], [c, c+l-1]}  \\
= \De_{[r, r+l], [c, c+l]} \De_{[r+1, r+l-1], [c+1, c+l-1]}, 
\end{multline*}
for all $r \in [1, m-l]$, $c \in [1, n-l]$, $l \in [2, \min(m,n)]$. This identity can be obtained from the equation
$$
t_{rc} t_{r+l, c+l} - q t_{r, c+l} t_{r+l, c} = \De_{ \{r, r+l\}, \{c, c+l \}}
$$
by applying the antipode in a copy of $\OO_q(GL_{l+1}(\KK))$. (See \cite[Lemma 4.1]{KLR} for how the antipode applies to quantum minors.)
\eex

The next result, which is the main result in this section, 
describes the relationship between the normal elements in 
\coref{u-elem}. It will play a key role in normalizing 
the scalars $\zeta$ in \thref{3}.

\bth{u-prod} Let $R$ be a symmetric CGL extension of length $N$.
For all $i \in [1,N]$ and $m \in \Zset_{>0}$ 
such that $s^{m+1}(i) \in [1,N]$,
\begin{equation}
\label{u-ident}
\lt( u_{[s(i),s^m(i)]} u_{[i,s^{m+1}(i)]})  = 
\theta_m \lt( u_{[i,s^m(i)]} u_{[s(i),s^{m+1}(i)]} )
\end{equation}
where $\theta_m = (\la^\sy_i)^{-1} \Om_\lab(e_{[s^2(i), s^m(i)]}, 2 e_i + 2 e_{s(i)})
\Om_\lab(e_{s^m(i)}, e_i)^{-1} \Om_\lab(e_{s(i)}, e_i)$, 
recall Section {\rm\ref{detailCGL}}
for the definition of leading terms $\lt(-)$.
\eth

For $m=1$, Eq. \eqref{u-ident} simplifies to 
\begin{equation}
\label{u-ident1}
\lt(u_{[i,s^2(i)]}) = (\la_i^*)^{-1} \lt( u_{[i,s(i)]} u_{[s(i), s^2(i)]}).
\end{equation}

\begin{proof} For a subspace $W$ of a $\KK$-vector space $V$ and 
$v_1, v_2 \in V$, we will write 
$$
v_1 \equiv v_2 \mod W \quad
\mbox{iff}
\quad
v_1 - v_2 \in W.
$$
Denote 
$$
W_m := \big( \oplus_{n_1,n_2 \in \Znn} x_i^{n_1} y_{[s(i), s^m(i)]} 
R_{[i+1, s^{m+1}(i)-1]} x_{s^{m+1}(i)}^{n_2} \big) \subset R_{[i,s^{m+1}(i)]}.
$$
By \cite[Proposition 4.7(b)]{GY} applied to $R_{[i+1,s^{m+1}(i)]}$,
$$
\delta_{s^{m+1}(i)} (y_{[s(i),s^m(i)]}) = \Om_\nu(e_{s^{m+1}(i)}, e_{[s(i),s^m(i)]})
(\la_{s^{m+1}(i)}-1) c_{[s(i),s^{m+1}(i)-1]},
$$
cf. also \thref{y-int} (d). Thus,
\begin{align*}
&y_{[i,s^{m+1}(i)]}  y_{[s(i), s^m(i)]} \equiv
y_{[i,s^m(i)]} x_{s^{m+1}(i)}  y_{[s(i), s^m(i)]}
\equiv 
y_{[i,s^m(i)]} \delta_{s^{m+1}(i)}(y_{[s(i), s^m(i)]})
\\
&\qquad=
\Om_\nu(e_{s^{m+1}(i)}, e_{[s(i),s^m(i)]})
(\la_{s^{m+1}(i)}-1) y_{[i,s^m(i)]} c_{[s(i),s^{m+1}(i)-1]}
\\
& \qquad\equiv
\Om_\nu(e_{s^{m+1}(i)}, e_{[s(i),s^m(i)]})
(1-\la_{s^{m+1}(i)}) y_{[i,s^m(i)]} y_{[s(i),s^{m+1}(i)]}
\quad
\mod W_m.
\end{align*}
It follows from \coref{u-elem} and \prref{la-equal} that
\begin{equation}
\label{in-equi}
\begin{aligned}
u_{[i,s^{m+1}(i)]} &=y_{[i, s^m(i)]} y_{[s(i), s^{m+1}(i)]} \,-  \\
&\qquad\qquad\qquad \Om_\lab(e_{s^{m+1}(i)}, e_{[s(i), s^m(i)]})^{-1}
y_{[i,s^{m+1}(i)]}  y_{[s(i), s^m(i)]}   \\
&\equiv 
(\la^\sy_i)^{-1}y_{[i, s^m(i)]} y_{[s(i), s^{m+1}(i)]} \quad \mod W_m.
\end{aligned}
\end{equation}

First, we consider the case $m>1$. 
Using Eq. \eqref{in-equi}, \coref{int-commut} and the fact that 
$y_{[s(i), s^m(i)]}$ is a prime element of $R_{[i+1, s^{m+1}(i)-1]}$, 
which follows from \coref{primeRjk},
we obtain 
\begin{align*}
&u_{[s(i),s^m(i)]} u_{[i, s^{m+1}(i)]}  
\\
&\qquad\equiv 
\big( 
y_{[s(i),s^{m-1}(i)]} y_{[s^2(i), s^m(i)]} - \Om_\lab(e_{s(i)}, e_{[s^2(i), s^{m-1}(i)]})
y_{[s^2(i),s^{m-1}(i)]} y_{[s(i), s^m(i)]} \big)
\\
&\qquad\qquad\times 
(\la^\sy_i)^{-1} y_{[i,s^m(i)]} y_{[s(i), s^{m+1}(i)]} 
\\
&\qquad\equiv 
(\la^\sy_i)^{-1}
y_{[s(i),s^{m-1}(i)]} y_{[s^2(i), s^m(i)]}
y_{[i,s^m(i)]} y_{[s(i), s^{m+1}(i)]} 
\\
&\qquad = (\la^\sy_i)^{-1}\Om_\lab(e_{[s^2(i), s^m(i)]}, e_{[i, s^m(i)]})
y_{[s(i),s^{m-1}(i)]} y_{[i,s^m(i)]} 
y_{[s^2(i), s^m(i)]} y_{[s(i), s^{m+1}(i)]} 
\\ 
&\qquad\equiv \theta_m u_{[i,s^m(i)]} u_{[s(i),s^{m+1}(i)]}
\quad \mod W_m.
\end{align*}
Since the first and last products above belong to $R_{[i+1, s^{m+1}(i)-1]}$,
\begin{equation}
\label{3terms}
u_{[s(i), s^m(i)]} u_{[i,s^{m+1}(i)]}  = \theta_m u_{[i,s^m(i)]} u_{[s(i),s^{m+1}(i)]}
+ y_{[s(i), s^m(i)]} r 
\end{equation}
for some $r \in R_{[i+1, s^{m+1}(i)-1]}$. The leading term of 
$y_{[s(i), s^m(i)]} r$ 
has the form $\xi x^f$ where $\xi \in \kx$ 
and $f = n_{i+1} e_{i+1} + \cdots + n_{s^{m+1}(i) -1} e_{s^{m+1}(i) -1}$ are 
such that $n_{s(i)}, \ldots, n_{s^m(i)} \geq 1$. Combining this with 
the last statement of \coref{u-elem} implies that the exponent of the leading term 
of the third product in \eqref{3terms} is different from the 
exponents of the leading terms of the first two products. Now the 
the validity of \eqref{u-ident} for $m>1$ follows from \eqref{3terms}.

To verify \eqref{u-ident1}, we use Eq. \eqref{in-equi} and 
the fact that $x_{s(i)}$ is a prime element of  
$R_{[i+1, s^2(i)-1]}$ to obtain 
$$
u_{[i, s^2(i)]}
\equiv (\la^\sy_i)^{-1} y_{[i,s(i)]} y_{[s(i), s^2(i)]} \equiv (\la^\sy_i)^{-1} u_{[i,s(i)]} u_{[s(i), s^2(i)]}
\quad \mod W_1.
$$
Hence, 
$$
u_{[i, s^2(i)]} = (\la^\sy_i)^{-1} u_{[i,s(i)]} u_{[s(i), s^2(i)]} + x_{s(i)} r
$$
for some $r \in R_{[i+1, s^2(i)-1]}$. Now \eqref{u-ident1} can be deduced 
analogously to the case of $m>1$.
\end{proof}

\chapter{Chains of mutations in symmetric CGL extensions}
\label{mut-sym}

For a given CGL extension $R$ one has the freedom of 
rescaling the generators $x_j$ by elements of the base field 
$\KK$. The prime elements $y_j$ and $\ol{y}_j$ from \thref{CGL} 
and Eq. \eqref{new-y} obviously depend (again up to rescaling) 
on the choice of $x_j$. 
In this chapter we prove that for each symmetric CGL extension $R$ 
its generators can be rescaled in such a way that all
scalars $\zeta$ in \thref{3} become equal to one. 
This implies a mutation theorem, proved in Chapter \ref{main}, 
for toric frames for $\Fract(R)$ associated 
to the elements of $\Xi_N$ via the sequences of prime elements 
from \thref{sym-prime}.

\section{The leading coefficients of $u_{[i,s^m(i)]}$}
\label{4b.1}
Throughout this chapter we will assume that $R$ is a symmetric CGL extension of length $N$ 
and that there exist square roots $\nu_{lj} = \sqrt{\la_{lj}} \in \KK$
such that the condition \eqref{-1} is satisfied. 
Recall \eqref{uuu}.
For $i \in [1,N]$ and $m \in \Zset_{\ge0}$ such that $s^m(i) \in [1,N]$, 
let
$$
\pi_{[i,s^m(i)]} \in \kx \quad \mbox{and} 
\quad
f_{[i, s^m(i)]} \in \sum_{j =i+1}^{s^m(i)-1} \Zset_{\geq 0}\, e_j
\subset \Zset^N
$$
\index{pismi@$\pi_{[i,s^m(i)]}$}  \index{fismi@$f_{[i, s^m(i)]}$}
be given by
\begin{equation}
\label{pi-f}
\lt(u_{[i,s^m(i)]})= \pi_{[i,s^m(i)]} x^{f_{[i,s^m(i)]}}.
\end{equation}
Note that $\pi_{[i,i]} = 1$ and $f_{[i,i]} = 0$, because $u_{[i,i]} = 1$.

A quick induction using \eqref{u-ident} shows that
\begin{equation}
\label{f-sum}
f_{[i,s^m(i)]} = f_{[i, s(i)]} + \cdots + f_{[s^{m-1}(i), s^m(i)]}
\end{equation}
for $m > 0$, and thus
\begin{equation}
\label{xf-prod}
x^{f_{[i, s^m(i)]}} = x^{f_{[i,s(i)]}} x^{f_{[s(i),s^m(i)]}} = x^{f_{[i,s^{m-1}(i)]}} x^{f_{[s^{m-1}(i), s^m(i)]}}.
\end{equation}
Since all terms in \eqref{uuu} are homogeneous 
of the same $\xh$-degree,
\begin{align}
\label{xf-deg}
\xh\mbox{-deg} ( x^{f_{[i, s^m(i)]}}) &= 
\xh\mbox{-deg}(u_{[i, s^m(i)]}) =\xh\mbox{-deg}(y_{[s(i), s^{m-1}(i)]} y_{[i,s^m(i)]})  \\
 &= \xh\mbox{-deg}(x^{e_i + 2 e_{[s(i), s^{m-1}(i)]} + e_{s^m(i)}})
\nn
\end{align}
for $m > 0$. The next result describes recursively the leading coefficients $\pi_{[i,s^m(i)]}$. Recall \eqref{Scr}.

\bpr{pi} Let $R$ be a symmetric CGL extension for which there exist
$\nu_{kj} = \sqrt{\la_{kj}} \in \kx$, $0 \leq j < k \leq N$ satisfying
condition \eqref{-1}. Then the following hold:

{\rm(a)} The scalars $\pi_{[i,s^m(i)]}$ satisfy the recursive relation
$$
\pi_{[s(i), s^m(i)]} \pi_{[i, s^{m+1}(i)]} =
(\la_i^\sy)^{-\delta_{m,1}} \Om_{\lab}(e_{s^m(i)}, e_{s(i)}) 
\pi_{[i,s^m(i)]} \pi_{[s(i), s^{m+1}(i)]}
$$
for all $i \in [1,N]$ and $m \in \Zset_{> 0}$ such that $s^{m+1}(i) \in [1,N]$.

{\rm(b)} If 
$$
\pi_{[i,s(i)]} = \Scr_\nub( - e_i + f_{[i,s(i)]} )
$$
for all $i \in [1,N]$ such that $s(i) \neq \infty$, then 
$$
\pi_{[i,s^m(i)]} = 
\Scr_\nub( e_{[s(i),s^m(i)]})^{-2} \Scr_\nub( - e_i + f_{[i, s^m(i)]})
$$
for $i \in [1,N]$ and $m \in \Zset_{\geq 0 }$ with $s^m(i) \in [1,N]$.
\epr

Before we proceed with the proof of \prref{pi} we derive several important properties 
of the form $\Om_\lab$.

\ble{pi-le} Let $R$ be a symmetric CGL extension of rank $N$
and $i \in [1,N]$, $m \in \Zset_{> 0}$ be such that $s^m(i) \in [1,N]$.
Then the following hold:

{\rm(a)} For all  $0 \leq k < n \leq m$,
$$
\Om_\lab (e_{[i,s^m(i)]}, f_{[s^k(i), s^n(i)]}) = 
\Om_\lab (e_{[i,s^m(i)]}, e_{s^k(i)} + 2 e_{[s^{k+1}(i), s^{n-1}(i)]} + e_{s^n(i)}). 
$$

{\rm(b)} If $j \in [i,s^m(i)] \backslash \{ i , s(i), \ldots, s^m(i)\}$ and $k \in \Zset_{\geq 0}$ 
are such that $p(j) < i$ and $s^{k+1}(j) > s^m(i) > s^k(j)$, then 
$$
\Om_\lab (f_{[i,s^m(i)]}, e_{[j,s^k(j)]} ) =
\Om_\lab (e_i + 2 e_{[s(i), s^{m-1}(i)]} + e_{s^m(i)}, e_{[j,s^k(j)]}). 
$$

{\rm(c)} For all  $0 \leq k < n \leq m$ and $j \in \{i,s(i),\dots,s^m(i)\} \cup [1,s^k(i)] \cup [s^n(i),N]$,
$$
\Om_\lab (e_j, f_{[s^k(i),s^n(i)]}) = (\la_i^*)^{\delta_{j,s^k(i)} - \delta_{j,s^n(i)}} \Om_\lab (e_j, e_{s^k(i)}+ 2e_{[s^{k+1}(i),s^{n-1}(i)]} + e_{s^n(i)}).
$$

{\rm(d)} For all  $0 \leq k < n \leq m$,
\begin{multline*}
\Om_\lab(f_{[i, s^m(i)]}, f_{[s^k(i), s^n(i)]}) 
\\
= (\la_i^\sy)^{- \delta_{k,0} + \delta_{n,m}}
\Om_\lab(e_i + 2 e_{[s(i), s^{m-1}(i)]}+ e_{s^m(i)}, e_{s^k(i)}+ 2 
e_{[s^{k+1}(i), s^{n-1}(i)]}+ e_{s^n(i)} ).  
\end{multline*}
\ele

\begin{proof} (a) By \cite[Theorem 5.3]{GY} $y_{[i,s^m(i)]}$ is an $\HH$-normal element 
of $R_{[i, s^m(i)]}$ when $\HH$ is replaced with the universal maximal torus 
from \thref{univ-torus}. This means that there exists $h \in \HH$ such that 
$y_{[i,s^m(i)]} u = (h . u) y_{[i,s^m(i)]}$ for all
$u \in R_{[i, s^m(i)]}$. In particular, $y_{[i,s^m(i)]} u = \chi_u(h) u y_{[i,s^m(i)]}$ when $u$ is homogeneous. Taking account of \thref{y-int}(e) and the result of \eqref{xf-deg} applied to the interval $[s^k(i),s^n(i)]$, we obtain the statement of (a).

Part (b) is obtained just as in (a), since under the given conditions, $y_{[j,s^k(j)]}$ is an $\Hmax$-normal element of $R_{[i,s^m(i)]}$.

(c) In view of \eqref{f-sum}, it suffices to prove that
\begin{equation}  \label{62c1}
\Om_\lab (e_j, f_{[s^{l-1}(i),s^l(i)]}) = (\la_i^*)^{\delta_{j,s^{l-1}(i)} - \delta_{j,s^l(i)}} \Om_\lab (e_j, e_{s^{l-1}(i)} +e_{s^l(i)})
\end{equation}
for all $l\in [1,m]$ and $j \in [1,s^{l-1}(i)] \cup [s^l(i),N]$. Fix such $l$ and $j$, let $\chi$ denote the $\HH$-eigenvalue of $x^{f_{[s^{l-1}(i),s^l(i)]}}$, and recall from \eqref{xf-deg} that $\chi$ equals the $\HH$-eigenvalue of $x_{s^{l-1}(i)} x_{s^l(i)}$.

Assume first that $j \ge s^l(i)$, and recall that
$$
\chi_{x_{i'}}(h_j) = \begin{cases} \la_{ji'} &\quad(1 \le i'<j)\\  \la_j &\quad(i'=j). \end{cases}
$$
Since $x^{f_{[s^{l-1}(i),s^l(i)]}} \in R_{[1,j-1]}$, it follows that
\begin{equation}  \label{62c2}
\chi(h_j) = \Om_\lab (e_j, f_{[s^{l-1}(i),s^l(i)]}) .
\end{equation}
Since $\chi$ is also the $\HH$-eigenvalue of $x_{s^{l-1}(i)} x_{s^l(i)}$, we have
\begin{equation}  \label{62c3}
\chi(h_j) = \begin{cases} \Om_\lab (e_j, e_{s^{l-1}(i)} + e_{s^l(i)}) &\quad(j > s^l(i))\\  \Om_\lab (e_j, e_{s^{l-1}(i)}) \la_{s^l(i)} &\quad(j= s^l(i))  \end{cases}
\end{equation}
as well. In the second case of \eqref{62c3}, we invoke \eqref{la-eq} and the skewsymmetry of $\Om_\lab$ to obtain
\begin{equation}  \label{62c4}
\chi(h_j) = (\la_i^*)^{-1} \Om_\lab (e_j, e_{s^{l-1}(i)} + e_{s^l(i)} ) \qquad(j = s^l(i)).
\end{equation}
Combining \eqref{62c2}--\eqref{62c4} yields \eqref{62c1} for $j \in [s^l(i),N]$.

We obtain \eqref{62c1} for $j \in [1,s^{l-1}(i)]$ in the same manner, by working with $\chi(h_j^*)$.

(d) The case $k=0$, $n=m$ is obvious since $\Om_\lab$ is 
multiplicatively skew-symmetric. 

In general, since $\Om_\lab (f_{[s^k(i),s^n(i)]}, f_{[s^k(i),s^n(i)]}) = 1$ and $f_{[i,s^k(i)]}  + f_{[s^n(i),s^m(i)]}$ is a $\Zset_{\ge0}$-linear combination of $e_j\,$s with $j < s^k(i)$ or $j > s^n(i)$, we see from \eqref{f-sum} and part (c) that
\begin{multline}  \label{62d1}
\Om_\lab (f_{[i,s^m(i)]},\, f_{[s^k(i),s^n(i)]}) = \Om_\lab (f_{[i,s^k(i)]}  + f_{[s^n(i),s^m(i)]}, f_{s^k(i),s^n(i)]})  \\
 = \Om_\lab (f_{[i,s^k(i)]}  + f_{[s^n(i),s^m(i)]}, e_{s^k(i)}+ 2 
e_{[s^{k+1}(i), s^{n-1}(i)]}+ e_{s^n(i)} ). 
\end{multline}
We now apply part (c) again, with $j= s^k(i), \dots, s^n(i)$. If $0 < k < n < m$, then \eqref{62d1} yields
\begin{multline*}
\Om_\lab (f_{[i,s^m(i)]}, f_{[s^k(i),s^n(i)]}) \\
 = \Om_\lab( e_i+ 2e_{[s(i),s^{k-1}(i)]}+ e_{s^k(i)} + e_{s^n(i)}+ 2e_{[s^{n+1}(i),  s^{m-1}(i)]}+ e_{s^m(i)},  \\
 e_{s^k(i)}+ 2 
e_{[s^{k+1}(i), s^{n-1}(i)]}+ e_{s^n(i)} ),
\end{multline*}
and the desired result follows because
$$
\Om_\lab (e_{s^k(i)}+ 2 
e_{[s^{k+1}(i), s^{n-1}(i)]}+ e_{s^n(i)} , e_{s^k(i)}+ 2 
e_{[s^{k+1}(i), s^{n-1}(i)]}+ e_{s^n(i)} ) = 1.
$$

If $k=0$ and $n<m$, then $f_{[i,s^k(i)]} = 0$ and \eqref{62d1} yields
\begin{multline*}
\Om_\lab (f_{[i,s^m(i)]}, f_{[s^k(i),s^n(i)]})  \\
= (\la_i^*)^{-1} \Om_\lab (e_{s^n(i)}+ 2 
e_{[s^{n+1}(i), s^{m-1}(i)]}+ e_{s^m(i)} , e_i+ 2 
e_{[s(i), s^{n-1}(i)]}+ e_{s^n(i)} ),
\end{multline*}
from which the desired result follows. The case $k>0$, $n=m$ is covered in the same fashion.
\end{proof}

Note from \leref{pi-le}(c) that
$$
\Om_\nub( e_i, f_{[i,s(i)]} - e_{s(i)})^2 = \Om_\lab( e_i, f_{[i,s(i)]} - e_{s(i)}) = \la_i^*, \quad \text{when\;} i, s(i) \in [1,N],
$$
and thus $\la_i^*$ has a square root in the subgroup $\langle \nu_{kj} \mid k,j \in [1,N] \rangle$ of $\kx$. All values of $\Scr_\nub$, $\Om_\nub$, and $\Om_\lab$ lie in this group.

\begin{proof}[Proof of \prref{pi}] (a) \thref{u-prod} and Eq. \eqref{xf-prod} imply
\begin{multline*}
\pi_{[s(i), s^m(i)]} \pi_{[i, s^{m+1}(i)]} 
\lt( x^{f_{[s(i),s^m(i)]}} x^{f_{[i,s(i)]}} x^{f_{[s(i),s^{m+1}(i)]}} )
\\
= \theta_m
\pi_{[i,s^m(i)]} \pi_{[s(i), s^{m+1}(i)]}
\lt( x^{f_{[i,s(i)]}} x^{f_{[s(i),s^m(i)]}} x^{f_{[s(i),s^{m+1}(i)]}})
\end{multline*}
for the scalar $\theta_m \in \kx$ from \thref{u-prod}.
It follows from Eq. \eqref{prod-lt}, \leref{pi-le} (d) and Eq. \eqref{f-sum} that
\begin{align*}
\frac{\pi_{[s(i), s^m(i)]} \pi_{[i, s^{m+1}(i)]}}
{\pi_{[i,s^m(i)]} \pi_{[s(i), s^{m+1}(i)]}} &=
\theta_m \Om_\lab( f_{[i,s(i)]}, f_{[s(i),s^m(i)]} )
= \theta_m \Om_\lab( f_{[i,s(i)]}, f_{[i,s^m(i)]} )
\\
&=\theta_m (\la_i^\sy)^{1-\delta_{m,1}} \Om_\lab( e_i + e_{s(i)}, e_i + 2 e_{[s(i),s^{m-1}(i)]} + e_{s^m(i)} ).
\end{align*}
Simplifying the last expression using the definition of $\theta_m$ 
from \thref{u-prod} proves part (a). 

(b) Denote
$$
\vp_{[i,s^m(i)]} := \Scr_\nub( e_{[s(i),s^m(i)]})^{-2} \Scr_\nub( - e_i + f_{[i, s^m(i)]}) .
$$
We will prove that
\begin{equation}
\label{vph}
\vp_{[s(i), s^m(i)]} \vp_{[i, s^{m+1}(i)]} =
(\la_i^\sy)^{-\delta_{m,1}} \Om_{\lab}(e_{s^m(i)}, e_{s(i)}) 
\vp_{[i,s^m(i)]} \vp_{[s(i), s^{m+1}(i)]}
\end{equation}
for all $i \in [1,N]$ and $m \in \Zset_{> 0}$ such that $s^{m+1}(i) \in [1,N]$.
Part (b) trivially holds for $m=0,1$ and
follows from the identity \eqref{vph} and part (a) by induction on $m$.

For all $j \in [2, N]$ and $f \in \Zset e_1 + \cdots + \Zset e_j$, 
$g \in \Zset e_{j+1} + \cdots + \Zset e_N$,
\begin{equation}
\label{fg}
\Scr_\nub(f+g) = \Om_\nub(f,g)^{-1} \Scr_\nub(f) \Scr_\nub(g).
\end{equation}
Hence, taking also \eqref{f-sum} into account,
$$
\frac{\vp_{[i, s^{m+1}(i)]}}{\vp_{[i,s^m(i)]}} = 
\frac{ \Om_\nub(e_{[s(i),s^m(i)]}, e_{s^{m+1}(i)})^2 \Scr_\nub(f_{[s^m(i),s^{m+1}(i)]}) } { \Om_\nub(-e_i + f_{[i, s^{m+1}(i)]}, f_{[s^m(i), s^{m+1}(i)]}) } .
$$
Similarly,
$$
\frac{ \vp_{[s(i), s^m(i)]} }{ \vp_{[s(i),s^{m+1}(i)]} } = 
\frac{ \Om_\nub(-e_{s(i)} + f_{[s(i), s^{m+1}(i)]}, f_{[s^m(i), s^{m+1}(i)]}) } 
{ \Om_\nub(e_{[s^2(i),s^m(i)]}, e_{s^{m+1}(i)})^2 \Scr_\nub(f_{[s^m(i),s^{m+1}(i)]}) } .
$$
Multiplying these two identities and simplifying 
leads to 
$$
\frac{\vp_{[s(i), s^m(i)]} \vp_{[i, s^{m+1}(i)]}}{\vp_{[i,s^m(i)]}\vp_{[s(i),s^{m+1}(i)]}}
= \frac{ \Om_\nub(e_{s(i)}, e_{s^{m+1}(i)})^2 \Om_\nub(e_i - e_{s(i)}, f_{[s^m(i),s^{m+1}(i)]}) } { \Om_\nub( f_{[i,s^{m+1}(i)]} - f_{[s(i),s^{m+1}(i)]}, f_{[s^m(i), s^{m+1}(i)]}) } .
$$
The square of the last term in the numerator simplifies, using \leref{pi-le}(c), to
$$
(\la_i^*)^{-\delta_{m,1}} \Om_\lab( e_i - e_{s(i)}, e_{s^m(i)} + e_{s^{m+1}(i)} ).
$$
Using \leref{pi-le}(d), the square of the denominator simplifies to
$$
(\la_i^*)^{\delta_{m,1}} \Om_\lab( e_i + e_{s(i)}, e_{s^m(i)} + e_{s^{m+1}(i)} ).
$$
Consequently, after some further simplifications we find that
$$
\left( \frac{\vp_{[s(i), s^m(i)]} \vp_{[i, s^{m+1}(i)]}}{\vp_{[i,s^m(i)]}\vp_{[s(i),s^{m+1}(i)]}} \right)^2 = 
(\la_i^*)^{-2\delta_{m,1}} \Om_\lab( e_{s^m(i)}, e_{s(i)} )^2.
$$
Invoking condition  \eqref{-1}, we obtain \eqref{vph}, as desired.
\end{proof}

\section{Rescaling of the generators of a symmetric CGL extension}
\label{4b.2}
For $\gab:=(\ga_1, \ldots, \ga_N)  \in (\kx)^N$ and $f:=(n_1, \ldots, n_N) \in \Zset^N$ set
$$
\gab^f := \ga_1^{n_1} \ldots \ga_N^{n_N} \in \kx.
$$

Given a symmetric CGL extension $R$ of length $N$ and $\gab=(\ga_1, \ldots, \ga_N) \in (\kx)^N$,
one can rescale the generators $x_j$ of $R$,
\begin{equation}
\label{rescale}
x_j \mt \ga_j x_j, \; \; j \in [1,N],
\end{equation}
by which we mean that one can use $\ga_1 x_1, \ldots, \ga_N x_N$ 
as a new sequence of generators of $R$. This obviously does not effect the $\HH$-action 
and the matrix $\lab$, but one obtains a new set of elements $y_k$, $y_{[i,s^m(i)]}$,
$u_{[i, s^m(i)]}$ by applying Theorems \ref{tCGL}, \ref{ty-int} and \coref{u-elem} for the new set 
of generators. (Note that this is not the same as substituting \eqref{rescale} 
in the formulas for $y_k$ and $y_{[i,s^m(i)]}$; those other kind of transformed elements 
will not be even prime because \eqref{rescale} does not determine an algebra automorphism.)
The uniqueness part of \thref{y-int} implies that 
the effect of \eqref{rescale} on the elements $y_{[i, s^m(i)]}$ 
is that they are rescaled by the rule
$$
y_{[i,s^m(i)]} \mt \gab^{e_{[i, s^m(i)]}} y_{[i,s^m(i)]}=(\ga_i \ldots \ga_{s^m(i)}) y_{[i,s^m(i)]}
$$
for all $i \in [1,N]$, $m \in \Zset_{\geq 0}$ such that $s^m(i) \in [1,N]$.
Hence, the effect of \eqref{rescale} on the elements $u_{[i,s^m(i)]}$ is that they are rescaled by
$$
u_{[i,s^m(i)]} \mt 
\gab^{e_{[i, s^m(i)]} + e_{[s(i), s^{m-1}(i)]}} u_{[i,s^m(i)]}= 
(\ga_i \ga_{s(i)}^2 \ldots \ga_{s^{m-1}(i)}^2 \ga_{s^m(i)}) u_{[i, s^m(i)]}
$$
for all $i \in [1,N]$, $m \in \Zset_{> 0}$ 
such that $s^m(i) \in [1,N]$. It follows from \eqref{pi-f} that
the effect of \eqref{rescale} on the scalars $\pi_{[i,s^m(i)]}$ is that 
they are rescaled by 
\begin{multline*}
\pi_{[i,s^m(i)]} \mt 
\gab^{e_{[i, s^m(i)]} + e_{[s(i), s^{m-1}(i)]} - f_{[i,s^m(i)]}} \pi_{[i,s^m(i)]}= 
\\
\ga_i \ga_{s^m(i)} \gab^{2 e_{[s(i), s^{m-1}(i)]} - f_{[i,s^m(i)]}} \pi_{[i, s^m(i)]}.
\end{multline*}
(Note that the rescaling \eqref{rescale} has no effect on the integer vector $f_{[i,s^m(i)]}$.)
This implies at once the following fact.

\bpr{resc} Let $R$ be a symmetric CGL extension of length $N$ and rank $\rk(R)$
for which there exist $\nu_{kj} = \sqrt{\la_{kj}} \in \kx$, $0 \le j < k \le N$. 
Then there exist $N$-tuples $\gab \in (\kx)^N$ such that after the rescaling 
\eqref{rescale} we have 
\begin{equation}
\label{pi-cond}
\pi_{[i,s(i)]} = \Scr_\nub(-e_i + f_{[i,s(i)]}), \; \;  
\forall i \in [1,N] \; \; \mbox{such that} \; \; 
s(i) \neq + \infty.
\end{equation}
The set of those $N$-tuples is parametrized by $(\kx)^{\rk(R)}$ and 
the entries $\ga_1, \ldots, \ga_N$ of all such $N$-tuples 
are recursively determined by
$$
\ga_i \; \; \mbox{is arbitrary if} \; \; p(i) = - \infty 
$$
and 
\begin{equation}
\label{ga-recursive}
\ga_i = \ga_{p(i)}^{-1} \gab^{f_{[p(i),i]}} \pi_{[p(i), i]}^{-1} \Scr_\nu(-e_{p(i)} + f_{[p(i), i]}), \; \; \text{if} \; \; p(i) \ne -\infty,
\end{equation}
where in the right hand side the $\pi$-scalars are the ones 
for the original generators $x_1, \ldots, x_N$ of $R$.
\epr

Note that the product of the first two terms of the right hand side of \eqref{ga-recursive} is a product of 
powers of $\ga_{p(i)}, \ldots, \ga_{i-1}$ since 
$f_{[p(i), i]} \in \Zset_{\geq 0} e_{p(i)+1} + \cdots + \Zset_{\geq 0} e_{i-1}$.  

We refer the reader to Chapter \ref{q-Schu} for an explicit example on 
how the rescaling in \prref{resc} works for the canonical generators of 
the quantum Schubert cell algebras. As for quantum matrices, we have the following

\bex{OqMmn8}
In the case of $R= \OqM$, condition \eqref{pi-cond} is already satisfied, and so no rescaling of the generators is necessary. To see this, let $i \in [1,N]$ with $s(i) \ne +\infty$. Then $i= (r-1)n+c$ for some $r \in [1,m-1]$ and $c\in [1,n-1]$, and $s(i) = rn+c+1$. From \eqref{OqMuisi}, we see that $\pi_{[i, s(i)]} = q$ and $f_{[i, s(i)]} =  e_{i+1} + e_{i+n}$. Recalling the choice of $\nub$ from \exref{OqMmn5}, we find that
$$
\Scr_\nub(-e_i + f_{[i,s(i)]}) = \nu_{i, i+1} \nu_{i, i+n} \nu^{-1}_{i+1, i+n} = (\sqrt q) (\sqrt q) (1) = q = \pi_{[i, s(i)]} .
$$
\eex

\section{Toric frames and mutations in symmetric CGL extensions}
\label{4b.3}
Let $R$ be a symmetric CGL extension of length $N$ for which $\sqrt{\la_{kj}} \in \kx$ 
for all $k,j \in [1,N]$ and let $\nub = (\nu_{kj}) \in M_N(\kx)$ be a 
multiplicatively skew-symmetric matrix such that $\nub^{\cdot 2} = \lab$. 

The procedure from Section \ref{4.3} defines a toric frame for $\Fract(R_{[i,s^m(i)]})$
for all $i \in [1,N]$ and $m \in \Zset_{> 0}$ such that $s^m(i) \in [1,N]$. We 
apply the construction from Section \ref{4.3} to the CGL extension presentation 
of $R_{[i,s^m(i)]}$ associated to the following order of adjoining 
the generators:
\begin{equation}
\label{gene-int}
x_{i+1}, \ldots, x_{s^m(i)-1}, x_i, x_{s^m(i)}.
\end{equation}
Extending \eqref{new-y}, set 
\begin{multline}
\label{ybarpnll}
\ol{y}_{[p^n(l),l]} :=  \\
 \Scr_{\nub}(e_{[p^n(l),l]}) y_{[p^n(l),l]}, \; \; 
\forall j \in [1,N], \; n \in \Zset_{\geq 0} \; \; \text{with} \; \; p^n(j) \in [1,N].
\end{multline}
\index{ypnllbar@$\ol{y}_{[p^n(l),l]}$}

Let $1\le j \le k \le N$, and recall the definition \eqref{genO} of the function $O_-^j$. 
By \coref{primeRjk},
$$
\Bigl\{ \ol{y}_{[p^{O_-^j(l)}(l), l]} \Bigm| l \in [j,k], \; s(l) >k \Bigr\}
$$
is a list of the the homogeneous prime elements of $R_{[j,k]}$ 
up to scalar multiples. The intermediate subalgebras for the 
CGL extension presentation of $R_{[i,s^m(i)]}$ associated to \eqref{gene-int} are 
$$
R_{[i+1, i+1]}, \ldots, R_{[i+1, s^m(i)-1]}, R_{[i,s^m(i)-1]}, R_{[i, s^m(i)]}.
$$
We identify $\Zset^{s^m(i)-i+1} \cong \Zset e_i + \cdots + \Zset e_{s^m(i)} \subseteq \Zset^N$ 
and enumerate the rows and columns of square matrices of size $(s^m(i)-i+1)$ by $i, \ldots, s^m(i)$. 
Define the multiplicatively skew-symmetric matrix $\rbf_{[i,s^m(i)]} \in M_{s^m(i)-i+1}(\kx)$ 
by  \index{rismi@$\rbf_{[i,s^m(i)]}$}
\begin{equation}
\label{rint-tor}
(\rbf_{[i,s^m(i)]})_{k,j} := 
\Om_{\nub} ( e'_k, e'_j) , \; \; 
j,k \in [i,s^m(i)],
\end{equation}
where
\begin{equation}
\label{ekprime}
e'_k =
\begin{cases}
e_{[p^{O_-^{i+1}(k)}(k),k]}, & \mbox{if} \; \; k \in [i+1, s^m(i)-1]
\\
e_{[i,s^{m-1}(i)]}, & \mbox{if} \; \; k = i
\\
e_{[i, s^m(i)]}, & \mbox{if} \; \; k = s^m(i).
\end{cases}  
\end{equation}
\prref{tframe} implies that there is a unique
toric frame $M_{[i,s^m(i)]} : \Zset^{s^m(i)-i+1} \to \Fract(R_{[i,s^m(i)]})$
with matrix $\rbf_{[i,s^m(i)]}$ satisfying  \index{Mismi@$M_{[i,s^m(i)]}$}
$$
M_{[i,s^m(i)]}(e_k) =
\begin{cases}
\ol{y}_{[p^{O_-^{i+1}(k)}(k),k]}, & \mbox{if} \; \; k \in [i+1, s^m(i)-1]
\\
\ol{y}_{[i,s^{m-1}(i)]}, & \mbox{if} \; \; k = i
\\
\ol{y}_{[i, s^m(i)]}, & \mbox{if} \; \; k = s^m(i).
\end{cases}
$$
Recall the definition \eqref{pi-f} of the vectors $f_{[i,s^m(i)]}$. It follows 
from \coref{u-elem} that there exists 
a unique vector 
$$
g_{[i, s^m(i)]} = \sum \{ m_k e_k \mid k \in P_{[i,s^m(i)]} \}    
$$
\index{gismi@$g_{[i, s^m(i)]}$}
such that 
$$
f_{[i,s^m(i)]} = \sum \Bigl\{ m_k e_{[p^{O_-^{i+1}(k)}(k), k]} \Bigm| 
k \in P_{[i,s^m(i)]} \Bigr\},
$$
where the integers $m_k$ are those from Theorems \ref{t2}, \ref{t3}.

Set $t := s^m(i)-i+1$. As above, we identify $\Zset^t$ with $\Zset e_i \oplus \cdots \oplus \Zset e_{s^m(i)}$, that is, with the sublattice of $\Zset^N$ 
with basis $e_j$, $j \in [i,s^m(i)]$. 
Define $\sigma \in M_t(\Zset)$ by $\sigma(e_k) = e'_k$ for all $k \in [i,s^m(i)]$, where the $e'_k$ are as in \eqref{ekprime}. It is easy to check that 
$$
\sigma \; \; {\mbox{is invertible and}} \; \; \sigma(\Zset_{\geq 0}^t) \subseteq \Zset_{\geq 0}^t.
$$
This choice of $\sigma$ ensures that
\begin{equation}
\label{gfsig}      
\sigma(g_{\ismi}) = f_{\ismi}
\end{equation}                  
and
\begin{equation}
\label{rsig-ident}
\rbf_{\ismi} = \bigl( \nub |_{\ismi \times \ismi} \bigr)_\sigma.
\end{equation}

The definition of the toric frame $M_{\ismi}$ also implies that 
\begin{equation}
\label{ltM}
\lt \left( M_{\ismi}(e_k) \right) = \Scr_\nub(\sigma(e_k)) x^{\sigma(e_k)}, \quad
\forall k \in [i, s^m(i)].
\end{equation}

\ble{leadterm} For all $g \in \Zset_{\geq 0}^t$,
$$
\lt \left( M_{\ismi}(g) \right) = \Scr_\nub(\sigma(g)) x^{\sigma(g)} .
$$
\ele

\begin{proof} We prove the lemma by using the inductive argument 
from the proof of Proposition 2.1.
By Eq. \eqref{ltM}, the lemma is valid for $g = e_k$. 
The validity of the lemma for $g$ and $g' \in \Zset_{\geq 0}^t$ 
implies its validity for $g + g' \in \Zset_{\geq 0}^t$ 
because of Eq. \eqref{mult-id} (applied to the toric frame 
$M_{\ismi}$), Eq. \eqref{rsig-ident} above, 
and the fact that
$$
\lt \left( \Scr_\nub(f) x^f \Scr_\nub(f') x^{f'}  \right)= \Om_\nub(f, f') 
\Scr_\nub(f+f') x^{f + f'}, \quad \forall f, f' \in \Zset_{\geq 0}^N.
$$
The last identity is an immediate consequence of Eqs. \eqref{prod-lt} and \eqref{mult-id}.
\end{proof}

\bth{CGLmuta} Let $R$ be a symmetric CGL extension of length $N$.
Assume that the base field $\KK$ contains square roots $\nu_{kj}$ of all scalars 
$\la_{kj}$, such that the subgroup of $\kx$ generated by $\nu_{kj}$, 
$1 \leq j < k \leq N$ does not contain elements of order {\rm2}. Assume also that 
the generators $x_1, \ldots, x_N$ of $R$ are normalized so that 
the condition \eqref{pi-cond} is satisfied. Then for all 
$i \in [1,N]$ and $m \in \Zset_{>0}$ such that $s^m(i) \in [1,N]$,
\begin{equation}
\label{ismi.mut}
\ol{y}_{[s(i), s^m(i)]} = \begin{cases}
M_{[i,s(i)]}(e_{s(i)} - e_i) + M_{[i,s(i)]}(g_{[i,s(i)]} - e_i),  &\text{if} \; \; m=1\\
M_{[i,s^m(i)]} (e_{s^{m-1}(i)} + e_{s^m(i)} - e_i)\, +\\
\qquad\qquad\qquad M_{[i,s^m(i)]} ( g_{[i,s^m(i)]} - e_i),  &\text{if} \; \; m>1. \end{cases}
\end{equation}
\eth

\begin{proof} \leref{leadterm} and Eq. \eqref{gfsig} imply 
$$
\lt \left( M_{\ismi}(g_{\ismi} ) \right) = \Scr_\nub(f_{\ismi}) x^{f_{\ismi}}.
$$
By Eq. \eqref{uformula}, $u_{\ismi}$ is a scalar multiple of $M_{\ismi}(g_{\ismi})$, and so we find that
\begin{equation}
\label{qu-scalar}
u_{[i,s^m(i)]} = 
\pi_{[i,s^m(i)]} \Scr_\nub(f_{[i,s^m(i)]})^{-1} 
M_{[i,s^m(i)]}(g_{[i,s^m(i)]}). 
\end{equation}

In case $m>1$, we observe that
\begin{multline*}
\ol{y}_{[s(i),s^m(i)]} - M_{\ismi}(e_{s^{m-1}(i)}+ e_{s^m(i)}- e_i) =   \\
  \Scr_\nub( e_{[s(i),s^m(i)]} ) y_{[s(i),s^m(i)]} - \xi y^{-1}_{[i,s^{m-1}(i)]} y_{[s(i), s^{m-1}(i)]} y_{\ismi} 
\end{multline*}
where, taking account of \eqref{fg} and \eqref{rint-tor} and simplifying,  
\begin{align*}
\xi &= \Scr_{\rbf_{\ismi}}(e_{s^{m-1}(i)}+ e_{s^m(i)}- e_i) \Scr_\nub(e_{[i,s^{m-1}(i)]})^{-1} \Scr_\nub(e_{[s(i),s^{m-1}(i)]}) \Scr_\nub(e_{\ismi})  \\
 &= \Scr_{\rbf_{\ismi}}(e_{s^{m-1}(i)}+ e_{s^m(i)}- e_i) \times  \\
&\qquad\qquad\times  \Om_\nub(e_i, e_{[s(i),s^{m-1}(i)]}) \Om_\nub(e_i, e_{[s(i),s^m(i)]})^{-1} \Scr_\nub(e_{[s(i),s^m(i)]})  \\
 &= \Om_\lab(e_i, e_{[s(i),s^{m-1}(i)]}) \Scr_\nub(e_{[s(i),s^m(i)]}).
\end{align*}
Thus, we have
\begin{multline}
\label{ybar-M}
\ol{y}_{[s(i),s^m(i)]} - M_{\ismi}(e_{s^{m-1}(i)}+ e_{s^m(i)}- e_i) =  \\
 \Scr_\nub(e_{[s(i),s^m(i)]}) y^{-1}_{[i,s^{m-1}(i)]} u_{\ismi}, \quad\text{if} \; \; m>1.
\end{multline}
A similar argument shows that
\begin{equation}
\label{ybar-M'}
\ol{y}_{[s(i),s(i)]} - M_{[i,s(i)]}(e_{s(i)} - e_i) = x_i^{-1} u_{[i,s(i)]}, \quad\text{if} \; \; m=1.
\end{equation}

From now until the last sentence of the proof, we combine the cases $m=1$ and $m>1$. Using the definition of $M_{\ismi}$ together with \eqref{fg} and simplifying, we obtain
\begin{align*}
M_{\ismi}&(g_{\ismi} - e_i)  \\
 &= \Scr_{\rbf_{\ismi}}(g_{\ismi} - e_i) M_{\ismi}(e_i)^{-1} \times  \\
 &\qquad\qquad \Scr_{\rbf_{\ismi}}(g_{\ismi})^{-1} M_{\ismi}(g_{\ismi})  \\
 &= \Om_{\rbf_{\ismi}}(e_i, g_{\ismi}) \Scr_\nub(e_{[i,s^{m-1}(i)]})^{-1} y^{-1}_{[i,s^{m-1}(i)]} M_{\ismi}(g_{\ismi}) .
\end{align*}
We also have
\begin{align*}
\Om_{\rbf_{[i,s^m(i)]}}( e_i , g_{[i,s^m(i)]}) &= \Om_\nub(e_{[i,s^{m-1}(i)]}, f_{[i,s^m(i)]})  \\
&=
\Om_\nub(e_i + e_{[s(i),s^{m-1}(i)]}, f_{[i,s^m(i)]})  \\
&= \Om_\nub(e_i, f_{[i,s^m(i)]}) \times  \\
&\qquad\qquad \Om_\nub(e_{[s(i),s^{m-1}(i)]}, e_i + 2 e_{[s(i),s^{m-1}(i)]} + e_{s^m(i)})  \\
 &= \Om_\nub(e_i, f_{[i,s^m(i)]}) 
\Om_\nub(e_{[s(i),s^{m-1}(i)]}, e_i + e_{s^m(i)}),
\end{align*}
where the first equality follows from the definition of $\rbf_{[i,s^m(i)]}$ and the third 
from \leref{pi-le} (c). 
It follows from Eq. \eqref{fg} that
\begin{align*} 
\Scr_\nub(e_{[i,s^{m-1}(i)]}) &= 
\Om_\nub(e_i, e_{[s(i), s^{m-1}(i)]})^{-1} \Scr_\nub(e_{[s(i),s^{m-1}(i)]})
\\
&=\Om_\nub(e_i, e_{[s(i), s^{m-1}(i)]})^{-1} \Om_\nub(e_{[s(i),s^{m-1}(i)]},e_{s^m(i)})
\Scr_\nub(e_{[s(i),s^m(i)]})
\end{align*}
and
\begin{equation}
\label{Snuvec}
\Scr_\nub(-e_i + f_{[i,s^m(i)]}) = \Om_\nub(e_i, f_{[i,s^m(i)]}) \Scr_\nub(f_{[i,s^m(i)]}).
\end{equation}
Combining the first three of the previous four equations yields
\begin{multline}
\label{Mg-e}
M_{\ismi}(g_{\ismi} - e_i) =  \\
 \Om_\nub(e_i, f_{\ismi}) \Scr_\nub(e_{[s(i),s^m(i)]})^{-1} y^{-1}_{[i,s^{m-1}(i)]} M_{\ismi}(g_{\ismi}) .
\end{multline}

Finally, applying \prref{pi} (b) and \eqref{fg} to \eqref{qu-scalar} yields
$$
u_{\ismi} = \Scr_\nub(e_{[s(i),s^m(i)]})^{-2} \Om_\nub(e_i, f_{\ismi}) M_{\ismi}(g_{\ismi}).
$$
Substituting this identity into \eqref{ybar-M}, \eqref{ybar-M'} and taking account of \eqref{Mg-e} yields \eqref{ismi.mut} and completes the proof of the theorem.
\end{proof}

\bre{appl} Because of \thref{y-int}, all instances when \thref{3} is applicable 
to a symmetric CGL extension with respect to the CGL extension presentations 
associated to the elements of the set $\Xi_N$ are covered by \thref{CGLmuta}. 
We refer the reader to \prref{cluster-tau-ind} for details.
\ere

Excerpting several steps from the proof of \thref{CGLmuta} yields the following:

\bco{new6.7} Let $R$ be a symmetric CGL extension of length $N$.
Assume that the base field $\KK$ contains square roots $\nu_{kj}$ of all scalars 
$\la_{kj}$, such that the subgroup of $\kx$ generated by $\nu_{kj}$, 
$1 \leq j < k \leq N$ does not contain elements of order {\rm2}.

{\rm(a)} For all 
$i \in [1,N]$ and $m \in \Zset_{>0}$ such that $s^m(i) \in [1,N]$, eq. \eqref{qu-scalar} holds:
$$
u_{[i,s^m(i)]} = 
\pi_{[i,s^m(i)]} \Scr_\nub(f_{[i,s^m(i)]})^{-1} 
M_{[i,s^m(i)]}(g_{[i,s^m(i)]}).
$$

{\rm(b)} Condition \eqref{pi-cond} is satisfied if and only if
\begin{equation}
\label{6.13alt}
M_{[i,s(i)]}( g_{[i,s(i)]} - e_i) = x_i^{-1} u_{[i,s(i)]}, \: \: \forall i \in [1,N], \; s(i) \ne +\infty.
\end{equation}
\eco

\begin{proof} Observe that \eqref{pi-cond} is not used in the proof of \thref{CGLmuta} until the very last paragraph, where \prref{pi} (b) is applied. Thus, all equations in the proof up through \eqref{Mg-e} hold without assuming \eqref{pi-cond}. This gives part (a).

(b) In case \eqref{pi-cond} holds, we obtain \eqref{6.13alt} by substituting \eqref{ybar-M'} into the case $m=1$ of \eqref{ismi.mut}.

Conversely, if \eqref{6.13alt} holds for some $i \in [1,N]$ with $s(i) \ne +\infty$, then \eqref{Mg-e} implies
$$
u_{[i,s(i)]} = \Om_\nub(e_i, f_{[i,s(i)]}) M_{[i,s(i)]}( g_{[i,s(i)]} ).
$$
Comparing this with \eqref{qu-scalar} yields
$$
\pi_{[i,s(i)]} = \Om_\nub(e_i, f_{[i,s(i)]}) \Scr_\nub(f_{[i,s(i)]}) = \Scr_\nub( -e_i + f_{[i,s(i)]}),
$$
recall \eqref{Snuvec}, and verifies \eqref{pi-cond}.
\end{proof}

\chapter[Division properties of mutations between CGL presentations]{Division properties of mutations between CGL extension presentations}
\label{Integr}
\section{Main result}
\label{5.1}
In this chapter we describe the intersection of the localizations of a 
CGL extension $R$ by the two sets of cluster variables ($y$-elements) 
in the setting of \thref{1} (b). This plays a key role in the next chapter 
where we prove that every symmetric CGL extension is a quantum cluster algebra which 
equals the corresponding upper quantum cluster algebra. 

Versions of the following two lemmas are well known. We give a
proof of the first for completeness.

\ble{primloc} Assume that $R$ is a $\KK$-algebra domain and $E$ a left or right 
Ore set in $R$. Let $p \in R$ be a prime element such that $\kx pE = \kx Ep$ and $Rp \cap E = \varnothing$.
Then $p$ is a prime element of $R[E^{-1}]$. If $s \in R$ and 
$p \nmid s$ as elements of $R$, then $p \nmid s$ as elements of $R[E^{-1}]$.
\ele

\begin{proof} If $e\in E$, then $pe = \al fp$ for some $\al \in \kx$ and $f \in E$, whence $pe^{-1} = \al^{-1} f^{-1} p$ in $R[E^{-1}]$. Similarly, $e^{-1}p = \be^{-1} p g^{-1}$ for some $\be \in \kx$ and $g \in E$. These observations, together with the normality of $p$ in $R$, imply that $p$ is normal in $R[E^{-1}]$. In particular, $R[E^{-1}]p$ is an ideal of $R[E^{-1}]$.

It follows from the assumption $Rp \cap E = \varnothing$ and the complete primeness of $Rp$ that $R[E^{-1}]p \cap R = Rp$. Consequently, $R[E^{-1}]/ R[E^{-1}] p$ is an Ore localization of the domain $R/Rp$, so it too is a domain, and thus $p$ is a prime element of $R[E^{-1}]$.

The last statement of the lemma follows from the equality $R[E^{-1}]p \cap R = Rp$. 
\end{proof}

\ble{loc} Let $B[x;\sigma,\delta]$ be an Ore extension {\rm(}recall Convention {\rm\ref{skewpoly})} and $E$ a  {\rm(}left or right\/{\rm)} Ore set of regular elements 
in $B$. If $\kx E$ is $\sigma$-stable, then $E$ is a {\rm(}left or right\/{\rm)} Ore set of regular elements 
in $B[x;\sigma,\delta]$.
\ele

\begin{proof} Replace $E$ by $\kx E$ and apply \cite[Lemma 1.4]{Go}.
\end{proof}

Recall that given a (noncommutative) domain $A$ and $a, b \in A$, we write $a \mid_l b$  \index{zzz1@$\mid_l$} to denote that $b \in a R$. \index{zzz1@$\mid_l$}

Assume now that $R$ is a CGL extension of length $N$ as in \eqref{itOre}.
Let $y_1, \ldots, y_N \in R$ be the sequence of 
elements from \thref{CGL}. Then $y_j$ is a normal element 
of $R_j$ and by \leref{loc}, $\{y_j^m \mid m \in \Zset_{\geq 0} \}$ 
is an Ore set in $R$, $\forall j \in [1,N]$. 
It follows from \eqref{ycomm} that for all $I \subseteq [1,N]$, 
\index{EI@$E_I$}
\begin{equation}
\label{Sloc}
E_I:=
\kx \prod_{j \in I} \{y_j^m \mid m \in \Zset_{\geq 0} \}
\end{equation}
is a multiplicative set and thus is an Ore set in $R$. Set $E:= E_{[1,N]}$.
We will say that $y_1^{m_1} \ldots y_N^{m_N}$ is a   \index{minimal denominator}
\emph{minimal denominator} of a nonzero element $v \in R[E^{-1}]$ if 
$$
v = y_1^{-m_1} \ldots y_N^{-m_N} s
$$
for some $s \in R$ such that $y_j \nmid_l s$ for all $j \in [1,N]$ 
with $m_j > 0$.

For the rest of this chapter, we assume the setting of \thref{1} (b). 
Then for all $I' \subseteq [1,N]$,
\begin{equation}
\label{E'set}
E'_{I'}:= \kx \prod_{j \in I'} \{(y'_j)^m \mid m \in \Zset_{\geq 0} \}
\end{equation}
is another Ore set in $R$. Recall from \thref{1} (b) that $y'_j$ is a scalar multiple of 
$y_j$ for $j \neq k$. Hence, $E'_{I'} = E_{I'}$ if $k \notin I'$. Denote $E':= E'_{[1,N]}$.   

\ble{EIEI'} For all $I, I' \subseteq [1,N]$,
\begin{equation}
\label{EcapE}
E_I \cap E'_{I'} = 
E_{(I \cap I') \backslash \{ k \}} = E'_{(I \cap I') \backslash \{ k \} }.
\end{equation}
\ele

\begin{proof} Let $e \in E_I \cap E'_{I'}$, and write
$$
e = \be \prod_{j=1}^N y_j^{m_j} = \be' \prod_{j=1}^N (y'_j)^{m'_j}
$$
for some $\be, \be' \in \kx$ and $m_j, m'_j \in \Znn$ such that $m_j=0$ for $j \notin I$ and $m'_j=0$ for $j \notin I'$. We will show that
\begin{equation}
\label{mmatchm'}
m_k = m'_k = 0 \quad \text{and} \quad m_j = m'_j, \; \; \forall j \in [1,N], \; j \ne k,
\end{equation}
which then implies $e \in E_{(I \cap I') \backslash \{ k \}} = E'_{(I \cap I') \backslash \{ k \} }$.

Since $y_N$ and $y'_N$ are prime elements and associates in $R$, and they do not divide $y_1,\dots,y_{N-1}$, unique factorization implies that $m_N = m'_N$. Then, after replacing $e$ by $e y_N^{-m_N}$ and adjusting $\be'$, we may reduce to the case in which $m_N = m'_N = 0$. Continuing in this manner, we reduce to the case in which $m_j = m'_j = 0$ for all $j > k$.

If $m_k > 0$, we observe using \thref{1} (b) that
$$
e y'_k = \be \biggl( \prod_{j<k} y_j^{m_j} \biggr) y_k^{m_k-1} (u + \al_{kp(k)} y_{p(k)} y_{k+1})
$$
for some $u \in R_{k-1}$. This implies that $ey'_k$ is an element of $R_{k+1}$ which has degree $1$ with respect to $x_{k+1}$. On the other hand, the equation
$$
e y'_k = \be' \biggl( \prod_{j<k} (y'_j)^{m'_j} \biggr) (y'_k)^{m'_k+1}
$$
shows that $e y'_k$ has degree $m'_k + 1$ with respect to $x_{k+1}$, in view of \eqref{ykyk'}. Hence, $m'_k = 0$ and $e \in R_{k-1}$, contradicting our assumption that $m_k > 0$.

Thus $m_k = 0$, and similarly $m'_k = 0$. Proceeding as above, we conclude that $m_j = m'_j$ for all $j < k$, verifying \eqref{mmatchm'}.
\end{proof}

The next theorem contains the main result of this chapter.

\bth{division} Assume the setting of Theorem {\rm\ref{t1} (b)}. Let 
$$
y_1^{m_1} \ldots y_N^{m_N} \quad
\mbox{and} \quad 
y_1^{m'_1} \ldots y_{k-1}^{m'_{k-1}} (y'_k)^{m'_k} y_{k+1}^{m'_{k+1}} \ldots y_N^{m'_N}
$$
be two minimal denominators of a nonzero element $v \in R[E^{-1}] \cap R[(E')^{-1}]$
{\rm(}with respect to the two different localizations\/{\rm)}. Then
$$
m_k = m'_k =0 \quad \mbox{and} 
\quad
m_j= m'_j, \; \; \forall j \in [1,N], \; j \neq k.
$$ 
In particular, 
$$
R[(E_I)^{-1}] \cap R[(E'_{I'})^{-1}] = R[ (E_I \cap E'_{I'})^{-1}] =
R[(E_{(I\cap I')\backslash\{k \}})^{-1}]
$$
for all $I,I' \subseteq [1,N]$.
\eth

\thref{division} implies that each nonzero element of $R[E^{-1}]$ has a unique minimal 
denominator. The theorem was inspired by \cite[Theorem 4.1]{GLSh2}.

\section{Proof of the main result}
\label{5.2}
For $s \in R$ and $j \in [1,N]$, denote
\begin{equation}
\label{aexp}
s = \sum_{l_j, \ldots, l_N \in \Zset_{\geq 0} } s_{l_j, \ldots, l_N} x_j^{l_j} 
\ldots x_N^{l_N}, \quad s_{l_j, \ldots, l_N} \in R_{j-1}.
\end{equation}
For the proof of \thref{division} we will need the following 
two lemmas. 

\ble{relprim} Let $R$ be a CGL extension of length $N$ and $y_1, \ldots, y_N$ 
the sequence of elements from Theorem {\rm\ref{tCGL}}. If 
$$
y_j \mid_l y_1^{n_1} \ldots y_{j-1}^{n_{j-1}} s
$$
for some $j \in [1,N]$, $s \in R$, and $n_1, \ldots, n_{j-1} \in \Zset_{\geq 0}$, then 
$y_j \mid_l s$.
\ele 

\begin{proof} There exists $s' \in R$ such that $y_j s' = y_1^{n_1} \ldots y_{j-1}^{n_{j-1}} s$.
Comparing the coefficients of $x_{j+1}^{l_{j+1}} \ldots x_N^{l_N}$ leads to 
$$
y_j \mid_l y_1^{n_1} \ldots y_{j-1}^{n_{j-1}} s_{l_{j+1}, \ldots, l_N}, \quad
\forall l_{j+1}, \ldots, l_N \in \Zset_{\geq 0}.
$$
By \thref{CGL}, $y_j$ is a prime element of $R_j$ and $y_j \nmid y_1, \ldots, y_{j-1}$. 
So $y_j \mid s_{l_{j+1}, \ldots, l_N}$, for all $l_{j+1}, \ldots, l_N \in \Zset_{\geq 0}$,
and thus $y_j \mid_l s$.
\end{proof}

\ble{relprim2} If, in the setting of Theorem {\rm\ref{t1} (b)}, $y_k \mid_l s (y'_k)^n$ for some 
$s \in R_{k+1}$ and $n \in \Zset_{\geq 0}$, then $y_k \mid_l s$.
\ele

\begin{proof} Write $s = \sum_{l \in \Zset_{ \geq 0}} s_l x_{k+1}^l$
for some $s_l \in R_k$. 
Similarly to the proof of \leref{relprim}, the conclusion  
of this lemma is equivalent to $y_k \mid s_l$, $\forall l \in \Zset_{\geq 0}$.
Assuming this is not the case, we replace $s$ with 
$$
s - \sum \{ s_l x_{k+1}^l \mid l \in \Zset_{\geq 0 }, y_k \nmid s_l \}. 
$$
Note that this new version of $s$ satisfies the assumption of the lemma. 
Moreover, it has the form 
$$
s = s_L x_{k+1}^L + \sum_{l=0}^{L-1} s_l x_{k+1}^l 
$$
for some $L \in \Zset_{\geq 0}$ such that $y_k \nmid s_L$. 
Then 
$$
s (y'_k)^n - \xi s_L y_{p(k)}^n x_{k+1}^{L+n} \in \oplus_{l=0}^{L+n-1} R_k x_{k+1}^l  
$$ 
for some $\xi \in \kx$. The assumption $y_k \mid_l s (y'_k)^n$ implies 
that $y_k \mid s_L y_{p(k)}^n$ which is a contradiction 
since $y_k$ is a prime element of $R_k$ which does not divide $s_L$ or $y_{p(k)}$.
\end{proof}

\begin{proof}[Proof of \thref{division}]
After multiplying $v$ 
by $\prod_{j \neq k} y_j^{\min (m_j, m'_j)}$ 
and using \eqref{ycomm} for both sets of $y$-elements,
we can assume that 
\begin{equation}
\label{min}
\min( m_j, m'_j ) = 0, \quad
\forall j \neq k. 
\end{equation}
(Under this condition we have to prove that $m_1 = m'_1 = \ldots = m_N = m'_N=0$.)
Write
$$
v = y_1^{-m_1} \ldots y_N^{-m_N} s = 
y_1^{-m'_1} \ldots y_{k-1}^{-m'_{k-1}} (y'_k)^{-m'_k} y_{k+1}^{-m'_{k+1}} \ldots y_N^{- m'_N}
s'   
$$
for some nonzero elements $s, s' \in R$. 

First, we prove that 
\begin{equation}
\label{a1}
m_{k+1}= m'_{k+1} = \ldots = m_N = m'_N =0.
\end{equation}
Assume that this is not the case. Then 
$$
i = \max \{ j \mid \max(m_j, m'_j) > 0 \} \geq k +1 .
$$
Without loss of generality, we can assume that $m'_i>0$. 
Applying \eqref{ycomm} to both sets of $y$-elements we obtain
$$
y_i^{m'_i} \Big( \prod_{j \in [1,i-1] \backslash \{k\} } 
y_j^{m'_j} \Big) s = \xi
\Big( \prod_{j \in [1,i-1]} y_j^{m_j} \Big) (y'_k)^{-m'_k} s' 
$$
for some $\xi \in \kx$. Comparing the coefficients of 
$x_{i+1}^{l_{i+1}} \ldots x_N^{l_N}$ gives
\begin{equation}
\label{yy'eq1}
y_i^{m'_i} \Big( \prod_{j \in [1,i-1] \backslash \{k\} } 
y_j^{m'_j} \Big) s_{l_{i+1}, \ldots, l_N} = \xi
\Big( \prod_{j \in [1,i-1]} y_j^{m_j} \Big) (y'_k)^{-m'_k} 
s'_{l_{i+1}, \ldots, l_N} \in R_i.
\end{equation}
Since $y_i$ is a prime element of $R_i$, by \leref{primloc}, 
it is a prime element of $R_i[(y'_k)^{-1}]$. Furthermore,
by \thref{CGL} and \leref{primloc}, $y_i \nmid y_j$ in $R_i[ (y'_k)^{-1}]$
for $j \in [1,i-1]$. Therefore for all $l_{i+1}, \ldots, l_N \in \Zset_{\geq 0}$,
$y_i \mid s'_{l_{i+1}, \ldots, l_N}$ 
in $R_i[(y'_k)^{-1}]$ and thus in $R_i$ (by \leref{primloc}).
Hence, $y_i \mid_l s'$ which contradicts 
the minimality assumption on the denominators.

Next, we prove that
\begin{equation}
\label{a2}
m_k = m'_k =0.
\end{equation}
Assume the opposite, e.g., $m_k >0$. (The other case is analogous because
of the symmetric nature of $y_k$ and $y'_k$.)
Similarly to \eqref{yy'eq1}, one obtains that there exists $\xi \in \kx$ such that 
$$
\Big( \prod_{j=1}^{k-1}
y_j^{m'_j} \Big) s_{l_{k+2}, \ldots, l_N} = \xi (y_k)^{m_k}
\Big( \prod_{j =1}^{k-1} y_j^{m_j} \Big) (y'_k)^{-m'_k} 
s'_{l_{k+2}, \ldots, l_N} \in R_{k+1}$$
for all $l_{k+2}, \ldots, l_N \in \Zset_{\geq 0}$.
Since $\{ (y'_k)^m \mid m \in \Zset_{\geq 0} \}$ 
is an Ore set in $R_{k+1}$,
$$
(y'_k)^{-m'_k} 
s'_{l_{k+2}, \ldots, l_N} = 
s''_{l_{k+2}, \ldots, l_N} (y'_k)^{-m''_k}
$$
for some $s''_{l_{k+2}, \ldots, l_N} \in R_{k+1}$ and $m''_k \in \Zset_{\geq 0}$
(depending on $l_{k+2}, \ldots, l_N$).
Therefore 
$$
y_k \mid_l \Big( \prod_{j =1}^{k-1} y_j^{m'_j} \Big) 
s_{l_{k+2}, \ldots, l_N} (y'_k)^{m''_k}.
$$
It follows from Lemmas \ref{lrelprim} and \ref{lrelprim2} that
$y_k \mid s_{l_{k+2}, \ldots, l_N}$, 
$\forall l_{k+2}, \ldots, l_N \in \Zset_{\geq 0}$ and thus $y_k \mid_l s$.
This contradicts the minimality assumption on denominators
and proves \eqref{a2}. 

Finally, analogously to \eqref{a1}
one shows that
$$
m_1 = m'_1 = \ldots = m_{k-1} = m'_{k-1} = 0.
$$
This completes the proof of \thref{division}. 
\end{proof}

\chapter{Symmetric CGL extensions and quantum cluster algebras}
\label{main}
In this chapter we prove that every symmetric CGL extension
possesses a quantum cluster algebra structure under a couple of very mild additional
assumptions. Those assumptions are satisfied by all known 
examples of symmetric CGL extensions.
The quantum cluster algebra structure is constructed
in an explicit fashion. 
We furthermore prove that each of these quantum 
cluster algebras equals the corresponding upper quantum cluster algebra. 
The proofs work for base fields of arbitrary characteristic.

\section{General setting}
\label{6.1}
Fix a symmetric CGL extension $R$ of length $N$ and rank $\rk(R)$ 
as in Section \ref{2.5} over a base field $\KK$ which contains the square roots 
of all scalars $\la_{lj}$, $1 \leq j < l \leq N$. 
Choose $\nu_{lj} = \sqrt{\la_{lj}} \in \kx$ and extend (as in \eqref{nu2}) to a multiplicatively skew-symmetric matrix $\nub = (\nu_{lj}) \in M_N(\kx)$ such that $\nub^{\cdot2} = \lab$. Define a multiplicatively 
skew-symmetric matrix $\rbf \in M_N(\kx)$ by \eqref{r} and 
a sequence of normalized homogeneous prime elements 
$\ol{y}_1, \ldots, \ol{y}_N$ by \eqref{new-y}. (We recall that generally 
each of those is a prime element of some of 
the subalgebras $R_l$, not of the full algebra $R=R_N$.) By \prref{tframe}
we obtain a toric frame $M : \Zset^N \to \Fract(R)$ whose matrix is $\rbf(M) := \rbf$ 
and such that $M(e_k) := \ol{y}_k$, for all $k \in [1,N]$.

Next, consider an arbitrary element $\tau \in \Xi_N \subset S_N$, recall \eqref{tau}.
We will associate to it a toric frame $M_\tau :\Zset^N \to \Fract(R)$ such that
$M_\id = M$. By \prref{tauOre} we have the CGL extension presentation
\begin{equation}
\label{tau-pres}
R = \KK [x_{\tau(1)}] [x_{\tau(2)}; \sigma''_{\tau(2)}, \delta''_{\tau(2)}] 
\cdots [x_{\tau(N)}; \sigma''_{\tau(N)}, \delta''_{\tau(N)}]
\end{equation}
with $\sigma''_{\tau(k)} = ( h''_{\tau(k)} \cdot)$ where 
$h''_{\tau(k)} \in \HH$, $\forall k \in [1,N]$ and  
\begin{align*}
&h''_{\tau(k)} = h_{\tau(k)}, \; \; \delta''_{\tau(k)} = \delta_{\tau(k)}, 
\quad \mbox{if} \; \;  \tau(k) = \max \, \tau( [1,k-1]) +1
\\  
&h''_{\tau(k)} = h^\sy_{\tau(k)}, \; \; \delta''_{\tau(k)} = \delta^*_{\tau(k)},
\quad \mbox{if} \; \; \tau(k) = \min \, \tau( [1,k-1]) -1.
\end{align*}
The $\lab$-matrix of the presentation \eqref{tau-pres} is $\lab_\tau:=\tau^{-1} \lab \tau$
\index{lambdatau@$\lab_\tau$}
(where as before we use the canonical embedding $S_N \hra GL_N(\Zset)$ 
via permutation matrices). In other words,
the entries of $\lab_\tau$ are given by $(\lab_\tau)_{lj}:= \la_{\tau(l)\tau(j)}$.
For this presentation we choose the multiplicatively skew-symmetric matrix
of square roots $\nub_\tau:= \tau^{-1} \nub \tau \in M_N(\kx)$ \index{nutau@$\nub_\tau$}
(with entries $(\nub_\tau)_{lj}:=\nu_{\tau(l)\tau(j)}$). Denote 
by $\wh{\rbf}_\tau$ \index{rtauhat@$\wh{\rbf}_\tau$} the corresponding multiplicatively skew-symmetric 
matrix derived from $\nub_\tau$ by applying \eqref{r}. Let $\ol{y}_{\tau,1}, \ldots, \ol{y}_{\tau,N}$ 
\index{ytauj@$\ol{y}_{\tau,j}$}
be the sequence of normalized prime elements given by \eqref{new-y}
applied for the presentation \eqref{tau-pres}.
By \prref{tframe} there exists a unique 
toric frame $\wh{M}_\tau : \Zset^N \to \Fract(R)$ \index{Mtauhat@$\wh{M}_\tau$} whose matrix is 
$\rbf(\wh{M}_\tau) := \wh{\rbf}_\tau$ and such that $\wh{M}_\tau(e_k) := \ol{y}_{\tau, k}$, 
$\forall k \in [1,N]$.    

{\em{The toric frames $\wh{M}_\tau$ will be the main toric frames that will be used 
to construct a quantum cluster algebra structure on the symmetric CGL extension 
$R$. In order to connect those frames with mutations {\rm(}without the need of additional 
permutations\/{\rm)} we need to permute the order of the cluster variables in each frame $\wh{M}_\tau$.}} 
This is done by applying a composition of $\tau$ 
and a permutation $\tau_\bu \in \prod_{a \in \Zset} S_{\eta^{-1}(a)}$,
where $\eta \colon [1,N] \to \Zset$ is a function satisfying the 
conditions of \thref{CGL} for the original CGL extension presentation of $R$.
Note that all terms in the above product are trivial except for the 
terms coming from the range of $\eta$.
For $a$ in the range of $\eta$, denote for brevity $|a|:= |\eta^{-1}(a)|$. Consider the 
set
\begin{equation}
\label{eta-1}
\eta^{-1}(a)= \{ \tau(k_1), \ldots, \tau(k_{|a|}) \mid k_1 < \cdots  < k_{|a|} \}
\end{equation}
and order its elements in an increasing order. Define $\tau_\bu \in S_N$ \index{taubullet@$\tau_\bu$} by setting
$\tau_\bu(\tau(k_i))$ to be equal to the $i$-th element in the list 
(for all choices of $a$ and $i$). 

\bex{tau-bu} Let $N=6$ and $\eta : [1,6] \to \Zset$ be given 
by 
$$
\eta(1)=\eta(4) = \eta(6)=1, \quad
\eta(2)=\eta(5)= 2, \quad
\eta(3) =3.
$$
(This is the $\eta$-function for the standard CGL extension presentation 
of the quantum Schubert cell algebra $\UU^-[w_\ci] \subset \UU_q(\slfrak_4)$ 
from Chapter \ref{q-gr}, where 
$w_\ci$ is the longest element of $S_4$.)
Let
$$
\tau := \tau_{(4,5)} = [3,4,2,5,1,6]
$$ 
in the notation \eqref{tau-seq} from Section \ref{4a.3}. Then
$$
\eta^{-1}(1) = \{ \tau(5)=1 < \tau(2) = 4 < \tau(6) = 6 \}
$$
and
$$
\tau_\bu(4) =1, \quad \tau_\bu(1) = 4, \quad \tau_\bu(6) =6.
$$
Similarly, one computes that $\tau_\bu(j)=j$ for $j =2$, $3$ and $5$.
\eex

Define the toric frame 
\begin{equation}
\label{Mtau}
M_\tau := \wh{M}_\tau (\tau_\bu \tau)^{-1} : \Zset^N \to \Fract(R).
\end{equation}
\index{Mtau@$M_\tau$}
It satisfies $M_\tau(e_k) = \ol{y}_{\tau, (\tau_\bu \tau)^{-1} (k)}$, 
$\forall k \in [1,N]$ and its matrix equals 
$$
\rbf_\tau: = \rbf(M_\tau) = (\tau_\bu \tau) \wh{\rbf}_\tau (\tau_\bu \tau )^{-1}.
$$
\index{rtau@$\rbf_\tau$}
The point of applying $\tau$ in this normalization is
to match the indexing of the $\ol{y}$-elements with the one of the 
$x$-elements in \eqref{tau-pres}. (Note that the 
order of the $x$-elements in \eqref{tau-pres} is 
$x_{\tau(1)}, \ldots, x_{\tau(N)}$.) The application 
of $\tau_\bu$ then rearranges the $\eta$-preimages 
$\tau(k_1), \ldots, \tau(k_{|a|})$ from \eqref{eta-1}
in increasing order. This is needed because in the setting
of \thref{1} (b) the element $y_k$ (not $y_{k+1}$) 
gets mutated. Clearly, $\tau_\bu$ preserves the level sets of 
$\eta$.

Recall that $P(N) = \{ k \in [1,N] \mid s(k) = + \infty \}$
parametrizes the set of homogeneous prime elements of $R$,
i.e.,  
\begin{equation}
\label{RPr}
\{y_k \mid k \in P(N) \} \; \; \mbox{is a list of the homogeneous prime
elements of $R$}
\end{equation}
up to associates. Define
\begin{equation}
\label{ex-prim}
\ex := [1,N] \setminus P(N) = \{ k \in [1,N] \mid s(k) \neq + \infty \}.
\end{equation}
Since $|P(N)| = \rk(R)$,
the cardinality of this set is $|\ex| = N - \rk(R)$. Finally, recall 
that for a homogeneous element $u \in R$, $\chi_u \in \xh$ denotes 
its $\HH$-eigenvalue.

\section{Statement of the main result}
\label{6.2}
The next theorem contains the main result of the paper. Keep $M$, $\tau_\bu$, $M_\tau$, $\rbf_\tau$, $\ex$, etc. as in Section \ref{6.1}.

\bth{cluster} Let $R$ be a symmetric CGL extension of length $N$ and rank $\rk(R)$
as in Definition {\rm\ref{dsymmetric}}. Define $\ex \subset [1,N]$ by \eqref{ex-prim}.
Assume that the base field $\KK$ contains square roots $\nu_{lj}$ of all scalars 
$\la_{lj}$, such that $\nub:= (\nu_{lj})$ is a multiplicatively 
skew-symmetric matrix and the subgroup of $\kx$ generated by $\{\nu_{lj} \mid 
1 \leq j < l \leq N\}$ does not contain elements of order {\rm$2$}. Assume also that there exist
positive integers $d_i$, $i \in {\mathrm{range}}(\eta)$ such that 
\begin{equation}
\label{d-prop}
(\la_l^\sy)^{d_{\eta(j)}} = (\la_j^\sy)^{d_{\eta(l)}}, \quad
\forall j, l \in \ex,
\end{equation}
recall the equality \eqref{la-eq}. Let the sequence of 
generators $x_1, \ldots, x_N$ of $R$ be normalized {\rm(}rescaled\/{\rm)}
so that \eqref{pi-cond} is satisfied {\rm(}recall Proposition {\rm\ref{presc})}. 

Then the following hold:

{\rm(a)} For all $\tau \in \Xi_N$ {\rm(}see \eqref{tau}{\rm)} and $l \in \ex$, 
there exists a unique vector $b_\tau^l \in \Zset^N$ such that $\chi_{M_\tau(b_\tau^l)}=1$ and 
\begin{equation}
\label{linear-eq}
\Om_{\rbf_\tau} ( b^l_\tau, e_j) = 1, \; \; \forall j \in [1,N], \; j \neq l
\quad \mbox{and} \quad
\Om_{\rbf_\tau} (b^l_\tau, e_l)^2 = \la_l^*.
\end{equation}
Denote by $\wt{B}_\tau \in M_{N \times |\ex| }(\Zset)$ the matrix with columns 
$b^l_\tau$, $l \in \ex$. Let $\wt{B}:= \wt{B}_\id$.
\index{Btautilde@$\wt{B}_\tau$}
 
{\rm(b)} For all $\tau \in \Xi_N$, the pair $(M_\tau, \wt{B}_\tau)$ is a 
quantum seed for $\Fract(R)$. The principal part of $\wt{B}_\tau$ 
is skew-symmetrizable via the integers $d_{\eta(k)}$, $k \in \ex$.

{\rm(c)} All such quantum seeds are mutation-equivalent to each other.
More precisely, they are linked by the following one-step mutations.
Let $\tau, \tau' \in \Xi_N$ be such that
$$
\tau' = ( \tau(k), \tau(k+1)) \tau = \tau (k, k+1)
$$
for some $k \in [1,N-1]$.
If $\eta(\tau(k)) \neq \eta (\tau(k+1))$, then $M_{\tau'} = M_\tau$.
If $\eta(\tau(k)) = \eta (\tau(k+1))$, 
then $M_{\tau'} = \mu_{k_\bu}(M_\tau)$, where $k_\bu = \tau_\bu \tau(k)$.

{\rm(d)} We have the following equality between the CGL extension $R$ and the quantum cluster and upper cluster algebras associated to $M$, $\wt{B}$, $\varnothing$:
$$
R = \Abb(M, \wt{B}, \varnothing)_\KK = \UU(M, \wt{B}, \varnothing)_\KK.
$$
In particular, $\Abb(M, \wt{B}, \varnothing)_\KK$ and $\UU(M, \wt{B}, \varnothing)_\KK$ 
are affine and noetherian, and more precisely $\Abb(M, \wt{B}, \varnothing)_\KK$ is 
generated by the cluster variables in the seeds parametrized by the subset 
$\Ga_N$ of $\Xi_N$, recall \eqref{sequence}.

{\rm(e)} Let $\inv$ be any subset of the set $P(N)$ of frozen variables, cf. \eqref{RPr}.
Then 
$$
R[y_k^{-1} \mid k \in \inv] = \Abb(M, \wt{B}, \inv)_\KK = \UU(M, \wt{B}, \inv)_\KK.
$$
\eth 

\thref{cluster} is proved in Section \ref{6.5}--\ref{6.7}. The strategy 
of the proof is summarized in Section \ref{6.5a}.
In Section \ref{6.3} 
we derive an explicit formula for the cluster variables 
of the quantum seeds that appear in the statement of \thref{cluster}.
\prref{pair-mut} (b) and \leref{1eig} imply at once the following properties of 
all quantum seeds of the quantum cluster algebras in \thref{cluster}.

\bpr{allclust} All quantum seeds $(M_\star, \wt{B}_\star)$ of the quantum cluster algebras 
in Theorem {\rm\ref{tcluster}} have the properties that 
$$
\chi_{M_\star(b^l)}=1, \quad
\Om_{\rbf(M_\star)} ( b^l_\star, e_j) = 1, \; \; \forall j \in [1,N], \; j \neq l,
\quad \mbox{and} 
\quad
\Om_{\rbf(M_\star)} (b^l_\star, e_l)^2 = \la_l^*
$$
for all $l \in \ex$,
where the $b^l_\star$, $l \in \ex$ are the columns of $\wt{B}_\star$.
\epr

In Section \ref{new8.9}, we give a ring theoretic interpretation of some columns of the initial exchange matrix $\wt{B}$, in terms of the exponent vectors $f_{[i,s(i)]}$ for certain of the elements $u_{[i,s(i)]}$.

\bex{OqMmn9}
Let $R = \OqM$ with the CGL extension presentation given in \exref{OqMmn}, and assume there exists a square root of $q$ in $\kx$. Fix a choice of $\sqrt q$. Define $\eta$ for $R$ as in \exref{OqMmn2} and $\ex$ as in \eqref{ex-prim}. Then
$$
\ex = \{ (r-1)n+c \mid r\in [1,m-1], \; c\in [1,n-1] \}.
$$
Next, define $\nu_{lj}$ for $l,j \in [1,N]$ by replacing $q$ with $\sqrt q$ in \eqref{OqMlakj}. Then the matrix $(\nu_{lj})$ is multiplicatively skew-symmetric, and the subgroup $\langle \nu_{lj} \mid 1\le j < l \le N \rangle$ of $\kx$ contains no elements of order $2$. We saw in \exref{OqMmn4}
 that the CGL extension presentation of $R$ is symmetric and that all $\lambda^*_j = q^2$. Hence, \eqref{d-prop} is trivially satisfied with $d_i := 1$ for all $i \in \range(\eta)$. Finally, \exref{OqMmn8} shows that condition \eqref{pi-cond} is satisfied.

Therefore all hypotheses of \thref{cluster} are satisfied, and the theorem provides a quantum cluster algebra structure on $R$. The mutation matrix $\wt{B} = (b_{ij})$ for $R$ is easily computed:
$$
b_{(r-1)n +c, \, (r'-1)n +c'} = \begin{cases}
\pm1, &\text{if} \; \; r=r', \; c'=c\pm1  \\
 &\text{or} \; \; c=c', \; r'=r\pm1  \\
 &\text{or} \; \; r=r'\pm1, \; c=c'\pm1,  \\
0, &\text{otherwise}
\end{cases} \qquad
\begin{matrix}
\forall r,r' \in [1,m],  \\
c,c' \in [1,n] .
\end{matrix}
$$
\eex

\bre{torsion-free}
As noted in \reref{sqroots}, after a field extension 
every symmetric CGL extension can be brought to one 
that satisfies $\sqrt{\la_{lj}} \in \kx$ for all $1 \leq j < l \leq N$.
Recall that a CGL extension $R$ is called torsionfree if the subgroup of $\kx$ generated
by $\{ \la_{lj} \mid 1 \leq j < l \leq N \}$ is torsionfree.
The torsionfree symmetric CGL extensions form 
the largest and most interesting class of known symmetric CGL extensions.
All CGL extensions coming from the theory of quantum groups are in this class. 
For every torsionfree symmetric CGL extension $R$, 
one can always make a choice of the square roots 
$\nu_{lj}= \sqrt{\la_{lj}} \in \KK$, $1 \leq j < k \leq N$  
such that the subgroup of $\kx$ generated by all of them 
is torsionfree and thus does not contain $-1$ 
(i.e., condition \eqref{-1} is satisfied for all of them).
\ere

\bre{d-pr} All of the symmetric CGL extensions that we are aware of satisfy 
$$
\la_k^\sy = q^{m_k}, \; \forall k \in [1,N]
$$ 
for some non-root of unity $q \in \kx$ and $m_1, \ldots, m_N \in \Zset_{> 0}$.
It follows from \prref{la-equal} that $m_j = m_k$ for all $j,k \in \ex$
with $\eta(j) = \eta(k)$. This implies that all such CGL extensions
have the property \eqref{d-prop} where the integers $d_i$ are simply chosen
as $d_{\eta(k)} : = m_k$, $\forall k \in \ex$.
\ere

\bre{rescaleBtil} For applications, it is useful to determine the exchange matrix $\wt{B}$ \emph{before} the generators $x_k$ have been rescaled to satisfy \eqref{pi-cond}. This is possible because the rescaling does not change $\wt{B}$, as we next note.

Assume $R$ is a symmetric CGL extension satisfying all the hypotheses of \thref{cluster} \emph{except} \eqref{pi-cond}. Let the elements $y_k$, $\ol{y}_k$, the skew-symmetric matrix $\rbf$, and the toric frame $M$ be as in Section \ref{6.1}. Now suppose that we rescale the $x_k$ according to \prref{resc}, say with new generators $x'_1,\dots,x'_N$, to make \eqref{pi-cond} hold. Build the toric frame and its matrix for the new setting as in Section \ref{6.1}. Since the scalars $\la_{kj}$ and $\nu_{kj}$ do not change, the matrix $\rbf$ does not change. The new versions of the $\ol{y}_k$, call them $\ol{y}'_k$, are scalar multiples of the old ones, and the new toric frame, call it $M'$, satisfies $M'(e_k) = \ol{y}'_k$ for all $k$. Hence, for any vector $b \in \Zset^N$, the element $M'(b)$ is a scalar multiple of $M(b)$, and thus $\chi_{M'(b)} = \chi_{M(b)}$. Finally, the same elements $h_k, h^*_k \in \HH$ which enter into the symmetric CGL conditions for the original generators are used with respect to the new generators, which means that the scalars $\la_k$ and $\la^*_k$ do not change under the rescaling. 

Thus, the conditions in \thref{cluster} (a) which uniquely determine the columns of $\wt{B} = \wt{B}_{\id}$ are the same before and after rescaling. Therefore $\wt{B}$ does not change under the rescaling.
\ere
\bre{2questions} 
\thref{cluster} raises several questions and problems.

(a) For a symmetric quantum nilpotent algebra $R$, the elements of $\Xi_N$
give rise to seeds for the quantum cluster algebras structure on $R$ 
whose cluster variables are prime elements in some subalgebras of $R$.
Do the remaining seeds of $R$ have any similar ring-theoretic meaning? 

(b) In general \thref{cluster} constructs quantum cluster algebras whose 
quantum tori are multiparameter in the sense of Chapter \ref{qClust} rather
than one-parameter in the sense of \cite{BZ}. However, all examples of 
such multiparameter quantum cluster algebras that we are aware of 
come from 2-cocycle twists of one-parameter CGL extensions (e.g., quantum 
Schubert cell algebras). Are there any families of symmetric 
CGL extensions that are genuinely multiparameter, in the sense that 
they are not obtained from each other by 2-cocycle twists?
\ere
\section{Cluster variables}
\label{6.3}
The next result gives an explicit formula for the cluster variables 
that appear in \thref{cluster}.

\bpr{int-tau} Assume the setting of Theorem {\rm\ref{tcluster}}. 
Let $\tau \in \Xi_N$ and $k \in [1,N]$.
  
If $\tau(k) \geq \tau(1)$, then $\ol{y}_{\tau, k} = \ol{y}_{[p^m(\tau(k)), \tau(k)]}$, 
where 
$$
m = \max \{ n \in \Znn \mid p^n( \tau(k)) \in \tau([1, k]) \}.
$$

If $\tau(k) \leq \tau(1)$, then $\ol{y}_{\tau, k} = \ol{y}_{[\tau(k), s^m(\tau(k))]}$, where 
$$
m = \max \{ n \in \Znn \mid s^n( \tau(k)) \in \tau([1, k]) \}.
$$
Here the predecessor and successor functions are computed
with respect to the original CGL extension presentation 
\eqref{itOre} of $R$.
\epr

\bre{gen} \thref{cluster} (d) and \prref{int-tau} imply that 
the quantum cluster algebra $\Abb(M, \wt{B}, \varnothing)_\KK =R$ 
coincides with the subalgebra of $\Fract(R)$ generated by the cluster 
variables in the (finite!) set of toric frames $\{ M_\tau \mid \tau \in \Ga_N\}$.
\ere

\begin{proof}[Proof of \prref{int-tau}]
We prove the two cases of the proposition in parallel, by induction on $k$. If $k=1$, both cases are clear because $\ol{y}_{\tau,1} = x_{\tau(1)} = \ol{y}_{[\tau(1),\tau(1)]}$. Assume now that $k > 1$. For the induction step, we will restrict to the second of the 
two cases, i.e., the case where $\tau(k) < \tau(1)$. The first case is similar; it is left to the reader.

By \thref{sym-prime}, $\ol{y}_{\tau, k} = \xi\, \ol{y}_{[\tau(k), s^m(\tau(k))]}$
for some $\xi \in \kx$, so all we need to show is that the leading terms 
of the two prime elements are equal. From the fact that $\tau([1,j])$ is an 
interval for all $j \leq k$ and the assumption that $\tau(k) \leq \tau(1)$
it follows that 
$$
\tau([1,k]) = [\tau(k), \tau(i)] \quad 
\mbox{for some} \quad i \in [1,k].
$$
Hence, $\tau([1,k-1]) = \tau([1,k]) \backslash \{ \tau(k) \} = [\tau(k)+1, \tau(i)]$ and 
thus
$$
R_{\tau, k} = R_{[\tau(k), \tau(i)]} 
\quad
\mbox{and}
\quad
R_{\tau, k-1} = R_{[\tau(k)+1, \tau(i)]}.
$$
All leading terms in what follows will be computed 
with respect to the standard CGL extension presentation of  
$R_{[\tau(k), \tau(i)]}$ obtained by adjoining 
$x_{\tau(k)}, x_{\tau(k)+1}, \ldots,\allowbreak x_{\tau(i)-1}, x_{\tau(i)}$ 
in this order. By \thref{CGL} and \coref{steps} (b),
$$
y_{\tau, k} = y_{\tau, p_\tau(k) } x_{\tau(k)} - d_{k-1}
$$
for some $d_{k-1} \in R_{\tau, k-1}$, where $p_\tau$ 
is the predecessor function for the level sets of $\eta_\tau = \eta \tau : [1,N]\to \Zset$.

Suppose first that $\tau(p_\tau(k)) \le \tau(1)$, and note that $\tau([1,p_\tau(k)]) = [\tau(p_\tau(k)), \tau(i')]$ for some $i' \in [1, p_\tau(k)]$. By definition of $p_\tau(k)$, we have $\eta_\tau(j) \ne \eta_\tau(k)$, $\forall j \in [p_\tau(k)+1, k-1]$, and so $\tau(p_\tau(k)) = s(\tau(k))$. Hence,
$$
\max \{ n \in \Znn \mid s^n( \tau(p_\tau(k))) \in \tau([1, p_\tau(k)]) \} = m-1.
$$
The induction hypothesis thus implies that $\ol{y}_{\tau, p_\tau(k) } = \ol{y}_{[s(\tau(k)), s^m(\tau(k))]}$. On the other hand, if $\tau(p_\tau(k)) \ge \tau(1)$, we find that $\tau(p_\tau(k)) = s^m(\tau(k))$, and the induction hypothesis again implies $\ol{y}_{\tau, p_\tau(k) } = \ol{y}_{[s(\tau(k)), s^m(\tau(k))]}$.

Applying Eq. \eqref{fg} for $\Scr_{\nub_\tau}$ and $\Scr_\nub$ gives 
\begin{align*}
\lt(\ol{y}_{\tau, k}) &= \Scr_{\nub_\tau}( e^\tau_{[p^m_\tau(k),k]} ) \lt( y_{\tau, p_\tau(k)} x_{\tau(k)} )  \\
 &= \Scr_{\nub_\tau}( e^\tau_{[p^m_\tau(k),k]} ) \Om_\nub (e_{\tau(k)}, e_{[s(\tau(k)), s^m(\tau(k))]} )^{-2} \Scr_{\nub_\tau}( e^\tau_{[p^m_\tau(k),p_\tau(k)]} )^{-1} \times  \\
 &\qquad\qquad\qquad x_{\tau(k)} \lt( \ol{y}_{\tau, p_\tau(k)} )  \\
 &= \Om_\nub (e_{\tau(k)}, e_{[s(\tau(k)), s^m(\tau(k))]} )^{-1} x_{\tau(k)} \lt( \ol{y}_{[s(\tau(k)), s^m(\tau(k))]} )  \\
 &=  \lt( \ol{y}_{[\tau(k), s^m(\tau(k))]} ),
\end{align*}
where $e^\tau_{[p_\tau^m(k), p_\tau^l(k)]} = e_{p^m_\tau(k)} + e_{p^{m-1}_\tau(k)} +\cdots+ e_{p^l_\tau(k)}$ for $m \ge l \ge 0$. This completes the proof of the proposition.
\end{proof}

\bco{yint=clvar}
Assume the setting of Theorem {\rm\ref{tcluster}}. The cluster variables of the seeds $(M_\tau, \wt{B}_\tau)$ for $\tau$ in $\Xi_N$ or in $\Gamma_N$ are exactly the homogeneous prime elements $\ol{y}_{[i,j]}$. More precisely,
\begin{align*}
\{ M_\tau(e_k) \mid \tau \in \Xi_N, \; k \in [1,N] \} &= \{ M_\tau(e_k) \mid \tau \in \Ga_N, \; k \in [1,N] \} \\
&= \{ \ol{y}_{[i,j]} \mid 1 \le i \le j \le N, \; \eta(i) = \eta(j) \}.
\end{align*}
\eco

\begin{proof}
The second of the displayed sets is contained in the first a priori, and the first is contained in the third by \prref{int-tau} and the definition of the $M_\tau$. It remains to show that the third set is contained in the second. Thus, let $1 \le i \le j \le N$ with $\eta(i) = \eta(j)$. Then $i = p^m(j)$ where
$$
m = \max \{ n \in \Znn \mid p^n(j) \in [i,j] \}.
$$
Set $k' := j-i+1$, and choose $\tau \in \Ga_N$ as follows:
$$
\tau = \begin{cases}
\id = \tau_{1,1}  &(\text{if} \; i = 1)  \\
\tau_{i-1,j} = [i,\dots,j, i-1, j+1, \dots, N, i-1, \dots, 1]  &(\text{if} \; i > 1).
\end{cases}
$$
Then $\tau(1) = i \le j = \tau(k')$ and $\tau([1,k']) = [i,j]$. Consequently,
$$
m = \max \{ n \in \Znn \mid p^n(\tau(k')) \in \tau([1,k']) \},
$$
and \prref{int-tau} shows that
$$
\ol{y}_{\tau,k'} = \ol{y}_{[ p^m(\tau(k')), \tau(k')]} = \ol{y}_{[ i,j ]}.
$$
Therefore $\ol{y}_{[i,j]} = \ol{y}_{\tau, (\tau_\bu \tau)^{-1}(k)} = M_\tau(e_k)$ where $k := \tau_\bu \tau(k')$.
\end{proof}

\bex{OqMmn10}
Let $R = \OqM$, assume there exists $\sqrt q \in \kx$, and let $R$ have the quantum cluster algebra structure coming from \thref{cluster} as in \exref{OqMmn9}. The cluster variables of the seeds $(M_\tau, \wt{B}_\tau)$ in \thref{cluster} are exactly the solid quantum minors within $[1,m] \times [1,n]$. To see this, first recall from \exref{OqMmn5} that $\nu_{ij} = 1$ for all $i,j \in [1,N]$ with $\eta(i) = \eta(j)$. Consequently, \eqref{ybarpnll} yields $\ol{y}_{[i, s^l(i)]} = y_{[i, s^l(i)]}$ for all $i \in [1,N]$ and $l \in \Znn$ with $s^l(i) \ne +\infty$. These elements, by \coref{yint=clvar}, are exactly the cluster variables $M_\tau(e_k)$ for $\tau \in \Xi_N$ and $k \in [1,N]$. On the other hand, each $y_{[i, s^l(i)]}$ is a solid quantum minor by \exref{OqMmn6}. Conversely, for any $l \in \Znn$ and $r,r+l \in [1,m]$, $c,c+l \in [1,n]$, we have
$$
\De_{[r, r+l], [c, c+l]} = y_{[(r-1)n+c, \, s^l((r-1)n+c) ]} .
$$
\eex

\section{Auxiliary results}
\label{6.4}
In this section we establish two results that will 
be needed for the proof of \thref{cluster}. 
The first one uses \thref{CGLmuta} and \prref{int-tau} to construct 
mutations between pairs of the toric frames $M_\tau$ 
for $\tau \in \Xi_N$. The corresponding mutations of quantum seeds 
(\thref{cluster} (c)) are constructed in Section \ref{6.5}--\ref{6.6}.

For $g = \sum_j g_j e_j \in \Zset^N$ set
$$
\supp(g) := \{ j \in [1,N] \mid g_j \neq 0 \}.
$$
\index{suppg@$\supp(g)$}

\bpr{cluster-tau-ind} 
Let $R$ be a symmetric CGL extension of length $N$ and $\nu_{kj} = \sqrt{\la_{kj}} \in \kx$ for $1 \leq j < k \leq N$ such that \eqref{-1} is satisfied. Assume that the generators 
of $R$ are rescaled so that the condition \eqref{pi-cond} is satisfied.

Let $\tau, \tau' \in \Xi_N$ be such that
$$
\tau' = ( \tau(k), \tau(k+1)) \tau = \tau (k, k+1)
$$
for some $k \in [1,N-1]$ such that $\tau(k) < \tau(k +1)$.

{\rm(a)} If $\eta(\tau(k)) \neq \eta (\tau(k+1))$, then $M_{\tau'} = M_\tau$.

{\rm(b)} Let $\eta(\tau(k)) = \eta (\tau(k+1))$. Set $k_\bu: = \tau_\bu \tau(k)$. Then 
$k_\bu = (\tau')_\bu \tau'(k)$ and 
\begin{equation}
\label{Mtautau'}
M_{\tau'}(e_j)= 
\begin{cases}
M_\tau(e_j), &\mbox{if} \; \; j \neq k_\bu
\\
M_\tau(e_{p(k_\bu)} + e_{s(k_\bu)} - e_{k_\bu}) + 
M_\tau(g - e_{k_\bu}), &\mbox{if} \; \; j =k_\bu
\end{cases}
\end{equation}
for some $g \in \Znn^N$ 
such that $\supp(g) \cap \eta^{-1} \eta(k_\bu) = \varnothing$ and 
$|\supp (g) \cap \eta^{-1}(a)| \leq 1$ for all $a \in \Zset$. Furthermore, 
the vector $e_{p(k_\bu)} + e_{s(k_\bu)} -g \in \Zset^N$ satisfies the identities
\begin{align}
\label{Omtautau'1}
&\Om_{\rbf_\tau}( e_{p(k_\bu)} + e_{s(k_\bu)} -g , e_j) = 1, \; \; \forall j \ne k_\bu,
\\
\label{Omtautau'2}
&\Om_{\rbf_\tau}( e_{p(k_\bu)} + e_{s(k_\bu)} -g , e_{k_\bu})^2
= \la^\sy_{k_\bu},
\end{align}
and
\begin{equation}
\label{chiMub}
\chi_{M_\tau(e_{p(k_\bu)} + e_{s(k_\bu)} -g)}=1. 
\end{equation}
\epr

Note that the condition $\tau(k) < \tau(k+1)$ is not essential 
since, if $\tau(k) > \tau(k+1)$ and all other conditions are satisfied, 
then one can interchange the roles of $\tau$ and $\tau'$.

\begin{proof} Part (a) follows from \thref{1} (a) applied to the CGL extension 
presentation \eqref{tau-pres} of $R$ associated to $\tau$.

(b) Since $\tau([1,j])$ is an interval for all $j \leq k+1$ and $\tau(k) < \tau(k+1)$, 
we have $\tau'([1,k+1]) = \tau([1,k+1]) = [\tau(i), \tau(k+1)]$ and 
$\tau([1,k]) = [\tau(i), \tau(k+1) -1]$ for some $i \in [1,k]$.
On the other hand, the set
$$
\tau'([1,k]) = \tau'([1,k+1]) \setminus \{ \tau'(k+1) \} 
= \tau([1,k+1]) \setminus \{ \tau(k) \}
$$
must be 
also an interval, so $\tau(k)  = \tau(i)$ and $i=k$. Therefore,
\begin{multline}
\label{tau-ident}
\tau'([1,k+1]) = \tau([1,k+1]) = [\tau(k), \tau(k+1)], \\
\tau([1,k]) = [\tau(k), \tau(k+1) -1], \quad
\mbox{and} \quad
\tau'([1,k]) = [\tau(k)+1, \tau(k+1)].
\end{multline}
This implies that 
$$
\tau(k+1) = s^m(\tau(k)) \; \; \mbox{for some} \; \;  m \in \Zset_{>0}
$$ 
and 
\begin{equation}
\label{tau's}
\eta^{-1} (\eta \tau(k)) \cap \tau([1,k+1]) 
= \{\tau(k), s(\tau(k)), \ldots, \tau(k+1) = s^m(\tau(k)) \}.
\end{equation}
From the last identity and the fact that $\tau(j) = \tau'(j)$ for $j \neq k, k+1$
we infer
$$
\tau^{-1} \eta^{-1} (\eta \tau(k)) = 
(\tau')^{-1} \eta^{-1} (\eta \tau(k)).
$$
By the definition of the permutations $\tau_\bu$ and $(\tau')_\bu$,
$$
k_\bu = \tau_\bu \tau(k) = (\tau')_\bu \tau' (k) = s^{m-1}(\tau(k)).  
$$

From \thref{3} (b) we have $\wh{M}_{\tau'} (e_l) = \wh{M}_\tau(e_l)$ for $l \neq k$. 
Eq. \eqref{Mtautau'} for $j \neq k_\bu$ follows from this.

The identities in \eqref{tau-ident} and \prref{int-tau} give
\begin{equation}
\label{2oly}
\begin{aligned}
\ol{y}_{\tau,k} &= \ol{y}_{[\tau(k), s^{m-1}(\tau (k))]},  &\qquad\ol{y}_{\tau, k+1} &= \ol{y}_{[\tau(k), s^m(\tau(k))]},  \\
\ol{y}_{\tau', k} &= \ol{y}_{[s(\tau(k)),s^m(\tau(k))]}.
\end{aligned}
\end{equation}
Recall the definition of the toric frames $M_{[i,s^m(i)]}$ from Section \ref{4b.3}. We will construct an isomorphism $\dot w$ from $\Zset^{k+1}$ to $\Zset e_{\tau(k)} +\cdots+ \Zset e_{s^m(\tau(k))}$ such that the frames $M_{[\tau(k), s^m(\tau(k))]} \dot w$ and $\wh{M}_\tau$ agree on suitable $e_j$.

Recall the definition of the set $P_{[\tau(k), s^m(\tau(k))]}$ from \eqref{Pset}, set
\begin{align*}
A &:= \eta(P_{[\tau(k), s^m(\tau(k))]}) = \eta( \tau(1,k+1] ) \setminus \{ \eta(\tau(k)) \}  \\
Q &:= \{ j \in [1,k-1] \mid \eta(\tau(j)) \ne \eta(\tau(l)), \; \; \forall l \in [j+1,k+1] \},
\end{align*}
and note that $\eta \tau$ restricts to a bijection of $Q$ onto $A$. Thus,
\begin{equation}
\label{|QAP|}
|Q| = |A| = | P_{[\tau(k), s^m(\tau(k))]} |.
\end{equation}
The definition of $Q$ also ensures that
\begin{equation}
\label{etatauQ}
\{ t \in [\tau(k), s^m(\tau(k))] \mid \eta(t) = \eta(\tau(j)) \} \subseteq \tau([1,j]), \; \; \forall j \in Q.
\end{equation}

If $j \in Q$ and $\tau(j) \ge \tau(1)$, then $\tau([1,j]) = [\tau(i_j), \tau(j)]$ for some $i_j \in [1,j]$, and we observe that $\tau(j) \in P_{[\tau(k), s^m(\tau(k))]}$. Moreover, the integer $m_j$ corresponding to $j$ in \prref{int-tau} equals $O_-^{\tau(k)+1}(\tau(j))$, and hence we obtain
\begin{equation}
\label{tauj1}
\ol{y}_{\tau,j} = \ol{y}_{[p^{m_j}(\tau(j)), \tau(j)]} = M_{[\tau(k), s^m(\tau(k))]}(e_{\tau(j)}), \; \; \forall j \in Q \; \; \text{with} \; \; \tau(j) \ge \tau(1).
\end{equation}
On the other hand, if $\tau(j) < \tau(1)$, then $\tau([1,j]) = [\tau(j), \tau(i_j)]$ for some $i_j \in [1,j-1]$. Let $m_j$ denote the integer corresponding to $j$ in \prref{int-tau}, and observe that $s^{m_j}(\tau(j)) = \tau(j^-) \in P_{[\tau(k), s^m(\tau(k))]}$ for some $j^- \in [1,j]$. Moreover, $m_j = O_-^{\tau(k)+1}(\tau(j^-))$, and so
\begin{multline}
\label{tauj2}
\ol{y}_{\tau,j} = \ol{y}_{[p^{m_j}(\tau(j^-), \tau(j^-)]} = M_{[\tau(k), s^m(\tau(k))]}(e_{\tau(j^-)}),  \\\forall j \in Q \; \; \text{with} \; \; \tau(j) < \tau(1).
\end{multline}
In case $m > 1$, we set
$$
t := \max ( \tau^{-1} \{ s (\tau(k)), \ldots, s^{m-1}(\tau(k)) \}) = \max \{ j \in [1,k-1] \mid \eta(\tau(j)) = \eta(\tau(k)) \}
$$
(see \eqref{tau's}). Then either $\tau(t) = s^{m-1}(\tau(k))$ or $\tau(t) = s(\tau(k))$, and \prref{int-tau} yields
\begin{equation}
\label{ytaul}
\ol{y}_{\tau,t} = \ol{y}_{[s(\tau(k)), s^{m-1}(\tau(k))]} = M_{[\tau(k), s^m(\tau(k))]}( e_{s^{m-1}(\tau(k))} ), \; \; \text{if} \; \; m > 1.
\end{equation}

Now choose a bijection $w : [1,k+1] \rightarrow [\tau(k), s^m(\tau(k))]$ such that
\begin{align*}
w(k+1) &= s^m(\tau(k))  \\
w(k) &= \tau(k)  \\
w(t) &= s^{m-1}(\tau(k)), \; \; \text{if} \; \; m > 1  \\
w(j) &= \tau(j), \; \; \forall j \in Q \; \; \text{with} \; \; \tau(j) \ge \tau(1)  \\
w(j) &= \tau(j^-), \; \; \forall j \in Q \; \; \text{with} \; \; \tau(j) < \tau(1),
\end{align*}
and let $\dot w$ denote the isomorphism $\Zset^{k+1} \to \Zset e_{\tau(k)} +\cdots+ \Zset e_{s^m(\tau(k))}$ such that $\dot w(e_j) = e_{w(j)}$ for $j \in [1,k+1]$.
In particular, $\eta w(j) = \eta(\tau(j))$, $\forall j \in Q$, so $\eta w|_Q$ is injective. By construction, $w(Q) \subseteq P_{[\tau(k), s^m(\tau(k))]}$, and so we conclude from \eqref{|QAP|} that
\begin{equation}
\label{wQP}
w|_Q : Q \longrightarrow P_{[\tau(k), s^m(\tau(k))]} \quad\text{is a bijection.}
\end{equation}
Combining Eqs. \eqref{2oly} and \eqref{tauj1}--\eqref{ytaul} with the definition of $w$, we see that
\begin{equation}
\label{MhatMw}
\wh{M}_\tau(e_j) = \ol{y}_{\tau,j} = M_{[\tau(k), s^m(\tau(k))]}\dot w(e_j), \; \; \begin{cases}
\forall j \in Q \cup \{t,k,k+1\} &(m > 1) \\
\forall j \in Q \cup \{k,k+1\} &(m = 1).
\end{cases}
\end{equation}
Comparing the matrices of the toric frames $\wh{M}_\tau|_{\Zset^{k+1}}$ and $M_{[\tau(k), s^m(\tau(k))]}\dot w$ (recall \eqref{r}, \eqref{rint-tor}, and \eqref{Msig}), we find that
\begin{multline}
\label{rtaurother}
\bigl( \rbf(\wh{M}_\tau|_{\Zset^{k+1}}) \bigr)_{lj} = \bigl( ( \rbf_{[\tau(k), s^m(\tau(k))]} )_{\dot w} \bigr)_{lj},  \\
\begin{cases}
\forall l,j \in Q \cup \{t,k,k+1\} &(m > 1) \\
\forall l,j \in Q \cup \{k,k+1\} &(m = 1).
\end{cases}
\end{multline}

The next step is to apply \thref{CGLmuta}. We do the case $m > 1$ and leave the case $m = 1$ to the reader. (In the latter case, $p(k_\bu) = -\infty$ and $e_{p(k_\bu)} = 0$.) Observe that \eqref{MhatMw} and \eqref{rtaurother} together imply that
\begin{equation}
\label{Mhatf}
\wh{M}_\tau(f) = M_{[\tau(k), s^m(\tau(k))]}\dot w(f), \; \; \forall f \in \Zset^N \; \; \text{with} \; \; \supp(f) \subseteq Q \cup \{t,k,k+1\}.
\end{equation}
Thus, taking account of \eqref{2oly} and \eqref{wQP}, \thref{CGLmuta} implies that
$$
\ol{y}_{\tau',k} = 
\wh{M}_\tau(e_t + e_{k+1} - e_k) + \wh{M}_\tau(g'-e_k),
$$
where $g' \in \Zset_{\geq 0} e_1 + \cdots + \Zset_{\geq 0} e_{k-1}$ is such that 
$\supp(g') \subseteq Q$.
By the definition of $\tau_\bu$ and Eq. \eqref{tau's},
$$
\tau_\bu \tau(t) = p(k_\bu) \; \; 
\mbox{and} \; \; \tau_\bu \tau(k+1) = s(k_\bu).
$$
Therefore,
$$
M_{\tau'}(e_{k_\bu}) = \wh{M}_{\tau'}(e_k)
= \ol{y}_{\tau', k}
= M_\tau(e_{p(k_\bu)} + e_{s(k_\bu)} - e_{k_\bu}) + M_\tau(\tau_\bu \tau(g')-e_{k_\bu}),
$$
which implies the validity of \eqref{Mtautau'} for $j = k_\bu$.

Finally, the identities \eqref{Omtautau'1}--\eqref{Omtautau'2} follow from
\thref{scalars}, Eq. \eqref{Mhatf} and the fact that 
$R_{\tau,k+1} = R_{[\tau(k),s^m(\tau(k))]}$, 
see \eqref{tau-ident}. We note that $\eta(k_\bu) = \eta( \tau(k) )$ and $s(\tau(k)) \neq + \infty$, 
$s(k_\bu) \neq + \infty$,
which follow from the definition of $\tau_\bu$ and Eq. \eqref{tau's}. Because of 
this and \prref{la-equal}, $\chi_{x_{\tau(k)}}(h_{\tau(k)}^\sy)= 
\la_{\tau(k)}^\sy = \la_{k_\bu}^\sy$. The identity \eqref{chiMub} follows 
from the fact that $\ol{y}_{\tau, j}$ and $\ol{y}_{\tau',j}$ are 
$\HH$-eigenvectors for all $j \in [1,N]$ and Eq. \eqref{Mtautau'}.
\end{proof}

The next lemma proves uniqueness of integral vectors satisfying bilinear identities 
of the form \eqref{Omtautau'1}--\eqref{Omtautau'2} from strong $\HH$-rationality of CGL extensions.

\ble{unique} Assume that $R$ is a symmetric CGL extension of length $N$ and
$( \nu_{kj} ) \in M_N(\kx)$ is a multiplicatively skew-symmetric matrix with $\nu_{kj}^2 = \la_{kj}$, $\forall j,k \in [1,N]$. Then for all $\tau \in \Xi_N$, 
$\theta \in \xh$, and $\xi_1, \ldots, \xi_N \in \kx$, there exists 
at most one vector $b \in \Zset^N$ such that $\chi_{M_\tau(b)}= \theta$ and 
$\Om_{\rbf_\tau}(b, e_j)^2 = \xi_j$, $\forall j \in [1,N]$.
\ele
\begin{proof}
Let $b_1, b_2 \in \Zset^N$ be such that $\chi_{M_\tau(b_1)} = \chi_{M_\tau}(b_2)= \theta$ and 
$\Om_{\rbf_\tau}(b_1, e_j)^2 = \Om_{\rbf_\tau}(b_2, e_j)^2$  $= \xi_j$, $\forall j \in [1,N]$.
Then $M_\tau(b_1) M_\tau(b_2)^{-1}$ commutes with $M_\tau(e_j)$ for all $j \in [1,N]$. 
This implies that $M_\tau(b_1) M_\tau(b_2)^{-1}$ belongs to the center of $\Fract(R)$, because
by \prref{CGLcluster} $\Fract(R)$ is generated (as a division algebra) by 
$M_\tau(e_1)^{\pm1}, \ldots, M_\tau(e_N)^{\pm1}$.
Furthermore, 
$$
\chi_{M_\tau(b_1) M_\tau(b_2)^{-1}} =1.
$$
By the strong $\HH$-rationality 
of the $0$ ideal of a CGL extension \cite[Theorem II.6.4]{BG},
$$
Z(\Fract(R))^\HH = \KK,
$$
where $Z(.)$  \index{Z@$Z(.)$}  stands for the center of an algebra and $(.)^\HH$ for the subalgebra fixed by $\HH$.
Hence, $M_\tau(b_1) M_\tau(b_2)^{-1} \in \KK$, which is only possible if $b_1=b_2$.
\end{proof}

\section{An overview of the proof of \thref{cluster}}
\label{6.5a}
In this section we give a summary of the strategy of 
our proof of \thref{cluster}.

In Section \ref{6.1} we constructed quantum frames $M_\tau : \Zset^N \to \Fract(R)$ 
associated to the elements 
of the set $\Xi_N$. In order to extend them to quantum seeds of $\Fract(R)$, 
one needs to construct a compatible matrix $\wt{B}_\tau \in M_{N \times |\ex|}(\Zset)$
for each of them.
This will be first done for the subset $\Ga_N$ of $\Xi_N$ 
in an iterative fashion with respect to the linear ordering \eqref{sequence}.
If $\tau$ and $\tau'$ are two consecutive 
elements of $\Ga_N$ in that linear ordering, then $\tau' = \tau(k,k+1)$ 
for some $k \in [1,N]$ such that $\tau(k) < \tau(k+1)$. 
If $\eta(\tau(k)) \neq \eta(\tau(k+1))$ then $M_{\tau'} = M_\tau$ 
by \prref{cluster-tau-ind} (a)
and nothing happens at that step. If $\eta(\tau(k)) = \eta(\tau(k+1))$,  
then we use \prref{cluster-tau-ind} (b) to construct $b_\tau^{k_\bu}$ and $b_{\tau'}^{k_\bu}$ 
where $k_\bu := (\tau_\bu \tau)(k)$. Up to $\pm$ sign these vectors 
are equal to $e_{p(k_\bu)} + e_{s(k_\bu)} -g$, where $g \in \Znn^N$ 
is the vector from \prref{cluster-tau-ind} (b). Then we use ``reverse'' mutation to construct $b_\sigma^{k_\bu}$ for 
$\sigma \in \Ga_N$, $\sigma \prec \tau$ in the linear order \eqref{sequence}. 
Effectively this amounts to starting with a quantum cluster algebra in which all 
variables are frozen and then recursively adding more exchangeable variables.

There are two things that could go wrong with this.
Firstly, the reverse mutations from 
different stages might not be synchronized. Secondly, there are many pairs 
of consecutive elements $\tau, \tau'$ for which $k_\bu$ is the same. 
So we need to prove that $b_\sigma^{k_\bu}$ is not {\em{overdetermined}}.
We use strong rationality of CGL extensions to handle both via \leref{unique}.
This part of the proof (of parts (a) and (b) of \thref{cluster}) 
is carried out in Section \ref{6.5}.

Once $\wt{B}_\id$ is (fully) constructed then the $\wt{B}_\tau$ are constructed 
inductively by applying \prref{cluster-tau-ind} and using the sequences of elements 
of $\Xi_N$ from \coref{steps} (a). At each step \leref{unique} is applied 
to match columns of mutation matrices. This proves parts (a) and (b) 
of \thref{cluster}. Part (c) of the theorem is obtained in a somewhat 
similar manner from
\prref{cluster-tau-ind}. This is done in Section \ref{6.6}.

The last two parts (d)--(e) of \thref{cluster} are proved in Section \ref{6.7}. For each  
$\tau \in \Xi_N$ we denote by $E_\tau$ the multiplicative subset of $R$ 
generated by $\kx$ and $M_\tau(e_j)$ for $j \in \ex$. An application 
of Section \ref{5.1} gives that it is an Ore subset of $R$. The idea of the proof 
of \thref{cluster} (d) is to obtain the following chain of embeddings
$$
R \subseteq \Abb(M, \wt{B}, \varnothing)_\KK \subseteq \UU (M, \wt{B}, \varnothing)_\KK
\subseteq \bigcap_{\tau \in \Ga_N} R[E_\tau^{-1}] = R,
$$
which clearly implies the desired result. The first inclusion follows from the fact that 
each generator $x_j$ of $R$ is a cluster variable for the toric frame $M_\tau$ associated 
to some $\tau \in \Ga_N$. The hardest is the last equality. It is derived 
from \thref{division} by using that the consecutive toric frames associated to the sequence 
of elements \eqref{sequence} are obtained from each other by one-step mutations and 
each $j \in \ex$ gets mutated at least one time along the way. Part (e) of 
\thref{cluster} is proved in a similar fashion.

Both parts (d) and (e) follow from results about intersections of localizations of symmetric CGL extensions. These intersection results, which do not require any rescaling of variables, are proved in Section \ref{prequantumcluster}.

\section{Recursive construction of quantum seeds for $\tau \in \Ga_N$}
\label{6.5}
Recall the linear ordering \eqref{sequence} on $\Ga_N \subset \Xi_N$. 
We start by constructing a chain of subsets $\ex_\tau \subseteq \ex$ \index{extau@$\ex_\tau$}
indexed by the elements of $\Ga_N$ such that 
$$
\ex_\id = \varnothing, \quad \ex_{w_\ci} = \ex \quad
\mbox{and} \quad
\ex_\sigma \subseteq \ex_\tau, \; \; \forall \sigma, \tau \in \Ga_N, \;
\sigma \prec \tau.
$$
This is constructed inductively by starting with $\ex_\id = \varnothing$.
If $\tau \prec \tau'$ are two consecutive elements in the linear ordering, 
then for some $1 \leq i < j \leq N$,
$$
\tau = \tau_{i,j-1}  \; \; 
\mbox{and} \; \; 
\tau'=\tau_{i,j}
$$
(recall \eqref{tauij} and the equalities in \eqref{sequence}).
Thus,
$$
\tau' = (ij) \tau = \tau(j-i, j-i+1).
$$
Assuming that $\ex_\tau$ has been constructed, we define
$$
\ex_{\tau'} := \begin{cases}
\ex_\tau \cup \{ p(j) \},  &\mbox{if} \; p(i) = -\infty, \; \eta(i) = \eta(j)  \\
\ex_\tau,  &\mbox{otherwise}.
\end{cases}
$$
It is clear that this process ends with $\ex_{w_\ci} = \ex$.

For a subset $X \subseteq \ex$, by an $N \times X$ matrix we will mean 
a matrix of size $N \times |X|$ whose columns are indexed by the set
$X$. The set of such matrices with integral entries will be denoted 
by $M_{N \times X} (\Zset)$. \index{MNXZ@$M_{N \times X} (\Zset)$}
(Recall from Section \ref{3.2} that the columns of all of our $N \times |\ex|$ matrices are indexed by 
$\ex \subset [1,N]$.) The following lemma provides an inductive 
procedure for establishing \thref{cluster} (a), (b) for $\tau \in \Ga_N$.

\ble{ind-ab} Assume that $R$ is a symmetric CGL extension of length $N$ 
satisfying \eqref{d-prop} and
$\{ \nu_{kj} \mid 1 \leq j < k \leq N \} \subset \kx$ is a set of square roots 
of $\la_{kj}$ which satisfies \eqref{-1}. Assume also that the generators 
of $R$ are rescaled so that the condition \eqref{pi-cond} is satisfied.

Let $\tau \in \Ga_N$. For all $\sigma \in \Ga_N$, $\sigma \preceq \tau$,
there exists a unique matrix $\wt{B}_{\sigma,\tau} \in M_{N \times \ex_\tau}(\Zset)$ 
whose columns $b_{\sigma, \tau}^l \in \Zset^N$, $l \in \ex_\tau$ satisfy 
\begin{equation}
\label{sigtau}
\begin{aligned}
\chi_{M_\sigma(b_{\sigma,\tau}^l)} &=1,  
&\quad\Om_{\rbf_\sigma} ( b^l_{\sigma,\tau}, e_j) = 1, \; \; \forall j \in [1,N], \; j \neq l  \\
\Om_{\rbf_\sigma} (b^l_{\sigma,\tau}, e_l)^2 &= \la_l^*
\end{aligned}
\end{equation}
for all $l \in \ex_\tau$.
The principal part of the matrix $\wt{B}_{\sigma, \tau}$ is skew-symmetrizable 
via the integers $\{ d_{\eta(k)} \mid k \in \ex_\tau\}$.
\ele  

\begin{proof} The uniqueness statement follows at once from \leref{unique}. 
If a matrix $\wt{B}_{\sigma, \tau}$ with such properties exists, 
then its principal part is skew-symmetrizable 
by \leref{alternative} and the condition \eqref{d-prop}, since the conditions \eqref{sigtau} imply that $(\rbf_\sigma, \wt{B}_{\sigma,\tau})$ is a compatible pair.

What remains to be proved is the existence statement in the 
lemma. It trivially holds for $\tau =\id$ since 
$\ex_\id = \varnothing$.

Let $\tau \prec \tau'$ be two consecutive elements of $\Ga_N$ 
in the linear ordering \eqref{sequence}. Assuming that the 
existence statement in the lemma holds for $\tau$, we will show 
that it holds for $\tau'$. The lemma will then follow by induction.

From our inductive assumption, $(\rbf_\tau, \wt{B}_{\tau,\tau})$ is a compatible pair, and consequently $(M_\tau, \wt{B}_{\tau,\tau})$ is a quantum seed.

As noted above, for some $1 \leq i < j \leq N$ we have
$\tau = \tau_{i,j-1}$ and $\tau'=\tau_{i,j}$. In particular, 
$\tau' = (ij) \tau = \tau(j-i, j-i+1)$ and $\tau(j-i) = i < j = \tau(j-i+1)$,
so \prref{cluster-tau-ind} is applicable to the pair $(\tau, \tau')$, with $k := j-i$.

If $\eta(i) \neq \eta(j)$, then $\ex_{\tau'} = \ex_\tau$ 
and $M_{\tau'} = M_\tau$ (by \prref{cluster-tau-ind} (a)).
So, $\Om_{\rbf_{\tau'}} = \Om_{\rbf_\tau}$.
These identities imply that the following matrices have the properties 
\eqref{sigtau} for the element $\tau' \in \Ga_N$: $\wt{B}_{\sigma,\tau'} := \wt{B}_{\sigma, \tau}$ 
for $\sigma \preceq \tau$ and $\wt{B}_{\tau',\tau'} := \wt{B}_{\tau, \tau}$.

Next, we consider the case $\eta(i) = \eta(j)$. This  implies that 
$j = s^m(i)$ for some $m \in \Zset_{>0}$. This fact and the definition of $\tau_\bu$
give that the element $\tau_\bu \tau(j-i)$ equals the $m$-th element 
of $\eta^{-1}(\eta(i))$ when the elements in the preimage are ordered from least to greatest.
Therefore this element is explicitly given by
\begin{equation}
\label{j-i}
\tau_\bu \tau(j-i) = s^{m-1} p^{O_-(i)}(i).
\end{equation}
Now set
$$
k_\bu:= \tau_\bu \tau(j-i)
$$
as in \prref{cluster-tau-ind} 
(b).
There are two subcases: (1) $p(i) \ne -\infty$ and (2) $p(i) = -\infty$.

Subcase (1). In this situation $\ex_{\tau'} = \ex_{\tau}$, so we do not need to generate an ``extra column'' for each matrix.
Set $\wt{B}_{\sigma, \tau'}:= \wt{B}_{\sigma, \tau}$ for $\sigma \in \Ga_N$, 
$\sigma \preceq \tau$. Eq. \eqref{sigtau} for the pairs $(\sigma, \tau')$ 
with $\sigma \preceq \tau$ follows from the equality $\ex_{\tau'} = \ex_\tau$.

Next we deal with the pair $(\sigma=\tau', \tau')$.
Applying the inductive assumption \eqref{sigtau} for $\wt{B}_{\tau, \tau}$ and \prref{cluster-tau-ind} 
(b) shows that the vector $g \in \Znn^N$
has the properties 
$$
\chi_{M_\tau(b_{\tau,\tau}^{k_\bu})}
= \chi_{M_\tau(e_{p(k_\bu)}+e_{s(k_\bu)} -g)}  
$$
and 
$$
\Om_{\rbf_\tau}(b_{\tau,\tau}^{k_\bu}, e_t)^2 = 
\Om_{\rbf_\tau}( e_{p(k_\bu)}+e_{s(k_\bu)} -g, e_t)^2, \; \; \forall t \in [1,N].
$$
\leref{unique} implies that $e_{p(k_\bu)}+e_{s(k_\bu)} -g = b_{\tau,\tau}^{k_\bu}$.
It follows from this, \coref{inv} (a), and Eq. \eqref{Mtautau'} that $\mu_{k_\bu}(M_\tau)(e_l) = M_{\tau'}(e_l)$ for all $l \in [1,N]$. Consequently,
$$
\Om_{\rbf_{\tau'}}(e_t, e_l)^2 M_{\tau'}(e_l) M_{\tau'}(e_t) = M_{\tau'}(e_t) M_{\tau'}(e_l) = 
\Om_{\mu_{k_\bu}(\rbf_\tau)}(e_t, e_l)^2 M_{\tau'}(e_l) M_{\tau'}(e_t)$$
for all $t, l \in [1,N]$, so $\Om_{\rbf_{\tau'}}(e_t, e_l)^2 = \Om_{\mu_{k_\bu}(\rbf_\tau)}(e_t, e_l)^2$, for all $j, l \in [1,N]$.
This is an equality in the subgroup of $\kx$ generated by $\{ \nu_{lt} \mid 1 \leq t < l \leq N \}$
and the condition \eqref{-1} implies that 
\begin{equation}
\label{Omtt}
\Om_{\rbf_{\tau'}}(e_t, e_l) = \Om_{\mu_{k_\bu}(\rbf_\tau)}(e_t, e_l), \; \; \forall t, l \in [1,N].
\end{equation}
Therefore $\rbf_{\tau'} = \mu_{k_\bu}(\rbf_\tau)$ and $M_{\tau'} = \mu_{k_\bu}(M_\tau)$.
We set $\wt{B}_{\tau', \tau'}: = \mu_{k_\bu}(\wt{B}_{\tau, \tau})$. \leref{1eig}
implies that the columns of $\wt{B}_{\tau', \tau'}$ have the property $\chi_{M_{\tau'}(b^l_{\tau', \tau'})}=1$
for all $l \in \ex_\tau = \ex_{\tau'}$. \prref{pair-mut} (b) implies that $(\rbf_{\tau'}, \wt{B}_{\tau',\tau'})$ is a compatible pair and
$$
\Om_{\rbf_{\tau'}} ( b^l_{\tau',\tau'}, e_t) = 1, \; \; \forall t \in [1,N], \; t \neq l
\quad \mbox{and} \quad
\Om_{\rbf_{\tau'}} (b^l_{\tau',\tau'}, e_l)^2 = \la_l^*
$$
for all $l \in \ex_{\tau'}$,
which completes the proof of the inductive step of the lemma in this subcase. 

Subcase (2). In this case, $\ex_{\tau'} = \ex_\tau \sqcup \{k_\bu\}$ and $k_\bu = s^{m-1}(i) = p(j)$. Define 
the matrix $\wt{B}_{\tau, \tau'}$ by 
$$
b^l_{\tau, \tau'} = 
\begin{cases}
b^l_{\tau, \tau}, & \mbox{if} \; \; l \neq k_\bu
\\
e_{p(k_\bu)} + e_{s(k_\bu)} - g, & \mbox{if} \; \; 
l = k_\bu,
\end{cases}
$$
where $g \in \Znn^N$ is the vector from \prref{cluster-tau-ind} (b).
Applying the assumption \eqref{sigtau} for $\wt{B}_{\tau, \tau}$ and \prref{cluster-tau-ind} (b), we obtain 
that the matrix $\wt{B}_{\tau, \tau'}$ has the properties \eqref{sigtau}. 
We set $\wt{B}_{\tau', \tau'}:= \mu_{k_\bu} (\wt{B}_{\tau, \tau'})$. As in subcase (1),
using \leref{1eig} and \prref{pair-mut} (b),  
one derives that $\wt{B}_{\tau', \tau'}$ satisfies the properties \eqref{sigtau}.

We are left with constructing $\wt{B}_{\sigma, \tau'} \in M_{N \times \ex_{\tau'}}(\Zset)$ 
for $\sigma \in \Ga_N$, $\sigma \prec \tau$.
We do this by a downward induction on the linear ordering \eqref{sequence} in a fashion that 
is similar to the proof of the lemma in the subcase (1). Assume that $\sigma \prec \sigma'$ is a pair of consecutive 
elements of $\Ga_N$ such that $\sigma' \preceq \tau$. As in the beginning of the section, 
we have that for some $1 \leq i\spcheck < j\spcheck \leq N$,
$$
\sigma = \tau_{i\spcheck,j\spcheck-1}  \; \; 
\mbox{and} \; \; 
\sigma'=\tau_{i\spcheck,j\spcheck},
$$
so
$$
\tau' = (i\spcheck j\spcheck) \tau = \tau(j\spcheck-i\spcheck, j\spcheck-i\spcheck+1).
$$
Assume that there exists a matrix $\wt{B}_{\sigma', \tau'} \in M_{N \times \ex_{\tau'}} (\Zset)$ 
that satisfies \eqref{sigtau}. We define the matrix $\wt{B}_{\sigma, \tau'} \in M_{N \times \ex_{\tau'}} (\Zset)$
by 
$$
\wt{B}_{\sigma, \tau'}:= 
\begin{cases}
\wt{B}_{\sigma', \tau'}, & \mbox{if} \; \; \eta(i\spcheck) \neq \eta(j\spcheck)
\\
\mu_{k\spcheck_\bu}(\wt{B}_{\sigma, \tau'}), & \mbox{if} \; \; 
\eta(i\spcheck) = \eta(j\spcheck),
\end{cases}
$$
where
$$
k\spcheck_\bu := \sigma_\bu \sigma(j\spcheck-i\spcheck).
$$
Analogously to the proof of the lemma in the subcase (1), 
using Propositions \ref{ppair-mut} (b) and \ref{pcluster-tau-ind}
and Lemmas \ref{l1eig} and \ref{lunique}, 
one proves that the matrix $\wt{B}_{\sigma, \tau'}$ has the properties \eqref{sigtau}.
This completes the proof of the lemma.
\end{proof}

\bre{restr} It follows from \leref{unique} that the matrices $\wt{B}_{\sigma, \tau}$ 
in \leref{ind-ab} have the following restriction property:

{\em{For all triples $\sigma \preceq \tau \prec \tau'$ of elements of $\Ga_N$,
$$
b^l_{\sigma, \tau} = b^l_{\sigma, \tau'}, \; \; \forall l \in \ex_\tau.
$$
In other words, $\wt{B}_{\sigma, \tau}$ is obtained from $\wt{B}_{\sigma, \tau'}$ by removing 
all columns indexed by the set $\ex_{\tau'} \backslash \ex_\tau$.}}

This justifies that \leref{ind-ab} gradually enlarges 
a matrix $\wt{B}_{\sigma, \sigma} \in M_{N \times \ex_\sigma}(\Zset)$ 
to a matrix $\wt{B}_{\sigma, w_0} \in M_{N \times \ex}(\Zset)$, 
for all $\sigma \in \Ga_N$. In the case of $\sigma = \id$, we start 
with an empty matrix ($\ex_\id = \varnothing$) and obtain a matrix $\wt{B}_{\id, w_\ci} \in M_{N \times \ex}(\Zset)$
which will be the needed mutation matrix for the initial toric frame $M_\id$.
\ere

\begin{proof}[Proof of \thref{cluster} {\rm(a), (b)} for $\tau \in \Ga_N$] Change $\tau$ to $\sigma$ in these statements.
These parts of the theorem for the elements of $\Ga_N$ follow from
\leref{ind-ab} applied to $(\sigma,\tau) = (\sigma,w_\ci)$. For all $\sigma \in \Ga_N$ we set 
$\wt{B}_\sigma := \wt{B}_{\sigma, w_\ci}$ and use that $\ex_{w_\ci} = \ex$. 
\end{proof}

\newcommand{\tcluster}{tcluster}
\section{Proofs of parts (a), (b) and (c) of Theorem \ref{\tcluster}}
\label{6.6}
Next we establish \thref{cluster} (a)--(b) in full generality. 
This will be done by using the result of Section \ref{6.5} for $\tau=\id$ 
and iteratively applying the following proposition.

\bpr{tau-seed-muta} Let $R$ be a symmetric CGL extension of length $N$ 
satisfying \eqref{d-prop} and
$\{ \nu_{kj} \mid 1 \leq j < k \leq N \} \subset \kx$ a set of square roots 
of $\la_{kj}$ which satisfies \eqref{-1}. Assume also that the generators 
of $R$ are rescaled so that the condition \eqref{pi-cond} is satisfied.
Let $\tau, \tau' \in \Xi_N$ be such that
$$
\tau' = ( \tau(k), \tau(k+1)) \tau = \tau (k, k+1)
$$
for some $k \in [1,N-1]$ and $\tau'(k) < \tau'(k +1)$,
$\eta(\tau(k)) = \eta(\tau(k+1))$. Set $k_\bu := \tau_\bu \tau(k)$.

Assume that there exists an $N \times |\ex|$ matrix $\wt{B}_\tau$ with integral 
entries whose columns 
$b_\tau^l \in \Zset^N$, $l \in \ex$ satisfy 
\begin{equation}
\label{B-Om}
\begin{aligned}
\chi_{M_\tau(b_\tau^l)} &=1,
&\quad\Om_{\rbf_\tau} ( b^l_\tau, e_j) &= 1, \; \; \forall j \in [1,N], \; j \neq l  \\
\Om_{\rbf_\tau} (b^l_\tau, e_l)^2 &= \la_l^*
\end{aligned}
\end{equation}
for all $l \in \ex$. Then its principal part is skew-symmetrizable 
and the columns $b_{\tau'}^j \in \Zset^N$, $j \in \ex$ 
of the matrix $\mu_{k_\bu}(\wt{B_\tau})$ satisfy  
\begin{equation}
\label{B-Om'}
\begin{aligned}
\chi_{M_{\tau'}(b_{\tau'}^l)} &=1, 
&\quad \Om_{\rbf_{\tau'}} ( b^l_{\tau'}, e_j) &= 1, \; \; \forall j \in [1,N], \; j \neq l  \\
\Om_{\rbf_{\tau'}} (b^l_{\tau'}, e_l)^2 &= \la_l^*
\end{aligned}
\end{equation}
for all $l \in \ex$. Furthermore,
\begin{equation}
\label{rtautau'}
\rbf_{\tau'} = \mu_k(\rbf_\tau)
\end{equation}
and
\begin{equation}
\label{bbg}
b_\tau^{k_\bu} = -b_{\tau'}^{k_\bu} = e_{p(k_\bu)} + e_{s(k_\bu)} - g,
\end{equation}
where $g \in \Znn^N$
is the vector from Proposition {\rm\ref{pcluster-tau-ind}}. 
\epr

The conditions \eqref{B-Om}--\eqref{B-Om'} are stronger than saying that 
the pairs $(\rbf_\tau, \wt{B}_\tau)$ and $(\rbf_{\tau'},\mu_{k_\bu}(\wt{B_\tau}))$ 
are compatible. The point is to recursively establish a 
stronger condition which matches the setting of \leref{unique} 
in order to use the uniqueness conclusion of the lemma.

\bre{diff}
The statements of \leref{ind-ab} and \prref{tau-seed-muta} have many similarities 
and their proofs use similar ideas. However, we note that there is a
principal difference between the two results. In the former case we have no 
mutation matrices to start with and we use \prref{cluster-tau-ind} (b) to gradually add
columns. In the latter case we already have a mutation matrix for one 
toric frame and construct a mutation matrix for another toric frame. 
\ere

\begin{proof}[Proof of \prref{tau-seed-muta}]
The fact that the principal part of $\wt{B}_\tau$ is skew-symmetrizable 
follows from \leref{alternative} and the condition \eqref{d-prop}. The assumptions 
on $\wt{B}_\tau$ and \prref{cluster-tau-ind} (b) imply
$$
\chi_{M_\tau(b_\tau^{k_\bu})} = \chi_{M_\tau(e_{p(k_\bu)} + e_{s(k_\bu)} - g)} 
$$
and
$$
\Om_{\rbf_\tau} (b^{k_\bu}_\tau, e_j)^2 = \Om_{\rbf_\tau}(e_{p(k_\bu)} + e_{s(k_\bu)} - g, e_j)^2, 
\; \forall j \in [1,N].
$$ 

By \leref{unique}, $b^{k_\bu}_\tau = e_{p(k_\bu)} + e_{s(k_\bu)} - g$. The 
mutation formula for $\mu_{k_\bu}(\wt{B}_\tau)$ also gives that $b^{k_\bu}_{\tau'} = - b^{k_\bu}_\tau$,
so we obtain \eqref{bbg}. Analogously to the proof of \eqref{Omtt}, one shows that
$$
\Om_{\rbf_{\tau'}}(e_j, e_l) = \Om_{\mu_{k_\bu}(\rbf_\tau)}(e_j, e_l), \; \; \forall j, l \in [1,N],
$$
which is equivalent to \eqref{rtautau'}. Finally, all identities in \eqref{B-Om'}
follow from the general mutation facts in \prref{pair-mut} (b) and \leref{1eig}.
\end{proof}

\begin{proof}[Proof of \thref{cluster} {\rm(a), (b)} for all $\tau \in \Xi_N$]
Similarly to the proof of \leref{ind-ab}, 
the uniqueness statement in part (a) follows from \leref{unique}.
We will prove the existence statement in part (a) by an inductive argument on $\tau$. 
Once the existence of the matrix $\wt{B}_\tau$ with the stated properties 
is established, the fact that the principal part of $\wt{B}_\tau$ is 
skew-symmetrizable follows from \leref{alternative} and 
the condition \eqref{d-prop}. Hence, $(M_\tau, \wt{B}_\tau)$ 
is a quantum seed and this yields part (b) of the theorem
for the given $\tau \in \Xi_N$.

For the existence statement in part (a)
we fix $\tau \in \Xi_N$. By \coref{steps} (a), there exists a sequence
$\tau_0 = \id, \tau_1, \ldots, \tau_n = \tau$ in $\Xi_N$ with the property that
for all $l \in [1,n]$,
$$
\tau_l = ( \tau_{l-1}(k_l), \tau_{l-1} (k_l+1)) \tau_{l-1} = \tau_{l-1}(k_l, k_l+1)
$$
for some $k_l \in [1,N-1]$ such that $\tau_{l-1}(k_l) < \tau_{l-1}(k_l +1)$.
In Section \ref{6.5} we established the validity of \thref{cluster} (a) for the identity element 
of $S_N$. By induction on $l$ we prove the validity of \thref{cluster} (a) for $\tau_l$. 
If $\eta(\tau_{l-1}(k_l)) \neq \eta(\tau_l(k_l))$, then \prref{cluster-tau-ind} (a) implies 
that $M_{\tau_l} = M_{\tau_{l-1}}$ and we can choose $\wt{B}_{\tau_l} = \wt{B}_{\tau_{l-1}}$.
If $\eta(\tau_{l-1}(k_l)) = \eta(\tau_l(k_l))$, then \prref{tau-seed-muta} proves that 
the validity of \thref{cluster} (a) 
for $\tau_{l-1}$ implies the validity of \thref{cluster} (a) 
for $\tau_l$. In this case $\wt{B}_{\tau_l} := \mu_{(k_l)_\bu} (\wt{B}_{\tau_{l-1}})$, 
where $(k_l)_\bu = ( (\tau_{k_l})_\bu  \tau_{k_l} )(k_l)$.
This completes the proof of \thref{cluster} (a) and (b). 
\end{proof}

\begin{proof}[Proof of \thref{cluster} {\rm(c)}] The one-step mutation statement 
in part (c) of \thref{cluster} and \coref{steps} (a) imply that all quantum seeds 
associated to the elements of $\Xi_N$ are mutation-equivalent 
to each other.

In the rest we prove the one-step mutation statement in part (c) 
of the theorem.
If $\eta(\tau(k)) \neq \eta(\tau(k+1))$, then the 
statement follows from \prref{cluster-tau-ind} (a).

Now let $\eta(\tau(k)) = \eta(\tau(k+1))$. We have that either $\tau(k) < \tau(k+1)$ or
$\tau'(k) = \tau(k+1) < \tau'(k+1) = \tau(k)$. In the first case we 
apply \prref{tau-seed-muta} to the pair $(\tau, \tau')$ and 
in the second case to the pair $(\tau', \tau)$. The  one-step mutation statement in \thref{cluster} (c) 
follows from this, the uniqueness statement in part (a) of the theorem
and the involutivity of mutations of quantum seeds (\coref{inv} (b)).
\end{proof}

\section{Intersections of localizations}
\label{prequantumcluster}

Parts (d) and (e) of \thref{cluster} involve showing that $R$ and $R[y_k^{-1} \mid k \in \inv]$ are equal to intersections of appropriate localizations. Proving this does not require either the existence of square roots of the $\la_{kj}$ in $\kx$ or the condition \eqref{pi-cond}, and leads to a general result of independent interest, showing that all symmetric CGL extensions are intersections of partially localized quantum affine space algebras. 

Throughout this section, we will assume that $R$ is a symmetric CGL extension of rank $N$ 
as in \deref{CGL}, and we fix a function $\eta : [1,N] \to \Zset$ satisfying the conditions of \thref{CGL}. For each $\tau \in \Xi_N$, there is a CGL presentation
$$
R = \KK [x_{\tau(1)}] [x_{\tau(2)}; \sigma''_{\tau(2)}, \delta''_{\tau(2)}] 
\cdots [x_{\tau(N)}; \sigma''_{\tau(N)}, \delta''_{\tau(N)}]
$$
as in \eqref{tauOre}. Let $p_\tau$ and $s_\tau$ denote the predecessor and successor functions for the level sets of $\eta_\tau := \eta \tau$, which by \coref{steps} (b) can be chosen as the $\eta$-function for the presentation \eqref{tauOre}. Let $y_{\tau,1}, \dots, y_{\tau, N}$ be the corresponding sequence of homogeneous prime elements of $R$ from \thref{CGL}, and denote
\begin{align*}
\Abb_\tau &:= \text{the} \; \KK\text{-subalgebra of} \; R \; \text{generated by} \; \{ y_{\tau,k} \mid k \in [1,N] \}  \\
\Tbb_\tau &:= \text{the} \; \KK\text{-subalgebra of} \; \Fract(R) \; \text{generated by} \; \{ y_{\tau,k}^{\pm1} \mid k \in [1,N] \}  \\
E_\tau &:= \text{the multiplicative subset of} \; \Abb_\tau \; \text{generated by}  \\
 & \qquad\qquad\qquad \qquad\qquad\qquad \kx \sqcup \{ y_{\tau,k} \mid k \in [1,N], \; s_\tau(k) \ne +\infty \}.
\end{align*}
\index{Atau@$\Abb_\tau$}  \index{Ttau@$\Tbb_\tau$}  \index{Etau@$E_\tau$}
By \prref{CGLcluster}, $\Tbb_\tau$ is a quantum torus with corresponding quantum affine space algebra $\Abb_\tau$. In particular, $E_\tau$ consists of normal elements of $\Abb_\tau$, so it is an Ore set in $\Abb_\tau$.

\bth{CGLalmostcluster}
Let $R$ be a symmetric CGL extension of length $N$.

{\rm (a)} $\Abb_\tau \subseteq R \subseteq \Abb_\tau[E_\tau^{-1}] \subseteq \Tbb_\tau \subseteq \Fract(R)$, for all $\tau \in \Xi_N$.

{\rm (b)} $R$ is generated as a $\KK$-algebra by $\{ y_{\tau,k} \mid \tau \in \Ga_N, \; k \in [1,N] \}$.

{\rm (c)} Each $E_\tau$ is an Ore set in $R$, and $R[E_\tau^{-1}] = \Abb_\tau[E_\tau^{-1}]$.

{\rm (d)} $R = \bigcap_{\tau \in \Ga_N} R[E_\tau^{-1}] = \bigcap_{\tau \in \Ga_N} \Abb_\tau[E_\tau^{-1}]$.

{\rm (e)} Let $\inv$ be any subset of $\{ k \in [1,N] \mid s(k) = +\infty \}$. Then
$$
R[y_k^{-1} \mid k \in \inv ] = \bigcap_{\tau \in \Ga_N} R[E_\tau^{-1}] [y_k^{-1} \mid k \in \inv ] .
$$
\eth

\begin{proof}[Proof of parts {\rm(a), (b), (c), (e)}] For part (a), the first, third, and fourth inclusions are clear. For the second, it suffices to show that $x_{\tau(k)} \in \Abb_\tau[E_\tau^{-1}]$ for all $k \in [1,N]$, which we do by induction on $k$. The case $k = 1$ is immediate from the fact that $x_{\tau(1)} = y_{\tau,1}$.

Now let $k \in [2,N]$. If $p_\tau(k) = -\infty$, then $x_{\tau(k)} = y_{\tau,k} \in \Abb_\tau[E_\tau^{-1}]$. If $p_\tau(k) = l \ne -\infty$, then $l < k$ and $y_{\tau,k} = y_{\tau,l} x_{\tau(k)} - c_{\tau,k}$ for some element $c_{\tau,k}$ in the $\KK$-subalgebra of $R$ generated by $x_{\tau(1)}, \dots, x_{\tau(k-1)}$. By induction, $c_{\tau, k} \in \Abb_\tau[E_\tau^{-1}]$. Further, $y_{\tau, l} \in E_\tau$ because $s_\tau(l) = k \ne +\infty$, and thus
$$
x_{\tau(k)} = y_{\tau,l}^{-1} \bigl( y_{\tau, k} + c_{\tau, k} \bigr) \in \Abb_\tau[E_\tau^{-1}]
$$
in this case also.

(b) For each $j \in [1,N]$, we have $\tau_{j,j} \in \Ga_N$ with $\tau_{j,j}(1) = j$. \thref{sym-prime} implies that $y_{\tau,1}$ is a scalar multiple of $y_{[j,j]} = x_j$. Thus, in fact, $R$ is generated by $\{ y_{\tau, 1} \mid \tau \in \Ga_N \}$.

(c) By \eqref{Sloc}, for each $\tau \in \Xi_N$ the set $E_\tau$ is an Ore set in $R$. It is clear from part (a) that $R[E_\tau^{-1}] = \Abb_\tau[E_\tau^{-1}]$.

(e) Let $v \in \bigcap_{\tau \in \Ga_N} R[E_\tau^{-1}] [y_k^{-1} \mid k \in \inv ]$. For each $\tau \in \Ga_N$, we can write $v$ as a fraction with numerator from $R[E_\tau^{-1}]$ and right-hand denominator from the multiplicative set $Y$ generated by $\{ y_k \mid k \in \inv \}$. (Note that $Y$ is generated by normal elements of $R$, so it is an Ore set in $R$ and in $R[E_\tau^{-1}]$.) Hence, choosing a common denominator, we obtain $y \in Y$ such that $vy \in R[E_\tau^{-1}]$ for all $\tau \in \Ga_N$. Once part (d) is proved, we can conclude that $vy \in R$, and thus $v \in R[y_k^{-1} \mid k \in \inv ]$, as required.
\end{proof}

In proving \thref{CGLalmostcluster} (d), we need to compare fraction expressions for an element $v$ of $\bigcap_{\tau \in \Ga_N} R[E_\tau^{-1}]$ across the given localizations, which is done via the reenumerations $\tau_\bu \tau$ from Section \ref{6.1}. Some properties of these permutations are given next. We keep the notation $
\ex := \{ j \in [1,N] \mid s(j) \ne +\infty \}$.

\ble{taubutau} {\rm (a)} For any $\tau \in \Xi_N$, the permutation $(\tau_\bu \tau)^{-1}$ maps $\ex$ bijectively onto the set $\{ k \in [1,N] \mid s_\tau(k) \ne +\infty \}$.

{\rm (b)} Suppose that $\tau, \tau' \in \Xi_N$ and $\tau' = \tau (k,k+1)$ for some $k \in [1,N-1]$. If $\eta \tau(k) = \eta \tau(k+1)$, then $\tau'_\bu \tau' = \tau_\bu \tau$.

{\rm (c)} Suppose that $\tau, \tau' \in \Xi_N$ and $\tau' = \tau (k,k+1)$ for some $k \in [1,N-1]$. If $\eta \tau(k) \ne \eta \tau(k+1)$, then 
$$
\tau'_\bu \tau'(k) = \tau_\bu \tau(k+1) \qquad\quad \text{and} \qquad\quad \tau'_\bu \tau'(k+1) = \tau_\bu \tau(k),
$$
while $\tau'_\bu \tau'(l) = \tau_\bu \tau(l)$ for all $l \in [1,N] \setminus \{k, k+1\}$.

{\rm (d)} Let $j \in \ex$. Then there exist $\tau, \tau' \in \Ga_N$ such that $\tau' = \tau (k,k+1)$ for some $k \in [1,N-1]$ with $\eta \tau(k) = \eta \tau(k+1)$ and $\tau_\bu \tau(k) = j$.
\ele

\begin{proof} (a) Let $j \in [1,N]$ and $k = (\tau_\bu \tau)^{-1}(j)$, and note that $\eta(j) = \eta \tau(k)$. Then $j \in \ex$ if and only if $j$ is not the largest element of the set $\eta^{-1}( \eta(j) )$, if and only if $k$ is not the largest element of $\tau^{-1} \eta^{-1} (\eta(j)) = \eta_\tau^{-1}( \eta_\tau(k))$, if and only if $s_\tau(k) \ne +\infty$.

(b) This follows from the observation that $\tau^{-1}(L) = (\tau')^{-1}(L)$ for all level sets $L$ of $\eta$.

(c) If $L$ is a level set of $\eta$ not containing $k$ or $k+1$, then $\tau^{-1}(L) = (\tau')^{-1}(L)$  and so $\tau'_\bu \tau'(l) = \tau_\bu \tau(l)$ for all $l \in L$. For the case $L = \eta^{-1}( \eta \tau(k) )$, we have
\begin{multline*}
L = \{ \tau(i_1) = \tau'(i_1), \dots, \tau(i_{r-1}) = \tau'(i_{r-1}), \, \tau(k) = \tau'(k+1),  \\
  \tau(i_{r+1}) = \tau'(i_{r+1}), \dots, \tau(i_t) = \tau'(i_t) \}
\end{multline*}
for some $1 \le i_1 < \cdots < i_{r_1} < k < k+1 < i_{r+1} < \cdots < i_t \le N$, from which it is clear that $\tau'_\bu \tau'(k+1) = \tau_\bu \tau(k)$ and $\tau'_\bu \tau'(l) = \tau_\bu \tau(l)$ for all $l \in \{ i_1, \dots, i_{r-1}, i_{r+1}, \dots, i_t \}$. Similarly, we see that  $\tau'_\bu \tau'(k) = \tau_\bu \tau(k+1)$
and $\tau'_\bu \tau'(l) = \tau_\bu \tau(l)$ for all $l \in \tau^{-1} \eta^{-1} ( \tau(k+1) ) \setminus \{ k+1 \}$.

(d) Set $m = O_-(j) \in \Znn$ and $i = p^m(j) \in [1,N]$, so that $p(i) = -\infty$ and $s^m(i) = j$. Since $s(j) \in [i+1,N]$, we may take $\tau := \tau_{i,s(j)-1}$ and $\tau' := \tau_{i,s(j)}$ in $\Ga_N$. Then $\tau' = \tau (k,k+1)$ where $k = s(j) - i$, and $\eta \tau(k) = \eta(i) = \eta(s(j)) = \eta \tau(k+1)$. Moreover, the level set $L := \eta^{-1} (\eta \tau(k)) = \eta^{-1} (\eta(i))$ equals $\eta^{-1} (\eta(j))$ and $\tau$ maps the elements of $[1,k] \cap \tau^{-1}(L)$ to $s(i), \dots, s^m(i), i$ in that order. Hence, $\tau_\bu \tau(k)$ equals the $(m+1)$-st element of $L$ when $L$ is written in ascending order. That element is $j$.
\end{proof}

\begin{proof}[Proof of \thref{CGLalmostcluster} {\rm (d)}] The second equality follows from part (c), and one inclusion of the first equality is obvious.

Let $v \in \bigcap_{\tau \in \Ga_N} R[E_\tau^{-1}]$ be a nonzero 
element. For each $\tau \in \Ga_N$, let 
$$
\prod_{l \in [1,N]} y_{\tau, l}^{m_{\tau, l}} 
$$
be a minimal denominator of $v$ with respect to the localization $R[E_\tau^{-1}]$, where all $m_{\tau, l} \in \Znn$ and $m_{\tau, l} = 0$ when $s_\tau(l) = +\infty$. We first verify the following

{\bf Claim}. Let $\tau, \tau' \in \Ga_N$ such that $\tau' = \tau (k, k+1)$ for some $k \in [1, N-1]$.
\begin{enumerate}
\item $m_{\tau, (\tau_\bu \tau)^{-1}(j)} = m_{\tau', (\tau'_\bu \tau')^{-1}(j)}$, for all $j \in [1,N]$.
\item If $\eta \tau(k)  = \eta \tau(k+1)$, then $m_{\tau, k} = 0$.
\end{enumerate}

If $\eta \tau(k) = \eta \tau(k+1)$, then (1) and (2) follow from \leref{taubutau} (b), \thref{division}, and \thref{1} (b). If $\eta \tau(k) \ne \eta \tau(k+1)$, we obtain from \thref{1} (a) that $y_{\tau', k} = y_{\tau, k+1}$ and $y_{\tau', k+1} = y_{\tau, k}$, while $y_{\tau', l} = y_{\tau, l}$ for all $l \ne k, k+1$. As a result, $E_{\tau'} = E_\tau$ and we see that $m_{\tau', k} = m_{\tau, k+1}$ and $m_{\tau', k+1} = m_{\tau, k}$, while $m_{\tau', l} = m_{\tau, l}$ for all $l \ne k, k+1$. In this case, (1) follows from \leref{taubutau} (c).

Since all the permutations in $\Ga_N$ appear in the chain \eqref{sequence}, part (1) of the claim implies that 
\begin{equation}
\label{matchm}
m_{\sigma, (\sigma_\bu \sigma)^{-1}(j)} = m_{\tau, (\tau_\bu \tau)^{-1}(j)}, \; \; \forall \sigma, \tau \in \Ga_N, \; j \in [1,N].
\end{equation}
For any $j \in \ex$, \leref{taubutau} (d) shows that there exist $\tau, \tau' \in \Ga_N$ such that $\tau' = \tau (k,k+1)$ for some $k \in [1,N-1]$ with $\eta \tau(k) = \eta \tau(k+1)$ and $\tau_\bu \tau(k) = j$. Part (2) of the claim above then implies that $m_{\tau, (\tau_\bu \tau)^{-1}}(j) = 0$. From \eqref{matchm}, we thus get
$$
m_{\sigma, (\sigma_\bu \sigma)^{-1}(j)} = 0, \; \; \forall \sigma \in \Ga_N, \; j \in \ex.
$$
In particular, $m_{\id, j} = 0$ for all $j \in \ex$, whence $m_{\id, j} = 0$ for all $j \in [1,N]$ and therefore $v \in R$, which completes the proof of \thref{CGLalmostcluster} (d).
\end{proof}

\section{Completion of the proof of Theorem \ref{\tcluster}} 
\label{6.7}
We finally prove the last and most important part of \thref{cluster} which 
establishes an equality between the CGL extension $R$, the quantum cluster algebra $\Abb(M, \wt{B}, \varnothing)_\KK$ 
and the corresponding upper quantum cluster algebra $\UU(M, \wt{B}, \varnothing)_\KK$.

Recall the setting of Section \ref{5.1}. For all $\tau \in \Xi_N$, define the multiplicative subsets
$$
E_\tau := \kx \Big\{ M_\tau(f) \mid f \in \sum_{j \in \ex} \Zset_{\geq 0} e_j \Big\}.
$$
\index{Etau@$E_\tau$}
By \eqref{Sloc}, $E_\tau$ is an Ore subset of $R$ for all $\tau \in \Xi_N$. In view of \leref{taubutau} (a) and the definition of $M_\tau$, we see that $E_\tau$ is generated (as a multiplicative set) by
$$
\kx \sqcup \{ y_{\tau,k} \mid k \in [1,N], \; s_\tau(k) \ne +\infty \},
$$
matching the definition used in Section \ref{prequantumcluster}.

\begin{proof}[Proof of \thref{cluster} {\rm(d)--(e)}] By \thref{cluster} (c), for all $\tau \in \Xi_N$ the quantum seeds $(M_\tau, \wt{B}_\tau)$ 
of $\Fract(R)$ are mutation-equivalent to each other. 
For each $j \in [1,N]$, we have $\tau_{j,j} \in \Ga_N$ with $\tau_{j,j}(1) = j$, and so $\ol{y}_{\tau_{j,j},1} = \ol{y}_{[j,j]} = x_j$ by \prref{int-tau}. Thus,
each generator $x_j$ of $R$ is a cluster 
variable for a quantum seed associated to some $\tau \in \Ga_N$. Hence,
\begin{equation}
\label{inc1}
R \subseteq \Abb(M, \wt{B}, \varnothing)_\KK.
\end{equation}
The quantum Laurent phenomenon (\thref{q-l}) implies that 
\begin{equation}
\label{inc2}
\Abb(M, \wt{B}, \varnothing)_\KK \subseteq \UU (M, \wt{B}, \varnothing)_\KK.
\end{equation}
Since $M_\tau(e_j) \in R$ for all $\tau \in \Xi_N$, $j \in [1,N]$,
\begin{equation}
\label{inc3}
\UU (M, \wt{B}, \varnothing)_\KK \subseteq \bigcap_{\tau \in \Ga_N} R[E_\tau^{-1}].
\end{equation} 
Combining the embeddings \eqref{inc1}, \eqref{inc2}, \eqref{inc3} and \thref{CGLalmostcluster}  (d)
leads to 
$$
R \subseteq \Abb(M, \wt{B}, \varnothing)_\KK \subseteq \UU (M, \wt{B}, \varnothing)_\KK
\subseteq \bigcap_{\tau \in \Ga_N} R[E_\tau^{-1}] = R,
$$
which establishes all equalities in \thref{cluster} (d).

For part (e) we have the embeddings
\begin{multline}
\label{chain}
R[y_k^{-1} \mid k \in \inv] \subseteq \Abb(M, \wt{B}, \inv)_\KK \subseteq 
\\
\subseteq
\UU (M, \wt{B}, \inv)_\KK
\subseteq \bigcap_{\tau \in \Ga_N} R[y_k^{-1} \mid k \in \inv][E_\tau^{-1}],
\end{multline}
which follow from the quantum Laurent phenomenon and
the fact that each generator of $R$ is a cluster variable in 
one of the seeds indexed by $\Ga_N$. \thref{CGLalmostcluster} (e) and the chain of embeddings \eqref{chain} imply the validity of 
part (e) of the theorem. 
\end{proof}

The chain of embeddings \eqref{chain} and \thref{CGLalmostcluster} (e) also imply the 
following description of the upper quantum cluster algebra in \thref{cluster}
as a finite intersection of mixed quantum tori--quantum affine 
space algebras 
of the form \eqref{tor}.

\bco{up-cl} In the setting of Theorem {\rm\ref{tcluster}},
$$
\UU (M, \wt{B}, \inv)_\KK =  \bigcap_{\tau \in \De} \Tbb\Abb_\tau(\inv)
$$
for every subset $\De$ of $\Ga_N$ which is an interval with respect to the 
linear ordering \eqref{sequence} and has the property that for each $k \in \ex$ 
there exist two consecutive elements $\tau=\tau_{i, j-1} \prec \tau' = \tau_{i,j}$ 
of $\De$ such that $\eta(i) = \eta(j)$ and $k = \tau_\bu \tau(j-i)$, recall Section {\rm\ref{6.5}}.
Here $\Tbb\Abb_\tau(\inv)$ denotes the mixed quantum torus--quantum affine space algebra 
from \eqref{tor} associated to the permutation $\tau$, $D = \KK$, and the set of inverted 
frozen variables $\inv$. 

In particular, this property holds for $\De = \Ga_N$. 
\eco

\section{Some vectors $f_{[i,s(i)]}$}
\label{new8.9}
The normalization condition \eqref{pi-cond} in \thref{cluster} depends on the leading terms of the elements $u_{[i,s(i)]}$ for $i \in [1,N]$ with $s(i) \ne +\infty$, and thus on the exponent vectors $f_{[i,s(i)]}$. The entries of these vectors actually arise from the exchange matrix $\wt{B} = (b_{lk})$, as we  now show in the case when $[i,s(i)] = [1,N]$. This will be used in Chapter \ref{q-Schu}.

\bpr{f1N} Let $R$ be a symmetric CGL extension of length $N$ which satisfies the hypotheses of Theorem {\rm\ref{tcluster}}. Assume that $N \ge 3$ and $s(1) = N$, and recall from \eqref{Pset} that $P_{[1,N]} = \{ l \in [2,N-1] \mid s(l) = +\infty \}$ in this case. Then
\begin{equation}
\label{f1Neqn}
f_{[1,N]} = - \sum_{l \in P_{[1,N]}} b_{l1} e_{[p^{O_-(l)}(l), l]}.
\end{equation}
\epr

\begin{proof} We will apply \prref{cluster-tau-ind} using the permutations
$$
\tau := \tau_{1,N-1} = [2,\dots,N-1,1,N] \qquad \text{and} \qquad \tau' := \tau_{1,N} = [2,\dots,N,1]
$$
from the subset $\Ga_N$ of $\Xi_N$.
In this case, $\tau' = \tau (k,k+1)$ for $k = N-1$ and $\tau_\bu \tau = \tau$, so that $k_\bu = 1$. We first show that
\begin{equation}
\label{MtauM}
M_\tau = M_\id = M  \qquad \text{and} \qquad \wt{B}_\tau = \wt{B}_\id = \wt{B}.
\end{equation}

We have $\tau_{1,1} = \id \prec \tau_{1,2} \prec \cdots \prec \tau_{1,N-2} \prec \tau$ in the linear ordering \eqref{sequence}. Note that $\tau_{1,j+1} = \tau_{1,j} (j, j+1)$ for all $j \in [1,N-2]$, and that $\eta(\tau_{1,j}(j)) = \eta(1) \ne \eta(j+1) = \eta(\tau_{1,j}(j+1))$ for such $j$. \thref{cluster} (c) implies that $M_{\tau_{1,j+1}} = M_{\tau_{1,j}}$, and hence also $\rbf_{\tau_{1,j+1}} = \rbf_{\tau_{1,j}}$, for all $j \in [1,N-2]$. It follows that the conditions determining the columns of $\wt{B}_{\tau_{1,j}}$ and $\wt{B}_{\tau_{1,j+1}}$ in \thref{cluster} (a) coincide, and thus $\wt{B}_{\tau_{1,j+1}} = \wt{B}_{\tau_{1,j}}$ by uniqueness. This verifies \eqref{MtauM}.

Recall the toric frame $M_{[1,N]}$ from Section \ref{4b.3}, and apply \prref{int-tau} to see that 
$$
M_{[1,N]}(e_l) = \ol{y}_{\tau, \tau^{-1}(l)} = \wh{M}_\tau \tau^{-1}(e_l) = M_\tau(e_l), \; \; \forall l \in [1,N].
$$
Thus, $M_{[1,N]} = M_\tau = M$.
From Section \ref{4b.3} we also recall the vector $g_{[1,N]} = \sum_{l \in P_{[1,N]}} m_l e_l$ and the relation
$$
f_{[1,N]} = \sum_{l \in P_{[1,N]}} m_l e_{[p^{O^2_-(l)}(l), l]} = \sum_{l \in P_{[1,N]}} m_l e_{[p^{O_-(l)}(l), l]},
$$
with the last equality following from the fact that $\eta(l) \ne \eta(1)$ for all $l \in P_{[1,N]}$. \coref{new6.7} (a) implies that
\begin{equation}
\label{u1N}
u_{[1,N]} \; \text{ is a scalar multiple of} \; M_{[1,N]}(g_{[1,N]}) = M(g_{[1,N]}).
\end{equation}

Set $g' := \sum_{l \in P_{[1,N]}} m_l e_{l-1}$, and observe from the proof of \prref{cluster-tau-ind} that the vector $g$ appearing in the statement of the proposition is $\tau_\bu \tau(g') = g_{[1,N]}$. Thus, equations \eqref{Omtautau'1}--\eqref{chiMub} imply, recalling \eqref{MtauM}, that
\begin{align*}
&\Om_{\rbf}( e_N -g_{[1,N]} , e_j) = 1, \; \; \forall j \ne 1,
\\
&\Om_{\rbf}( e_N -g_{[1,N]} , e_1)^2
= \la^\sy_1,
\end{align*}
and
$$
\chi_{M(e_N -g_{[1,N]})} = 1. 
$$
\thref{cluster} (a) now implies that $e_N - g_{[1,N]} = b^1$. From this we obtain $b_{N1} = 1$ and 
$$b_{l1} = \begin{cases} -m_l , &l \in P_{[1,N]} \\
0  & l \in [1,N-1] \setminus P_{[1,N]}. \end{cases}
$$
Consequently,
$$
M(g_{[1,N]}) = \prod_{l \in P_{[1,N]}} \ol{y}_l^{\, -b_{l1}}.
$$
Eq. \eqref{f1Neqn} follows from this and \eqref{u1N}.
\end{proof}

\chapter{Quantum groups and quantum Schubert cell algebras}
\label{q-gr}
In this chapter we first set up notation and review material on 
quantum groups and quantum Schubert cell algebras.
Then we derive an explicit description of the sequences of prime elements 
of the latter algebras constructed in the previous chapters.

\section{Quantized universal enveloping algebras}
\label{9.1}
Fix a finite dimensional 
complex simple Lie algebra $\g$ of rank $r$ with
Weyl group $W$ and set of simple roots $\Pi= \{ \al_1, \ldots, \al_r \}$.
\index{W@$W$}  \index{Pi@$\Pi$}
Let $\lcor.,. \rcor$  \index{zzz2@$\lcor.,. \rcor$}  be the invariant bilinear form on $\Rset \Pi$
normalized by $\lcor \al_i, \al_i \rcor = 2$ for short roots 
$\al_i$. For $\ga \in \Rset \Pi$, set
$$
\| \ga \|^2 = \lcor \ga, \ga \rcor.
$$
\index{zzzz@$\Vert \ga \Vert$}
Denote by $\{s_i\}$, $\{\al_i\spcheck \}$ and $\{\vpi_i\}$ 
\index{si@$s_i$}  \index{alphaicheck@$\al_i\spcheck$}  \index{pii@$\vpi_i$}
the corresponding sets of simple reflections, coroots and fundamental weights. 
Let $\QQ$ and $\PP$ be  \index{Q@$\QQ$}  \index{P@$\PP$}
the root and weight lattices of $\g$, 
and $\PP^+ = \sum_i \Zset_{\geq 0} \vpi_i$  \index{Pplus@$\PP^+$}
the set of dominant integral weights of $\g$. 
Denote the Cartan matrix of $\g$ by
\begin{equation}
\label{Cartan}
(c_{ij}) := ( \lcor \al_i\spcheck, \al_j \rcor ) \in M_r(\Zset).
\end{equation}
\index{cij@$(c_{ij})$}

As in the previous chapters, we will work over a base field
$\KK$ of arbitrary characteristic. Choose a non-root of unity $q \in \KK^*$ and
denote by $\UU_q(\g)$ the quantized universal enveloping algebra of $\g$ 
over $\KK$ with deformation parameter $q$. \index{Uqg@$\UU_q(\g)$} We will mostly follow
the notation of Jantzen \cite{Ja}, except for denoting 
the standard generators of $\UU_q(\g)$ by $K_i^{\pm 1}$, $E_i$, $F_i$ instead of 
\index{Ki@$K_i$}  \index{Ei@$E_i$}  \index{Fi@$F_i$}  
$K_{\pm \al_i}$, $E_{\al_i}$, $F_{\al_i}$ which better fits with the combinatorial notation 
from Chapters \ref{review}--\ref{main}. We will use the form 
of the Hopf algebra $\UU_q(\g)$ with relations 
listed in \cite[\S 4.3]{Ja}, and
comultiplication, counit, and antipode given by 
\begin{align*}
\De(K_i)   &= K_i \otimes K_i,  &\ep(K_i) &=1,  &S(K_i) &= K^{-1}_i,  \\
\De(E_i)   &= E_i \otimes 1 + K_i \otimes E_i,  &\ep(E_i) &= 0,  &S(E_i) &= - K^{-1}_i E_i,  \\
\De(F_i)   &= F_i \otimes K_i^{-1} + 1 \otimes F_i,  &\ep(F_i) &=0,  &S(F_i) &= - F_i K_i,
\end{align*}
for all $i \in [1,r]$ \cite[Proposition 4.11]{Ja}.

The algebra $\UU_q(\g)$ is $\QQ$-graded 
with $\deg K_i = 0$, $\deg E_i = \al_i$, and $\deg F_i = - \al_i$
for all $i \in [1,r]$ \cite[\S 4.7]{Ja}. The corresponding graded components will 
be denoted by $\UU_q(\g)_\ga$, $\ga \in \QQ$.  \index{Uqggamma@$\UU_q(\g)_\ga$}
Define the torus  \index{H@$\HH$}
\begin{equation}
\label{Lie_tor}
\HH:= ( \KK^*)^r.
\end{equation}
Its rational character lattice is isomorphic to $\QQ$, 
where the simple root $\al_j$ is mapped to the character
\begin{equation}
\label{Hchar}
h \mt 
h^{\al_j} := t_j, \; \; \forall h =(t_1, \ldots, t_r) \in \HH.
\end{equation}
The $\QQ$-grading of $\UU_q(\g)$ gives rise to the
rational $\HH$-action on $\UU_q(\g)$ 
by algebra automorphisms such that
\begin{equation}
\label{torus-act}
h \cdot u = h^\ga u, \; \; \forall \ga \in \QQ, \; u \in \UU_q(\g)_\ga.
\end{equation}

\section{Quantum Schubert cell algebras}
\label{9.2}
Let $\B_\g$  \index{Bg@$\B_\g$}  denote the braid group of $\g$ and $\{T_i\}_{i=1}^r$
its standard generating set.  \index{Ti@$T_i$}  We will use Lusztig's action 
of $\B_\g$ on $\UU_q(\g)$ by algebra automorphisms 
in the version given in \cite[\S 8.14, Eqs. 
8.14 (2), (3), (7), (8)]{Ja}. The canonical section $W \to \B_\g$ will be denoted 
by $w \mt T_w$.  \index{Tw@$T_w$}  It follows from \cite[\S 4.7, Eq. (1) and \S 8.14, Eq. (2)]{Ja} that
\begin{equation}
\label{TiUga}
T_i \bigl( \UU_q(\g)_\ga \bigr) = \UU_q(\g)_{s_i\ga}, \; \; \forall i \in [1,r], \; \ga \in \QQ.
\end{equation}

The quantum Schubert cell algebras $\UU^\pm[w]$, $w \in W$  
\index{quantum Schubert cell algebra}
were defined by De Concini, Kac, and Procesi 
\cite{DKP}, and Lusztig \cite[\S 40.2]{L} as follows. 
Fix a reduced expression
\begin{equation}
\label{reduced}
w = s_{i_1} \ldots s_{i_N}
\end{equation}
and denote
\begin{align*}
w_{\leq k} &:= s_{i_1} \ldots s_{i_k}, \qquad w_{\leq 0 }:=1, 
\\
w_{[j,k]} &:= 
\begin{cases}
s_{i_j} \ldots s_{i_k} ,& \mbox{if} \; \; j \leq k
\\
1, &\mbox{if} \; \; j>k
\end{cases}
\end{align*}
\index{wlek@$w_{\leq k}$}  \index{wjk@$w_{[j,k]}$}
for all $j,k \in [1,N]$.
(The above notation depends on the choice of reduced expression, 
but this dependence will not be displayed explicitly for simplicity of the notation.) 
Define the roots
\begin{equation}
\label{beta}
\beta_k := w_{\leq (k-1)} \al_{i_k}, \; \; \forall k \in [1, N]
\end{equation}
of $\g$ and the Lusztig root vectors  \index{Lusztig root vector}
\begin{multline}
E_{\be_k} := T_{w_{\leq k-1}}(E_{i_k})=
T_{i_1} \ldots T_{i_{k-1}}
(E_{i_k}), \\   
F_{\be_k} := T_{w_{\leq k-1}} (F_{i_k}) =
T_{i_1} \ldots T_{i_{k-1}} 
(F_{i_k}) \in \UU_q(\g), \; \; 
\forall k \in [1, N],
\label{rootv}
\end{multline}
\index{Ebetak@$E_{\be_k}$}  \index{Fbetak@$F_{\be_k}$}  
see \cite[\S 39.3]{L}. Note from \eqref{TiUga} that $E_{\be_k} \in \UU_q(\g)_{\be_k}$ and $F_{\be_k} \in \UU_q(\g)_{-\be_k}$.

By \cite[Proposition 2.2]{DKP} and 
\cite[Proposition 40.2.1]{L}, the subalgebras 
$\UU^\pm[w]$  \index{Uplisminusw@$\UU^\pm[w]$}
of $\UU_q(\g)$ generated by 
$E_{\be_k}$, $k \in [1,N]$ and $F_{\be_k}$, $k \in [1,N]$, respectively,
do not depend on the choice of a 
reduced expression for $w$ and have the $\KK$-bases
\begin{align}
\nn
&\{(E_{\be_N})^{m_N} \ldots (E_{\be_1})^{m_1} \mid 
m_N, \ldots, m_1 \in \Zset_{\geq 0}
\} 
\; \; \mbox{and}
\\
&\{(F_{\be_N})^{m_N} \ldots (F_{\be_1})^{m_1} \mid 
m_1, \ldots, m_N \in \Zset_{\geq 0} \}.
\label{PBW}
\end{align}
The algebras $\UU^\pm[w]$ are $\QQ$-graded subalgebras 
of $\UU_q(\g)$ and are thus stable under the action \eqref{torus-act}.

There is a unique algebra automorphism $\omega$  \index{omega@$\omega$} 
of $\UU_q(\g)$ such that 
\[
\om(E_i) = F_i, \quad
\om(F_i) = E_i, \quad
\om(K_i) = K_i^{-1}, \quad
\forall i \in [1,r].
\]
By \cite[eq. 8.14(9)]{Ja}, $\om( T_i(u)) = (-1)^{\lcor \al_i\spcheck, \ga \rcor } 
q^{ \lcor \al_i, \ga \rcor }  T_i( \om (u))$ for all
$i \in [1,r]$, $\ga \in \QQ$, and $u \in \UU_q(\g)_\ga$. 
This implies that the restrictions of $\om$ induce 
isomorphisms
$$
\om : \UU^+[w]  \stackrel{\cong}{\longrightarrow} \UU^-[w], 
\quad \om(E_{\be_k}) = 
(-1)^{ \lcor \be_k - \al_{i_k}, \rho\spcheck \rcor } 
q^{ - \lcor \be_k - \al_{i_k}, \rho \rcor } 
F_{\be_k}, 
\; \; \forall k \in [1, N],
$$
where $\rho$ and $\rho\spcheck$ are the half sums of the positive 
roots and coroots of $\g$, respectively.

We will restrict ourselves to $\UU^-[w]$ since these algebras are naturally realized 
in terms of highest weight vectors for $\UU_q(\g)$-modules (see Section \ref{9.3} for details),
while the $\UU^+[w]$ are realized in terms of lowest weight vectors \cite[Theorem 2.6]{Y-sqg}. 
The above isomorphisms can be used to translate all results to $\UU^+[w]$.

The Levendorskii--Soibelman straightening law  \index{Levendorskii--Soibelman straightening law}  is the following commutation relation 
in $\UU^-[w]$:
\begin{multline}
\label{LS}
F_{\be_k} F_{\be_j} - 
q^{ - \lcor \be_k, \be_j \rcor }
F_{\be_j} F_{\be_k}  \\
= \sum_{ {\bf{m}} = (m_{j+1}, \ldots, m_{k-1}) \, \in \, \Znn^{k-j-1} }
\xi_{\bf{m}} (F_{\be_{k-1}})^{m_{k-1}} \ldots (F_{\be_{j+1}})^{m_{j+1}},
\; \; \xi_{\bf{m}} \in \KK,
\end{multline}
for all $1 \le j < k \le N$ (see e.g. \cite[Proposition I.6.10]{BG} and apply $\omega$). 
This identity and \eqref{PBW} easily imply, as recorded in the following lemma, that $\UU^-[w]$ is a symmetric CGL extension 
for the obvious order of its canonical generators: 
\begin{equation}
\label{cano-gen}
x_1:=F_{\be_1}, \ldots, x_N:=F_{\be_N}.
\end{equation}
Note that the unital subalgebra of $\UU^-[w]$ 
generated by $F_{\be_1}, \ldots, F_{\be_k}$ is equal to $\UU^-[w_{\leq k}]$.
It follows from \eqref{Hchar} 
that for all $k \in [1, N]$ there exist $h_k, h^\sy_k \in \HH$ 
such that 
\begin{equation}
\label{hj}
h_k^{\be_j} = 
q^{\lcor \be_k, \be_j \rcor }, \; \; 
\forall j \in [1, k]
\quad \mbox{and} \quad
(h^\sy_k)^{\be_l} = 
q^{- \lcor \be_k, \be_l \rcor }, \; \; 
\forall l \in [k, N].
\end{equation}
For $i \in [1,r]$, 
set $q_i := q^{ \| \al_i \|^2/ 2 } \in \kx$. Recall that $\|\al_i \|^2/2 \in \{1,2,3\}$
for all $i \in [1,r]$. 

\ble{Uw} For all Weyl group elements $w \in W$ of length $N$ and reduced expressions
\eqref{reduced} of $w$ we have:

{\rm(a)} For all $k \in [1,N]$, the algebra $\UU^-[w_{\leq k}]$ is an Ore extension of the form $ \UU^-[w_{\leq (k-1)}][F_{\be_k}; \sigma_k, \delta_k]$,
where $\sigma_k = (h_k \cdot) \in \Aut(\UU^-[w_{\leq (k-1)}])$ and
$\delta_k$ is the locally nilpotent $\sigma_k$-derivation
of $\UU^-[w_{\leq (k-1)}]$ given by
$$
\delta_k(u) := F_{\be_k} u - q^{\lcor \be_k, \ga \rcor} u F_{\be_k}, \; \; 
\forall
u \in (\UU^-[w_{\leq (k-1)}])_\ga, \; \ga \in \QQ.
$$
The $h_k$-eigenvalue of $F_{\be_k}$ is $q_{i_k}^{-2}$ and is not a root of unity. 

{\rm(b)} The algebra 
\begin{equation}
\label{UwCGL}
\UU^-[w] = \KK [F_{\be_1}] [F_{\be_2}; \sigma_2, \delta_2] \cdots [F_{\be_N}; \sigma_N, \delta_N]
\end{equation}
is a symmetric CGL extension 
for the choice of generators \eqref{cano-gen} and 
the choice of elements $h_k, h^\sy_k \in \HH$.

{\rm(c)} The corresponding ``interval'' subalgebras are given by 
$$
\UU^-[w]_{[j,k]} = T_{w_{\leq j-1}} \bigl( \UU^-[w_{[j,k]}] \bigr).
$$ 
\ele

For parts (a)-(b) of the lemma see e.g. \cite[Lemma 2.1]{GeY}.
Part (c) of the lemma follows at once from the definition of the 
root vectors $F_{\be_k}$. The multiplicatively skew-symmetric matrix 
$\lab \in M_N(\kx)$ associated to the above CGL extension presentation of $\UU^-[w]$ 
is given by 
\begin{equation}
\label{Uwlakj}
\la_{kj} = q^{-\lcor \be_k , \be_j \rcor }, \; \; \forall 1 \leq j < k \leq N.
\end{equation}  
Furthermore,
\begin{equation}
\label{lalastar}
\la_k = q_{i_k}^{-2} \; \; \mbox{and} \; \; 
\la^\sy_k = q_{i_k}^2, \; \; \forall k \in [1,N],
\end{equation}
because $\| \be_k \|^2 = \| \al_{i_k} \|^2$.

\section{Quantum function algebras and homomorphisms}
\label{9.3}
Recall that a $\UU_q(\g)$-mod\-ule $V$ is a type one module  \index{type one module}  if it is a direct sum 
of its weight spaces  \index{weight space}  defined by 
\[
V_\mu := \{ v \in V \mid K_i v = q^{ \lcor \mu, \al_i \rcor} v, \; \; 
\forall i \in [1,r] \}, \; \mu \in \PP.
\]
\index{Vmu@$V_\mu$}
The category of  finite dimensional type one (left) 
$\UU_q(\g)$-modules is semisimple 
(see  \cite[Theorem 5.17]{Ja} and the remark on p. 85 of 
\cite{Ja}), and is closed under taking tensor products and duals.
The irreducible modules in this category are classified by 
the dominant integral weights of $\g$
(\cite[Theorem 5.10]{Ja}). For each $\mu \in \PP^+$, denote 
by $V(\mu)$ the corresponding irreducible module. Let $v_\mu$ 
be a highest weight vector of $V(\mu)$.  \index{V(mu)@$V(\mu)$}  \index{vmu@$v_\mu$}

Let $G$ denote the connected, simply connected, complex simple algebraic group with Lie algebra $\g$.
The quantum function algebra  \index{quantum function algebra}  $R_q[G]$  \index{RqG@$R_q[G]$}  is the Hopf subalgebra 
of the restricted dual $\UU_q(\g)^\circ$ 
spanned by the matrix coefficients  \index{matrix coefficient}  $c^\mu_{\xi, v}$  \index{cmuxiv@$c^\mu_{\xi, v}$}  of the modules $V(\mu)$, 
$\mu \in \PP^+$, defined by
\begin{equation} 
\label{c-notation}
c^\mu_{\xi, v}(u) := \xi ( u v ), \; \; \forall v \in V(\mu), \;
\xi \in V(\mu)^*, \; u \in \UU_q(\g).
\end{equation}
Set for brevity
\[
c_{\xi}^\mu := c^\mu_{\xi, v_\mu}, \; \; \forall \mu \in \PP^+, \;
\xi \in V(\mu)^*.
\]
The space
\[
R^+ := \Span \{ c^\mu_\xi \mid \mu \in \PP^+, \; \xi \in V(\mu)^* \}
\]
\index{Rplus@$R^+$}
is a subalgebra of $R_q[G]$ \cite[\S 9.1.6]{Jbook}. The braid group $\B_\g$ acts on 
the finite dimensional type one $\UU_q(\g)$-modules $V$ (see  
\cite[\S 8.6]{Ja}) in a compatible way with its action on
$\UU_q(\g)$:
$$
T_w ( u . v ) := (T_w u) . (T_w v), \; \; \forall 
w \in W, \; u \in \UU_q(\g), \; v \in V 
$$
(cf. \cite[eq. 8.14 (1)]{Ja}). Given $\mu \in \PP^+$ and $w \in W$,
there exists a unique vector $\xi_{w, \mu} \in (V(\mu)^*)_{- w\mu}$ 
such that
\begin{equation}
\label{xi-w}
\lcor \xi_{w, \mu}, T^{-1}_{w^{-1}} v_\mu \rcor =1. 
\end{equation}
For $y, w \in W$ and $\mu \in \PP^+$ define the (generalized) quantum minors  \index{generalized quantum minor}
\begin{equation}
\label{e}
\De_{y \mu, w \mu} := c^\mu_{\xi_{y,\mu}, T_{w^{-1}}^{-1} v_\mu} \in R_q[G] \; \; 
\mbox{and} \; \; 
\De_{w \mu} := \De_{w \mu, \mu} = c^\mu_{\xi_{w,\mu}} \in R^+.
\end{equation}
\index{Deltaymuwmu@$\De_{y \mu, w \mu}$}
For all $\mu \in \PP^+$ and $y,w \in W$, the quantum minor $\De_{y \mu, w \mu}$ above
equals the quantum minor $\De_{y \mu, w \mu}$ of Berenstein--Zelevinsky 
from \cite[Eq. (9.10)]{BZ}. One shows this by a repetitive application of 
the formulas 
$$
T_i^{-1} v_{ m \vpi_i} = \frac{1}{[m]_{q_i} !} F_i^m v_{m \vpi_i}
\quad
\mbox{and} 
\quad 
E^m_i F^m_i v_{m \vpi_i} 
= ([m]_{q_i}!)^2 v_{m \vpi_i}
$$
for $m \in \Zset_{\geq 0}$, $i \in [1,r]$. 
(Here and below we use the 
notation of \cite{BG,Ja} for $q$-integers and factorials.) 
This is why everywhere in this and the next chapter
we use $T_{w^{-1}}^{-1}$ instead of $T_w$. Using the latter would bring 
some unwanted extra scalars coming from the formula
$$
T_i v_{ m \vpi_i} = \frac{(-q_i)^m}{[m]_{q_i} !} F_i^m v_{m \vpi_i}.
$$

For $w \in W$, the set $E_w := \{ \De_{w \mu} \mid \mu\in \PP^+ \} \subset R^+$  \index{Ew@$E_w$}
is a multiplicative subset of $R^+$ because 
$\De_{w \mu_1}  \De_{w \mu_2} = \De_{w \mu_1 + w \mu_2 }$ for all
$\mu_1, \mu_2 \in \PP^+$, $w \in W$.
Joseph proved that it is an Ore set \cite[Lemma 9.1.10]{Jbook}
for $\charr \KK =0$ and $q$ transcendental. His proof works for all base fields 
$\KK$, $q \in \KK^*$ not a root of unity, \cite[\S 2.2]{Y-multi}.
Following \cite[\S 10.4.8]{Jbook}, define the algebras
$$
R^w := R^+ [E_w^{-1}] \; \; 
\mbox{and} \; \; 
R^w_0 := \{ c^\mu_\xi \De_{w \mu}^{-1} \mid \mu \in \PP^+, \; \xi \in V(\mu)^* \}.
$$
\index{Rw@$R^w$}  \index{Rwzero@$R^w_0$}
One does not need to take a span in the definition of the second algebra, but note that the elements of $R^w_0$ do not have 
unique presentations in the form $c^\mu_\xi \De_{w \mu}^{-1}$.
The algebra $R^w_0$ is $\QQ$-graded by 
\[
(R_0^w)_\ga := \{  c^\mu_\xi \De_{w \mu}^{-1} \mid \mu \in \PP^+, \;
\xi \in (V(\mu)^*)_{\ga - w \mu } \}, \; \; 
\forall\, \ga \in \QQ.
\]
In particular,
$$
\De_{y \mu} \De_{w \mu}^{-1} =
c^\mu_{\xi_{y,\mu}} \De_{w \mu}^{-1} \in (R^w_0)_{(w-y)\mu} , \; \; \forall \mu \in \PP^+, \; y,w \in W.
$$

The algebra $\UU^-[w]$ is realized as a quotient of $R_0^w$ as follows. 
Denote $\QQ^+ := \sum_i \Zset_{\geq 0} \al_i$.  \index{Qaplus@$\QQ^+$}
For $\ga \in \QQ^+$, set
$m_w(\ga) := \dim (\UU^+[w])_\ga= \dim (\UU^-[w])_{-\ga}$.  \index{mwgamma@$m_w(\ga)$}  Let
$\{u_{\ga, n} \}_{n=1}^{m_w(\ga)}$ and 
$\{u_{-\ga, n} \}_{n=1}^{m_w(\ga)}$ be 
dual bases of $(\UU^+[w])_\ga$ and 
$(\UU^-[w])_{-\ga}$ with respect to the Rosso--Tanisaki form,
see \cite[Ch. 6]{Ja}. The quantum $R$-matrix corresponding 
to $w$ is the element
$$
\RR^w := 1 \otimes 1 + \sum_{\ga \in \QQ^+, \, \ga \neq 0} \sum_{n=1}^{m_w(\ga)} 
u_{\ga, n} \otimes u_{- \ga, n}
$$
\index{Rwscript@$\RR_w$}
of the completion $\UU^+ \wh{\otimes} \, \UU^-$ of $\UU^+ \otimes \UU^-$ 
with respect to the descending filtration \cite[\S 4.1.1]{L}.

There is a unique graded algebra antiautomorphism 
$\tau$  \index{tau@$\tau$}  of $\UU_q(\g)$ given by 
\begin{equation}
\label{iota}
\tau(E_i) = E_i,
\; \;
\tau(F_i) = F_i, 
\; \;
\tau(K_i) = K_i^{-1}, \; \; 
\forall i \in [1,r],
\end{equation}
see \cite[Lemma 4.6(b)]{Ja} for details.
It is compatible with the braid group action:
\begin{equation}
\label{iota-ident}
\tau (T_w u) = T_{w^{-1}}^{-1} ( \tau (u)), \; \; 
\forall u \in \UU_q(\g), \; w \in W,
\end{equation}
see \cite[Eq. 8.18(6)]{Ja}. 

\bth{antihom} \cite[Theorem 2.6]{Y-sqg}
For all Weyl group elements $w \in W$,
the map $\vp_w :  R^w_0 \to \UU^-[w]$  \index{phiw@$\vp_w$}  given by
$$
\vp_w \big( c^\mu_\xi \De_{w \mu}^{-1} \big) 
:= \big( c^\mu_{\xi, T^{-1}_{w^{-1}} v_\mu} \otimes \id \big) (\tau \otimes \id) 
\RR^w, \; \; 
\forall \mu \in \PP^+, \; \xi \in V(\mu)^*
$$
is a well defined surjective $\QQ$-graded algebra antihomomorphism.
\eth

The kernel of $\vp_w$ has an explicit 
description in terms of Demazure modules, \cite[Theorem 2.6]{Y-sqg}.
Later we will need the following explicit formula for 
$\vp_w$ which follows at once from the standard formula
for the inner product of pairs of monomials 
with respect to the Rosso--Tanisaki form 
\cite[Eqs. 8.30 (1) and (2)]{Ja}:
\begin{multline}
\label{phi1} 
\vp_w(c_\xi^\mu \De_{w \mu}^{-1} ) = \sum_{m_1, \ldots, m_N \, \in \, \Zset_{\geq 0}}
\left( \prod_{j=1}^N
\frac{ (q_{i_j}^{-1} - q_{i_j})^{m_j}}
{q_{i_j}^{m_j (m_j-1)/2} [m_j]_{q_{i_j}}! } \right) 
\\
\times 
\lcor \xi, (\tau E_{\be_1})^{m_1} \ldots (\tau E_{\be_N})^{m_N} T^{-1}_{w^{-1}} v_\mu \rcor 
F_{\be_N}^{m_N} \ldots F_{\be_1}^{m_1},
\end{multline}
for all $\mu \in \PP^+$, $\xi \in V(\mu)^*$. 

\section{Quantum minors and sequences of prime elements}
\label{9.4}
For $y, w \in W$ and $i \in [1,r]$, 
set 
\begin{equation}
\label{defDet}
\Det_{y\vpi_i, w\vpi_i} := \vp_w( \De_{y \vpi_i} \De_{w \vpi_i}^{-1} )
= (\De_{y \vpi_i, w \vpi_i} \tau \otimes \id) \RR^w \in \UU^-[w]_{(w-y) \vpi_i}.
\end{equation}
\index{Deltaytildevarpi@$\Det_{y\vpi_i, w\vpi_i}$}
The elements $\De_\mu \De_{w \mu}^{-1} \in R_0^w$ are normal 
modulo $\ker \vp_w$ for all $\mu \in \PP^+$, 
which implies that $\Det_{\vpi_i, w \vpi_i}$  
are (nonzero) normal elements of $\UU^-[w]$ for all $i \in [1,r]$.
More precisely,
\begin{equation}
\label{commute}
\Det_{\vpi_i,w \vpi_i} u
= 
q^{ - \lcor (1+w)\vpi_i, \ga \rcor }
u \Det_{\vpi_i,w \vpi_i}, \; \; 
\forall 
i \in [1,r], \;
u \in \UU^-[w]_\ga, \; \ga \in \QQ,
\end{equation}
see \cite[Eq. (3.30)]{Y-sqg}. 

Denote the support  \index{support}  of $w \in W$ by
$$
\Sbb(w) := \{ i \in [1,r] \mid s_i \leq w \} = \{ i_1,\dots,i_N \},
$$
\index{Sw@$\Sbb(w)$}
where $\leq$ refers to the Bruhat order on $W$. The homogeneous 
prime elements of the algebras $\UU^-[w]$ are given by the following theorem.

\bth{U-prim} \cite[Theorem 6.2]{Y-sqg} For all  $w \in W$,
$$
\{ \Det_{\vpi_i, w \vpi_i} \mid i\in \Sbb(w) \}
$$ 
is a list of the homogeneous prime elements of $\UU^-[w]$ up to scalar multiples.
\eth 

\bco{rkUw} In the above setting,
$$
\rk(\UU^-[w]) = | \Sbb(w)|.
$$
\eco

\leref{Uw} (c) and \thref{U-prim} imply that all prime elements $y_{[j,s^m(j)]}$ 
are scalar multiples of elements of the form $T_{w'} (\Det_{\vpi_i, w'' \vpi_i})$ where 
$w'$ and $w''$ are subwords of \eqref{reduced}, $i \in \Sbb(w)$.
In the remaining part of this chapter we obtain a more explicit form of this fact.

From now on we fix a reduced expression \eqref{reduced} of $w \in W$.
Consider the function
\begin{equation}
\label{Uweta}
\eta : [1,N] \to [1,r], \; \; 
\eta(k) = i_k.
\end{equation}
The associated functions $p : [1,N] \to [1,N] \sqcup \{ - \infty \}$ and 
$s : [1,N] \to [1,N] \sqcup \{ + \infty \}$ are the functions $k \mt k^-$ and 
$k \mt k^+$ from \cite{BFZ}:
\begin{equation}
\label{Uw-p}
p(k) = 
\begin{cases} 
\max \{ j < k \mid i_j = i_k \}, & \mbox{if such $j$ exists} \\
-\infty, & \mbox{otherwise}  
\end{cases}
\end{equation}
and 
\begin{equation}
\label{Uw-s}
s(k) = 
\begin{cases} 
\min \{ j > k \mid i_j = i_k \}, & \mbox{if such $j$ exists} \\
+\infty, & \mbox{otherwise}.  
\end{cases}
\end{equation}
(One should note that \cite{FZ} uses integers instead of $\pm \infty$.)
Recall from \eqref{tau-ci} that $w_\circ$ denotes the longest element of $S_N$. The CGL extension presentation of $\UU^-[w]$ corresponding to this permutation has the form
\begin{equation}
\label{UwCGLw0}
\UU^-[w] = \KK [F_{\be_N}] [F_{\be_{N-1}}; \sigma^*_{N-1}, \delta^*_{N-1}] \cdots [F_{\be_1}; \sigma^*_1, \delta^*_1].
\end{equation}
We write $y_{w_\ci,1}, \dots, y_{w_\ci,N}$  \index{ywcirclek@$y_{w_\ci,k}$}  for the sequence of $y$-elements given by \thref{CGL} for the above presentation.

\bth{Uw-prim} For all Weyl group elements $w \in W$ of length $N$, reduced expressions
\eqref{reduced} of $w$, and $k \in [1,N]$, we have
\begin{equation}
\label{y-De}
y_{w_\ci,k} =  
\bigl( q_{i_{N-k+1}}^{-1} - q_{i_{N-k+1}} \bigr)^{ -O_+(N-k+1) -1} \,
\Det_{w_{\leq N-k} \vpi_{i_{N-k+1}}, w \vpi_{i_{N-k+1}}}.
\end{equation}
The function $\eta$ from \eqref{Uweta} satisfies the conditions of Theorem {\rm\ref{tCGL}}
for the CGL extension presentation \eqref{UwCGL} of $\UU^-[w]$, and 
the corresponding functions $p$ and $s$ are given by the formulas \eqref{Uw-p} and \eqref{Uw-s}. 

Furthermore, for all $1\leq j < k \leq N$,
\begin{multline}
\Det_{w_{\leq j-1} \vpi_{i_j}, w \vpi_{i_j}} \Det_{w_{\leq k-1} \vpi_{i_k}, w \vpi_{i_k}} 
=  \\
 = q^{ \lcor (w_{\leq j-1} + w) \vpi_{i_j}, ( w_{\leq_{k-1}} - w) \vpi_{i_k}  \rcor } 
\Det_{w_{\leq k-1} \vpi_{i_k}, w \vpi_{i_k}} \Det_{w_{\leq j-1} \vpi_{i_j}, w \vpi_{i_j}}. 
\label{Decomm}
\end{multline}
\eth

Using the definition of the braid group action, one easily sees that
$$
\Det_{\vpi_i, w \vpi_i} = \Det_{w_{\leq j - 1 } \vpi_{i_j}, w \vpi_{i_j}}, \; \; 
\mbox{for} \; j = \min \{ k \in [1,N] \mid i_k = i \}.
$$

\begin{proof} First we prove \eqref{Decomm}. Applying \cite[Eq. (10.2)]{BZ}, we obtain 
\begin{multline*}
\De_{w_{\leq k-1} \vpi_{i_k}, \vpi_{i_k}} \De_{w_{\leq j-1} \vpi_{i_j}, \vpi_{i_j}} =  \\ 
 q^{ \lcor w_{\leq k-1} \vpi_{i_k}, w_{\leq j-1} \vpi_{i_j} \rcor - 
\lcor \vpi_{i_k}, \vpi_{i_j} \rcor } 
\De_{w_{\leq j-1} \vpi_{i_j}, \vpi_{i_j}} \De_{w_{\leq k-1} \vpi_{i_k}, \vpi_{i_k}}  
\end{multline*}
for all $1\leq j < k \leq N$.
Eq. (10.2) in \cite{BZ} was stated for $\KK = \Qset(q)$, but its 
proof works for all fields $\KK$ since it only uses
the standard left and right actions of $\UU_q(\g)$ on $R_q[G]$. In addition, 
we have
$$
\De_{w \mu, \mu} u = q^{ - \lcor w \mu, \ga \rcor} u \De_{w \mu, \mu} 
\mod R^+[E_w^{-1}] \ker \vp_w, \; \;
\forall u \in (R_0^w)_\ga, \; \ga \in \QQ, \; \mu \in \PP^+.
$$
see \cite[Eq. (2.22) and Theorem 2.6]{Y-sqg}. 
Eq. \eqref{Decomm} follows from these 
two identities using that $\vp_w : R^w_0 \to \UU^-[w]$ is a 
graded antihomomorphism by \thref{antihom}.

Proposition 3.3 from \cite{GeY} implies that 
$\Det_{w_{\leq j-1} \vpi_{i_j}, w \vpi_{i_j}} \in \UU^-[w]_{[j,N]}$ for all $j \in [1,N]$ 
and that $\UU^-[w]_{[j,N]}$ sits inside the 
subalgebra of $\Fract(\UU^-[w])$ generated by the elements
$\Det_{w_{\leq l-1} \vpi_{i_l}, w \vpi_{i_l}}^{ \pm 1}$, $l \in [j,N]$. Since 
$\Det_{w_{\leq k-1} \vpi_{i_k}, w \vpi_{i_k}}$ belongs to $(\UU^-[w])_{-(w_{\leq k-1}-w) \vpi_{i_k}}$, 
it follows from \eqref{Decomm} that
$\Det_{w_{\leq j-1} \vpi_{i_j}, w \vpi_{i_j}}$ is a normal element of $\UU^-[w]_{[j,N]}$:
$$
\Det_{w_{\leq j-1} \vpi_{i_j}, w \vpi_{i_j}} u = q^{- \lcor (w_{\leq j-1} + w) \vpi_{i_j}, \ga \rcor }
u \Det_{w_{\leq j-1} \vpi_{i_j}, w \vpi_{i_j}}
$$
for all $u \in \left( \UU^-[w]_{[j,N]} \right)_\ga$, $\ga \in \QQ$. 

Invoking again \cite[Proposition 3.3]{GeY}, we have
\begin{multline}
\label{firsts}
\Det_{w_{\leq j-1} \vpi_{i_j}, w \vpi_{i_j}} \equiv 
(q_{i_j}^{-1} - q_{i_j}) \Det_{w_{\leq s(j)-1} \vpi_{i_{s(j)}}, w \vpi_{i_{s(j)}}} F_{\be_j}  \\
\mod \UU^-[w]_{[j+1, N]},
\; \; \mbox{if} \; \; s(j) \neq + \infty
\end{multline}
and 
\begin{equation}
\label{seconds}
\Det_{w_{\leq j-1} \vpi_{i_j}, w \vpi_{i_j}} \equiv 
(q_{i_j}^{-1} - q_{i_j}) F_{\be_j}
\mod \UU^-[w]_{[j+1, N]},
\; \; \mbox{if} \; \; s(j) = + \infty.
\end{equation}
We apply \prref{constr-y} for the presentation \eqref{UwCGLw0} of $\UU^-[w]$. The previous arguments verify conditions (i)--(iii) of the proposition, where $y'_k$ is the element on the right  hand side of \eqref{y-De} for $k \in [1,N]$ and $c'_1,\dots,c'_N$ are obtained from \eqref{firsts} and \eqref{seconds}.
We are left with showing that condition (iv) is satisfied. Denote 
by $\wt{s} : [1,N] \to [1,N] \sqcup \{ + \infty \}$ the successor 
function from \thref{CGL} for the presentation \eqref{UwCGLw0} of $\UU^-[w]$. We need to 
prove that, if $\wt{s}(j) = + \infty$, 
then $s(j) = + \infty$, i.e we are in the case of Eq. \eqref{seconds}.
If $\wt{s}(j) = + \infty$, then $\rk \UU^-[w_{[j,N]}]= \rk \UU^-[w_{[j+1,N]}]+1$.
Applying \coref{rkUw}, we obtain that $|\Sbb(w_{[j,N]})| = | \Sbb(w_{[j+1,N]})| +1$, 
and thus $s(j) = + \infty$. Eq. \eqref{y-De} and the statements 
for $\eta$, $p$ and $s$ now follow from \prref{constr-y}.
\end{proof}

Combining \thref{Uw-prim} (applied to $\UU^-[w_{\le s^m(j)}]$) and \prref{int-tau} (applied to the case when $\tau$ equals the longest element of $S_{s^m(j)}$) leads to the first part of the following:

\bco{Uw-prim-int} In the setting of Theorem {\rm\ref{tUw-prim}}, assume there exists $\sqrt{q} \in \kx$. Then
for all $j \in [1,N]$ and $m \in \Zset_{\geq 0}$ 
such that $s^m(j) \in [1,N]$,
\begin{equation}
\label{yjsmjDet}
y_{[j,s^m(j)]} = \bigl( q_{i_j}^{-1} - q_{i_j} \bigr)^{-m -1}
\biggl( \prod_{0 \leq l < n \le m} q^{ \lcor \be_{s^l(j)}, \be_{s^n(j)} \rcor } \biggr)
\Det_{w_{\leq j-1} \vpi_{i_j}, w_{\le s^m(j)} \vpi_{i_j}}.
\end{equation}
and
\begin{equation}
\label{betasum}
\be_j + \be_{s(j)} + \cdots + \be_{s^m(j)} = ( w_{\le j-1} - w_{\le s^m(j)} ) \vpi_{i_j}.
\end{equation}
\eco

Eq. \eqref{betasum} follows from \eqref{yjsmjDet} by inspecting degrees with respect to the $\QQ$-grading. On one hand, $\Det_{w_{\leq j-1} \vpi_{i_j}, w_{\le s^m(j)} \vpi_{i_j}}$ has degree $( w_{s^m(j)} - w_{\le j-1} ) \vpi_{i_j}$ by \eqref{defDet}. On the other hand,
$$\deg( y_{[j, s^m(j)]} ) = \sum_{l=0}^m \deg( F_{\be_{s^l(j)}} ) = - \sum_{l=0}^m \be_{s^l(j)}.$$

\chapter[Quantum cluster structures on quantum Schubert cells]{Quantum cluster algebra structures on quantum Schubert cell algebras}
\label{q-Schu}
As an application of the results of Chapter \ref{main}, we
prove that for all finite dimensional simple Lie algebras $\g$, 
the quantum Schubert cell algebras $\UU^-[w]$ have canonical structures 
of quantum cluster algebras. Previously this was known for symmetric Kac--Moody
algebras $\g$ due to Gei\ss, Leclerc and Schr\"oer \cite{GLSh}. Our proof 
works under very general assumptions on the base field $\KK$ and
the deformation parameter $q \in \kx$. The field can have arbitrary 
characteristic, it does not need to be algebraically closed, 
and $q$ is only assumed to be a non-root of unity.  

The existence of \emph{some} quantum cluster algebra structure on $U^-[w]$ is guaranteed by \thref{cluster}, once we enlarge the base field to include a square root of $q$ and rescale the generators $F_{\be_1}, \dots, F_{\be_N}$ to satisfy condition \eqref{pi-cond}. We actually put these generators in reverse order, $F_{\be_N}, \dots, F_{\be_1}$, before arranging to apply \thref{cluster}. What then remains is to identify the initial quantum seed and the cluster variables of this structure.

\section{Statement of the main result}
\label{10.1}
Fix a finite dimensional complex simple Lie algebra $\g$, a Weyl group element $w \in W$, and a reduced expression 
\eqref{reduced} of $w$. Throughout this chapter, 
the predecessor and successor functions $p$ and $s$ will 
refer to the ones given by \eqref{Uw-p}--\eqref{Uw-s}.
Recall that the multiplicatively skew-symmetric matrix $\lab \in M_N(\kx)$ associated to 
the CGL extension presentation of $\UU^-[w]$ for the sequence of 
generators \eqref{cano-gen} is given by 
$$
\la_{jk} = q^{\lcor \be_j , \be_k \rcor }, \; \; \forall 1 \leq j < k \leq N.
$$  

From now on we will assume that the base field $\KK$ contains a square root of $q$
and fix such a root $\sqrt{q} \in \kx$. Let $\nub \in M_N(\kx)$ be the unique 
multiplicatively skew-symmetric matrix given by 
$$
\nu_{jk} = {\sqrt{q}}^{\,\lcor \be_j, \be_k \rcor}, \; \; \forall 1 \leq j < k \leq N.
$$
The results of Chapters \ref{mCGL}--\ref{main} will be applied to the 
algebra $\UU^-[w]$ for this choice of the matrix $\nub$.

\thref{Uw-prim} implies that there is a unique toric frame 
$M^w \colon \Zset^N \to \Fract(\UU^-[w])$  \index{Mw1@$M^w$}  such that 
$$
M^w(e_k) = \sqrt{q}^{\, \| (w - w_{\leq k-1}) \vpi_{i_k} \|^2/2} \, \Det_{w_{\leq k-1} \vpi_{i_k}, w \vpi_{i_k}}, 
\; \; \forall k \in [1,N]
$$
whose matrix is given by 
\begin{equation}
\label{rM-q}
\rbf(M^w)_{jk} = \sqrt{q}^{\, \lcor (w_{\leq j-1} + w)\vpi_{i_j}, (w_{\leq k-1} - w) \vpi_{i_k} \rcor },
\; \; \forall 1 \leq j < k \leq N. 
\end{equation}
The cluster variables $M^w(e_k)$ are also given by the expressions
\begin{align*}
M^w(e_k) &= 
\sqrt{q}^{\, \| (w_{\leq O_+(k)} - w_{\leq k-1}) \vpi_{i_k} \|^2/2} \, \Det_{w_{\leq k-1} \vpi_{i_k}, w \vpi_{i_k}}  \\
 &= \sqrt{q}^{\, \| (w_{[k,N]} - 1) \vpi_{i_k} \|^2/2} \, \Det_{w_{\leq k-1} \vpi_{i_k}, w  \vpi_{i_k}}.
\end{align*}

The frozen variables of the quantum cluster algebra structure that we will define will 
be indexed by 
\begin{equation}
\label{frozen-u}
\{ k \in [1,N] \mid p(k) = - \infty \}. 
\end{equation}
Note that this set has the same cardinality as $\Sbb(w)$. We will use the convention that
\begin{multline}
\label{convent2}
{\mbox{\em{the columns of all}}} \; \; n \times (N - |\Sbb(w)|) \; \; 
{\mbox{\em{matrices}}}
\\
{\mbox{\em{are indexed by}}} \; \;  
\{ k \in [1,N] \mid p(k) \neq - \infty \}. 
\end{multline}
Under this convention,
define the $N \times ( N - |\Sbb(w)|)$-matrix $\wt{B}^w$  \index{Btildew1@$\wt{B}^w$}  with entries
$$
b_{jk} = 
\begin{cases}
1, &\mbox{if} \; \; j = p(k)
\\
-1, &\mbox{if} \; \; j = s(k)
\\
c_{i_j i_k}, & \mbox{if} \; \; 
p(j) < p(k) < j < k
\\ 
- c_{i_j i_k}, & \mbox{if} \; \; 
p(k) < p(j) < k < j
\\
0, & \mbox{otherwise}.
\end{cases}
$$

\bth{Uw-main} Consider an arbitrary finite dimensional complex simple Lie algebra $\g$, a Weyl group element $w \in W$,
a reduced expression \eqref{reduced} of $w$, an arbitrary base field $\KK$ and a non-root of 
unity $q \in \kx$ such that $\sqrt{q} \in \KK$. Then $(M^w, \wt{B}^w)$ is a 
quantum seed and the quantum Schubert cell algebra $\UU^-[w]$ 
equals the quantum cluster algebra $\Abb(M^w, \wt{B}^w, \varnothing)_\KK$
with set of frozen variables indexed by \eqref{frozen-u}. Furthermore, this 
quantum cluster algebra equals the upper quantum cluster algebra 
$\UU(M^w, \wt{B}^w, \varnothing)_\KK$. For all $j \in [1,N]$ and
$m \in \Zset_{\geq 0}$ such that $s^m(j) \in [1,N]$, the elements
\begin{equation}
\label{Uw-int-prime}
\ol{y}_{[j, s^m(j)]} = 
{\sqrt{q}}^{\, \| (w_{\leq s^m(j)} - w_{\leq j-1})\vpi_{i_j} \|^2/2 }
\, \Det_{w_{\leq j-1} \vpi_{i_j}, w_{\leq s^m(j)} \vpi_{i_j}}
\end{equation}
are cluster variables of $\UU^-[w]$.
\eth

As mentioned, Gei\ss, Leclerc and Schr\"oer \cite{GLSh} constructed quantum cluster algebra structures 
on the algebras $\UU^-[w]$ in the case of symmetric Kac--Moody algebras $\g$. We note that the normalization scalars for the cluster variables in \thref{Uw-main} match the ones used in \cite{GLSh}. In the case when $w$ 
is the square of a Coxeter element and $\g$ is an arbitrary Kac--Moody algebra, a quantum cluster algebra 
structure on $\UU^-[w]$ was constructed by Berenstein and Rupel \cite{BR} simultaneously 
to our work. 

Because of the invariance of the bilinear form 
$\lcor.,.\rcor$ with respect to the Weyl group $W$, the scalar in \eqref{Uw-int-prime}
equals 
$$
{\sqrt{q}}^{\, \| (w_{[j, s^m(j)]} - 1)\vpi_{i_j} \|^2/2 }.
$$ 
\thref{Uw-main} is proved in Section \ref{10.3}. In Section \ref{10.2} we prove that the
toric frame $M^w$ and the matrix $\wt{B}^w$ are compatible.
Most of the conditions needed to apply \thref{cluster} are straightforward to 
verify for the algebra $\UU^-[w]$, except for the condition \eqref{pi-cond}. Set for brevity
$$
\sqrt{q_i} := \sqrt{q}^{\, \| \al_i \|^2/2 }, \; \; \forall i \in [1,r].
$$
\index{qiroot@$\sqrt{q_i}$}
We prove in Section \ref{10.3} that the condition \eqref{pi-cond} is satisfied after the rescaling
$$
x_k \mt \sqrt{q_{i_k}} (q_{i_k}^{-1} - q_{i_k}) F_{\be_k}, \; \; \forall k \in [1,N]
$$
of the standard generators \eqref{cano-gen} of $\, \UU^-[w]$ and the reversal of the order of these generators.

By setting $m=0$ in the last statement of \thref{Uw-main}, one obtains that 
the rescaled generators 
$$
\sqrt{q_{i_k}} (q_{i_k}^{-1} - q_{i_k}) F_{\be_k}
$$
of $\UU^-[w]$ are cluster variables for all $k \in [1,N]$.

\bre{other-cluster} The quantum seed of $\UU^-[w]$ in \thref{Uw-main} comes from the one in 
\thref{cluster} associated to the longest element $w_\ci \in \Xi_N$. More precisely, $M^w = M_{w_\ci} (w_\ci)_\bu$ if the generators are rescaled appropriately. 
Using \eqref{Uw-int-prime} and \prref{int-tau} it is straightforward to write down explicitly the 
toric frames of the quantum seeds of $\UU^-[w]$ associated 
to all elements of $\Xi_N$ via the construction of \thref{cluster}.
To compute the matrices for all those other seeds 
one needs to solve explicitly the system of linear equations \eqref{linear-eq}.

There is one additional difference between Theorems \ref{tcluster} and \ref{tUw-main}.
In the first case 
the cluster variables for the seed corresponding to $w_\ci \in \Xi_N$ 
(and as a matter of fact for all clusters corresponding to elements of $\Xi_N$)
are reenumerated according to the rule of Section \ref{6.1}. This is needed in order 
to match the combinatorics for the different seeds. In the case of 
\thref{Uw-main}, we do not perform this reenumeration, in order 
to match our results to the conventions in the existing literature. 
Because of this difference, the roles of the predecessor and successor functions 
in Theorems \ref{tcluster} and \ref{tUw-main} are interchanged. Details are given next.
\ere 

For the remainder of Chapter \ref{q-Schu}, we fix  the following choice of generators of $\UU^-[w]$: 
\begin{equation}
\label{new-x}
x_k = \sqrt{q_{i_k}} (q_{i_k}^{-1} - q_{i_k}) F_{\be_k}, \; \; \forall k \in [1,N],
\end{equation}
which as previously noted are rescalings of the canonical 
generators \eqref{cano-gen}. Since $\UU^-[w]$ is a symmetric CGL extension with respect to these generators, it  is also a symmetric CGL extension with respect to the presentation with the order of generators reversed:
\begin{equation}
\label{revpres}
\UU^-[w] = \KK[x_N] [x_{N-1}; \sigma^*_{N-1}, \delta^*_{N-1}] \cdots [x_1; \sigma^*_1, \delta^*_1].
\end{equation}
We will apply \thref{cluster} to this presentation as the \emph{starting} presentation of $\UU^-[w]$. To do so, we must verify that conditions \eqref{d-prop} and \eqref{pi-cond} hold for the presentation \eqref{revpres}.

We will need to identify the  data involved when \thref{cluster} is applied as above. (These data are not the same as those appearing at the stage $\tau = w_\ci$ when the starting presentation is the one with the generators in the usual order $x_1,\dots,x_N$.) First, we label the generators of \eqref{revpres} in ascending order in the form $x_{w_\ci,k} = x_{w_\ci(k)}$. By \coref{steps}, we can use $\eta_{w_\ci} := \eta w_\ci$ for the  $\eta$-function going with this indexing. The corresponding predecessor and successor functions are given by $p_{w_\ci} = w_\ci s w_\ci$ and $s_{w_\ci} = w_\ci p w_\ci$. The multiplicatively skew-symmetric $\lab$-matrix associated to \eqref{revpres} is $w_\ci \lab w_\ci$, where $\lab$ is the one associated to the presentation \eqref{UwCGL}, recall \eqref{Uwlakj}. Given the choice of $\nub$ above, we choose $w_\ci \nub w_\ci$ as the $\nub$-matrix associated to \eqref{revpres}.

We have already labelled the sequence of normalized prime elements for the presentation \eqref{revpres} coming from \thref{CGL} as $\ol{y}_{w_\ci,k}$, $k \in [1,N]$, and we will show later that $M^w(e_k) = \ol{y}_{w_\ci, w_\ci(k)}$ for all $k \in [1,N]$.

\prref{tauOre} shows that the elements of $\HH$ needed for the CGL conditions on \eqref{revpres} can be chosen to be $h^*_N, \dots, h^*_1$. Since $h^*_j . x_j = \la^*_j x_j$ for all $j \in [1,N]$, the singly-indexed $\la$-elements for \eqref{revpres}, which we denote $\la_{w_\ci,k}$, have the form 
$$
\la_{w_\ci,k} = \la^*_{w_\ci(k)} = q^2_{i_{w_\ci(k)}}, \; \; \forall k \in [1,N],
$$
recall \eqref{lalastar}.  \index{lambdawcirclek@$\la_{w_\ci,k}$}
The set of exchangeable indices for the presentation \eqref{revpres} is
\begin{align*}
\ex_\ci &:= \{ k \in [1,N] \mid s_{w_\ci}(k) \ne +\infty \} = \{ k \in [1,N] \mid p w_\ci (k) \ne -\infty \}  \\
 &\phantom{:}= \{ k \in [1,N] \mid s (w_\ci)_\bu w_\ci(k) \ne +\infty \} = ((w_\ci)_\bu w_\ci)^{-1} (\ex).
\end{align*}
\index{excircle@$\ex_\ci$}
In the third equality, we have used the observation that $(w_\ci)_\bu$ acts by reversing the order of the elements in each level set of $\eta$.
Taking account of \prref{la-equal}, we see that
\begin{equation}
\label{lawostar}
\la^*_{w_\ci,k} := (\la_{w_\ci,k})^* = \la^{-1}_{w_\ci, w_\ci pw_\ci(k)} = q^{-2}_{i_{pw_\ci(k)}} = q^{-2}_{i_{w_\ci(k)}}, \; \; \forall k \in \ex_\ci .
\end{equation}
In particular, it follows that
$$
\bigl( \la^*_{w_\ci, l} \bigr)^{ \| \al_{i_{w_\ci(j)}} \|^2} = q^{ -  \| \al_{i_{w_\ci(j)}} \|^2 \| \al_{i_{w_\ci(l)}} \|^2} = \bigl( \la^*_{w_\ci, j} \bigr)^{ \| \al_{i_{w_\ci(l)}} \|^2}, \; \; \forall j,l \in \ex_\ci ,
$$
which verifies \eqref{d-prop} for the presentation \eqref{revpres}.

\ble{new10.5}
Let $j \in [1,N]$.

{\rm(a)} If $m \in \Znn$ such that $s^m(j) \in [1,N]$, then \eqref{Uw-int-prime} holds:
$$
\ol{y}_{[j, s^m(j)]} = 
{\sqrt{q}}^{\, \| (w_{\leq s^m(j)} - w_{\leq j-1})\vpi_{i_j} \|^2/2 }
\, \Det_{w_{\leq j-1} \vpi_{i_j}, w_{\leq s^m(j)} \vpi_{i_j}} .
$$

{\rm(b)} If $k = w_\ci(j)$, then
$$
\ol{y}_{w_\ci,k} = \sqrt{q}^{ \, \| (w-w_{\le j-1}) \vpi_{i_j} \|^2/2 } \, \Det_{ w_{\le j-1} \vpi_{i_j}, w \vpi_{i_j} } .
$$

{\rm(c)} If the presentation \eqref{revpres} satisfies condition \eqref{pi-cond}, then
$$
M^w = M_{w_\ci} (w_\ci)_\bu .
$$
\ele

\begin{proof} (a) Using \coref{Uw-prim-int}, the form of the rescaling, 
and the equality 
$\sqrt{q_{i_k}} = \sqrt{q}^{\, \| \al_{i_k} \|^2/2}= \sqrt{q}^{\, \|\be_k\|^2/2} $, we obtain
$\ol{y}_{[j, s^m(j)]} = \Theta \Det_{w_{\leq j-1} \vpi_{i_j}, w_{\leq s^m(j)} \vpi_{i_j}}$ where
\begin{align*}
\Theta &=
\biggl( \, \prod_{0 \leq l < n \leq m} {\sqrt{q}}^{ \, - \lcor \be_{s^l(j)} , \be_{s^n(j)} \rcor} \biggr) \sqrt{q_{i_j}}^{ \, m+1} \biggl( \, \prod_{0 \leq l < n \leq m} q^{ \, \lcor \be_{s^l(j)} , \be_{s^n(j)} \rcor} \biggr)   \\
 &= \sqrt{q}^{ \, (m+1) \| \be_j \|^2/2 } \biggl( \, \prod_{0 \leq l < n \leq m} {\sqrt{q}}^{ \, \lcor \be_{s^l(j)} , \be_{s^n(j)} \rcor} \biggr) = {\sqrt{q}}^{ \, \| \be_j + \be_{s(j)} + \cdots + \be_{s^m(j)} \|^2/2}  \\
 &=
{\sqrt{q}}^{ \, \| (w_{\leq s^m(j)} - w_{\leq j-1})\vpi_{i_j} \|^2/2} .
\end{align*}

(b) \prref{int-tau} shows that $\ol{y}_{w_\ci,k} = \ol{y}_{[j, s^m(j)]}$ where
$$
m = \max \{ n \in \Znn \mid s^n(j) \in [j,N] \} = O_+(j).
$$
Thus $w_{\le s^m(j)} \vpi_{i_j} = w \vpi_{i_j}$, so the desired equation follows from part (a).

(c) Under the given assumptions, it follows from part (b) and the definition of $M^w$ that
\begin{equation}
\label{MwMwo}
M^w(e_k) = \ol{y}_{w_\ci, w_\ci(k)} = \wh{M}_{w_\ci}(e_{w_\ci(k)}) = M_{w_\ci} (w_\ci)_\bu (e_k), \; \; \forall k \in [1,N] .  \qedhere
\end{equation}
\end{proof}

\section{Compatibility of the toric frame $M^w$ and the matrix $\wt{B}^w$}
\label{10.2}
Recall that the rational character lattice $\xh$ of the torus
$\HH = (\kx)^r$ can be identified with $\QQ$ via \eqref{Hchar}. Thus, we now view $\xh$ as an additive group with identity $0$.

\bpr{Uw-compat} In the setting of Theorem {\rm\ref{tUw-prim}}, the multiplicatively skew-symmetric matrix $\rbf(M^w)$ and the matrix $\wt{B}^w$ are compatible, and more precisely
\begin{equation}
\label{Om-be}
\Om_{\rbf(M^w)}(b^k, e_l) = q_{i_k}^{\; -\delta_{kl}}, 
\; \; 
\forall k ,l \in [1,N], \; p(k) \neq - \infty.
\end{equation}

Moreover, the columns 
$b^k$ of the matrix $\wt{B}^w$ satisfy
\begin{equation}
\label{b-hom}
\chi_{M^w(b^k)} = 0, \; \; \forall k \in [1,N], \; p(k) \neq - \infty.
\end{equation}
\epr

\begin{proof} We derive the proposition from results of \cite{BFZ,BZ}.
This appears to be the shortest way to prove the proposition, although 
some constructions that match our setting to the setting of \cite{BFZ,BZ} 
might appear to be a bit artificial. To avoid this, one could prove the 
proposition directly following the idea of \cite{BFZ,BZ}. We leave this to the 
interested reader.

Define the function $\eta\spcheck : [1,N+r] \to [1,r]$ by
$$
\eta\spcheck(j) = i_j, \; \; \forall j \in [1,N]
\quad \mbox{and} 
\quad \eta\spcheck(N+i) = i, \; \; \forall i \in [1,r]. 
$$
The predecessor map $p \spcheck : [1,N+r] \to [1,N+r] \sqcup \{ - \infty \}$ 
for the level sets of $\eta\spcheck$ coincides with $p$ on $[1,N]$ and satisfies 
$$
p\spcheck(N+i) := \begin{cases} 
\max \{j \in [1,N] \mid i_j = i \}, & \mbox{if such exists}
\\
- \infty, & \mbox{otherwise}
\end{cases}
$$
for $i \in [1,r]$.
Similarly, the successor map $s\spcheck$ 
for the level sets of $\eta\spcheck$ is given by
$$
s\spcheck(j) = 
\begin{cases} 
s(j), & \mbox{if} \; \; 
s(j) \neq + \infty
\\
N+i_j, & \mbox{if} \; \; s(j) = + \infty
\end{cases}
$$
for $j \in [1,N]$ and $s\spcheck(N+i)= + \infty$ for $i \in [1,r]$. Our setting differs 
from \cite{BFZ,BZ} in that we place additional indices at the end of a word for a 
Weyl group element, not at the beginning. The reason for this is that the results 
of \cite{BFZ,BZ} will be applied to the inverse of a certain reduced expression. 

We extend $\wt{B}^w$ to an $(N+r) \times (N - |\Sbb(w)|)$-matrix
by setting 
$$
b_{N+i,k}:= 
\begin{cases}
-1, & \mbox{if} \; \; s\spcheck(k)=N+i
\\
-c_{i,i_k}, & \mbox{if} \; \; 
p(k) < p\spcheck(N+i) < k
\\
0, & \mbox{otherwise}
\end{cases}
$$
for all $i \in [1,r]$ and $k \in [1,N]$, $p(k) \neq - \infty$,
using the convention \eqref{convent2}.
Applying \cite[Theorem 8.3]{BZ} for the Weyl
group elements $u := 1$, $v :=w^{-1}$ and
the double word corresponding to the
reduced expression $s_N \ldots s_1$ which is reverse to \eqref{reduced}
leads to
\begin{multline}
\label{long-ident}
\sum_{j=1}^N b_{jk} \sign(j-l) 
\bigl( \lcor \vpi_{i_j}, \vpi_{i_l} \rcor
- \lcor w_{\leq j-1} \vpi_{i_j}, w_{\leq l-1}\vpi_{i_l} \rcor \bigr)
\\
- \sum_{i=1}^r b_{N+i,k} \lcor w \vpi_i, (w_{\leq l-1}-w)\vpi_{i_l} \rcor 
= - \delta_{kl} \lcor \al_{i_l}, \al_{i_l} \rcor 
\end{multline}
for $k,l \in [1,N]$, $p(k) \ne -\infty$.

Using the fact that the one-step mutations in \cite{BZ}
are eigenfunctions of both the left and right regular actions 
of the maximal torus of the connected, simply connected 
complex algebraic group with Lie algebra $\g$, we obtain that 
\begin{align}
&\sum_{j=1}^N b_{jk} \vpi_{i_j} + \sum_{i=1}^r b_{N+i,k} \vpi_i = 0
\label{bz1}
\\
&\sum_{j=1}^N b_{jk} w_{\leq j-1} \vpi_{i_j} + \sum_{i=1}^r b_{N+i,k} w \vpi_i = 0
\label{bz2}
\end{align}
for all $k \in [1,N]$, $p(k) \neq - \infty$.

Acting by $w$ on \eqref{bz1} and subtracting \eqref{bz2} gives
$$
\chi_{M(b^k)} = \sum_{j=1}^N b_{jk} (w - w_{\leq j-1}) \vpi_{i_j}= 0,
\; \; \forall k \in [1,N], \; p(k) \neq - \infty,
$$
which proves \eqref{b-hom}. Similarly, \eqref{bz1} and \eqref{bz2} 
imply
\begin{align*}
& \sum_{j=1}^N b_{jk} \lcor ( w_{\leq j-1} - w) \vpi_{i_j}, w \vpi_{i_l} \rcor = 0 \quad 
\mbox{and}
\\
& \sum_{j=1}^N b_{jk} \lcor w \vpi_{i_j}, (w_{\leq l-1} - w) \vpi_{i_l} \rcor 
+ \sum_{i=1}^r b_{N+i, k} \lcor w \vpi_i, (w_{\leq l-1} - w) \vpi_{i_l} \rcor = 0
\end{align*}
for the same values of $k$ and all $l \in [1,N]$.
Combining these two identities, \eqref{long-ident}, and the fact that 
$\Om_{\rbf(M^w)}( e_j, e_l) = \sqrt{q}^{\, m_{jl}}$ where
\begin{multline*}
m_{jl}:= \sign(j-l) \bigl(  \lcor \vpi_{i_j}, \vpi_{i_l} \rcor -
\lcor w_{\leq j-1} \vpi_{i_j}, w_{\leq l-1}\vpi_{i_l} \rcor \bigr)
\\
- \lcor (w_{\leq j-1} - w) \vpi_{i_j}, w \vpi_{i_l} \rcor +
\lcor w \vpi_{i_j}, (w_{\leq l-1} - w) \vpi_{i_l} \rcor, \; \; \forall j,l \in [1,N]
\end{multline*}
(which follows from \eqref{rM-q}), we obtain $\Om_{\rbf(M^w)}( b^k, e_l) = \sqrt{q}^{\, -\delta_{kl} \| \al_{i_k} \|^2}$ and thus \eqref{Om-be}.
\end{proof}

Up to a possible additional rescaling, \thref{cluster} applied to \eqref{revpres} produces an initial quantum seed that we will label $(M_\ci, \wt{B}_\ci)$.  \index{Mcircle@$M_\ci$}  \index{Bcircletilde@$\wt{B}_\ci$}  For $k \in [1,N]$, the toric frame $M_\ci$ sends $e_k$ to $\ol{y}_{w_\ci, k}$ (possibly rescaled), so we see from \leref{new10.5} (b) and the definition of $M^w$ that up to rescaling, $M_\ci$ agrees with $M^w w_\ci$. Since $\wt{B}_\ci$ is not affected by any rescaling, recall \reref{rescaleBtil}, we may use $M^w w_\ci$ in place of $M_\ci$ to determine $\wt{B}_\ci$, as follows.

\bco{Botil} $\wt{B}_\ci = w_\ci \wt{B}^w \bigl( w_\ci |_{\ex_\ci} \bigr).$
\eco

\begin{proof} Recalling \eqref{Msig}, we have
$$
\rbf(M^w w_\ci) = {}^{w_\ci^T} \rbf(M^w)^{w_\ci} = w_\ci \rbf(M^w) w_\ci .
$$
Hence, \prref{Uw-compat} implies that 
\begin{multline*}
\Om_{\rbf(M^w w_\ci)} ( w_\ci b^{w_\ci(k)}, e_l) = \Om_{\rbf(M^w)} ( b^{w_\ci(k)}, e_{w_\ci(l)}) = q^{-\delta_{kl}}_{i_{w_\ci(k)}} = \sqrt{\la^*_{w_\ci,k}}^{ \, \delta_{kl}},  \\
 \; \; \forall k \in \ex_\ci, \; l \in [1,N],
\end{multline*}
and also
$$
\chi_{(M^w w_\ci)( w_\ci b^{w_\ci(k)})} = 0, \; \; \forall k \in \ex_\ci .
$$
Recall from \reref{rescaleBtil} that neither $\rbf(M_\ci)$ nor $\chi_{M_\ci(b)}$, for $b\in \Zset^N$, changes under a rescaling of the generators $x_{w_\ci, k}$. Consequently, the two previous equations also hold when $M^w w_\ci$ is replaced by $M_\ci$. \thref{cluster} (a) thus implies that $w_\ci b^{w_\ci(k)}$ equals the $k$-th column of $\wt{B}_\ci$, for all $k \in \ex_\ci$, as desired.
\end{proof}

\newcommand{\tuwmain}{tUw-main}
\section{Proof of Theorem {\ref{\tuwmain}}}
\label{10.3}
Let us write $y_{w_\ci, [i,s^m(i)]}$ and $u_{w_\ci, [i,s^m(i)]}$ for the ``interval prime elements" and corresponding $u$-elements obtained from \thref{y-int} and \coref{u-elem} applied to the CGL extension \eqref{revpres}. In particular,
$$
u_{w_\ci, [i,s(i)]} = x_{w_\ci, i} x_{w_\ci, s(i)} - y_{w_\ci, [i,s(i)]}, \; \; \forall i \in [1,N], \; s(i) \ne +\infty.
$$
As in eq. \eqref{pi-f}, write
$$
\lt( u_{w_\ci, [i,s(i)]} ) = \pi_{w_\ci, [i,s(i)]} (x_{w_\ci})^{ f_{w_\ci, [i,s(i)]} }, \; \; \forall i \in [1,N], \; s(i) \ne +\infty,
$$
for some $\pi_{w_\ci, [i,s(i)]} \in \kx$ and $f_{w_\ci, [i,s(i)]} \in \sum_{j=i+1}^{s(i)-1} \Znn e_j$, where the leading term and the monomial $( x_{w_\ci})^{ f_{w_\ci, [i,s(i)]} }$ are computed with respect to the order of generators $x_N,\dots,x_1$.

\ble{fwo1N} Assume that in the setting of Theorem {\rm\ref{tUw-main}},
$$
w = s_{i_1} \ldots s_{i_N}
$$
is a reduced expression such that $i_1 = i_N = i$ and $i_k \neq i$, 
$\forall k \in [2,N-1]$, for some $i \in [1,r]$. Then 
$$
f_{w_\ci, [1,N]} = - \sum_{j=2}^{N-1} c_{i_{w_\ci(j)} i} e_j ,
$$
so that
\begin{equation}
\label{ltuwo1N}
\lt( u_{w_\ci, [1,N]} ) = \pi_{w_\ci, [1,N]} x_{N-1}^{- c_{i_{N-1} i}} \cdots x_2^{ - c_{i_2 i}}.
\end{equation}
\ele

\begin{proof} In view of \coref{Botil}, the entries in the first column of $\wt{B}_\ci$ satisfy
$$
\bigl( \wt{B}_\ci \bigr)_{j1} = \bigl( \wt{B}^w \bigr)_{w_\ci(j),N} = c_{i_{w_\ci(j)} i} , \; \; \forall j \in [2,N-1].
$$
Since the vector $f_{w_\ci, [1,N]}$ does not change under a rescaling of the generators $x_{w_\ci, k}$, we obtain the desired conclusion by applying \prref{f1N} to the (possibly rescaled) presentation \eqref{revpres}.
\end{proof}

Note that in the setting of \leref{fwo1N}, we have $\be_1 = \al_i$, $\be_N = - w \al_i$ and 
$$
y_{w_\ci, N} = q_i \Det_{\vpi_i, w \vpi_i}
$$
because of \thref{Uw-prim} and the effect of the rescaling 
of the generators \eqref{new-x}. Hence,
\begin{equation}
\label{uwo1N}
u_{w_\ci, [1,N]} = x_N x_1 - q_i  \Det_{\vpi_i, w \vpi_i} .
\end{equation}

For use in the next proof, we note the following identities for the braid group action 
in the case of $\g = \slfrak_2$:
\begin{equation}
E_1^m T_1^{-1} v_{m \vpi_1} = 
[m]_q! \, v_{m \vpi_1}
\quad \mbox{and} \quad
T_1^2 v_{m \vpi_1} = \frac{(-q)^m}{([m]_q!)^2} v_{m \vpi_1} , \; \; \forall m \in \Znn.
\label{bg}
\end{equation}
They are easily deduced from the standard facts (\cite[\S 8.6 and Lemma 1.7]{Ja}) 
for the $\UU_q(\slfrak_2)$-braid group action. We leave this to the reader.

\ble{expan} Under the assumptions of Lemma {\rm\ref{lfwo1N}},
$$
\pi_{w_\ci, [1,N]} = \Scr_{w_\ci \nub w_\ci}( -e_1 + f_{w_\ci, [1,N]} ).
$$
\ele

\begin{proof} Using \eqref{fg}, we see that
\begin{align*}
\Scr_{w_\ci \nub w_\ci}( -e_1 + f_{w_\ci, [1,N]} ) &= \Om_{w_\ci \nub w_\ci}( e_1, f_{w_\ci, [1,N]}) \Scr_{w_\ci \nub w_\ci}( f_{w_\ci, [1,N]} )  \\
&= \Om_\nub( e_N, f ) \Scr_\nub( f )^{-1} =  \Om_\nub( f, -e_N ) \Scr_\nub( f )^{-1}  \\
&= \Scr_\nub( -e_N + f )^{-1},
\end{align*}
where
$$
f :=  w_\ci f_{w_\ci, [1,N]} = - \sum_{j=2}^{N-1} c_{i_j i} e_j .
$$
Hence, what must be proved is that
\begin{equation}
\label{pigoal}
\pi_{w_\ci, [1,N]} = \Scr_\nub( -e_N + f )^{-1}.
\end{equation}

Since $x_N x_1$, $y_{w_\ci, N}$, and $u_{w_\ci, [1,N]}$ are homogeneous of the same degree with respect to the $\QQ$-grading of $\UU^-[w]$, 
\begin{equation}
\label{al-be}
-\be_N -\be_1 = (w-1) \vpi_i = (-c_{i_{N-1} i}) (-\be_{N-1}) +\cdots+ (- c_{i_2 i}) (-\be_2) ,
\end{equation}
recall \eqref{defDet}.
Therefore, 
\begin{align*}
\Scr_\nub(-e_N + f)^{-1} &= 
\sqrt{q}^{ \, \lcor \be_N, c_{i_2 i} \be_2 + \cdots + c_{i_{N-1} i} \be_{N-1} \rcor }
\prod_{1<j<k<N} \sqrt{q}^{\,  c_{i_j i} c_{i_k i} \lcor \be_j, \be_k \rcor }
\\
&= 
\sqrt{q}^{\, \lcor \al_i, (w-1) \al_i \rcor}
\prod_{1<j<k<N} \sqrt{q}^{\, c_{i_j i} c_{i_k i} \lcor \be_j, \be_k \rcor } ,
\end{align*}
because $\lcor -w\al_i, (w-1) \al_i \rcor =  \lcor \al_i, (w-1) \al_i \rcor$.
Taking the square length of the vector in \eqref{al-be} and 
using that $\|\be_k \|^2 = \| \al_{i_k} \|^2$ leads to 
\begin{align*}
\sum_{k=2}^{N-1} c_{i_k i}^2 \| \al_{i_k} \|^2 
+ 2 \sum_{1 < j < k < N} c_{i_j i} c_{i_k i} \lcor \be_j, \be_k \rcor
&= \lcor (1-w) \al_i, (1-w) \al_i \rcor 
\\
&= 2 \lcor \al_i, (1-w) \al_i \rcor.
\end{align*}
Hence,
\begin{equation}
\label{Ste}
\Scr_\nub( - e_N + f)^{-1} = \prod_{k=2}^{N-1} \sqrt{q_{i_k}}^{ \; - c_{i_k i}^2} 
\end{equation}
because $\sqrt{q_{i_k}} = \sqrt{q}^{ \, \| \al_{i_k} \|^2/2}$ for all $k \in [1,N]$.

Denote for brevity $\xi_\mu := \xi_{1, \mu}$ 
for $\mu \in \PP^+$ (i.e., $\xi_\mu \in (V(\mu)^*)_{-\mu}$ is the 
unique vector such that $\lcor \xi_\mu, v_\mu \rcor =1$). 
The explicit form \eqref{phi1} of the antihomomorphism $\vp_w$ 
and the rescaling of the generators \eqref{new-x} imply
\begin{equation}
\label{piU}
\pi_{w_\ci, [1,N]} = - q_i
\frac{
\lcor \xi_{\vpi_i}, (\tau E_{\be_2})^{-c_{i_2 i}} \ldots (\tau E_{\be_{N-1}})^{-c_{ i_{N-1} i}} T^{-1}_{w^{-1}} 
v_{\vpi_i} \rcor}
{\prod_{j=2}^{N-1}
{\sqrt{q_{i_j}}^{ \; c_{i_j i}^2} [-c_{i_j i}]_{q_{i_j}}! }}
\cdot
\end{equation}
For all $i' \in [1,N]$ with $i' \neq i$, the element $T_i^{-1} v_{\vpi_i}$ is a 
highest weight vector for the $\UU_{q_{i'}}(\slfrak_2)$-sub\-al\-ge\-bra 
of $\UU_q(\g)$ generated by $\{ E_{i'}, F_{i'}, K_{i'}^{-1} \}$ 
with highest weight 
$$\lcor s_i \vpi_i, \al_{i'}\spcheck \rcor \vpi_{i'} = \lcor \vpi_i - \al_i, \al_{i'}\spcheck \rcor \vpi_{i'} 
= - c_{i' i} \vpi_{i'},$$
because $E_{i'} F_i^m v_{\vpi_i} = F_i^m E_{i'} v_{\vpi_i} = 0$ for all $m \in \Zset_{\geq 0}$.  
Applying \eqref{bg}, we obtain
$$
E_{i_j}^{-c_{i_j i}} T_{i_j}^{-1} (T_i^{-1} v_{\vpi_i}) = 
[-c_{i_j i}]_{q_{i_j}}! \, T_i^{-1} v_{\vpi_i}, \; \; \forall j \in [2,N-1].
$$
Using that $\tau(E_{\be_j}) = T_{i_1}^{-1} \ldots T_{i_{j-1}}^{-1}(E_{i_j})$ and
repeatedly applying the above identity gives
\begin{multline*}
\lcor \xi_{\vpi_i}, (\tau E_{\be_2})^{-c_{i_2 i}} \ldots (\tau E_{\be_{N-1}})^{-c_{ i_{N-1} i}} T^{-1}_{w^{-1}} 
v_{\vpi_i} \rcor = 
\\
= \lcor \xi_{\vpi_i}, T_i^{-2} v_{\vpi_i} \rcor 
\prod_{j=2}^{N-1} [-c_{i_j i}]_{q_{i_j}}! 
=
- q_i^{-1} \prod_{j=2}^{N-1} [-c_{i_j i}]_{q_{i_j}}!. 
\end{multline*}
For the last equality we use the second identity in \eqref{bg} for $m=1$.
Hence,
$$
\pi_{w_\ci, [1,N]} = 
\prod_{j=2}^{N-1}
\sqrt{q_{i_j}}^{ \; - c_{i_j i}^2}
= \Scr_\nub( - e_N + f)^{-1},
$$
which verifies \eqref{pigoal} and completes the proof of the proposition.
\end{proof}

\begin{proof}[Proof of \thref{Uw-main}] Applying \leref{expan} to each of the interval subalgebras
$$
\KK[x_{w_\ci, i}] [x_{w_\ci, i+1}; \sigma^*_{w_\ci(i+1)}, \delta^*_{w_\ci(i+1)}] \cdots [x_{w_\ci, s_{w_\ci}(i)}; \sigma^*_{w_\ci s_{w_\ci}(i)}, \delta^*_{w_\ci s_{w_\ci}(i)}],
$$
for $i \in [1,N]$ such that $s_{w_\ci}(i) \ne +\infty$, we conclude that the symmetric CGL extension \eqref{revpres} satisfies condition \eqref{pi-cond}. Thus, the hypotheses of \thref{cluster} hold for \eqref{revpres}. Eq. \eqref{Uw-int-prime} has been verified in \leref{new10.5} (b).

By \thref{cluster}, $(M_\ci, \wt{B}_\ci)$ is a quantum seed and
$$
\UU^-[w] = \Abb( M_\ci, \wt{B}_\ci, \varnothing)_\KK = \UU( M_\ci, \wt{B}_\ci, \varnothing)_\KK ,
$$
with set of frozen variables indexed by $\ex_\ci$. Applying this result to $\UU^-[w_{\le s^m(j)}]$, for $j \in [1,N]$ and $m \in \Znn$ such that $s^m(j) \in [1,N]$, we see that $\ol{y}_{[j,s^m(j)]}$ is a cluster variable of $\UU^-[w_{\le s^m(j)}]$ and hence a cluster variable of $\UU^-[w]$. 

Since $M_\ci(e_k) = \ol{y}_{w_\ci, k} = M^w(e_{w_\ci(k)})$ for all $k \in [1,N]$, we have $M_\ci = M^w w_\ci$. Taking \coref{Botil} into account, we therefore obtain \thref{Uw-main}.
\end{proof}

\backmatter

\bibliographystyle{amsalpha}


\printindex

\end{document}